\documentclass[a4paper,11pt]{amsart}
\usepackage{graphicx}
\usepackage{subcaption}
\usepackage{float}
\usepackage{amsmath,amsthm,amssymb}
\usepackage{slashed}
\usepackage[all]{xy}
\usepackage{color}
\usepackage{comment}
\newcommand{\mychoice}[2]{#1
}
\mychoice{
\newcommand{\plabel}[1]{ \label{#1}}
\newcommand{\gbibitem}[1]{ \bibitem{#1}}
\newcommand{\snewpage}{}

\specialcomment{commentx}{   }{   }\excludecomment{commentx}
}{
\newcommand{\plabel}[1]{ \label{#1}\rlap{\smash{${}^{^{[#1]}}$}}}
\newcommand{\gbibitem}[1]{ \bibitem{#1}\rlap{\smash{${}^{^{[#1]}}$}}}
\newcommand{\snewpage}{\newpage}

\newenvironment{commentx}{\color{magenta} }{\color{black} }
}
\DeclareMathOperator{\sgn}{sgn}
\DeclareMathOperator{\Id}{Id}
\DeclareMathOperator{\tr}{tr}
\DeclareMathOperator{\dist}{dist}
\DeclareMathOperator{\proj}{proj}
\DeclareMathOperator{\arcosh}{arcosh}
\DeclareMathOperator{\arsinh}{arsinh}

\DeclareMathOperator{\SO}{SO}
\DeclareMathOperator{\smc}{sc}
\DeclareMathOperator{\SL}{SL}
\DeclareMathOperator{\AN}{AN}
\DeclareMathOperator{\ATT}{AT}
\DeclareMathOperator{\AC}{AC}
\DeclareMathOperator{\PD}{PD}
\DeclareMathOperator{\CD}{CD}
\DeclareMathOperator{\MS}{MS}
\DeclareMathOperator{\MC}{MC}
\DeclareMathOperator{\ext}{ext}
\DeclareMathOperator{\modu}{\,mod\,}
\DeclareMathOperator{\Sinh}{\slashed{\mathrm S}inh}
\DeclareMathOperator{\Cosh}{\slashed{\mathrm C}osh}
\DeclareMathOperator{\Sin}{\slashed{\mathrm S}in}
\DeclareMathOperator{\Cos}{\slashed{\mathrm C}os}
\DeclareMathOperator{\intt}{int}
\DeclareMathOperator{\conv}{conv}
\DeclareMathOperator{\intD}{\mathring D}
\newcommand{\real}{\mathrm{real}}

\DeclareMathOperator{\cq}{cq}

\newcommand{\upper}{\overline{\mathbb C}^{+}}
\DeclareMathOperator{\Dbar}{D}
\DeclareMathOperator{\re}{Re}
\DeclareMathOperator{\ima}{Im}
\DeclareMathOperator{\Rea}{Re}
\DeclareMathOperator{\Ima}{Im}
\DeclareMathOperator{\Discr}{Discr}
\newcommand{\logst}{\log^{\mathrm{st}}}
\DeclareMathOperator{\Rexp}{exp_{R}}
\DeclareMathOperator{\Lexp}{exp_{L}}

\DeclareMathOperator{\CRext}{CR^{ext}}
\DeclareMathOperator{\CR}{CR}
\DeclareMathOperator{\ccCR}{\boldsymbol{CR}}
\newcommand{\ccCRext}{\ccCR^{\mathrm{ext}}}
\DeclareMathOperator{\DW}{DW}
\DeclareMathOperator{\spec}{sp}
\DeclareMathOperator{\Dual}{\mathcal D}
\DeclareMathOperator{\dett}{\mathbf{det}}

\theoremstyle{definition}
\newtheorem{point}{}[section]

\newtheorem{disc}[point]{Discussion}
\newtheorem{remark}[point]{Remark}
\newtheorem{example}[point]{Example}
\newtheorem{exaprop}[point]{Example (Proposition)}
\theoremstyle{plain}
\newtheorem{prop}[point]{Proposition}
\newtheorem{lemma}[point]{Lemma}
\newtheorem{cor}[point]{Corollary}
\newtheorem{theorem}[point]{Theorem}
\newcommand{\bem}{\begin{bmatrix}}
\newcommand{\eem}{\end{bmatrix}}

\newenvironment{bsmallmatrix}{\left[\begin{smallmatrix}}{\end{smallmatrix}\right]}
\newcommand{\leaveout}[1]{}

\newcommand{\qedexer}{  \renewcommand{\qedsymbol}{$\diamondsuit$} \qed \renewcommand{\qedsymbol}{$\Box$}}
\newcommand{\qedremark}{  \renewcommand{\qedsymbol}{$\triangle$} \qed \renewcommand{\qedsymbol}{$\Box$}}
\newcommand{\qedno}{\renewcommand{\qedsymbol}{}}
\newcommand{\proofremark}[1]{
\begin{proof}[Remark] #1
\renewcommand{\qedsymbol}{}
\end{proof}
}
\newcommand{\proofremarkqed}[1]{
\begin{proof}[Remark] #1
\renewcommand{\qedsymbol}{$\triangle$}
\end{proof}
\renewcommand{\qedsymbol}{$\Box$}
}
\newcommand{\eqed}{
\pushQED{\qed}
\qedhere
\popQED
}

\newcommand{\marginextend}[1]{ \addtolength{\oddsidemargin}{-#1}  \addtolength{\evensidemargin}{-#1}
  \addtolength{\textwidth}{#1}\addtolength{\textwidth}{#1}}
\newcommand{\updownextend}[1]{ \addtolength{\topmargin}{-#1}  \addtolength{\textheight}{#1}
\addtolength{\textheight}{#1}}
\marginextend{1cm}
\updownextend{0cm}
\allowdisplaybreaks[4]
\title{Convergence estimates for the Magnus expansion II. $C^*$-algebras}
\author{Gyula Lakos}
\address{Department of Geometry, Institute of Mathematics, E\"otv\"os University, P\'azm\'any P\'eter s.~1/C,  Budapest, H--1117, Hungary}
\email{lakos@cs.elte.hu}
\begin{document}
\begin{abstract}
We  review and provide simplified proofs related to the Magnus expansion, and improve convergence estimates.
Observations and improvements concerning the Baker--Campbell--Hausdorff expansion are also made.

In this Part II, we consider the case of $C^*$-algebras, i.~e.~essentially the case of operators on Hilbert spaces.
We present the spectral approach to the Magnus expansion in the context of the conformal range
(which is a projection of the Davis--Wielandt shell), allowing a more effective approach.
This makes possible to clarify certain convergence properties of the BCH expansion related to the critical cumulative norm $\pi$.
In particular, we prove that for finite dimensional matrices $A,B$,
the norm condition $\|A\|_2+\|B\|_2\leq\pi$ implies that the BCH expansion of $A$ and $B$ is convergent.
Several counterexamples regarding convergence of the Magnus and BCH expansions are presented.
In the rest, we prove growth estimates for the Magnus expansion in the setting of Hilbert space operators,
both in terms of the overall sum and the individuals terms.

\end{abstract}
\maketitle
\snewpage
\section*{Introduction to Part II}
This paper is a continuation of Part I, \cite{L1}.
We assume general familiarity with the results presented there, although the techniques we use here are quite different.
General sources for algebra, analysis and combinatorics should also be taken from there.

\textbf{Introduction to the Hilbert space operator setting.}
Practically, most applications of the Magnus expansion deal with matrices or Hilbert space operators.
A common feature of these cases is that certain spectral methods apply.
In fact, not only the available tools are more numerous in these cases, but stronger convergence theorems hold.
These are the subjects of this Part II.

Convergence in the case of finite matrices is more like a matter of complex analysis, as finite dimensional norms are equivalent.
Even so, it a very reasonable choice to use the operator norm in order to keep track of matters of convergence.
In the infinite-dimensional case, the choice of norm is, of course, essential.

Matrices with operator norm and Hilbert space operators have always been the principal subjects of investigations
of convergence of the Baker--Campbell--Hausdorff and Magnus expansions  (cf.~references in Part I).
However, the first result which is truly specific to these classes seems to be of Mityagin \cite{Mity} (1990), unfortunately unpublished
(cf. Day, So, Thompson \cite{DST}, Blanes, Casas \cite{BC}, and, especially, Casas \cite{Ca}).
It uses spectral arguments and establishes the convergence of BCH series with cumulative norm less than $\pi$.
The result was extended to the Magnus expansion by Moan, Niesen \cite{MN} (2008) for finite matrices,
and by Casas \cite{Ca}, ultimately, for operators on Hilbert spaces.

Divergence was considered earlier by Wei \cite{W} (who does not actually considers the norm)
and Michel \cite{Mi} (who uses different norms but his results and methods
have consequences in the Hilbert setting of $2\times2$ matrices); then  Vinokurov \cite{V} (1991)
gave a simple counterexample to the convergence of the BCH expansion with cumulative norm greater than $\pi$ in the BCH case.
Moan \cite{Ma} (2002) (cf.  Moan, Niesen \cite{MN}) gave a counterexample to the convergence of the Magnus expansion with cumulative norm $\pi$.
It is notable, however, that, regarding the Magnus expansion, all the crucial ideas and examples
were already available to Sch\"affer \cite{Sch} (1964), except he does not apply them to the convergence radius of the Magnus expansion as such.
In any case, the qualitative convergence radius $\pi$ is  well-established for Hilbert space operators.
\snewpage

\textbf{Outline of content.}
In this paper, we refine the results above.
In Section \ref{sec:SpectA} we present the basic idea of the spectral approach.
Section \ref{sec:ConformalRange} provides an introduction the conformal range of operators on Hilbert spaces, which is a reduced version (in fact, a projection) of the Davis--Wielandt shell.
We present the above mentioned convergence results in Section \ref{sec:MagnusHilbert}.
As we use the logarithmic variation instead of the angular variation, this allows us
a slightly more effective extension to the BCH expansion in the critical case.
In particular, we prove that if $A_1$, $A_2$ are linear operators on a finite dimensional Hilbert
space with $\|A_1\|_2+\|A\|_2=\pi$, then the converge radius of the BCH expansion is greater than $1$.
In Section \ref{sec:CounterEx}, we start presenting counterexamples.
but with a more complete set of examples, especially in the case of BCH expansion.
Section \ref{sec:LogA} collects some standard facts for $2\times2$ matrices for later use.
In Section \ref{sec:BCHsharp}, an extended discussion of counterexamples takes place,
especially in the regard of skew-Hermitian / unitary setting.
Up to this point, with some exceptions in Section  \ref{sec:LogA}, only material related to the question of convergence radius was presented.

In the rest, more particular but related problems are addressed.
In Section \ref{sec:ConformalRangeTwo} we include further facts regarding the conformal range.
In Section \ref{sec:MagnusHilbertTwo} we take the dual view regarding the conformal range of time-ordered exponentials.
In Section \ref{sec:MagnusGrowth} we give explicit growth estimates for the Magnus expansion,
both in terms of general growth and individual terms.

At this point, it seams reasonable to test the effectiveness of our results, mainly our range and norm estimates, against examples.
As it happens,
(a) there are several possible examples;
(b) the accessible ones deal much with $2\times2$ matrixes;
(c) in which case certain problems can be clarified relatively satisfactorily;
(d) nevertheless the computations are quite involved, and somewhat lengthy.
For this reason, the rest of the material dealing with $2\times2$ matrices has been separated into Part IIA.

Technically, most of the present paper deals with the conformal range.
While the conformal range is a simpler version of Davis--Wielandt shell, a comparison is very instructive  in many ways.
For this reason, two appendices are provided at the end.
In Appendix \ref{app:hyprev}, the hyperbolic space is reviewed.
Appendix \ref{app:DWrev} contains information about the Davis--Wielandt shell, provided for comparison to the conformal range.

\textbf{Acknowledgements.}
The author would like to thank Bal\'azs Csik\'os for some useful discussions.
\snewpage

\textbf{Notation and terminology.}
In this part we mainly work in the setting of Hilbert spaces.
$\mathfrak H$ will always be a real or complex Hilbert space.
(And it can be any of them, if the type is not specified explicitly.)
In the course of this paper, the identity element in a general Banach algebra will be denoted by $1$,
and the norm will be denoted by $|\cdot|$ (as an abbreviation of $1_{\mathfrak A}$ and $|\cdot|_{\mathfrak A}$).
However, in the case of a Hilbert space, the identity operator will be denoted by $\Id$,
and the operator norm will be denoted as $\|\cdot\|_2$
(as an abbreviation of $1_{\mathcal B(\mathfrak H)}=\Id_{\mathfrak H}$ and $|\cdot|_{\mathcal B(\mathfrak H)}$);
furthermore, the Hilbert space norm will be denoted by $|\cdot|_2$
(as an abbreviation of $|\cdot|_{\mathfrak H}$).
 But, for $2\times 2$ matrices we prefer use $\Id_2$ as the identity matrix.

Suppose that $A$ is an operator on $\mathfrak H$.
We define the co-norm $\|A\|_2^-$ of
$A$ as
\[\|A\|_2^-=\inf \{ |A\mathbf x|_2\,:\mathbf x\in\mathfrak H , |x|_2=1\}.\]
Thus $\|A\|_2^-=\|A^{-1}\|_2^{-1}$, where $A^{-1}$ is, in general, the graph inverse.
If $A$ is a finite dimensional matrix and $A$ is not invertible, then we can consider $\|A^{-1}\|_2=\infty$.

As a not entirely standard notation, we use
\begin{equation}
\mathrm G_n=(-1)^{n}\binom{-1/2}{n}=\frac{(2n)!}{2^{2n}(n!)^2}= \frac{1}{2^{2n}}\binom{2n}{n}\plabel{eq:gfunc}
\end{equation}
for the normalized central binomial coefficients ($n\geq0$).
It is well-known (an application of the Stirling formula,
or the other way around, the Wallis formula) that
\begin{equation}\mathrm G_n=\frac{1}{\sqrt{\pi n}}(1+o(1)),\plabel{eq:wallis}\end{equation}
as $n\rightarrow\infty$. It is also easy to see that $\mathrm G_n$ is strictly  decreasing.

We use notation
\[\upper=\{z\in\mathbb C\,:\,\Ima z\geq0\}\]
for the closed upper complex plane (which is a sort of inconsistent usage of overline, for typographical reasons).
Otherwise, overline is used for closure of sets, except when applied to individual complex numbers, where it denotes conjugation.

If $a$ and $b$ are points in a real affine space, then $[a,b]_e$ denotes the closed segment connecting them.
One an similarly define $[a,b)_e$, etc.

For us, $\sqrt A$ and $\log A$ are to be defined only for operators $A$ which spectrally avoid $(-\infty,0]$, i. e.
$(\spec A)\cap (-\infty,0]=\infty$.
However, for $\lambda\in\mathbb C$,  one typically defines the standard values $\sqrt \lambda^{\mathrm{st}}$
and $\log^{\mathrm{st}} \lambda$ as the values on the branch cut along $(-\infty,0)-\mathrm i\epsilon $.
As we will not really use square roots for operators, no conflict arises;
we keep the standard notation $\sqrt\lambda$ for $\sqrt \lambda^{\mathrm{st}}$.
(Henceforth, $\Rea\sqrt\lambda\geq0$.)
For the $\log$, it yields
\[\log^{\mathrm{st}} (-p)=
\begin{cases}
\pi\mathrm i+\log p& \text{ if } p>0,\\
\infty& \text{ if } p=0.
\end{cases}\]
In case of `$\log$', we use it only in operator sense, and consider it to be defined only if the spectrum is disjoint from $(-\infty,0]$;
and we keep the distinct notation $\log^{\mathrm{st}}$ for the whole $\mathbb C$.
If we need multivalued functions, then we write  $\sqrt \lambda^{\mathrm{mv}}$ and $\log^{\mathrm{mv}} \lambda$.

We will use $\backsimeq$ for similarity of matrices
(and sometimes $\simeq$ for being orthogonally conjugated), and $\sim$ for ``around'' (or, sometimes, for proportionality by scalars).
\begin{commentx}
And this should be corrected everywhere.
\end{commentx}
\snewpage
\begin{commentx}
\tableofcontents
\end{commentx}
\snewpage
\section{The spectral approach to the Magnus expansion}
\plabel{sec:SpectA}

The content of this section applies to general Banach algebras, but we have already some particular cases in mind.

Let us recall the definition of the logarithm in a Banach algebra $\mathfrak A$.
This is
\begin{align}
\log A&=\int_{\lambda=0}^1 \frac{A-1}{\lambda +(1-\lambda)A }\,d\lambda.\plabel{eq:logdef1}
\end{align}

According to the definition,
$\log A$ is well-defined, i. e. $A$ is log-able  if and only if the segment $(1-t)1+tA$ ($t\in[0,1]$) contains only invertible operators;
or, equivalently, if  $\spec(A)\cap(-\infty,0]=0$.
Thus,  in our sense, the domain of $\log$ is an open subset $\mathfrak A$.
Also, in that domain $\log$ is analytic and smooth in any possible sense; in particular,
\[\lim_{u\rightarrow0}\frac{\log(A+uB)-\log A}{u}=\int_{\lambda=0}^1
\frac{1}{\lambda +(1-\lambda)A }\,B\,\frac{1}{\lambda +(1-\lambda)A }\,d\lambda,\]
etc.
Higher derivatives can also be described, but with more occurrences of $\frac{1}{\lambda +(1-\lambda)A }$
under the integral sign.

Assume that $\sum_{k=1}^\infty s_k$ is a series in $\mathfrak A$, and we want to assign a value $S$ to it.
The series may be absolutely convergent, or simply just convergent to $S$.
The series maybe Abel summable to $S$, i. e. $\sum_{k=1}^\infty t^ks_k$
is (absolutely) convergent for any $t\in(0,1)$ such that $\lim_{t\nearrow 1} \sum_{k=1}^\infty t^ks_k=S$.
We say that the series in is radially convergent to $S$ if there is function
$s:(0,1)\rightarrow\mathfrak A$ such that
(i) $\sum_{k=1}^\infty t^ks_k$ (absolutely) converges to $s(t)$ for small $t$,
(ii) $s:(0,1)\rightarrow\mathfrak A$ is real-analytic,
(iii) $\lim_{t\nearrow 1} s(t)=S$.
These are, of course, gradually weaker notions of convergence.
If the series is not even radially convergent, we say that it is completely divergent.

The most general phenomenon related to Magnus expansion is

\begin{theorem}\plabel{th:MagnusExp}
If $\phi$ is a Banach-algebra $\mathfrak A$ valued ordered measure, and the Magnus expansion
\begin{equation}\sum_{k=1}^{\infty} \mu_{k,\mathrm R}(\phi)  \plabel{eq:MagnusExp}\end{equation}
absolutely converges to / converges to / Abel-summable to / radially convergent to the value $M$, then
\[\Rexp \phi=\exp M.\]
\begin{proof} See Part I for the absolute convergent case.
A trivial but important observation is that
$\mu_{k,\mathrm R}(t\cdot \phi)=t^k \mu_{k,\mathrm R}( \phi)$.
Now, the exponential $\Rexp(t\cdot \phi)$ is analytic in $t$, thus the statement extends even to the radially convergent case.
\end{proof}
\end{theorem}
This the reason of primary interest in the convergence of the Magnus expansion.
(In what follows we use right or left exponential and Magnus expansions somewhat eclectically, it does not really matters.)
\snewpage

With some abuse of terminology, we say that the convergence radius of the Magnus expansion \eqref{eq:MagnusExp} is the convergence radius of
\begin{equation}\sum_{k=1}^{\infty} \mu_{k,\mathrm R}(\phi) t^k \plabel{eq:MagnusExpT}\end{equation}
in terms of $t$.
If the convergence radius is greater than $1$, then the Magnus expansion converges absolutely.
However, if the convergence radius is the Magnus expansion is $1$, then
it might be absolute convergent or or not even Abel summable
(in that case Abel summability is equivalent to radial convergence).
If the radius of convergence is less than $1$, then it is not even Abel-summable, but radial convergence might happen.

The nontrivial absolute convergence statement connecting to $\log$ is the following:
\begin{theorem}\plabel{th:MPshort}
If $\phi$ is an $\mathfrak A$-valued ordered measure, and $\int |\phi|<2$, the Magnus expansion
\eqref{eq:MagnusExp} is  absolute convergent,  $\Rexp(\phi)$ is $\log$-able, and
\[\log \Rexp(\phi)=\sum_{k=1}^{\infty} \mu_{k,\mathrm R}(\phi).\]

In particular, if $\phi$ is an $\mathfrak A$-valued ordered measure, then the converges radius
of its Magnus expansion \eqref{eq:MagnusExp} is at least $2/\int |\phi|$ .
\begin{proof} See Part I, except for the last sentence.
The last sentence follows from the fact that \eqref{eq:MagnusExpT} will be absolute convergent for any
$t<2/\int |\phi|$.
\end{proof}
\end{theorem}

A very weak consequence is
\begin{theorem}\plabel{th:MagnusAbel}
If $\phi$ is an $\mathfrak A$-valued ordered measure, and the Magnus expansion
\eqref{eq:MagnusExp} has convergence radius at least $1$, and $\Rexp(t\cdot \phi)$ is $\log$-able for any $t\in[0,1]$,
then the Magnus expansion is Abel summable and
\[\log \Rexp(\phi)=\mathrm{Abel}\sum_{k=1}^{\infty} \mu_{k,\mathrm R}(\phi).\]

\begin{proof}  \eqref{eq:MagnusExpT} gives an analytical function on $t\in\intD(0,1)$, while the
$\log\Rexp(t\cdot \phi)$ extend to an analytic function in an open rectangle around $t\in[0,1]$.
According to Theorem \ref{th:MPshort} the two functions are the same for $t\sim0$.
Thus, by unicity, we have an common extension, leading to Abel summability directly.
\end{proof}
\end{theorem}
\snewpage
The statement, however, which uses the full power of analytic extension and may lead
to possibly stronger convergence results is
\begin{theorem} \plabel{th:specsh}
Suppose that  $\phi$ is an $\mathfrak A$-valued ordered measure.
Assume that for any $t\in\Dbar(0,1)$ the value $\Rexp(t\cdot\phi)$ is $\log$-able
(i.~e.~its spectrum is disjoint from $(-\infty,0]$).

Then, we claim, there is an $R>1$ ($R$ can be infinite), such that
$\log \Rexp(t\cdot\phi)$ is well-defined, and analytic for $t\in\intD(0,R)$.
On such a such a disk $\intD(0,R)$,
\[\log(\Rexp(t\cdot\phi)) =\sum_{k=1}^{\infty} \mu_{k,\mathrm R}(\phi) t^k  \]
holds. As the convergence radius of the series is larger than $1$, in particular,
the Magnus expansion converges absolutely.
\begin{proof} The elements
$\lambda+(1-\lambda)\Rexp(t\phi)$
are invertible for $(\lambda,t)\in[0,1]\times\Dbar(0,1)$, and, due to continuity, even
in a neighborhood of $[0,1]\times\Dbar(0,1)$.
This proves that $f(t)=\log \Rexp(t\phi)$ is well-defined, and analytic for $t$ in a neighborhood of
$\Dbar(0,1)$.
We know that the power series expansion of $f(t)$ is given by $f_k=\mu_{k,\mathrm R}(\phi)$ around $0$.
Then, a standard application of the generalized Cauchy formula shows that the growth of the
coefficients is limited by the analytic radius; which we know to be larger than $1$, etc.
\end{proof}
\end{theorem}

Thus Theorem \ref{th:specsh} offers a way to deal with the convergence problem using spectral arguments.
(Specifically: by testing $\log$-ability.)
Controlling spectral behaviour is difficult in general but quite doable in the case of Hilbert spaces,
or, more generally, $C^*$-algebras.
This will be the subject of the next two sections.

\begin{remark}\plabel{rem:logradius}
A trivial consequence of Theorem \ref{th:specsh} is that ``the $\log$-able radius of $\Rexp(t\cdot \phi)$'' is less or equal
as the convergence radius of (the analytic extension) of $\log \Rexp(t \cdot\phi)$.
In general, equality cannot be expected (nor it is easy to prove when it happens).
Indeed, taking the multiple concatenation $\phi\boldsymbol.^n=\phi\boldsymbol.\ldots\boldsymbol.\phi$ ($n$ times), we see that
$\Rexp(t\cdot \phi\boldsymbol.^n)=(\Rexp(t\cdot \phi))^n$.
Then the Magnus expansion is simply multiplied by $n$, thus its convergence radius remains invariant;
while (the radius of) $\log$-ability can easily be destroyed by taking powers.
Cf. Remark \ref{rem:logradius2}.
Theorem \ref{th:MMNC} may, though, grant optimal $\log$-ability in certain cases.
\qedremark
\end{remark}

\snewpage\section{The conformal range}\plabel{sec:ConformalRange}
This section introduces the conformal range of Hilbert space operators.
It is just a projection of the so-called Davis--Wielandt shell.
Its properties will be explained in greater detail in Section \ref{sec:ConformalRangeTwo}.
Here we present only the bare minimum needed to the convergence results.

For $\mathbf x,\mathbf y \in\mathfrak H\setminus\{0\}$ let $\sphericalangle(\mathbf x,\mathbf y)$ be denote their angle.
This can already be obtained from the underlying real scalar product
$\langle\mathbf x,\mathbf y\rangle_{\mathrm{real}}=\Rea\,\langle\mathbf x,\mathbf y\rangle$.
For $\mathbf x,\mathbf y\in\mathfrak H$, $\mathbf x\neq 0$,  let
\[\mathbf y:\mathbf x=\frac{
\langle \mathbf y,\mathbf x\rangle_{\real}}{|\mathbf x|_2^2}
+\mathrm i\left|\frac{\mathbf y}{|\mathbf x|_2}-  \frac{
\langle \mathbf y,\mathbf x\rangle_{\real}}{|\mathbf x|_2^2}\frac{\mathbf x}{|\mathbf x|_2}\right|_2.\]
(This is the metric information of the real orthogonal decomposition of $\mathbf y$ with respect to $\mathbf x$.
For $\mathbf y\neq 0$ it would be reasonable to define $\mathbf y:0=\infty$.)
Note that
\begin{equation}|\mathbf y:\mathbf x |= |\mathbf y|_2:|\mathbf x |_2.\plabel{eq:wquot}\end{equation}

For $A\in \mathcal B(\mathfrak H)$, we define the (extended) conformal range as
\[\CRext(A)=\{ A\mathbf x:\mathbf x, \,\overline{(A\mathbf x:\mathbf x)}
\,:\, \mathbf x\in\mathfrak H\setminus\{0\}\};\]
and the restricted conformal range as
\[\CR(A)=\{ A\mathbf x:\mathbf x,
\,:\, \mathbf x\in\mathfrak H\setminus\{0\}\}.\]

(If $A$ were but just a linear relation, derived from its graph, it would yield a
subset of the Riemann sphere $\mathbb C^\star=\mathbb C \cup \{\infty\}$.
Thus, in general, $\CR(A)=\CRext(A)\cap (\upper)^\star$ where $(\upper)^\star=\upper\cup \{\infty\}$.)

Let $A\in\mathcal B(\mathfrak H)$. From \eqref{eq:wquot}, it is immediate that
\begin{equation}\|A\|_2=\sup \{|\omega|\,:\,\omega\in\CRext(A)\},\plabel{eq:wrsup}\end{equation}
and
\begin{equation}\|A\|_2^-=\inf \{|\omega|\,:\,\omega\in\CRext(A)\},\plabel{eq:wrsupmin}\end{equation}

Assume, for now, that $\lambda\in\mathbb R$ or  $\mathfrak H$ is complex and $\lambda\in\mathbb C$.
Then,
\begin{equation}|(A-\lambda\Id)\mathbf x|_2\geq \dist(\lambda, \CRext(A)) |\mathbf x|_2.\plabel{eq:wrav}\end{equation}
\begin{lemma}\plabel{lem:spectrange}
 Suppose that $A$ is a linear operator on a complex Hilbert space.
Then, for its spectrum,
\begin{equation}\spec(A)\subset \overline{\CRext(A)}\cup \CRext(A^*).\plabel{eq:spectrange}\end{equation}
\begin{proof}
For $\lambda\in\mathbb C\setminus \overline{\CRext(A)}$, we have  $\dist(\lambda, \CRext(A))>0$,
hence, by \eqref{eq:wrav}, $A-\lambda\Id$ is invertible on its (closed) range.
This range is  $\mathfrak H$ if $\ker A^*-\bar\lambda\Id =0$, in which case $\lambda\notin \spec(A)$.
Otherise, if $\ker A^*-\bar\lambda\Id \neq0$, then $\bar\lambda\in\CRext(A^*)$, thus $\lambda\in\CRext(A^*)$.
\end{proof}
\end{lemma}
\snewpage
For us, it may also be practical to use the following temporary
\begin{lemma}\plabel{lem:spectrealtemp}
 Suppose that $A$ is a linear operator on a Hilbert space.

Then, $\overline{\CRext(A)}\cap (-\infty,0]=\emptyset$ implies $\spec(A)\cap\mathbb R=\emptyset$.
\begin{proof}
Assume, that $\overline{\CRext(A)}\cap (-\infty,0]=\emptyset$ and $\lambda=\min \spec(A)\cap\mathbb R$ exists.
Then $d=\dist(\lambda, \CRext(A))>0$.
Now $A-(\lambda-\frac13d)\Id$ invertible and $\|A-(\lambda-d/3)\Id\|_2^-\geq\frac 23d$.
Hence, $\|(A-(\lambda-\frac13d)\Id)^{-1}\|_2<(\frac 23d)^{-1}$.
Consequently, $(A-\lambda\Id)^{-1}$ is constructed by the Neumann series
$\sum_{n=0}^\infty (A-(\lambda-\frac13d)\Id)^{-1}\left((-\frac13d\Id   )(A-(\lambda-\frac13d)\Id)^{-1}\right)^n$.
This is a contradiction.
\end{proof}
\end{lemma}

The finite dimensional case is, in general, simpler:
\begin{lemma}\plabel{lem:CRFinDim}
If $\dim\mathfrak H<\infty$, then

(a)
\[\text{$\CRext(A)$ is compact};\]

(b)
\[\spec(A)\subset\CRext(A);\]

(c)
\[\spec(A)\cap\mathbb R=\CRext(A)\cap\mathbb R.\]
\begin{proof}
(a) $\CR(A)$ is a continuous image of the compact unit sphere, thus compact.
Hence, $\CRext(A)$ is too.
(b) and (c) hold as the spectrum is the point spectrum
(i. e. every spectrum point comes from an eigenvector).
\end{proof}
\end{lemma}

\snewpage\section{Conformal range and time-ordered exponentials}\plabel{sec:MagnusHilbert}
Suppose that $\mathfrak H$ is a finite dimensional real  Hilbert space.
Let us consider the map
\[\log_{\mathfrak H}:\mathfrak H\setminus\{0\}\rightarrow \mathbb R\times\mathbb S(\mathfrak H)\]
\[\mathbf v\mapsto\left(\log |\mathbf v|_2, \frac{\mathbf v}{|\mathbf v|_2}\right).\]
$\mathbb R$ and $\mathbb S(\mathfrak H)$ both possess standard Riemannian manifold structures with respect to the ordinary arc length.
So there is a standard Riemannian product structure on $\mathbb R\times\mathbb S(\mathfrak H)$, which is, in fact, complete.
This induces a Riemannian structure on $\mathfrak H\setminus\{0\}$ via $\log_{\mathfrak H}$,
which we simply call the logarithmic Riemannian structure on $\mathfrak H\setminus\{0\}$.

Distance  minimizing geodesic segments are easy to describe there:
If $\mathbf v\in\mathfrak H\setminus\{0\}$, and $\mathbf w\in\mathfrak H\setminus\{0\}$ is sufficiently close to $\mathbf v$,
then there is such a unique geodesic segment connecting $\mathbf v$ and  $\mathbf w$.
Due to symmetry, this segment must lay in $(\mathbb R\mathbf v +\mathbb R\mathbf w)\setminus \{0\}$.
In this case, it can be assumed that $\mathbb R\mathbf v +\mathbb R\mathbf w=\mathbb C$,
$\mathbf v=1$, and $\mathbf w\in \mathbb C\setminus (-\infty,0]$.
Here, the Riemannian structure and the geodesic segment is induced via the ordinary $\log$.
As the continuation of geodesics is unique, it is easy to see that this picture is valid as long as
$(\mathbf v,\mathbf w)\sphericalangle<\pi$, i. e. $\mathbf v$ and $\mathbf w$ not antipodal.
If $\mathbf v$ and $\mathbf w$ are antipodal, then any geodesic segment connecting them
should contain a point outside of $\mathbb R\mathbf v=\mathbb R\mathbf w$.
Thus, it is easy to see that the geodesic segment is in $(\mathbb R\mathbf v +\mathbb R\mathbf z)\setminus \{0\}$.
After making the appropriate identifications, it is induced by $\log$ again.
The distance minimizing geodesic segment is unique if a half-plane bounding $\mathbb R\mathbf v$ is chosen.
In any case, the geodesic distance is given by
\begin{equation}
d_{\log}(\mathbf v,\mathbf w)=\sqrt{\left(\log \frac{|\mathbf w|_2}{|\mathbf v|_2} \right)^2+\left((\mathbf v,\mathbf w)\sphericalangle\right)^2}.
\plabel{eq:logdist}\end{equation}

Let us now consider  any real or complex Hilbert space $\mathfrak H$.
In the complex case we can simply consider the underlying real Hilbert space structure.
Then $\mathfrak H\setminus\{0\}$ can be endowed by the metric \eqref{eq:logdist}.
We can keep the terminology `distance segment' for distance minimizing geodesic segments
coming from any finite (in fact: two) dimensional real Hilbert subspace of $\mathfrak H$.
\begin{lemma}\plabel{lem:logdist}
Let $\mathfrak H$ be a real or complex Hilbert space. Then

(a) $\mathfrak H\setminus \{0\}$ is a metric space with respect to \eqref{eq:logdist}.

(b) If $\mathbf v$ and $\mathbf w$ are not antipodal (i. e. $(\mathbf v,\mathbf w)\sphericalangle<\pi$), then there is a unique distance segment connecting $\mathbf v$ to $\mathbf w$,
and any distance-minimizing finite of infinite path connecting $\mathbf v$ to $\mathbf w$  is a (not necessary strictly) increasing
subpath of this distance segment.
(In particular, if $d_{\log}(\mathbf v,\mathbf z)+d_{\log}(\mathbf z,\mathbf w)=d_{\log}(\mathbf v,\mathbf w)$, then
$\mathbf z$ is a point of the segment).

(c)  If $\mathbf v$ and $\mathbf w$ are  antipodal (i. e. $(\mathbf v,\mathbf w)\sphericalangle=\pi$)
and
$\mathbf z$ is a unit vector perpendicular to $\mathbf v$ (and $\mathbf w$), then there is a unique distance segment connecting $\mathbf v$ to $\mathbf w$ contained in $\mathbb R\mathbf v+[0,\infty)\cdot\mathbf z$.
This exhibits all distance segments.
These distance segments are pairwise disjoint except at end points.
Any distance-minimizing finite of infinite path connecting $\mathbf v$ to $\mathbf w$  is a (not necessary strictly) increasing
subpath of such a distance segment.
\begin{proof}
It is sufficient to consider finite paths, which are, in turn, come from finite dimensional Hilbert subspaces,
where the finite dimensional picture is sufficient.
(Of course, one comfortable with Hilbert manifolds may approach the situation more directly.)
\end{proof}
\end{lemma}
For the rest of this section, $\mathfrak H$ will be a real or complex Hilbert space.
\snewpage
\begin{lemma}\plabel{lem:logvarest}
Assume that $\mathbf v,\mathbf w\in \mathfrak H\setminus\{0\}$, $|\mathbf w|_2<|\mathbf v|_2$.
Then
\[d_{\log}(\mathbf v,\mathbf v+\mathbf w)\leq\log \frac{|\mathbf v|_2}{|\mathbf v|_2-|\mathbf w|_2}.\]
\begin{proof}
We can assume that the Hilbert space is $\mathbb C$, $\mathbf v=1$.
For $t\in[0,1]$ the absolute value of the derivative of $\log (1+\mathbf wt)$ is less or equal the
$-1$ times the derivative of  $\log (1-|\mathbf w|t)$.
Integrated, this yields the statement.
\end{proof}
\end{lemma}
\begin{theorem}\plabel{th:cs}
Suppose that $\mathbf z:[a,b]\rightarrow\mathfrak H\setminus\{0\}$ is continuous. Then
\[d_{\log}(\mathbf z(a),\mathbf z(b) )\leq\int_{t\in[a,b]}\frac{|\mathrm d\mathbf z(t)|_2}{|\mathbf z(t)|_2}.\]
In case of equality, $\mathbf z$ is a (not necessarily strictly) monotone subpath of a
distance segment connecting $\mathbf z(a)$ to $\mathbf z(b)$
\begin{proof}
The statement is non-vacuous only if $\mathbf z$ is of finite variation.
For small $\mathbf w$ ($\mathbf v$ bounded in $(0,+ \infty)$),
\[\log \frac{|\mathbf v|_2}{|\mathbf v|_2-|\mathbf w|_2}=\frac{|\mathbf w|_2}{|\mathbf v|_2}+O\left(\left(\frac{|\mathbf w|_2}{|\mathbf v|_2}\right)^2\right),\]
thus appropriate refinements with $\mathbf v\sim\mathbf z(t_i)$, $\mathbf w\sim\mathbf z(t_{i+1})-\mathbf z(t_i)$
yield the estimate.

Equality implies that \[d_{\log}(\mathbf z(a),\mathbf z(b) )=d_{\log}(\mathbf z(a),\mathbf z(t_1) )+d_{\log}(\mathbf z(t_1),\mathbf z(t_2) )+d_{\log}(\mathbf z(t_2),\mathbf z(b) )\]
for and $t_1\leq t_2\in[a,b]$.
Compared to Lemma \ref{lem:logdist}, this implies the equality statement.
\end{proof}
\begin{proof}[Geometric proof] (Bal\'azs Csik\'os, \cite{Cs}.)
The following argument captures the geometric essence of the statement:
The statement is non-vacuous only if the logarithmic variation  $\int_{t\in[a,b]}\frac{|\mathrm d\mathbf z(t)|_2}{|\mathbf z(t)|_2}$ is finite.
This, however, implies that the (less or equal) angular variation $\int_{t\in[a,b]}\left|\mathrm d\frac{\mathbf z(t)}{|\mathbf z(t)|_2}\right|_2$ is finite.
This allows to define a continuous map $\tilde {\mathbf z}:[a,b]\rightarrow\widetilde{\mathbb C}$ by
\[\tilde {\mathbf z}(t)=\left(|\mathbf z(t)|_2,  \int_{s\in[a,t]}\left|\mathrm d\frac{\mathbf z(s)}{|\mathbf z(s)|_2}\right|_2 \right)_{\mathrm{polar}},\]
where $\widetilde{\mathbb C}$ is the universal covering space of $\mathbb C\setminus\{0\}$.
The intuitive idea is that one can consider the cone over $\mathbf z$, which is a developable surface,
which we unfold to $\widetilde{\mathbb C}$. The curves ${\mathbf z}$ and $\tilde {\mathbf z}$
look quite different but their (log)variations are the same because the
their (log)radial and angular variations are the same, and the (log)variations can be
assembled from them in the same manner. Then
 \[
d_{\log}(\mathbf z(a),\mathbf z(b) )
\leq |\log \tilde {\mathbf z}(b)  -\log\tilde {\mathbf z}(b) =
 \left| \int_{t\in[a,b]}\frac{\mathrm d\tilde{\mathbf z}(t)}{\tilde{\mathbf z}(t)}\right|
 \leq \int_{t\in[a,b]}\frac{|\mathrm d\tilde{\mathbf z}(t)|}{|\tilde{\mathbf z}(t)|}=
 \int_{t\in[a,b]}\frac{|\mathrm d\mathbf z(t)|_2}{|\mathbf z(t)|_2}\notag
\]
shows the inequality part.
If we map the developed surface into the Hilbert space, then any contraction
in the angular variation will yields a contraction with respect to the points of
starting and ending ray relative to each other.
Thus we deal along great circles of the unit sphere, hence the equality statement can be reduced to the planar case, which is simple.
\end{proof}
\end{theorem}

\begin{theorem}\plabel{th:CRrange} (Time ordered exponential mapping theorem.)

If $\phi$ is $\mathcal B(\mathfrak H)$-valued ordered measure, then
\begin{equation} \CRext(\Lexp \phi)\subset \exp \Dbar(0,\textstyle{\int \|\phi\|_2}),\plabel{eq:CRrange}\end{equation}
and
\begin{equation}  \spec(\Lexp \phi)\subset \exp \Dbar(0,\textstyle{\int \|\phi\|_2}).\plabel{eq:spect}\end{equation}

In particular, if $\int \|\phi\|_2<\pi$, then $\log \Lexp \phi$ is well-defined, and for its spectral radius
\begin{equation} \mathrm r(\log\Lexp \phi) \leq\textstyle{\int \|\phi\|_2}.\plabel{eq:spect2}\end{equation}
\begin{proof}

Let $\mathbf x\in\mathfrak H$, $|\mathbf x|_2=1$.
Let us define $\mathbf z:[a,b]\rightarrow\mathfrak H$ by
\begin{equation}
\mathbf z(t)=\Lexp(\phi|_{[a,t]}) \mathbf x.
\plabel{eq:sample}
\end{equation}
Apply Theorem \ref{th:cs}.
Due to $\mathbf z(a)= \mathbf x $, $\mathbf z(b)= \Lexp(\phi)\mathbf x$, and the estimate
\[\int_{t\in[a,b]}\frac{|\mathrm d\mathbf z(t)|_2}{|\mathbf z(t)|_2}\leq \int \|\phi\|_2,\]
we obtain \eqref{eq:CRrange} immediately.
If we replace $\phi$ by $(\phi^*)^\dag$, i.e. adjoined and order-reversed, then it yields
$\CRext((\Lexp \phi)^*)\subset \exp \Dbar(0,\textstyle{\int \|\phi\|_2})$.
For the rest of the argument we can assume that
in the case of a real Hilbert space, we have already considered the complexified setting from the beginning.
(It leaves the norms invariant.)
Then Lemma \ref{lem:spectrange} implies \eqref{eq:spect}.
If $\int \|\phi\|_2<\pi$, then $\spec(\Lexp \phi)$ is disjoint from $(-\infty,0]$, and \eqref{eq:spect2} follows from the spectral mapping theorem.
\end{proof}
\end{theorem}
An immediate consequence is
\begin{theorem}\plabel{th:MMNC} (Moan, Niesen \cite{MN}, Casas \cite{Ca}.
For the BCH case: Mityagin \cite{Mity}.)

If $\phi$ is a $\mathcal B(\mathfrak H)$-valued ordered measure and  $\int \|\phi\|_2<\pi$,
then the Magnus expansion $\sum_{k=1}^\infty \mu_{k,\mathrm L}(\phi)$ is absolute convergent.
In fact, $\log\Lexp(\phi)=\sum_{k=1}^\infty \mu_{k,\mathrm L}(\phi)$ also holds.

The statement also holds if $\mathcal B(\mathfrak H)$ is replaced by any $C^*$-algebra.
\begin{proof}
By Theorem \ref{th:CRrange}, the spectral conditions of  Theorem \ref{th:specsh} are satisfied, thus it applies.
Regarding $C^*$-algebras, the extension immediately follows about the Gelfand--Naimark isometric representation theorem.
\end{proof}
\end{theorem}
\begin{remark}\plabel{rem:MMNC}
The authors quoted above use the angular variation, thus they prove (in our terminology)
\[\CRext(\Lexp \phi)\subset \exp\{z\,:\,|\Ima z|\leq \textstyle{\int \|\phi\|_2}\}.\]
This, used in conjunction with
\[\CRext(\Lexp \phi)\subset \exp\{z\,:\,|\Rea z|\leq \textstyle{\int \|\phi\|_2}\}\]
(which is rather trivial), yields
\begin{equation} \CRext(\Lexp \phi)\subset \exp\{z\,:\,|\Rea z|,|\Ima z| \leq \textstyle{\int \|\phi\|_2}\};
\plabel{eq:mapcrude}\end{equation}
etc.
This is sufficient to establish Theorem \ref{th:MMNC}.
Otherwise, their ideas are basically the same, see, in particular, Casas  \cite{Ca}.

Thus, there are two crucial steps:
One is the observation that the convergence radius of the Magnus expansion of $\phi$ is
at least  the $\log$-able radius  around $0$ (the easy step)
and other is that actual $\log$-ability works out in the given setting (the hard step).
It is true, however, that the harder part has already been solved much earlier in great generality.
Not only that Kre\u{\i}n \cite{K} (1969) is most explicit about angular localization,
but the relevant $\log$-ability / existence of $\log$ problem was already solved sharply by Sch\"affer \cite{Sch} (1964).
\qedremark
\end{remark}

The following lemma  is a natural addition to Theorem  \ref{th:CRrange}:

\begin{lemma}\plabel{lem:CRrange}
Consider the situation of Theorem  \ref{th:CRrange}.
Assume that $\mathbf x\in\mathfrak H\setminus\{0\}$ such that
\[ \lambda= ((\Lexp \phi)\mathbf x:\mathbf x)\in \partial \exp \Dbar(0,\textstyle{\int \|\phi\|_2}).\]
Then the following is true:

(i) If $\lambda\in(0,\infty)$, then
\[(\Lexp \phi|_{[0,t)})\mathbf x= \mathrm e^{ \|\phi|_{[0,t)}\|_2\cdot \sgn\log\lambda}\mathbf x.\]

(ii) Otherwise,  there is a unique $\alpha\in(0,\pi)$ and $\mathbf y$ real-orthogonal to $\mathbf x$, $|\mathbf y|_2=|\mathbf x|_2$, such that
\[\cos\alpha+\mathrm i \sin\alpha=
\left(\logst \lambda\right)^0
=\frac{\logst \lambda}{\int\|\phi\|_2}\]
and
\[(\Lexp \phi|_{[0,t)})\mathbf x= \mathrm e^{ \|\phi|_{[0,t)}\|_2\cdot \cos\alpha}\left(
\cos\left(  \|\phi|_{[0,t)}\|_2\cdot  \sin\alpha\right)\cdot\mathbf x+\sin\left(  \|\phi|_{[0,t)}\|_2\cdot  \sin\alpha\right)\cdot\mathbf y
 \right).\]
\begin{proof}
Apply the equality part of Theorem \ref{th:cs} to \eqref{eq:sample}.
\end{proof}
\end{lemma}

We say that the ordered measure $\phi$ is a multiple Baker--Campbell--Hausdorff (mBCH) measure,
if, up to reparametrization,  $\phi$ is of form $A_1\mathbf 1\boldsymbol.\ldots\boldsymbol.A_n\mathbf 1$.
In this case, $\phi$ also allows a mass-normalized version
\begin{equation}
\psi=B_1\mathbf 1_{(0,t_1]}\boldsymbol.\ldots\boldsymbol.B_k\mathbf 1_{[t_{k-1},t_k]}
\plabel{eq:mmBCH}
\end{equation}
where $t_i<t_{i+1}$, $\|B_i\|_2=1$ and thus $t_k=\int\|\phi\|_2$.
It is constructed by replacing $A_i$ with $A_i/\|A_i\|_2$ if $\|A_i\|_2\neq0$, and eliminating the term if $A_i=0$.
As it is obtained by a kind a reparametrization, its Magnus expansion is not affected.

\begin{prop}\plabel{pr:CRrangeF}
Let $\psi$ be a mass-normalized mBCH-measure as in \eqref{eq:mmBCH}, $\int\|\psi\|_2>0$.
Assume that $\mathbf x\in\mathfrak H\setminus\{0\}$ such that
\[ \lambda= ((\Lexp \psi)\mathbf x:\mathbf x)\in \partial \exp \Dbar(0,\textstyle{\int \|\psi\|_2}).\]
Then the following is true:

(i) If $\lambda\in(0,\infty)$, then the real orthogonal decomposition $\mathbb R\mathbf x \oplus (\mathbb R\mathbf x)^\bot$ is invariant for all $B_i$ and
\[B_i\mathbf x=(\sgn\log\lambda)\mathbf x\]
(action independent from $i$).

(ii) Otherwise,  there is a unique $\alpha\in(0,\pi)$ and $\mathbf y$ real-orthogonal to $\mathbf x$, $|\mathbf y|_2=|\mathbf x|_2$,  such that
\[\cos\alpha+\mathrm i \sin\alpha=\left(\logst \lambda\right)^0=
\frac{\logst \lambda}{\int\|\psi\|_2},\]
the orthogonal decomposition $(\mathbb R\mathbf x +\mathbb R\mathbf y)\oplus (\mathbb R\mathbf x+\mathbb R\mathbf y)^\bot$ is invariant for all $B_i$, and
\[B_i\mathbf x=(\cos\alpha)\mathbf x+(\sin\alpha)\mathbf y,\]
\[B_i\mathbf y=-(\sin\alpha)\mathbf x+(\cos\alpha)\mathbf y\]
(action independent from $i$).

(ii)' If $\mathfrak H$ is a complex Hilbert space and case (ii) holds, then there
is a common eigenvector $\tilde {\mathbf x}$ with eigenvalue $\cos\alpha+\mathrm i\sin\alpha$,
thus case (ii) applies with $\tilde {\mathbf y}=\mathrm i\tilde {\mathbf x}$ or $\tilde {\mathbf y}=-\mathrm i\tilde {\mathbf x}$.
\begin{proof}
Apply Lemma \ref{lem:CRrange}.
In case (i) the restricted action of the $B_i $ is straightforward. Consider case (ii) now.
On any interval $(t_i,t_{i+1})$ we know the (infinitesimal) action (of the exponential) of $B_i$ at linearly independent places, which determines it.
As restricted actions are orthogonal, and the $B_i$ are orthogonal, we know that the decomposition of the Hilbert space is orthogonal.
Regarding (ii)', $\tilde{\mathbf x}=\mathbf x+\mathrm i\mathbf y $ or $\mathbf x-\mathrm i\mathbf y $ will do.
\end{proof}
\end{prop}
\begin{theorem}\plabel{th:CRrangeF}
Let $\psi$ be a mass-normalized mBCH-measure as in \eqref{eq:mmBCH}, $\int\|\psi\|_2>0$.
Assume that the underlying Hilbert space $\mathfrak H$ is finite-dimensional.
Then there is an orthogonal decomposition $\mathfrak H=\mathfrak H_1\oplus\mathfrak H_2$ such that

(i) the operators $B_i$ leave the decomposition invariant;

(ii) $\psi|_{\mathfrak H_1}$ is a constant measure $U\mathbf 1_{[0,t_k]}$ which
is orthogonal (unitary) in case of real (complex) Hilbert spaces, and its eigenvalues are from
\[\left\{w\in\partial\Dbar(0,1)\,:\,  \ima w\cdot \int\|\psi\|_2\leq\pi\right\};\]

(iii) for $\psi|_{\mathfrak H_2}$,
\[ \spec(\Lexp \psi|_{\mathfrak H_2})\subset\CRext(\Lexp \psi|_{\mathfrak H_2})\subset\mathrm{int}\,\exp \Dbar(0,\textstyle{\int \|\psi\|_2}).\]
\begin{proof}
We iterate Proposition \ref{pr:CRrangeF} in order to chip away constant orthogonal (unitary) parts.
The spectrum is contained in the conformal range by finite dimensionality.
\end{proof}
\end{theorem}
\begin{lemma}\plabel{lem:common}
Let $\psi$ be a mass-normalized mBCH-measure as in  \eqref{eq:mmBCH}, $\int\|\psi\|_2>0$.
Consider all the Hilbert subspaces $\mathfrak V$ of $\mathfrak H$ such that

(i) $\mathfrak H=\mathfrak V\oplus\mathfrak V^\bot$ is an invariant orthogonal decomposition for all $B_i$.

(ii) $B_1|_{\mathfrak V}=\ldots=B_k|_{\mathfrak V}$, and these are orthogonal (unitary).

Then there is a single maximal such $\mathfrak V$.

\begin{proof}
Take the closure of the unions, and the intersection of the orthogonal complements, respectively.
\end{proof}
\end{lemma}
With respect of such a maximal $\mathfrak V$ as in the previous lemma, we may phrase
$B_i|_{\mathfrak V}$ as the common part, and  $B_i|_{\mathfrak V^\bot}$ as the reduced part of the operators $B_i$.
If $\mathfrak V=0$, then we call $\psi$ reduced.
In particular, this applies if $\bigcap_{i<j}\ker(B_i-B_j)=0$.

\begin{theorem}\plabel{th:better2}
(Finite critical BCH convergence theorem.)

Let $\mathfrak H$ be a finite dimensional Hilbert space .
Consider the $\mathcal B(\mathfrak H)$ valued mBCH measure $\phi=A_1\mathbf 1\boldsymbol.\ldots\boldsymbol.A_k\mathbf 1$  with cumulative norm
$\|A_1\|_2+\ldots+\|A_k\|_2=\pi$.

Then, the convergence radius of the Magnus (mBCH) expansion of $\phi$ is greater than $1$.
In particular, finite dimensional mBCH expansions with cumulative norm $\pi$ converge.

\begin{proof}
First, we pass from the original $\phi$ to a mass-normalized version $\psi$ as in \eqref{eq:mmBCH}.

 Using Lemma \ref{lem:common}, we can decompose to common and reduced parts, $\mathfrak H=\mathfrak V\oplus\mathfrak V^\bot$.
On the common part the Magnus expansion is trivial (the higher Magnus brackets vanish).
Thus can assume $\mathfrak V=0$, $\mathfrak V^\bot=\mathfrak H$.
Taking any $t\in\Dbar(0,1)\setminus\{0\}$, the system $t\cdot B_1,\ldots,t\cdot B_k$ is still reduced.
Thus by Theorem \ref{th:CRrangeF} ($\mathfrak H_1$ must be $0$), $\spec(\Lexp(t\cdot\psi))\subset \intt\exp\Dbar(0,\pi|t|)\subset\intt\exp\Dbar(0,\pi)$.
This implies $\spec(\Lexp(t\cdot\psi))\cap(-\infty,0]=\emptyset$.
(This, of course, also holds for $t=0$.)
Thus Theorem \ref{th:specsh} can be applied to prove convergence.
\end{proof}
\end{theorem}

\begin{theorem}\plabel{th:better3}
(Finite logarithmic critical BCH convergence theorem.)

Let $\mathfrak H$ be a finite dimensional Hilbert space.
Consider the $\mathcal B(\mathfrak H)$ valued mass-normalized mBCH measure $\psi$ as in \eqref{eq:mmBCH}
with cumulative norm $\int\|\psi\|_2=\pi$.

(a) Unless the component operators $B_i$ have a common eigenvector for $\mathrm i$ or $-\mathrm i$ (complex case), or a common eigenblock
$\bem&-1\\1&\eem$ (real case), then  $\log\Lexp(\psi)=\sum_{k=1}^\infty \mu_{k,\mathrm L}(\psi)$ also holds.

(b) If $\psi$ is reduced, then
for any $t\in\Dbar(0,1)$, $\log\Lexp(t\cdot \psi)=\sum_{k=1}^\infty t^k\mu_{k,\mathrm L}(\psi)$  holds.
Thus, the $log$-able radius of the Magnus (BCH) expansion is also greater than $1$.
\begin{proof}
In continuation to the proof of the previous theorem,
 the logarithm formula holds on the reduced part, but on the common part we must avoid the indicated eigenvalues.
\end{proof}
\end{theorem}

Theorem \ref{th:better3} allows an infinite-dimensional version which is somewhat weaker.
The key observation is that Lemma \ref{lem:CRrange} allows a stable version:
\begin{lemma}\plabel{lem:CRrangeG}
For all $\delta>0$ and $M>0$ there is $\varepsilon>0$ such that the following holds:

Consider the situation of Theorem  \ref{th:CRrange} with $\int \|\phi\|_2\leq M$.
Assume that
\[ \lambda \in\partial \exp \Dbar(0,\textstyle{\int \|\phi\|_2})\cap \upper .\]

Assume that for the unit vector $\mathbf x\in\mathfrak H$,
\[\left| ((\Lexp \phi)\mathbf x:\mathbf x)-\lambda\right|_2<\delta.\]
Then the following is true:

(i) If $\lambda\in(0,\infty)$, then for all possible $t$,
\[\left|(\Lexp \phi|_{[0,t)})\mathbf x-\mathrm e^{ \|\phi|_{[0,t)}\|_2\cdot \sgn\log\lambda}\mathbf x\right|_2<\varepsilon.\]

(ii) Otherwise,  with $\alpha\in(0,\pi)$,
\[\cos\alpha+\mathrm i \sin\alpha=\left(\logst \lambda\right)^0=\frac{\logst \lambda}{\int\|\phi\|_2},\]
there is another unit vector $\mathbf y$ such that
\[\mathbf x_m\bot_{\mathrm{real}}\mathbf y_m;\]
and for all possible $t$,
\[\left|(\Lexp \phi|_{[0,t)})\mathbf x- \mathrm e^{ \|\phi|_{[0,t)}\|_2\cdot \cos\alpha}\left(
\cos\left(  \|\phi|_{[0,t)}\|_2\cdot  \sin\alpha\right)\cdot\mathbf x+\sin\left(  \|\phi|_{[0,t)}\|_2\cdot  \sin\alpha\right)\cdot\mathbf y
 \right)\right|_2<\varepsilon.\]
\begin{proof}
This is probably the best to explained in non-technical terms.
The idea is that paths which ``quite optimal'' are still closed to geodesics,
otherwise the ``triangle equality'' would compromised too much.
This affects not only the trajectory but the pace of the path, as it is also controlled by the limited variation.
Firstly, rigidity is best to be established in terms of $((\Lexp \phi|_{[0,t)})\mathbf x:\mathbf x)$.
If this is done, one candidate for $\mathbf y$ based on a middle value, if necessary, can be found.
As the situation is sufficiently compact, independence from $\lambda$ and $\int \|\phi\|_2$ can be achieved.
(Although it would not matter for us later.)
\end{proof}
\end{lemma}
\snewpage
\begin{prop}\plabel{pr:CRrangeG}
Let $\psi$ be a mass-normalized mBCH-measure as in \eqref{eq:mmBCH}, $\int\|\psi\|_2>0$.
Assume that
\[ \lambda\in   \overline{\CRext(\Lexp(\phi)) }\cap \partial \exp \Dbar(0,\textstyle{\int \|\psi\|_2})\cap\upper.\]
with $\Ima\lambda\geq0$.

(i) If $\lambda\in(0,\infty)$, then there is a sequence of unit vectors $\mathbf x_m$ such that for all $B_i$,
\[B_i\mathbf x_m - (\sgn\log\lambda)\mathbf x_m\rightarrow 0.\]

(ii) Otherwise,  with $\alpha\in(0,\pi)$
\[\cos\alpha+\mathrm i \sin\alpha=\left(\logst \lambda\right)^0=\frac{\logst \lambda}{\int\|\phi\|_2},\]
there are sequences unit vectors $\mathbf x_m$ and  $\mathbf y_m$ such that
\[\mathbf x_m\bot_{\mathrm{real}}\mathbf y_m\]
holds for all $m$;
and
\[B_i\mathbf x_m-\Bigl((\cos\alpha)\mathbf x_m+(\sin\alpha)\mathbf y_m\Bigr)\rightarrow0,\]
and
\[B_i\mathbf y_m-\Bigl(-(\sin\alpha)\mathbf x_m+(\cos\alpha)\mathbf y_m\Bigr)\rightarrow0.\]
holds for all $B_i$.

(ii)' If $\mathfrak H$ is a complex Hilbert space and case (ii) holds, then there
there is a sequence  of unit vectors $\tilde{\mathbf x}_m$ such that
\[B_i\tilde{\mathbf x}_m - (\cos\alpha+\mathrm i \sin\alpha )\tilde{\mathbf x}_m\rightarrow 0,\]
or
\[B_i\tilde{\mathbf x}_m - (\cos\alpha-\mathrm i \sin\alpha )\tilde{\mathbf x}_m\rightarrow 0.\]
holds for all $B_i$ (independently from $i$).
\begin{proof}
Let us approximate $\lambda$ by $((\Lexp \phi)\mathbf x_m:\mathbf x_m)$
For any $[u,v)\subset[t_{i-1},t_i)$,
\[ (\Lexp \phi|_{[0,v)})\mathbf x_m- (\Lexp \phi|_{[0,u)})\mathbf x_m =B_i\int_{t=u}^v (\Lexp \phi|_{[0,t)})\mathbf x_m \mathrm dt\]
holds.
Putting the approximating geodesics there, we obtain an appropriately quantified approximating \textit{linear} equation for $\mathbf x_m$ and $\mathbf y_m$
and   $B_i\mathbf x_m$ and $B_i\mathbf y_m$.
Taking another interval, we obtain another equation; thus ultimately we obtain an approximating linear expression for  $B_i\mathbf x_m$ and $B_i\mathbf y_m$
by $\mathbf x_m$ and $\mathbf y_m$.
As the approximations get better, the approximating terms vanish; leading, predictably, to the statement.
The complex case can be recovered by the standard trick.
\end{proof}
\end{prop}
\snewpage
The content of Proposition \ref{pr:CRrangeG} can phrased so that all indicated $\lambda$
come from common approximate eigenvectors or common approximate eigenblocks.
As a consequence,
\begin{theorem}\plabel{th:better1}
Consider the $\mathcal B(\mathfrak H)$ valued mass-normalized mBCH measure $\psi$ as in \eqref{eq:mmBCH}
with cumulative norm $\int\|\psi\|_2=\pi$.

(a) Unless the component operators $B_i$ have approximate common eigenvectors for $\mathrm i$ or $-\mathrm i$ (complex case), or a common approximate eigenblock
$\bem&-1\\1&\eem$ (real case), then the Magnus expansion is Abel summable to the logarithm of the exponential,
\[\log\Lexp(\psi)=\mathrm{Abel}\sum_{k=1}^\infty \mu_{k,\mathrm L}(\psi).\]

(b) Unless the component operators $B_i$ have approximate common eigenvectors for an unit complex number (complex case),
or a common approximate eigenvalue $1$, $-1$, or commmon approximate eigenblock
$\bem\cos\alpha&-\sin\alpha\\\sin\alpha&\cos\alpha\eem$ (real case)
then, for any $t\in\Dbar(0,1)$ the Magnus expansion is absolutely convergent to the logarithm of the exponential,
\[\log\Lexp(t\cdot \psi)=\sum_{k=1}^\infty \mu_{k,\mathrm L}(t\cdot \psi).\]
Thus, the $\log$-able radius of the Magnus (BCH) expansion is also greater than $1$.

In particular, this latter situation holds if the operators
$B_i-B_j$ $(i<j)$ has no common approximate $0$.
\begin{proof}
(a) If the assumption holds, then $-1\notin \overline{\CRext(\Lexp(\psi))} $, while $\overline{\CRext(\Lexp(\psi))}\subset\exp\Dbar(0,\pi)$ holds.
Hence, by Lemma \ref{lem:spectrealtemp},  $\spec(\Lexp(\psi))\cap(-\infty,0]=\emptyset$.
Thus, Theorem \ref{th:MagnusAbel} can be applied.

(b) In the complex case, the previous argument
applied to complex unit multiples of $\psi$ shows that $\log$-able radius is greater than $1$.
Thus, Theorem \ref{th:specsh} can be applied.
In the real case, one has show that complexification does not destroy the key assumption, but this is not difficult.
\end{proof}
\end{theorem}

\snewpage
\section{Counterexamples}
\plabel{sec:CounterEx}

The counterexamples presented in this section are historical but with improvements.
\\~

\subsection{Some restrictions on exponentials}
\plabel{ss:sre}
~\\

Recall that `$A$ is $\log$-able' simply means that the spectrum of $A$ is disjoint from $(-\infty,0]$.
If $A$ is $\log$-able, then it is the exponential of its logarithm.
For the sake of completeness, we include the following well-known

\begin{lemma}\plabel{lem:rc}
Concerning $2\times2$ matrices:

(a) In the real case, if $A$ is not $\log$-able and $A$ is not a negative scalar matrix,
then it does not occur as an exponential of another real matrix.

(b) In the complex case, however, if $A$ is not $\log$-able but invertible, then it does  occur as an exponential
of another complex matrix.
\begin{proof}[Note]In the real case:
If $\det A>0$ and $A$ is not $\log$-able and $A$ is not a negative scalar matrix, then
$A$ is the strictly parabolic case with negative trace,
i. e. \[\text{$A\backsimeq\bem\lambda&1\\&\lambda\eem$ with $\lambda<0$},\]
or $A$ is  strictly
hyperbolic case with two negative eigenvalues, i. e.
\[\text{$A\backsimeq\bem\lambda_1&\\&\lambda_2\eem$ with $\lambda_1<\lambda_2<0$}.\]
If $\det A\leq0$, then $A$ is not $\log$-able.
\qedno
\end{proof}

\begin{proof}
(a) Assume $A=\exp B$ where $B$ is a real matrix, and $A$ is not $\log$-able.
Due to the second assumption, $A$ must have an eigenvalue in $(-\infty,0]$.
Now, $\det A=\exp\tr A>0$  implies that $A$ should have two negative eigenvalues.
If $\beta$ is an eigenvalues of $B$, then $\exp \beta$ is an eigenvalue of $A$.
This implies that $\beta=r+(2k+1)\pi\mathrm i$ where $r\in\mathbb R$ and $k\in\mathbb Z$.
Due to $B$ being real, the other eigenvalue of $B$ should be the complex conjugate
$\bar\beta=r-(2k+1)\pi\mathrm i\neq \beta$.
Thus $B$ is diagonalizable with (distrinct) eigenvalues  $r\pm(2k+1)\pi\mathrm i$.
That makes $A=\exp B$ diagonalizable with equal eigenvalues $-\exp r$.
Therefore, $A$ is a negative scalar matrix.

(b) By inspection, one can check that every invertible complex Jordan block occurs as an exponential.
\end{proof}
\end{lemma}
\begin{remark}\plabel{rem:rc}
In general, regarding $n\times n$ matrices:

(a) In the real case, $A$ is an exponential if and only if it is invertible and its complex Jordan blocks with negative diagonals are ``doubled''.

(b)  In the complex case, $A$ is an exponential if and only if it is invertible.
\qedremark
\end{remark}

Wei \cite{W} (1963) already uses the observations above, systematically, to present counterexamples
 of the Magnus (BCH) expansion (but without optimization to the norm).

\snewpage

Another useful observation is the following:
Assume that $\phi$ is ordered measure of finite variation with values in $n\times n$ real or complex matrices.
Then there are the usual trace and determinant operations in those matrix algebras.
Now, the observation is that all higher term in the Magnus expansion of $\phi$ are of vanishing trace, i. e.
$\tr \mu_{k,\mathrm R}(\phi)=0$ for $k\geq2$.

Indeed, this follows from the fact those higher terms are of integrals of commutator expressions.
Or, alternatively one can also argue as follows:
One can see that $\exp\left(\int\tr\phi \right)=\det \Rexp(\phi)$, indeed, this follow from continuous
 extension from the case of mBCH measures, where the statement is trivial.
 Then, for $t\sim 0$,
\[ \exp\left(\int t\cdot\tr\phi \right)=\det\Rexp(t\cdot\phi)
=\det\exp\left(\sum_{k=1}^\infty t^k\cdot\mu_{k,\mathrm R}(\phi)\right)
=\exp\left(\sum_{k=1}^\infty t^k\cdot\tr\mu_{k,\mathrm R}(\phi)\right).\]
Taking logarithm is possible for $t\sim0$, leading to $t\cdot\int\tr\phi=\sum_{k=1}^\infty t^k\cdot\tr\mu_{k,\mathrm R}(\phi)$ for $t\sim 0$.
This implies the statement.

Now, if the complex matrix $A$ has Jordan form  $\bem\lambda&1\\&\lambda\eem$
 and $A=\exp C$, then $C$ has  Jordan form  $\bem\varkappa&1\\&\varkappa\eem$ with $\lambda=\exp\varkappa$.
Therefore, if $\lambda\neq 1$, this precludes $C$ to be a matrix of trace $0$.

The following counterexamples in Examples \ref{ex:Moan}, \ref{ex:MagCrit}, \ref{ex:vin}, \ref{ex:vinimp}, \ref{ex:vinimp15},
work similarly:

We will have a measure in terms of real $2\times2$ matrices, whose
time-ordered exponential $E$ is not an exponential of a real $2\times2$ matrix.
That precludes the absolute convergence / convergence / Abel-summability / radial convergence of the
Magnus-expansion, because the limit value should also be real  $2\times2$ matrix $M$,
whose exponential is $E$ (cf. Theorem \ref{th:MagnusExp}), which is a contradiction.
(Of, course, after we have established divergence, it does not matter if we consider
the matrices as complex matrices.)

Or, we will have a measure in terms of complex $2\times2$ matrices but whose cumulative trace is $0$, whose
time-ordered exponential $E$ is a matrix of Jordan form  $\bem-1&1\\&-1\eem$.
This also precludes even radial convergence, as radial convergence should go through matrices
of vanishing trace, which do not exponentiate to the indicated Jordan form.
\\
\snewpage

\subsection{Counterexamples to the Magnus expansion}
\plabel{ss:CountMagnus}
~\\

The following example due to Moan \cite{Ma}, cf. also Sch\"affer  \cite{Sch}, shows that the convergence bound $\pi$ cannot be improved
for the Magnus expansion, not even in the case of $2\times2$ real matrices:

\begin{example}\plabel{ex:Moan}
(Moan's / Sch\"affer's example.)
For $\theta\in[0,\pi]$, let
\[\MS(\theta)=\begin{bmatrix}
\cos \theta&\theta\cos \theta -\sin \theta\\\sin \theta&\theta\sin \theta+\cos \theta
\end{bmatrix}.\]
Let us consider the measure
\begin{equation}
\hat\Phi= \frac{\mathrm d\MS(\theta)}{\mathrm dt} (\MS(\theta))^{-1}\,\mathrm d\theta|_{[0,\pi]}\equiv
\frac12
\begin{bmatrix}
-\sin2\theta& -1+\cos2\theta\\1+\cos2\theta&\sin2\theta
\end{bmatrix}
\,\mathrm d\theta|_{[0,\pi]}.
\plabel{eq:moan}
\end{equation}
As
$\frac12\begin{bsmallmatrix}-\sin2\theta& -1+\cos2\theta\\1+\cos2\theta&\sin2\theta\end{bsmallmatrix}
=\begin{bsmallmatrix}\cos\theta& -\sin\theta\\\sin\theta&\cos\theta\end{bsmallmatrix}
\begin{bsmallmatrix}&0 \\1&\end{bsmallmatrix}
\begin{bsmallmatrix}\cos\theta& -\sin\theta\\\sin\theta&\cos\theta\end{bsmallmatrix}^{-1}
$ has norm $1$, we find
\[\int\|\hat\Phi\|_2=\pi.\]

On the  other hand, due to its construction,
\[\Lexp(\hat\Phi)=\MS(\pi)\MS(0)^{-1}=\bem-1&-\pi\\&-1\eem.\]
This is not the exponential of a real matrix, thus the Magnus expansion of $\hat\Phi$ cannot be convergent,
meanwhile the cumulative norm of $\hat\Phi$ is $\pi$.
\qedexer
\end{example}
\snewpage
Another notable example is
\begin{example}\plabel{ex:MagCrit}
(Magnus ``critical'' case.)
For $\theta\in[0,\pi]$, let
\[\MC(\theta)=\begin{bmatrix}
\cos \theta&2\theta\cos \theta -\sin \theta\\\sin \theta&2\theta\sin \theta+\cos \theta
\end{bmatrix}.\]
Let us consider the measure
\begin{equation}
\Phi= \frac{\mathrm d\MC(\theta)}{\mathrm dt} (\MC(\theta))^{-1}\,\mathrm d\theta|_{[0,\pi]}\equiv
\begin{bmatrix}
-\sin2\theta& \cos2\theta\\\cos2\theta&\sin2\theta
\end{bmatrix}
\,\mathrm d\theta|_{[0,\pi]}.
\plabel{eq:MagCrit}
\end{equation}
As we see orthogonal (in fact, reflection) matrices in \eqref{eq:MagCrit},
\[\int\|\Phi\|_2=\pi.\]

On the  other hand, due to its construction,
\[\Lexp(\Phi)=\MC(\pi)\MC(0)^{-1}=\bem-1&-2\pi\\&-1\eem.\]
This is not the exponential of a real matrix, thus the Magnus expansion of $\Phi$ cannot be convergent,
meanwhile the cumulative norm of $\Phi$ is $\pi$.
\qedexer
\end{example}
Examples \ref{ex:Moan} and \ref{ex:MagCrit} do not seem to be particularly different.
However, as we see later, Example \ref{ex:MagCrit} is not only more ``extremal'' than Example \ref{ex:Moan},
but also more manageable.
In later view (see Part IIA), we could term  Example \ref{ex:Moan} (Moan's example) as a Magnus elliptic development with
contraction factor  $h=\frac12$ and total ``mass'' $\pi$, and
Example \ref{ex:MagCrit} (Magnus critical case) as a Magnus parabolic development with total ``mass'' $\pi$.
\\
\snewpage

\subsection{Counterexamples to the BCH expansion}
\plabel{ss:CountBCH}
~\\

The ``Minimal Examples'' of \cite{L11} also apply to the case of $2\times2$ matrices with the operator norm.
However, we will redevelop counterexamples here.
We do this partly because of historical reasons, and
partly because we already have an eye toward the unitary / quaternionic cases and some other settings.
(Having skew-Hermitian matrices for the counterexamples is more ``physical'' and considered to  have more edge.)

\begin{example}\plabel{ex:vin}
Here we exhibit, for any $\varepsilon>0$, two $2\times2$ real matrices $V_1$, $V_2$
such that $\|V_1\|_2=\pi$, $\|V_2\|_2<\varepsilon$
but such that the Magnus (BCH) expansion of  the measure
$V_1\mathbf 1_{[0,1)}\boldsymbol.V_2\mathbf 1_{[1,2)}$ is divergent:

(a) Let
\[V_1^{[\delta]}=\pi\left[ \begin{matrix} &-1\\1&\end{matrix}\right],\qquad\text{and}
\qquad V_2^{[\delta]}=2\delta\left[\begin{matrix} 0&1\\&0\end{matrix}\right],\]
where $\delta\in\mathbb R\setminus \{0\}$.
The Magnus expansion of $V_1^{[\delta]}\mathbf 1_{[0,1)}\boldsymbol.V_2^{[\delta]}\mathbf 1_{[1,2)}$ is (completely) divergent,
 because otherwise it would be radially convergent to a
$2\times2$ matrix which exponentiates to
\[\exp_{\mathrm R}(V_1^{[\delta]}\mathbf 1_{[0,1)}\boldsymbol.V_2^{[\delta]}\mathbf 1_{[1,2)})
=(\exp V_1^{[\delta]})\cdot(\exp V_2^{[\delta]})=
\left[ \begin{matrix}-1&\\&-1\end{matrix}\right]\left[\begin{matrix} 1&2\delta\\&1\end{matrix}\right]=
\left[ \begin{matrix} -1&-2\delta\\&-1\end{matrix}\right];
\]
however, it is known that no $2\times 2$ real matrix exponentiates to the latter matrix.
As $\delta\rightarrow0$, it satisfies the indicated properties.
This is essentially the counterexample of Vinokurov \cite{V} (1991); it is a ``parabolic'' counterexample.

\begin{commentx}
[Warning: There is no $\delta\leftrightarrow-\delta$ symmetry here.]
\end{commentx}

(b) A similar, ``hyperbolic version'' can be made as follows: Let
\[\tilde V_1^{[\eta]}=\pi\left[ \begin{matrix} &-1\\1&\end{matrix}\right],
\qquad\text{and}\qquad \tilde V_2^{[\eta]}=\eta\left[\begin{matrix} 1&\\&-1\end{matrix}\right],\]
 where $\eta\in\mathbb R\setminus \{0\}$.
Here
\[\exp_{\mathrm R}(\tilde V_1^{[\eta]}\mathbf 1_{[0,1)}\boldsymbol.\tilde V_2^{[\eta]}\mathbf 1_{[1,2)})
=(\exp \tilde V_1^{[\eta]})\cdot(\exp \tilde V_2^{[\eta]})=
\left[ \begin{matrix}-1&\\&-1\end{matrix}\right]
\left[\begin{matrix} \mathrm e^{\eta}& \\&\!\!\mathrm e^{-\eta}\end{matrix}\right]=
\left[\begin{matrix} -\mathrm e^{\eta}& \\&\!\!-\mathrm e^{-\eta}\end{matrix}\right] \]
is of strictly hyperbolic type with (different) negative eigenvalues, which is also not an exponential of real $2\times2$ matrices.
Again, as $\eta\rightarrow0$ we have the desired counterexamples.

(c) Part (a) extends \textit{easily} to the complex case $\delta\in\mathbb C\setminus \{0\}$.
In this case we have to use that the time-ordered exponential is not an exponential of a complex matrix of trace $0$.
(We know $V_1^{[\delta]}\mathbf 1_{[0,1)}\boldsymbol.V_2^{[\delta]}\mathbf 1_{[1,2)}$ is of cumulative trace $0$,
so the (limiting) sum of its Magnus expansion should also be of trace $0$.)
\qedexer
\end{example}

The  examples above are ``unbalanced'' in the sense that the participating operators have norms of different scale.
There are, however, examples which are ``balanced'' in the sense that the participating operators have norms of equal scale.
Instead of modifying the previous example to a family
of balanced and unbalanced examples (which is possible) we will consider a balanced example on a slightly different ground:
\snewpage
\begin{example}\plabel{ex:vinimp}
Here we exhibit, for any $\varepsilon>0$, two $2\times2$ real matrices $A_1$ and $A_2$ such that $\|A_1\|_2=\|A_2\|_2<\frac\pi2+\varepsilon$
but such that the Magnus (BCH) expansion of  the measure $A_1\mathbf 1_{[0,1)}\boldsymbol.A_2\mathbf 1_{[1,2)}$ is divergent.

(a) Let
\[A_1^{[\delta]}=\log\left[\begin{matrix} 2\delta&-1\\1&\end{matrix}\right],
\qquad\text{and}\qquad A_2^{[\delta]}=\log\left[ \begin{matrix} &-1\\1&2\delta\end{matrix}\right],\]
where $\delta\in(-1,+\infty)$.
Indeed, the expressions above are meaningful.
The characteristic polynomials of the matrices $\left[\begin{matrix} 2\delta&-1\\1&\end{matrix}\right]$
and $\left[ \begin{matrix} &-1\\1&2\delta\end{matrix}\right]$
have characteristic polynomial $\varkappa^2-2\delta\varkappa+1$ in $\varkappa$.
This has either no real root (for $\delta\in(-1,1)$) or have positive roots $\delta\pm\sqrt{\delta^2-1}$  (for $\delta\in[1,+\infty)$).
Altogether, there are no spectral values in $(-\infty,0]$, showing these
matrices to be in the (open) logarithmic domain.

But let us restrict to the case $\delta\neq 0$ from now on. We find the following:

Firstly, the Magnus expansion of $A_1^{[\delta]}\mathbf 1_{[0,1)}\boldsymbol.A_2^{[\delta]}\mathbf 1_{[1,2)}$ is (completely) divergent,
because otherwise it would be radially convergent to a real
$2\times2$ matrix which exponentiates to
\[\exp_{\mathrm R}(A_1^{[\delta]}\mathbf 1_{[0,1)}\boldsymbol.A_2^{[\delta]}\mathbf 1_{[1,2)})
=(\exp A_1^{[\delta]})\cdot(\exp A_2^{[\delta]})=
\left[\begin{matrix} 2\delta&-1\\1&\end{matrix}\right]\left[ \begin{matrix} &-1\\1&2\delta\end{matrix}\right]=
\left[ \begin{matrix} -1&-4\delta\\&-1\end{matrix}\right];
\]
however, it is known that no $2\times 2$ real matrix exponentiates to the latter matrix.

Secondly, due to the continuity of $\log$ on its domain, as $\delta\rightarrow0$, we see that
\[A_1^{[\delta]},A_2^{[\delta]}\rightarrow\frac\pi2\left[\begin{matrix} &-1\\1&\end{matrix}\right]. \]
This shows that we can choose $\delta$ such that $\|A_1^{[\delta]}\|_2,\|A_2^{[\delta]}\|_2$ can be arbitrarily close to $\frac\pi2$.
Furthermore, $\|A_1^{[\delta]}\|_2=\|A_2^{[\delta]}\|_2$ because $A_2^{[\delta]}=\begin{bsmallmatrix}&-1\\1&\end{bsmallmatrix}
A_1^{[\delta]}\begin{bsmallmatrix}&-1\\1&\end{bsmallmatrix}^{-1}$.

This example which we prefer to call the case of balanced ``rigid pairs'' corresponds to the balanced subcase of the ``intermediately selected'' critical examples of Michel \cite{Mi}.
(Except there the Frobenius Banach--Lie norm on $\mathfrak{sl}_2(\mathbb R )$ is considered.)

\begin{commentx}
[Warning: There is no $\delta\leftrightarrow-\delta$ symmetry here.]
\end{commentx}

(b)
A similar hyperbolic example is given by
\[\tilde A_1^{[\eta]}=\log\left[\begin{matrix}&-\mathrm e^{-\eta} \\\mathrm e^{\eta}&\end{matrix}\right],
\qquad\text{and}\qquad \tilde A_2^{[\eta]}=\log\left[\begin{matrix}&-\mathrm e^{\eta} \\\mathrm e^{-\eta}&\end{matrix}\right],\]
where $\eta\in\mathbb R\setminus \{0\}$.
(The matrices, whose logarithm is taken have no  eigenvalues from $(-\infty,0]$ but $\pm\mathrm i$.)
Again, here $\exp_{\mathrm R}(\tilde A_1^{[\eta]} \mathbf 1_{[0,1)}\boldsymbol.\tilde A_2^{[\eta]} \mathbf 1_{[1,2)})$
 will be of strictly hyperbolic type with negative eigenvalues, which is also not an exponential of real $2\times2$ matrices.
(We remark that
\[\tilde A_1^{[\eta]}=\frac\pi2\left[\begin{matrix}&-\mathrm e^{-\eta} \\\mathrm e^{\eta}&\end{matrix}\right],
\qquad\text{and}\qquad \tilde A_2^{[\eta]}=\frac\pi2\left[\begin{matrix}&-\mathrm e^{\eta} \\\mathrm e^{-\eta}&\end{matrix}\right],\]
as
$\begin{bsmallmatrix}&-\mathrm e^{\mp\eta} \\\mathrm e^{\pm\eta}&\end{bsmallmatrix}$
are just the conjugates of
$\begin{bsmallmatrix}&-1 \\1&\end{bsmallmatrix}$
by
$\begin{bsmallmatrix}\mathrm e^{\mp\eta/2} &\\&\mathrm e^{\pm\eta/2}\end{bsmallmatrix}$.)
\snewpage

(c) Part (a) extends to the complex case $\delta\in\left(\mathbb C\setminus (-\infty,-1]\right)\setminus \{0\}$.
Indeed, if $\delta\in \mathbb C\setminus (-\infty,-1]$, then the logarithms of $\left[\begin{matrix} 2\delta&-1\\1&\end{matrix}\right]$
 and $\left[ \begin{matrix} &-1\\1&2\delta\end{matrix}\right]$ can be taken:
The additional cases are when $\delta\in\mathbb C\setminus\mathbb R$,
 but then $\varkappa^2+1=\delta\varkappa$ cannot happen with $\varkappa\in(-\infty,0])$.
Furthermore, from the determinant of the matrices $\left[\begin{matrix} 2\delta&-1\\1&\end{matrix}\right]$
 and $\left[ \begin{matrix} &-1\\1&2\delta\end{matrix}\right]$ (which is $1$), we can see that the
 trace of the logarithms is constant, an element of $2\pi\mathrm i\mathbb Z$;
 but then it must be $0$ as there are real matrices among the logarithms.
In the rest, the argument is as in (a), but
 we have to use that the time-ordered exponential is not an exponential of a matrix of trace $0$.
($A_1^{[\delta]}\mathbf 1_{[0,1)}\boldsymbol.A_2^{[\delta]}\mathbf 1_{[1,2)}$ is of cumulative trace $0$.)
\qedexer
\end{example}
\begin{commentx}
For $0<\delta<1$, we have $\log\left[\begin{matrix} 2\delta&-1\\1&\end{matrix}\right]=\dfrac{\arccos\delta}{\sqrt{1-\delta^2}}
\bem \delta&-1\\1&-\delta \eem$.

Or, with $\delta=\cos u$, we have
$\log\left[\begin{matrix} 2\cos u&-1\\1&\end{matrix}\right]=u\bem \cot u&-\csc u\\ \csc u&-\cot u \eem$.
\end{commentx}

At this point, however, the lack of ready-to-use computational tools regarding $\exp$, $\log$ and the operator norm becomes impeding.
In the following Section \ref{sec:LogA}, we will review these technical tools concerning $2\times2$ matrices.
Then, in Section \ref{sec:BCHsharp}, we will pick up the topic of counterexamples again.

\begin{commentx}
There, we will address:
\begin{itemize}
\item Detailed convergence analysis of classical cases.
\item The unitary / quaternionic cases.
\item The topic of arbitrary balance ratios.
\item ``Sharp'' counterexamples to the BCH case (in the infinite dimensional setting).
\end{itemize}
\end{commentx}

\snewpage
\section{Exponentials, logarithms, and norms of $2\times2$ matrices}
\plabel{sec:LogA}
Here we collect various tools from elementary real and complex analysis
concerning $2\times2$ matrices, which will be used in this paper, and also in Part IIA.
Everything is easy to check here but it might be distracting doing so later in particular examples.
\\

\subsection{The skew-quaternionic form}~\\

One can write the $2\times2$ matrix $A$ in skew-quaternionic form
\begin{equation}
A= \tilde a\Id_2 +\tilde b\tilde I+\tilde c\tilde J+\tilde d\tilde K\equiv
\tilde a\begin{bmatrix} 1&\\&1\end{bmatrix}
+\tilde b\begin{bmatrix} &-1\\1&\end{bmatrix}
+\tilde c\begin{bmatrix} 1&\\&-1\end{bmatrix}
+\tilde d\begin{bmatrix} &1\\1&\end{bmatrix}.
\plabel{eq:skqform}
\end{equation}

This is particularly, useful in the case of real matrices.
Conjugation by matrices of shape $(\cos\alpha)\Id_2+(\sin\alpha)\tilde I$
 leaves $\Id_2$ and $\tilde I$ invariant, but rotates the basis element $\tilde J,\tilde K$,
 while it is also isometric with respect to the usual $\ell^2$ operator norm.
The following observations are extremely elementary but they have some practical consequences
 regarding the presentation of some examples.
Let
\[\tilde L=\left(\cos\frac\pi4\right)\Id_2+\left(\sin\frac\pi4\right)\tilde I=
\frac{\sqrt2}2\bem1&-1\\1&1\eem.\]
Then conjugation by $\tilde L$, i. e. the map $X\mapsto \tilde LX\tilde L^{-1}$,
takes $\tilde J\mapsto\tilde K\mapsto-\tilde J \mapsto -\tilde K$, while leaves $\Id_2$ and $\tilde I$ invariant.
In particular, the seemingly different-looking matrices
\[\tilde I+\delta\tilde J=\bem\delta&-1\\1&-\delta\eem
,\qquad
\tilde I+\delta\tilde K=\bem& -(1-\delta)\\1+\delta&\eem
,
\]\[
\tilde I-\delta\tilde J=\bem-\delta&-1\\1&\delta\eem
,\qquad
\tilde I-\delta\tilde K=\bem& -(1+\delta)\\1-\delta&\eem
\]
are isometric to each other.
Thus, it is often a choice how to present counterexamples.
This, unfortunately, cannot be done in a completely uniformed way,
as different things look simpler or more complicated in logarithmic or exponential form.
Furthermore, the isometry concrete matrices above implies (without studying the computation of the operator norm in general)
that
\begin{multline}
\left\|\bem\delta&-1\\1&-\delta\eem\right\|_2=\left\|\bem&\!\!\!\!-(1-\delta)\\1+\delta&\eem\right\|_2
=\left\|\bem-\delta&-1\\1&\delta\eem\right\|_2=\left\|\bem&\!\!\!\!-(1+\delta)\\1-\delta&\eem\right\|_2=
\\\\
=\left\| \bem 1+\delta& \\&1-\delta\eem\right\|_2= {\max(|1-\delta|,|1+\delta|)}=:{\Delta(\delta)}.
\plabel{eq:uninorm}
\end{multline}
Thus, we can relatively easily compute with the norm of these matrices, although
one may have to discriminate on the cases $\Rea\delta\geq0$ and $ \Rea\delta\leq0$.

In form \eqref{eq:skqform},
\[\frac{\tr A}2=\tilde a, \]
and
\[\det A=\tilde a^2+\tilde b^2-\tilde c^2-\tilde d^2.\]
\subsection{Spectral type}~\\

Let us use the notation
\[D_A=\det\left( A-\frac{\tr A}2\Id_2\right)=(\det A)-\frac{(\tr A)^2}4.\]
It is essentially the discriminant of $A$, as the eigenvalues of $A$ are $\frac12\tr A\pm\sqrt{-D_A}$.

In form \eqref{eq:skqform},
\begin{equation}D_A=\tilde b^2-\tilde c^2-\tilde d^2.\plabel{nt:DA}\end{equation}

In the special case of real $2\times2$ matrices, we use the classification

$\bullet$ elliptic case: two conjugate strictly complex eigenvalues,

$\bullet$ parabolic case: two equal real eigenvalues,

$\bullet$ hyperbolic case: two distinct real eigenvalues.

Then, for real $2\times2$ matrices,  $D_A$ measures `ellipticity/parabolicity/hiperbolicity':
If $D_A>0$, then $A$ is elliptic;
if $D_A=0$, then $A$ is parabolic;
if $D_A<0$, then $A$ is hyperbolic.

In the general complex case, there are two main categories: parabolic $(D_A=0)$ and non-parabolic $(D_A\neq 0)$.
\\

\subsection{Principal and chiral disks}~\\

For $2\times2$ real matrices we can refine the spectral data as follows:
Assume that $A=\bem a&b\\c&d\eem =\tilde a\Id_2 +\tilde b\tilde I+\tilde c\tilde J+\tilde d\tilde K$.
Its principal disk is
\begin{align*}
\PD(A):=&\Dbar\left(\frac{a+d}2+\frac{|c-b|}2\mathrm i,
\sqrt{\left(\frac{a-d}2\right)^2+\left(\frac{b+c}2\right)^2} \right)
\\
=&\Dbar\left( \tilde a+|\tilde b|\mathrm i,\sqrt{\tilde c^2+\tilde d^2} \right)
.
\end{align*}
The principal disk is a point if $A$ has the effect of a complex multiplication (that $A$ is a quasicomplex matrix).
In general, matrices $A$ fall into three categories: elliptic, parabolic, hyperbolic;
such that the principal disk are disjoint, tangent or secant to the real axis, respectively.
This is refined further by the chiral disk
\begin{align*}
\CD(A):=&\Dbar\left(\frac{a+d}2+\frac{c-b}2\mathrm i,
\sqrt{\left(\frac{a-d}2\right)^2+\left(\frac{b+c}2\right)^2} \right)
\\
=&\Dbar\left( \tilde a+\tilde b\mathrm i,\sqrt{\tilde c^2+\tilde d^2} \right)
.
\end{align*}
The additional data in the chiral disk is the chirality, which is the sign of the twisted trace,
$\sgn\left( \tr\begin{bmatrix}&1\\-1&\end{bmatrix}A\right)=\sgn (c-b)=\sgn\tilde b.$
This chirality is, in fact, understood with respect to a fixed orientation of $\mathbb R^2$.
It does not change if we conjugate $A$ by a rotation, but it changes sign if
we conjugate $A$ by a reflection.
From the properties of the twisted trace, it is also easy too see that $\log$ respects chirality.

One can read off many data from the disks.
For example, if $\PD(A)=\Dbar((\tilde a,\tilde b),r)$, then $\det A=\tilde a^2+\tilde b^2-r^2$.
This is not surprising in the light of
\begin{lemma}\plabel{lem:PDCD}
$\CD$ makes a bijective correspondence between
possibly degenerated disks in $\mathbb C$ and
the orbits of $\mathrm M_2(\mathbb R)$ with respect to
conjugacy by special orthogonal matrices (i. e. rotations).

$\PD$ makes a bijective correspondence between
possibly degenerated disks with center in $\mathbb C^+$
and the orbits of $\mathrm M_2(\mathbb R)$ with respect to
conjugacy by  orthogonal matrices.
\begin{proof}
One can write $A\in\mathrm M_2(\mathbb R)$ in the skew-quaternionic form \eqref{eq:skqform}.
In that way, it is clear that every possibly degenerated disk occurs as conformal disk.
On the other hand, conjugation by $
(\cos\alpha)\Id_2+(\sin\alpha)\tilde I=
\begin{bmatrix} \cos\alpha&-\sin\alpha\\\sin\alpha&\cos\alpha\end{bmatrix}$
takes $A$ into $\tilde a\Id +\tilde b\tilde I+(\tilde c\cos2\alpha-\tilde d\sin2\alpha)\tilde J
+(\tilde c\sin2\alpha+\tilde d\cos2\alpha)\tilde K$.
This shows that the rotational orbit data is the same as the conformal disk data.
Conjugation by $\tilde J=\begin{bmatrix} 1&\\&-1\end{bmatrix}$
takes $A$ into $\tilde a\Id -\tilde b\tilde I+\tilde c\tilde J-\tilde d\tilde K$.
This shows the second part.
\end{proof}
\end{lemma}

\snewpage
\subsection{Exponentials}

\begin{lemma}\plabel{lem:expSep}
Consider a real $2\times2$ matrix $A$, written as in \eqref{eq:skqform}.

If $-\tilde b^2+\tilde c^2+\tilde d^2<0$, then
\[\exp A=(\exp \tilde a)\cdot\left(\left(\cos \sqrt{\tilde b^2-\tilde c^2-\tilde d^2}\right)\Id_2
+\left(\sin \sqrt{\tilde b^2-\tilde c^2-\tilde d^2}\right)\frac {\tilde b\tilde I+\tilde c\tilde J+\tilde d\tilde K}{\sqrt{\tilde b^2-\tilde c^2-\tilde d^2}} \right);\]
if $-\tilde b^2+\tilde c^2+\tilde d^2=0$, then
\[\exp A=(\exp \tilde a)\cdot\left(\Id_2
+\tilde b\tilde I+\tilde c\tilde J+\tilde d\tilde K \right);\]
if $-\tilde b^2+\tilde c^2+\tilde d^2>0$, then
\[\exp A=(\exp \tilde a)\cdot\left(\left(\cosh \sqrt{-\tilde b^2+\tilde c^2+\tilde d^2}\right)\Id_2
+\left(\sinh \sqrt{-\tilde b^2+\tilde c^2+\tilde d^2}\right)\frac {\tilde b\tilde I+\tilde c\tilde J+\tilde d\tilde K}{\sqrt{-\tilde b^2+\tilde c^2+\tilde d^2}} \right).\]
\begin{proof}
Here one just faces the situation $\exp(\tilde a+\beta H)$ with $H^2=-1,0,$ or $1$.
\end{proof}
\end{lemma}

We define the functions $\Cos$ and $\Sin$  by
\[\Cos(x)=\begin{cases}
\cos\sqrt{x}&\text{if }x>0 \\
1&\text{if }x=0\\
\cosh\sqrt{-x}&\text{if }x<0,
\end{cases}\]
and
\[\Sin(x)=\begin{cases}
\dfrac{\sin\sqrt{x}}{\sqrt{x}}&\text{if }x>0 \\
1&\text{if }x=0\\
\dfrac{\sinh\sqrt {-x}}{\sqrt {-x}}&\text{if }x<0,
\end{cases}\]
on the real domain.
However, it is easy to see, $\Cos$ and $\Sin$ extend to entire functions on the complex plane.
In this way,
\begin{lemma}\plabel{lem:expUni}
Assume that $A$ is a complex $2\times2$ matrix, written as in \eqref{eq:skqform}.
Then
\[\exp A=(\exp \tilde a)\cdot\left(\left(\Cos \left(\tilde b^2-\tilde c^2-\tilde d^2\right)\right)\Id_2
+\left(\Sin \left(\tilde b^2-\tilde c^2-\tilde d^2\right)\right)   \left(\tilde b\tilde I+\tilde c\tilde J+\tilde d\tilde K\right) \right).\]
\begin{proof}
The real case extends by analytic continuation.
\end{proof}
\end{lemma}

\begin{lemma}\plabel{lem:expQuasi}
Let $A$ be a complex $2\times2$ matrix. Then
\[\exp A=\left(\exp\frac{\tr A}2\right)\cdot
\left(
\Cos \left(D_A\right)\Id_2
+\Sin \left(D_A\right)\left(A-\frac{\tr A}2\Id_2\right)
\right).\]
\begin{proof}
This is just the transcription of the previous lemma using \eqref{nt:DA}.
\end{proof}
\end{lemma}~

\subsection{A simple differential equation}~\\

If $M$ and $W$ are matrices (or just elements of a Banach algebra), then the functions
 $A(\theta)=\exp(\theta W)\exp(\theta(M-W))$ and $B(\theta)=\exp(\theta(M+W))\exp(-\theta W)$
 (for $\theta\in\mathbb R$ or $\mathbb C$) satisfy the differential eqiuations
\[\frac{\mathrm dA(\theta)}{\mathrm d\theta}A(\theta)^{-1}=\exp(\theta W)M\exp(-\theta W)\]
 and
\[B(\theta)^{-1}\frac{\mathrm dB(\theta)}{\mathrm d\theta}=\exp(\theta W)M\exp(-\theta W)\]
respectively; with initial data $A(0)=1$ and $B(0)=1$, respectively.
This leads to
\begin{lemma}\plabel{lem:genprec}
Assume that $r \geq0$, $M$, $W$ are elements of a Banach algebra.
Then
\[\Lexp\Bigl(\theta\in[0,r)\mapsto \exp(\theta W)M\exp(-\theta W)\Bigr)=\exp(\theta W)\exp(\theta(M-W)),\]
\[\Rexp\Bigl(\theta\in[0,r)\mapsto \exp(\theta W)M\exp(-\theta W)\Bigr)= \exp(\theta(M+W))\exp(-\theta W). \]
\begin{proof}
The previous observation implies the statement in the case finite dimensional algebras,
which implies the formal case, which implies the general case.
\end{proof}
\end{lemma}

For us a special case will be relatively important; we spell it out as
\begin{lemma}\plabel{lem:fullprec}
The solution of the ordinary differential equation
\[\frac{\mathrm dA(\theta)}{\mathrm d\theta}A(\theta)^{-1}=
a\bem&-1\\1&\eem+
b\begin{bmatrix}
-\sin 2c\theta& \cos2c\theta\\\cos2c\theta&\sin2c\theta
\end{bmatrix}
\equiv \exp(c\theta\tilde I)(a\tilde I+ b\tilde K )\exp(-c\theta\tilde I),
\]
with initial data
\[A(0)=\begin{bmatrix}
1& \\&1
\end{bmatrix} \equiv\Id_2,\]
is given by
\[A(\theta)=F(a\theta,b\theta,c\theta);\]
where
\[
F(s,p,w)=  \exp(w\tilde I) \exp ( (s-w)\tilde I+p\tilde K) .
\]
\begin{commentx}
\begin{align*}
F(s,p,w)&=\exp(w\tilde I) \exp ( (s-w)\tilde I+p\tilde K)\\
&\equiv  \exp(w\tilde I)(\Cos((s-w)^2-p^2 ) \Id_2+\Sin((s-w)^2 -p^2) ((s-w)\tilde I+p\tilde K) ) \\
&\equiv  \begin{bmatrix}
\cos w&-\sin w\\\sin w&\cos w
\end{bmatrix}
\begin{bmatrix}
\Cos((s-w)^2-p^2)&(p-(s-w))\Sin((s-w)^2 -p^2) \\(p+(s-w))\Sin((s-w)^2 -p^2)&\Cos((s-w)^2-p^2)
\end{bmatrix}
\end{align*}
\end{commentx}
That, in particular, makes
\[\Lexp\left(\left(a\bem&-1\\1&\eem+
b\begin{bmatrix}
-\sin 2c\theta& \cos2c\theta\\\cos2c\theta&\sin2c\theta
\end{bmatrix}\right)d\theta|_{[0,r]}\right)=F(ar,br,cr).
\eqed\]
\begin{commentx}
\begin{lemma}
The solution of the ordinary differential equation
\[\frac{\mathrm dB(\theta)}{\mathrm d\theta}B(\theta)^{-1}=
a\bem&-1\\1&\eem+
b\begin{bmatrix}
-\sin 2c\theta& \cos2c\theta\\\cos2c\theta&\sin2c\theta
\end{bmatrix}
\equiv \exp(c\theta\tilde I)(a\tilde I+ b\tilde K )\exp(-c\theta\tilde I),
\]
with initial data
\[B(0)=\begin{bmatrix}
1& \\&1
\end{bmatrix} \equiv\Id_2,\]
is given by
\[B(\theta)=G(a\theta,b\theta,c\theta);\]
where
\[
G(s,p,w)=  \exp ( (s+w)\tilde I+p\tilde K)\exp(-w\tilde I) .
\]
That, in particular, makes
\[\Rexp\left(\left(a\bem&-1\\1&\eem+
b\begin{bmatrix}
-\sin 2c\theta& \cos2c\theta\\\cos2c\theta&\sin2c\theta
\end{bmatrix}\right)d\theta|_{[0,r]}\right)=G(ar,br,cr).
\eqed\]
\end{lemma}

\end{commentx}
\end{lemma}

We will also use the special notation
\[W(p,w)=F(0,p,w)\equiv  \exp(w\tilde I) \exp (-w\tilde I+p\tilde K) .\]
\begin{commentx}
Also, let
\[E(p,w)=F(p-w,w,p)=\begin{bmatrix}
\cos p&2w\cos p -\sin p\\\sin p&2w\sin p+\cos p
\end{bmatrix}=
(\cos p\Id_2+\sin p\tilde I)(\Id_2-w\tilde I+w\tilde K).\]
\end{commentx}

\begin{remark}
By conjugation, $(\exp -\frac\pi4 \tilde I)\tilde K(\exp \frac\pi4 \tilde I)=\tilde J$,
thus we could have used $\tilde J$ instead of $\tilde K$.
Here and later the choice and use of basis elements is either accidental or influenced by
practical considerations like that the nilpotent element
$\frac12(\tilde I+\tilde K)=\begin{bsmallmatrix}&0\\1&\end{bsmallmatrix}$ is simpler to the human eye than
$\frac12(\tilde I+\tilde J)=\frac12\begin{bsmallmatrix}1&-1\\1&-1\end{bsmallmatrix}$.
\qedremark
\end{remark}
~

\subsection{The quaternionic alternative}~\\

In the complex case, it is tempting to use
quaternionic base,
\begin{equation*}
A= \check a\Id_2 +\check b I+\check c J+\check d K\equiv
\check a\begin{bmatrix} 1&\\&1\end{bmatrix}
+\check b\begin{bmatrix} &-1\\1&\end{bmatrix}
+\check c\begin{bmatrix} \mathrm i&\\&-\mathrm i\end{bmatrix}
+\check d\begin{bmatrix} &\mathrm i\\\mathrm i&\end{bmatrix};
\end{equation*}
($I=\tilde I$, $ J=\mathrm i\tilde J$, $ K=\mathrm i\tilde K$),
or simply to put forward quaternionic examples,
like $F(a,\mathrm i b,c)$, instead of real matrices.
For us, the advantages would be limited, but this happens in physics, see Remark \ref{rem:physex}.
(In terms of Pauli matrices $I=-\mathrm i\sigma_2$, $J= \mathrm i\sigma_1$, $K=\mathrm i\sigma_3$.)
\\

\snewpage
\subsection{The differential calculus of $\Cos$ and $\Sin$}~\\

First of all, it is useful to notice that
\[z=\frac{1-\Cos(z)^2}{\Sin(z)^2}\]
(as entire analytic functions).
Then one can easily see that
\[\Cos'(z)=-\frac12\Sin(z)\]
and
\begin{align*}
\Sin'(z)&=\frac{\Cos(z)-\Sin(z)}{2z}\\
&=\frac12  \frac{\Sin(z)^2\cdot(\Cos(z)-\Sin(z))}{1-\Cos(z)^2}
\end{align*}
(as entire analytic functions).
In particular, differentiation will not lead out of the rational field
generated by $\Cos$ and $\Sin$.
\\

\snewpage
\subsection{Logarithms}~\\

Let us define the function $\AC$ by
\[\AC(x)=\begin{cases}
\dfrac{\arccos x}{\sqrt{1-x^2}}&\text{if }-1< x<1\\[3mm]
1&\text{if }x=1\\[1mm]
\dfrac{\arcosh x}{\sqrt{x^2-1}}\qquad&\text{if } 1<x.\\
\end{cases}\]
\begin{lemma} \plabel{lem:logreal}
Let $A$ be a  real $2\times 2$ matrix.

Then  $A$ is a $\log$-able if and only if
 $\det A>0$ and $\dfrac{\tr A}{2\sqrt{\det A}}>-1$.

In the $\log$-able case
\begin{equation}\log A= (\log\sqrt{\det A})\Id_2+\frac{\AC\left(\dfrac{\tr A}{2\sqrt{\det A}}\right)}{\sqrt{\det A}}
\left(A- \frac{\tr A}2\Id_2\right).\plabel{eq:log2}\end{equation}

\begin{proof}
If $\det A\leq0$ then $\det \lambda\Id_2 +(1-\lambda)A =0$ for some $\lambda\in[0,1]$, in which case
$\lambda\Id_2 +(1-\lambda)A$ is not invertible, thus this clearly falls to the the not  $\log$-able case,
with no more to prove.

So, we can assume that $\det A>0$. In this case $A$ is $\sqrt{\det A}$ times an element from $\SL_2(\mathbb R)$.
According to general functional calculus (valid even in general Banach algebras), multiplication by $\alpha>0$ leaves being
$\log$-able invariant, and in  the $\log$-able case $\log(\alpha A)=(\log\alpha)\Id_2+\log A$.
Thus, as the determinant can be factorized out, and
everything is conjugation invariant, it is sufficient to consider only a set of orbit type representatives with respect to  the
$\SL_2(\mathbb R)$ part.

In the $\log$-able case the corresponding orbit types are
$\begin{bmatrix}1&\\&1\end{bmatrix}$,
$\begin{bmatrix}\cos \alpha&-\sin\alpha\\\sin\alpha&\cos\alpha\end{bmatrix}$ ($\alpha\in (0,\pi/2]$),
$\begin{bmatrix}1&1\\&1\end{bmatrix}$,
$\begin{bmatrix}\mathrm e^\beta&\\& \mathrm e^{-\beta}\end{bmatrix}$ ($\beta>0$);
and the statement be checked separately in each case.
(Alternatively, standard diagonalizability and density tricks can be used for the formula.)

In the not $\log$-able case the corresponding orbit types are
$\begin{bmatrix}-1&\\&-1\end{bmatrix}$ ,
$\begin{bmatrix}-1&1\\&-1\end{bmatrix}$,
$\begin{bmatrix}-\mathrm e^\beta&\\& -\mathrm e^{-\beta}\end{bmatrix}$ ($\beta>0$) , in which cases
$\dfrac{\tr A}{2\sqrt{\det A}}\leq-1$.
\end{proof}
\end{lemma}
As the proof indicates, we compute $\AC$  by $\arccos$ for elliptic matrices, by $\arcosh$ for hyperbolic matrices, and as $1$
for parabolic matrices.
\snewpage
\begin{lemma}\plabel{lem:AC}
The function $\AC$ extends to $\mathbb C\setminus (-\infty,-1]$ analytically.
$\AC$ is monotone decreasing on $(-1,\infty)$ with range $(0,\infty)$.
Special values are $\AC(0)=\frac\pi2$, $\AC(1)=1$.
\begin{proof}
According to Lemma \ref{lem:logreal}, for $z\in(-1,\infty)$,
\begin{equation}
\AC(z)=\frac12\tr\left(\begin{bmatrix}&1\\-1&\end{bmatrix}\log\begin{bmatrix}z&z-1\\z+1&z\end{bmatrix}\right).\plabel{eq:ACext}
\end{equation}
hold. However, he RHS of the equation, is well-defined for any $z\in\mathbb C\setminus (-\infty,-1]$.
Indeed, the eigenvalues of the matrix under the $\log$ are $z\pm\sqrt{z^2-1}$.
The equation $z\pm\sqrt{z^2-1}=r\leq0$, however, solves to $z=\frac{r+1/r}{2}<0$, excluded by assumption.
Thus, the corresponding matrices are in the domain of $\log$.
The expression is analytic, and it extends the original definition.
Monotonicity can be obtained by elementary function calculus.
\end{proof}
\begin{proof}[Second proof]
We start again by the observation that \eqref{eq:ACext} holds for $z\in(-1,\infty)$.
Using the spectral formula \eqref{eq:logdef3}, this explicitly yields
\begin{equation}
\AC(z)
=\int_{\lambda=-\infty}^{0}\frac{\mathrm d\lambda}{\lambda^2-2\lambda z+1}
=\int_{\lambda=-1}^{0}\frac{2\,\mathrm d\lambda}{\lambda^2-2\lambda z+1}
=\int_{\lambda=-\infty}^{-1}\frac{2\,\mathrm d\lambda}{\lambda^2-2\lambda z+1}
\plabel{eq:ACext15}
.\end{equation}
Applying change of variable for $\lambda=\kappa-\sqrt{\kappa^2-1}$, $\kappa\in(-\infty,-1)$ to the
very last term of \eqref{eq:ACext15}, this yields
\begin{equation}
\AC(z)=-\int_{\kappa=-\infty}^{-1}\frac{1}{\sqrt{\kappa^2-1}}\,\frac{\mathrm d\kappa}{\kappa-z}.
\plabel{eq:ACext2}
\end{equation}
In this latter form, however, analytic extendibility to $z\in\mathbb C\setminus (-\infty,-1]$
is immediate.
The monotonicity of $\AC$ on $(-1,\infty)$ is particularly transparent from \eqref{eq:ACext2}.
(Ultimately, we have traded checking the spectral condition for matrices to change of variables
in integration.
In fact, the spectral integral \eqref{eq:ACext2} can also be obtained by other means.)
\end{proof}
\begin{proof}[Third proof]
For the simply connected domain  $\mathbb C\setminus \left(  (-\infty,-1]\cup[1,\infty) \right))$,
$\arccos z$ extends as a primitive function of $\frac1{\sqrt{1-z}\sqrt{1+z}}$.
Consequently, $\AC(z)$ also extends to this domain as $\frac{\arccos z}{\sqrt{1-z}\sqrt{1+z}} $.
If we can prove that $\AC$ is real-analytic on $(-1,\infty)$, then by unicity, $\AC(z)$
analytically extends to $\mathbb C\setminus(\infty,1]$.
Now, the real-analyticity of $\AC$ is obvious except at $z=1$.
There we can apply the following argument:
$\Cos(0)=1$, and $\Cos'(0)=-\frac12\neq0$.
Thus by the 1-dimensional complex inverse function theorem $\Cos$
can be inverted in a neighborhood of $0$.
Let us temporarily denote this local inverse by $\AN$.
Hence, $\AN(1)=0$, and $\AN$ is analytic near $1$.
Consequenty, $-\frac12\AN'$ is also analytic.
However, by explicit formulas it is easy to check that $-\frac12\AN'(z)=\AC(z)$
for in a pointed real neighbourhood of $1$ (cf.  \eqref{eq:ANreal}), thus, by the continuity of $\AC$,
$-\frac12\AN'(z)=\AC(z)$ is also true in a neighborhood of $1$.
Hence $\AC$ is real-analytic on $(-1,\infty)$. Etc.
\end{proof}
\end{lemma}

\begin{remark}
\plabel{rem:ATT}
In Part I, \cite{L1}, the function $\ATT$ was considered.
The simple relationship between $\ATT$ and $\AC$ is the following.
For $z\in (-1,+\infty)$ it is easy to check that
$\AC(z)=\dfrac2{1+z}\ATT\left(\dfrac{1-z}{1+z}\right)$
or, equivalently,
$\ATT(z)
=\dfrac1{1+z}{\AC\left(\dfrac{1-z}{1+z}\right)}$
holds. This can naturally be extended to $z\in\mathbb C\setminus(-\infty,-1]$.
\qedremark
\end{remark}

Suppose that $A$ is $n\times n$ complex matrix which is $\log$-able.
Then let $\sqrt{\dett A}$ denote the value  of the standard branch of the square root of the determinant on $\log$-able elements.
It can be realized as
\[\sqrt{\dett A}=\exp\frac{\tr \log A}{2}=\exp\frac12\int_{t=0}^1 \tr\frac{\mathrm d((1-t)\Id+tA )}{(1-t)\Id+tA},\]
or as $\sqrt{\varepsilon_1}\cdot\ldots\cdot\sqrt{\varepsilon_n}$, where $\varepsilon_i$ are the eigenvalues of $A$, and the square
roots are in $\mathbb C\setminus(-\infty,0]$. (A this point it would be very reasonable to use `$\det\sqrt A$' instead of the symbolic notation `$\sqrt{\dett A}$',
but we make the point the that we prefer to consider it as a single function, not as a composite function.)

\begin{lemma}\plabel{lem:logcomplex}
Suppose that $A$ is $2\times2$ complex matrix which is $\log$-able.

Then $\sqrt{\dett A}\in\mathbb C\setminus (\infty,0]$, $\dfrac{\tr A}{2\sqrt{\dett A}}\in\mathbb C\setminus (-\infty,-1]$, and
the extended form of \eqref{eq:log2} holds:
\begin{equation}\log A= (\log\sqrt{\dett A})\Id_2+\frac{\AC\left(\dfrac{\tr A}{2\sqrt{\dett A}}\right)}{\sqrt{\dett A}}
\left(A- \frac{\tr A}2\Id_2\right).\plabel{eq:log2c}\end{equation}
\begin{commentx}
(b) Assume that the eigenvalues of $A$ are on the unit circle.
Then, in particular, $\sqrt{\dett A}\in\partial\Dbar(0,1)\setminus \{-1\}$, $\dfrac{\tr A}{2\sqrt{\dett A}}\in(-1,1]$.
\end{commentx}
\begin{proof}Then $\varepsilon_i=\mathrm e^{\alpha_i}$, with $-\pi< \Ima \alpha_i <\pi$.
Hence,
$\det A=\mathrm e^{\frac{\alpha_1+\alpha_2}2}$, and $\left|\Ima\frac{\alpha_1+\alpha_2}2\right|<\pi$ is transparent.
Indirectly,
\[\frac{\tr A}{2\sqrt{\dett A}}=\frac{\mathrm e^{\frac{\alpha_1-\alpha_2}2}+ \mathrm e^{-\frac{\alpha_1-\alpha_2}2}  }2=r\leq-1\]
solves to
\[\mathrm e^{\pm\frac{\alpha_1-\alpha_2}2}=r\pm\sqrt{r^2-1}\leq 0.\]
But this contradicts  $\left|\Ima\frac{\alpha_1-\alpha_2}2\right|<\pi$.
The logarithm formula extends from \eqref{eq:log2} analytically.
\begin{commentx}
(b) is the case when $\alpha_i$ are purely imaginary.
\end{commentx}
\end{proof}
\end{lemma}
If $A$ is real, then
\[\log A=\left(\log \sqrt{\tilde a^2+\tilde b^2-\tilde c^2-\tilde d^2}\right)\Id_2+
\frac{\AC\left(  \dfrac{\tilde a}{\sqrt{\tilde a^2+\tilde b^2-\tilde c^2-\tilde d^2}}  \right)}{\sqrt{\tilde a^2+\tilde b^2-\tilde c^2-\tilde d^2}}\left(\tilde b\tilde I+\tilde c\tilde J+\tilde d\tilde K\right). \]
The formula can also be used in the complex case but the choice of
$\sqrt{\tilde a^2+\tilde b^2-\tilde c^2-\tilde d^2}^{\mathrm{mv}}$ requires special care.~\\

\snewpage
\subsection{The differential calculus of $\AC$}~\\

Extending analytically from $(-1,1)$, one finds that
\[\AC'(z)=\frac{z\AC(z)-1}{1-z^2}\]
as analytic functions on $z\in\mathbb C\setminus (-\infty,-1]$.

We define
\begin{equation}
\AN(z)=\AC(z)^2(1-z^2),
\plabel{eq:ANdef}
\end{equation}
an analytic function on $z\in\mathbb C\setminus (-\infty,-1]$.
We remark that it is easy to check that
\[\AN'(z)=-2\AC(z).\]
(Thus $-\frac12\AN(z)$ is a primitive function of $\AC(z)$.)

From the definition \eqref{eq:ANdef}, it is easy to see that for $x\in(-1,+\infty)$,
\begin{equation}
\AN(x)=\begin{cases}
(\arccos x)^2&\text{if }-1< x<1\\[3mm]
0&\text{if }x=1\\[1mm]
-(\arcosh x)^2\qquad&\text{if } 1<x.\\
\end{cases}
\plabel{eq:ANreal}
\end{equation}
This latter one is an inverse function to $\Cos|_{(-\infty,\pi^2)}$.
Hence, by analytic continuation,
\begin{equation}
z=\Cos\circ \AN(z)
\plabel{eq:xANcomp}
\end{equation}
holds for $z\in\mathbb C\setminus (-\infty,-1]$.
Thus, we safely consider $\AN$ as \textit{the} principal branch of $\Cos^{-1}$, understood on $\mathbb C\setminus (-\infty,-1]$.
The identity $\AC(x) \cdot \Sin(\AN(x))=1$ also extends from $x\in(-1,+\infty)$; thus we find
\begin{equation}
\AC(z)=\frac{1}{\Sin\circ \AN(z)}
\plabel{eq:ACANcomp}
\end{equation}
for $z\in\mathbb C\setminus (-\infty,-1]$.

It is useful to consider the set
\[\beth=\left\{w\in\mathbb C\,:\,(\Rea w)<\pi^2-\frac{(\Ima w)^2}{4\pi^2}\right\}.\]
Note that $\beth$ is an open set containing $(-\infty,\pi^2)$.
If $w\in \beth$, then $\sqrt w^{\mathrm{mv}}$ (two-valued) is in the set
$\left\{w\in\mathbb C\,:\,|\Rea w|<\pi\right\}$.
Consequently, $\Cos(w)$ is in $\mathbb C\setminus[-1,\infty)$, thus in particular, in the range of $\AC$.
Combining \eqref{eq:xANcomp} and \eqref{eq:ACANcomp} and analytic continuation, we find
\begin{equation}
\AC(\Cos(w))=\frac1{\Sin(w)}
\plabel{eq:ACCoscomp}
\end{equation}
for
$w\in\beth$; extending from $w\in(-\infty,\pi^2)$; yielding an analytic function on $\beth$.

(Alternatively, the composition properties \eqref{eq:xANcomp}--\eqref{eq:ACANcomp}  could have been derived from the fact that
$\exp(\log A)=A$ whenever   $A$ is $\log$-able;
and the composition property \eqref{eq:ACCoscomp}  could have been derived from the fact that
$\log(\exp A)=A$ whenever the spectrum of $A$ is a subset of $\{z\,:\,|\Ima z|<\pi\}$.)
~\\

\snewpage
\subsection{$\AC$ near $-1$}~\\

For $x\in(-1,1)$, one has $\AC(x)=\frac{\arccos x}{\sqrt{1-x^2}}$.
Rewriting the standard identity $\arccos x=\pi-\arccos (-x)$, one obtains
\[\AC(x)=\frac{\pi}{\sqrt{1-x^2}}-\AC(-x).\]
By analytic continuation, we find that
\[\AC(z)=\frac{\pi}{\sqrt{1-z}\sqrt{1+z}}-\AC(-z)\]
is valid for $x\in\mathbb C\setminus ( (-\infty,-1] \cup [1,\infty  ) )$.
(As a reminder, $\sqrt\cdot$ is understood as the standard branch on $\mathbb C\setminus (-\infty,0)$.)
Making the change of variable $x=y-1$, one obtains
\[\AC(y-1)=\frac{\pi}{\sqrt{y}\sqrt{2-y}}-\AC(1-y),\]
valid for $y\in\mathbb C\setminus ( (-\infty,0] \cup [2,\infty  ) )$.
Applying power series expansions, it yields
\begin{equation}
\AC(y-1)=\frac{\pi\sqrt2}2\,\frac1{\sqrt y}-1+\frac{\pi\sqrt2}8\sqrt{y}-\frac13\,{\sqrt y}^2+\ldots\,.
\plabel{eq:ACbranch}
\end{equation}
This is properly convergent for $y\in\intD(0,2)\setminus (-2,0]$.
Also, note that, apart from the first term, we have a power series in $\sqrt y$ with convergent for `$\sqrt y$'$\in \intD(0,\sqrt{2})$.
Thus, if we relax the meaning of the square roots to a $2$-valued function, then
\eqref{eq:ACbranch} yields a ``branched'' or fractional Laurent series expansion for $y\in \intD(0,2)$.
In particular, the ``first analytic continuation'' of $\AC$ near $-1$ is $2$-branching.
This behaviour is typical:~\\

\snewpage
\subsection{The general analytic continuation of $\AC$}~\\

More generally, the meaning of square roots relaxed,  $\sqrt{1-z^2}$ yields a $2$-valued function.
Then, the formula
\[(\arccos z)'=-\frac1{\sqrt{1-z^2}}\text{ \qquad(with ambiguity, locally)}\]
shows that we can surely continue $\arccos$, thus $\AC$ locally as long as $z\neq\pm1$.
(From this viewpoint, it is a luck that a nice branch of $\AC$ extends through $1$.)
In particular, for essential behaviour, it is sufficient to consider the branching
in a neighborhood of $[-1,1]$.
In fact, due to the composition properties, one can easily reconstruct the Riemann surface of $\AC$
by considering the map
\[w\in\mathbb C\mapsto \left(\Cos(w),\frac1{\Sin(w)} \right)\in\mathbb C^2 .\]
(Indeed, apart $z=\pm1$, $w=\AN(z)$ can be extended in the manner of \eqref{eq:ANdef},
which, by the extension of \eqref{eq:xANcomp} and \eqref{eq:ACANcomp}  acts as reparametrization.)
The possible singularities correspond to  $\Cos(w)=\pm1$, i. e. $w=k^2\pi^2$ with $k\in\mathbb N$.
In the case $k=0$, the singularity is removable, as we have seen at the complex definition of $\AC$.
Otherwise, for $k\geq1$, $\Cos(w)$ takes a double, while $1/\Sin(w)$ has a simple pole.
Therefore the singularity is a double branching around $\infty$.

Consequently, in order to recover the continuation between $1$ and $-1$,
 it is sufficient to go back and forth between $1$ and $-1$,
 it does not matter that in what direction we encircle  $1$ and $-1$.
The $0$th branch line is $\frac{\arccos x}{\sqrt{1-x^2}}=\AC(x)$ (that is $\AC$ proper).
After encircling $-1$, the $1$st branch line is  $\frac{2\pi-\arccos x}{-\sqrt{1-x^2}}=\AC(x)-\frac{2\pi}{\sqrt{1-x^2}}$.
Then, after encirling $1$,  the $2$nd branch line is  $\frac{2\pi+\arccos x}{\sqrt{1-x^2}}=\AC(x)+\frac{2\pi}{\sqrt{1-x^2}}$.
After that,  after encirling $-1$,  the $3$rd branch line is  $\frac{4\pi-\arccos x}{-\sqrt{1-x^2}}=\AC(x)-\frac{4\pi}{\sqrt{1-x^2}}$.
Etc.
In terms of singularities: The $0$th central sheet of the covering
(meaning over $\mathbb C\setminus (  (-\infty,-1] \cup  [1,+\infty)$) yields
$\AC(z)$ for $z\in\intD(1,2)$, and
$\frac\pi{\sqrt{1-z}\sqrt{1+z}}-\AC(-z)$  for  $z\in\intD(-1,2)\setminus [-1,-3)$.
The $1$st sheet yields
$-\frac\pi{\sqrt{1-z}\sqrt{1+z}}-\AC(-z)$  for  $z\in\intD(-1,2)\setminus [-1,-3)$ and
$\AC(z)-\frac{2\pi}{\sqrt{1-z}\sqrt{1+z}}$ for $z\in\intD(1,2)\setminus[1,3)$.
The $2$nd sheet yields
$\AC(z)+\frac{2\pi}{\sqrt{1-z}\sqrt{1+z}}$ for $z\in\intD(1,2)\setminus[1,3)$, and
$\frac{3\pi}{\sqrt{1-z}\sqrt{1+z}}-\AC(-z)$  for  $z\in\intD(-1,2)\setminus [-1,-3)$.
The $3$rd sheet yields
$-\frac{3\pi}{\sqrt{1-z}\sqrt{1+z}}-\AC(-z)$  for  $z\in\intD(-1,2)\setminus [-1,-3)$ and
$\AC(z)-\frac{4\pi}{\sqrt{1-z}\sqrt{1+z}}$ for $z\in\intD(1,2)\setminus[1,3)$.
Etc.
In particular, we find that the singularity corresponding to $w=k^2\pi^2$ is
\[\AC^{\ext}(z)=\frac{\pi k}{\sqrt{1-z^2}^{\mathrm{mv}}}+\underbrace{(-1)^k\AC((-1)^kz)}_\text{holomorphic}\qquad\text{ for }z\sim (-1)^k\]
(in multivalued sense).

We can see that the analytic continuation of $\AC$ and of $\AN$ are essentially equivalent.
Indeed, by taking (extensions according to) \eqref{eq:ANdef} and \eqref{eq:ACANcomp}, one obtains one from the other.
In this manner, regarding $\AN$,  the singularity corresponding to $w=k^2\pi^2$ is
\[\AN^{\ext}(z)=\pi^2k^2+2
{\pi k}{\sqrt{1-z^2}^{\mathrm{mv}}}{(-1)^k\AC((-1)^kz)}+\AN((-1)^kz)\qquad\text{ for }z\sim (-1)^k\]
(in multivalued sense), cf. $\AC(1)=1$, $\AN(1)=0$.
This is a slightly milder singularity compared to the $\AC$, and
self-indexing the sense that the value at the singularity informs about the location.
(This is not surprising as $\AN$ is the inverse of $\Cos$.)~\\

\subsection{The meromorphic continuation of $\AC$ along analytic curves}~\\

Assume we want to extend $\AC$ along an analytic curve $z=f(t)$  $(t\in[0,\alpha)$, $f(0)=1$.
This can be done in a neighborbood of $t=0$, and can be continued up until $z$ hits $1$ or $-1$.
Say, this happens at $\tilde t$, $\tilde z=f(\tilde t)$.
If the analytic extension is $g(t)$ up to that point, then $w=g(t)^2(1-f(t)^2)$ keeps track
that if we hit a proper singularity or not, by taking the limit $\tilde w$.
The only case when the singularity is not proper this happens is $\tilde w=0$; then the analytic continuation is unproblematic.
If the singularity is proper, then meromorphic continuation is possible if and only if $\tilde z$ is point of even multiplicity at $t$.
In that case, $\tilde w=\pi^2k^2$, $k\in \mathbb N^+$, $\tilde z=(-1)^k$.
Assume that locally $f(t)=\tilde z+(t-\tilde t)^{2s} p(t)$, $p(\tilde t)\neq 0$.
Locally, for the extension
\[g(t)= \frac{\pi k}{(t-\tilde t)^s q(t)  }+(-1)^k\AC( (-1)^kf(t) ),\]
where $q(\tilde t)\neq0$. Then
\[w=g(t)^2(1-f(t)^2)=\left(\frac{\pi k}{(t-\tilde t)^sq(t)  }+(-1)^k\AC( (-1)^kf(t) )\right)^2
\left( (t-\tilde t)^sq(t)  \right)^2 \]
\[=\pi^2k^2+2{\pi k}{(t-\tilde t)^s(-1)^kq(t)  }+\AC( (-1)^kf(t) )^2\left( (t-\tilde t)^sq(t)  \right)^2 . \]

In particular, we see that as we continue meromorphically, $w$ extends analytically.
After that we carry on the continuation as usual.
In short, meromorphic continuation is possible as long as proper singularies are taken at even multiplicities.
Notice, that for $(-1)^s=1$ the continuation is ``bounced back to the sheet'',
while  for $(-1)^s=-1$ the continuation ``crosses to another sheet''.

Note that if $g(t)$ is the meromorphic continuation of $\AC(f(t))$, then the function
$h(t)=g(t)^2(1-f(t)^2)$ is the analytic continuation of $\AN(f(t))$.
Conversely, if $h(t)$ is the analytic continuation of $\AN(f(t))$, then
taking (extension by) \eqref{eq:ACANcomp} we obtain $g(t)$ as the meromorphic continuation of $\AC(f(t))$.

\begin{example}
Assume that we want to continue $g(t)=\AC(f(t))$ for $f(t)=\cos t$, $t\in[0,+\infty)$.
The original expression is well-defined only in a neighborhood of $t=0$.
The first possible singularity occurs for $\tilde t=\pi$, $\tilde z=-1$.
Here the multiplicity of $\tilde z$ is $2s=2$.
Here meromorphic continuation is possible.
(Indeed, it is given by $g(t)=\frac{t}{\sin t}$ for $t\sim \pi$.)
At this point the continuation is crosses sheets, $w=\pi^2-$ turns into $w=\pi^2+$.
The next  possible singularity occurs for $\tilde t=2\pi$, $\tilde z=1$.
Here the multiplicity of $\tilde z$ is $2s=2$ again.
Here meromorphic continuation is possible; again we cross sheets etc.
Indeed, we can continue for arbitrary $t$, in accordance to the meromorphicity of
$g(t)=\frac{t}{\sin t}$.
Note that  in the process, $w=g(t)^2(1-f(t)^2)=t^2$;
which is the rather simple analytic continuation $\AN(f(t))$.
\end{example}

The previous example was particularly simple, because the  meromorphic and analytic continuations
can be given in explicit form. Note, however, that in the case of a purely real analytic continuation (as above)
the value of $w$ can be  tracked just by considering the multiplicities at the critical points.
\snewpage
\subsection{Analytic continuation of $\log$ of the $2\times2$ matrices along analytic curves}
\begin{disc}\plabel{disc:cont}
Let us assume that we analytically continue   $\log A(t)$
 along the matrix valued analytic curve $A(t)$  $(t\in[0,\alpha)$, $A(0)=\Id_2$.
In our examples the analytic extension of $\sqrt{\dett A(t)}$ will be particularly simple (mostly trivial).
Thus we assume that an analytic extension $d(t)$ of $\sqrt{\dett A(t)}$ is given.
We can take reduced matrices $A^0(t)=\frac{A(t)}{d(t)}$, which are of determinant $1$.
Then the extension problem of log is related to the meromorphic extension of
$\AC\left(\dfrac{\tr A^0(t)}2\right)\equiv \AC\left(\dfrac{\tr A}{2\sqrt{\dett A}}\right)$.
It may happen that $\dfrac{\tr A^0(t)}2\equiv 1$.
Then $A^0(t)$ is unipotent, and the analytical computation of $\log A(t)$ is possible.
Otherwise $\dfrac{\tr A^0(t)}2$ is not constant.

According to the previous discussion, for $2\times2$ complex matrices, starting from $\Id_2$,
we can surely continue $\log$ until $z=f(t)=\dfrac{\tr A^0(t)}2$ hits $+1$ or $-1$.
Applying the notation of the previous subsection, this
happens at $\tilde t$, corresponding $\tilde w=k^2\pi^2$, and $\tilde z=(-1)^k$.

If $k=0$, then the analytic extension is unproblematic.
In fact,

(a) $A^0(\tilde t)$ is unipotent and $\AC(f(t))$ extends continuously.

(The case $\dfrac{\tr A^0(t)}2\equiv 1$ can be also be considered  as such a $k=0$ case.)

For $k\in\mathbb N^+$, however, the extension may exist under very special circumstances:
It must happen that

(b) $\tr A^0(t)$ has multiplicity $2s$ at $t=\tilde t$, and each coordinate of $A^0(t)-\frac{\tr A^0(t)}2\Id_2$
has multiplicity at least $s$ at $t=\tilde t$.
In particular, in this case $A^0(\tilde t)=\tilde z\Id_2$.

In case (b), we may say that $A^0(t)$ has sharp multiplicity $s$ at $t=\tilde t$.
Notice, however, that the content of point (b) is not only necessary but sufficient:
As long as the value $\tilde z\Id_2$ is taken with a sharp multiplicity, the continuation is possible.

If the continuation is successful, then one can continue until
the extension of $f(t)=\dfrac{\tr A^0(t)}2$ hits $+1$ or $-1$ again; etc.
Notice, however, that if the continuation is unsuccessful at $t=\tilde t$, the on some level,
a square root type singularity must occur somehow.

A small observation regarding the critical points:
If we have only $\log$-able matrices $A(t)$ or $A^0(r)$ for $t\in[0,\tilde t)$, then only the case $\tilde z=-1$ is critical (with $k=1$),
because $(-\infty,-1]$ will not have been crossed.
(This can be practical, as the curves themselves may be images of radial segments from a disc $\Dbar(0,R)$ with a $\log$-able zone provided.)

A very special but important case is when when the matrices $A(t)$ are real.
In this case, the curve starts from  $A(0)=\Id_2$,  $A^0(0)=\Id_2$.
Until $\frac{\tr A^0(t)}2$ hits $-1$, the matrix curve is in the $\log$-able zone; thus the first proper singularity must happen only with $\tilde z=-1$ reached,
which, on the other hand, quite typically leads to non-continuability.
\end{disc}
\snewpage
\subsection{Analytic continuation of $\log$ of the $2\times2$ matrices, Magnus expansion}
\begin{disc}\plabel{disc:cont2}
A situation where the previous subsection applies is the following.
Assume that $\phi$ is an ordered measure of $2\times2$ matrices, $\int\|\phi\|_2<+\infty$.

Then $A(u)=\Rexp(u\cdot \phi)$ is an entire function in $u$.
Furthermore, $\log A(u)=\log \Rexp(u\cdot \phi)$ is well-defined and analytic for small $u$.
We want to extend to $B(u)$ to $u\in\intD(0,R)$.
For that reason it is sufficient  to consider only radial extensions, that is
extensions of $\log A(\omega t)$ for $t\in[0,R)$ where $\omega$ is a complex unit.

In this case $\sqrt{\dett A(\omega t)}$ extends as $\exp\left(\frac\omega2\int\tr\phi \right)$.
In fact, as a method of reduction, $\phi$ can replaced by $\phi^0=\phi-\frac{\tr \phi}2\Id_2$;
as $A^0(u)= \Rexp(u\cdot \phi^0)$.

Now it sufficient to consider the radial extensions only through the points $p$
where $\frac{\tr A^0(p)}{2}=1$ or $-1$.
Examining these, the qualitative picture is sufficient again: we have to separate the cases
$k=0$ and the cases $k\in\mathbb N^+$ and extendibility is decided qualitatively by multiplicities as in the previous subsection.
If we want to prove non-extendibility, a candidate with $\frac{\tr A^0(p)}{2}=-1$ is often successful.

In order to  express situation in compact terms, in what follows,
we call every point $p$ with $\frac{\tr A^0(p)}{2}=1$ or $-1$ critical.
We will call the value $k$ above as the $\mathrm A$-index of $p$.
In summary,  if the $\mathrm A$-index is $0$, then one can continue;
if the $\mathrm A$-index is $k>0$, then one can continue iff $A$ takes the value $(-1)^k\Id_2$ with a sharp multiplicity at $p$.

We can compare this picture to the qualitative one we have already used along the real axis for real matrices.
So assume that $\phi$ is a measure of real $2\times2$ matrices, and
we want to continue $\log A^0$ along the real axis.
As we have seen, $\frac{\tr A^0(p)}{2}=-1$ immediately leads to divergence, unless $A^0(p)=- \Id_2$.
(Otherwise the extension would also be real, which is impossible in the lack of any possible real logarithm, cf. Lemma \ref{lem:rc}.)
This, of course also follows from the previous discussion.
Another fact is that $\frac{\tr A^0(p)}{2}<-1$ also implies divergence.
(Again, in the lack of a possible real logarithm.)
In that case let $\hat p\in(0,p)_e$ be the point with the smallest absolute value where $\frac{\tr A^0(\check p)}{2}$ crosses under $-1$.
Then  $\frac{\tr A^0(\hat p)}{2}=-1$ and $\log A(u)$ cannot be continued through $u=\hat p$.
This, however, could have been argued from the fact that continuability requires that
the multiplicity of $\frac{\tr A^0(p)}{2}$ should be even.
Hence, the real case fits into our complex picture nicely.

In any case, in our examples later, extension along the real axis (or, in a related way, along the imaginary axis)
will dominate the discussion.
Nevertheless, having a general understanding as above is useful.
\end{disc}

\subsection{Alternative notations}
\begin{remark}
It is somewhat of a struggle to choose between our current notational conventions applying
$D_A$, $\Cos(x)$, $\Sin(x)$, $\AN(x)$ or the other one using
$\mathrm{Discr}_A=-D_A$, $\Cosh(x)=\Cos(-x)$,  $\Sinh(x)=\Sin(-x)$, $\mathrm{ANh}(x)=-\AN(x)$.
The latter conventions are very natural, but we will stick to the former ones.
The function $\AC$ is, however, defined very naturally.
\end{remark}
~
\snewpage
\subsection{Some observations connected to $\arccos$ and $\cos$}
~\\

Recall that $\arcsin$ and $\arccos$ and can be extended to  $\mathbb C\setminus\left( (-\infty,-1]\cup [1,\infty)\right)$
as
$\arcsin z=\mathrm i\log\left(-\mathrm iz+\sqrt{1-z^2}\right)$
and
$\arccos z=\frac\pi2-\arcsin z$.

The following lemma is useful taking limits in $\arccos$ (or $\AC$) at $-1$.
\begin{lemma}
\plabel{lem:arccos1}
Let us consider the function $\arccos$ on $\mathbb C\setminus\left( (-\infty,-1]\cup [1,\infty)\right)$
but extended by the classical values $\arccos(-1)=\pi$ and $\arccos1=0$, but even
possibly with the $2$-valued extensions
$\arccos \varrho=\pm\mathrm i\arcosh\varrho$ for $\varrho>1$ and
$\arccos \varrho=\pi\mp\mathrm i\arcosh(-\varrho)$ for $\varrho<-1$,
or even in the local multivalued branching.
Then $\arccos$ is continuous at $1$ and at $-1$ (despite cuts and multivaluedness).
\begin{proof}
Continuity at $1$ follows from either from of $\arccos z=\sqrt{\AN(z)}^{\mathrm{mv}}$ or
$\arccos z=\AC(z)\sqrt{1+z}\sqrt{1-z}^{\mathrm{mv}}$ at $z\sim1$ (in the local branching).
Continuity at $-1$ follows from the symmetry property $\arccos z=\pi-\arccos(-z)$.
Or, more laboriously, the statement also follows directly from the properties of $\log$.
\end{proof}
\end{lemma}
\begin{lemma}
\plabel{lem:arccos2}
For $\delta\in\mathbb C\setminus\left( (-\infty,-1]\cup [1,\infty)\right)$,
\[\left|\frac{\AC(\delta)}{\AC(-\delta)}\right|\equiv \left|\frac{\arccos(\delta)}{\arccos(-\delta)}\right|\lesseqqgtr1
\qquad\text{is according to}\qquad  0\lesseqqgtr \Rea \delta.\]
\begin{proof}
Using the logarithmic definition, one can show that, on the indicated domain,
$\Rea\arcsin \delta\in\left(-\frac\pi2,\frac\pi2\right)$, and  $\sgn\Rea\arcsin \delta=\sgn\Rea \delta$.
In that light the behaviour of
$\left|\frac{\arccos(\delta)}{\arccos(-\delta)}\right|
=\left|\frac{\frac\pi2-\arcsin \delta}{ \frac\pi2+\arcsin \delta}\right|$ is transparent.
\end{proof}
\end{lemma}
\begin{lemma}
\plabel{lem:ACir}
For $\varrho\in\mathbb R$,
\[\AC(\mathrm i\varrho)=\frac{\arccos(\mathrm i\varrho)}{\sqrt{1-(\mathrm i\varrho)^2}}
=\frac{\frac\pi2-\arcsin(\mathrm i\varrho)}{\sqrt{1+\varrho ^2}}=\frac{\frac\pi2-\mathrm i\arsinh\varrho}{\sqrt{1+\varrho^2}}.
\eqed\]
\end{lemma}
\begin{lemma}
\plabel{lem:cos1}
(a) $\cos s$ is real only $s\in\mathbb R$ and for $s\in\mathrm i\mathbb R+\pi\mathbb Z$; and

(b) $\cos s$ is purely imaginary only for $s\in\mathrm i\mathbb R+\pi\left(\mathbb Z+\frac12\right) $.
\begin{proof}
This follows from the formula $\cos(x+\mathrm iy)=\cos x\cosh y-\mathrm i\sin x\sinh y$.
\end{proof}
\end{lemma}
\begin{lemma}
\plabel{lem:cos2}
Let $r>0$. We claim:

(a) If $r\leq\frac\pi2$, then the smallest real value   $(\cos s)^2$ for $s\in\Dbar(0,r)$ is
$(\cos r)^2$, which is taken at the two places $s=\pm r$.

(b) If $r>\frac\pi2$, then the smallest real value   $(\cos s)^2$ for $s\in\Dbar(0,r)$ is
$-(\sinh u)^2 $, which is taken at the four places $s=\pm \frac\pi2\pm \mathrm iu $, where $u=\sqrt{r^2-\left(\frac\pi2\right)^2}$.
\begin{proof}
The previous lemma restricts the possibilities when $(\cos s)^2$ is real, which can be inspected directly.
\end{proof}
\end{lemma}

\subsection{Norms}
\begin{lemma} \plabel{lem:normcomputeC}
Let $A=\begin{bmatrix}a&b\\c&d\end{bmatrix}
=\tilde a\Id_2+\tilde b\tilde I+\tilde c\tilde J+\tilde d\tilde K$ be a real or complex matrix. Then
\begin{align*}
\left\|A \right\|_2
&=\sqrt{\frac{\tr(A^*A)}2+\sqrt{-D_{A^*A}}}\\
&=\sqrt{\frac{\tr(A^*A)}2+\sqrt{\frac{(\tr(A^*A))^2}4-|\det A|^2}}\\
&=\frac{\sqrt{\tr(A^*A)+2|\det A|}+\sqrt{\tr(A^*A)-2|\det A|} }2\\
&=\frac{\sqrt{|a|^2+|b|^2+|c|^2+|d|^2+2|ad-bc| }+\sqrt{|a|^2+|b|^2+|c|^2+|d|^2-2|ad-bc|}}2\\
&=\sqrt{\frac{|\tilde a|^2+|\tilde b|^2+|\tilde c|^2+|\tilde d|^2+ |\tilde a^2+\tilde b^2-\tilde c^2-\tilde d^2| }{2}}\\
&\qquad +
\sqrt{\frac{|\tilde a|^2+|\tilde b|^2+|\tilde c|^2+|\tilde d|^2- |\tilde a^2+\tilde b^2-\tilde c^2-\tilde d^2| }{2}}
;
\end{align*}
and
\begin{align*}
\|A\|_2^-&=\left\|A^{-1} \right\|_2^{-1}\\
&=\sqrt{\frac{\tr(A^*A)}2-\sqrt{-D_{A^*A}}}\\
&=\sqrt{\frac{\tr(A^*A)}2-\sqrt{\frac{(\tr(A^*A))^2}4-|\det A|^2}}\\
&=\frac{\sqrt{\tr(A^*A)+2|\det A|}-\sqrt{\tr(A^*A)-2|\det A|} }2\\
&=\frac{\sqrt{|a|^2+|b|^2+|c|^2+|d|^2+2|ad-bc| }-\sqrt{|a|^2+|b|^2+|c|^2+|d|^2-2|ad-bc|}}2\\
&=\sqrt{\frac{|\tilde a|^2+|\tilde b|^2+|\tilde c|^2+|\tilde d|^2+ |\tilde a^2+\tilde b^2-\tilde c^2-\tilde d^2| }{2}}\\
&\qquad -
\sqrt{\frac{|\tilde a|^2+|\tilde b|^2+|\tilde c|^2+|\tilde d|^2- |\tilde a^2+\tilde b^2-\tilde c^2-\tilde d^2| }{2}}
.
\end{align*}
In particular,
\begin{equation}\|A\|_2\cdot\|A\|_2^-=|\det A|.\plabel{eq:compmult}\end{equation}

In the case of real matrices, the results are the same for the Hilbert spaces  $\mathbb R^2$ and $\mathbb C^2$.
\begin{proof}
This follows from $\|A\|_2=\sqrt{\max \spec(A^*A)}$ and $\|A\|_2^-=\sqrt{\min \spec(A^*A)}$, computing the eigenvalues.
\end{proof}
\end{lemma}

Motivated by \eqref{eq:compmult}, for  $2\times2$ matrices, one can define the signed co-norm as
\begin{equation}
\lfloor A\rfloor_2=
\begin{cases}
0&\text{if }A=0\\\\
\dfrac{\det A}{\|A\|_2}&\text{if }A\neq 0.\\
\end{cases}
\plabel{eq:signedconorm}
\end{equation}

Then, from \eqref{eq:compmult},
\begin{equation}
\left| \left\lfloor A\right\rfloor_2 \right|=\|A\|_2^-;
\end{equation}
and
\begin{equation}
\|A\|_2\cdot \left\lfloor A\right\rfloor_2 =\det A.\plabel{eq:realmult}
\end{equation}
Nevertheless, we will consider the signed co-norm only for real matrices.
\begin{lemma} \plabel{lem:normcompute}
Let $A=\begin{bmatrix}a&b\\c&d\end{bmatrix}
=\tilde a\Id_2+\tilde b\tilde I+\tilde c\tilde J+\tilde d\tilde K$ be a real matrix. Then
\begin{align}
\left\|A \right\|_2
&=\frac{\sqrt{\tr(A^*A)+2\det A}+\sqrt{\tr(A^*A)-2\det A} }2
\plabel{eq:2norm}\\
&=\frac{\sqrt{(a+d)^2+(c-b)^2}+\sqrt{(a-d)^2+(b+c)^2}}2\notag\\
&=\sqrt{\tilde a^2+\tilde b^2}+\sqrt{\tilde c^2+\tilde d^2}.
\plabel{eq:norm22skew}
\end{align}
On the other hand,
\begin{align}
\|A\|_2^-=\left\|A^{-1} \right\|_2^{-1}
&=\left|\frac{\sqrt{\tr(A^*A)+2\det A}-\sqrt{\tr(A^*A)-2\det A} }2\right|
\plabel{eq:2conorm1}\\
&=\left|\frac{\sqrt{(a+d)^2+(c-b)^2}-\sqrt{(a-d)^2+(b+c)^2}}2\right|\notag\\
&=\left|\sqrt{\tilde a^2+\tilde b^2}-\sqrt{\tilde c^2+\tilde d^2}\right|
.
\notag
\end{align}
It is true that
\begin{align}
\sgn \det A
&=\sgn\frac{\sqrt{\tr(A^*A)+2\det A}-\sqrt{\tr(A^*A)-2\det A} }2
\plabel{eq:2conorm2}\\
&= \sgn \frac{\sqrt{(a+d)^2+(c-b)^2}-\sqrt{(a-d)^2+(b+c)^2}}2\notag\\
&=\sgn\left(\sqrt{\tilde a^2+\tilde b^2}-\sqrt{\tilde c^2+\tilde d^2}\right).
\notag
\end{align}
Furthermore,
\begin{align}
\left\lfloor A\right\rfloor_2=\sgn(\det A )\left\|A \right\|_2^-
&=\frac{\sqrt{\tr(A^*A)+2\det A}-\sqrt{\tr(A^*A)-2\det A} }2
\plabel{eq:2conorm}\\
&=\frac{\sqrt{(a+d)^2+(c-b)^2}-\sqrt{(a-d)^2+(b+c)^2}}2\notag\\
&=\sqrt{\tilde a^2+\tilde b^2}-\sqrt{\tilde c^2+\tilde d^2}
.
\notag
\end{align}

\begin{proof}
\eqref{eq:2norm} and \eqref{eq:2conorm1} are immediate from the complex case. \eqref{eq:2conorm2} is trivial.
\eqref{eq:2conorm} follows from \eqref{eq:2conorm2} and the definition \eqref{eq:signedconorm}.
\end{proof}
\end{lemma}

\snewpage
\section{Counterexamples (continued)}
\plabel{sec:BCHsharp}

\subsection{Cases of convergence analysis for the Magnus expansion}
\plabel{ss:CAMagnus}~\\

\begin{example}\plabel{ex:Moan2}
(Moan's example.)
 On the interval $[0,\pi]$, we consider the measure $\hat\Phi$, such that
\[\hat\Phi(\theta)=
\frac12\begin{bmatrix}
-\sin2\theta& -1+\cos2\theta\\1+\cos2\theta&\sin2\theta
\end{bmatrix}
\,\mathrm d\theta|_{[0,\pi]}.\]
Then,
\[\int\|\hat\Phi\|_2=\pi.\]

For $t\in\intD(0,1)$, therefore,
$\mu_{\mathrm L}(t\cdot\hat \Phi)=\log\Lexp (t\cdot\hat\Phi)$ holds.
Explicitly,
\begin{align}
\Lexp (t\cdot\hat\Phi)&=\Lexp\left(
\frac t2
\bem&-1\\1&\eem+
\frac t2\begin{bmatrix}
-\sin2\theta& \cos2\theta\\\cos2\theta&\sin2\theta
\end{bmatrix}\mathrm d\theta|_{[0,\pi]}\right)
\notag\\
&=F\left(\frac\pi2 t,\frac\pi2 t,\pi \right)
\notag\\
&=-\begin{bmatrix}\cos( \pi\sqrt{1-t})& \frac{\sin( \pi\sqrt{1-t})}{\sqrt{1-t}}\\
\frac{\sin( \pi\sqrt{1-t})}{\sqrt{1-t}}(t-1)&\cos( \pi\sqrt{1-t})\end{bmatrix}.
\notag
\end{align}
Thus,
\begin{align}
\mu_{\mathrm L}(t\cdot\hat\Phi)&=\log\Lexp (t\cdot\hat\Phi)
\notag\\
&=\frac{\AC(-\cos (\pi\sqrt{1-t}))\sin (\pi\sqrt{1-t})  }{\sqrt{1-t}}
\begin{bmatrix}& -1\\-t+1&\end{bmatrix}
\notag\\
&=\pi\left(\frac1{\sqrt{1-t}}-1\right)\begin{bmatrix}& -1\\-t+1&\end{bmatrix}.
\notag
\end{align}
Consequently, if $t\in[0,1)$, then
\[\|\mu_{\mathrm L}(t\cdot\hat\Phi)\|_2=\pi\left(\frac1{\sqrt{1-t}}-1\right).\]
This indicates the rate of divergence of the Magnus expansion as $t\nearrow1$.

Regarding the individual terms of the Magnus expansion, we see that
\[\mu_{n,\mathrm L}(\hat\Phi)=
\begin{cases}
0 & \text{if }n=1\\[2mm]
\pi\bem&-\frac12\\\frac12&\eem&\text{if  $n=2$}
\\[6mm]
\pi\bem&-(-1)^{n}\binom{-1/2}{n}\\(-1)^{n}\binom{1/2}{n}&\eem=\pi \mathrm G_n\bem&-1\\-\frac1{2n-1}&\eem&\text{if  $n\geq 3$,}
\end{cases}\]
where $\mathrm G_n$ is as in \eqref{eq:gfunc}.
Due to \eqref{eq:wallis}, we see explicitly how the Magnus expansion diverges.
In particular, $\|\mu_{n,\mathrm L}(\Phi)\|_2=\pi\mathrm G_{n}$ for $n\geq2$; thus
\begin{equation}
\|\mu_{n, \mathrm L}(\hat\Phi) \|_2=\sqrt{\frac{\pi}{n}}(1+o(1)),
\plabel{eq:critast}
\end{equation}
as $n\rightarrow\infty$.
\qedexer
\end{example}

\begin{example}\plabel{ex:critical}
(Magnus critical case.)
 On the interval $[0,\pi]$, we consider the measure $\Phi$, such that
\[\Phi(\theta)=
\begin{bmatrix}
-\sin2\theta& \cos2\theta\\\cos2\theta&\sin2\theta
\end{bmatrix}
\,\mathrm d\theta|_{[0,\pi]}.\]
Then,
\[\int\|\Phi\|_2=\pi.\]

For $t\in\intD(0,1)$, therefore,
$\mu_{\mathrm L}(t\cdot \Phi)=\log\Lexp (t\cdot\Phi)$ holds.
Explicitly,
\begin{align}
\Lexp (t\cdot\Phi)&=\Lexp\left( t\begin{bmatrix}
-\sin2\theta& \cos2\theta\\\cos2\theta&\sin2\theta
\end{bmatrix}\mathrm d\theta|_{[0,\pi]}\right)
\notag\\
&=W(\pi t,\pi)
\notag\\
&=-\begin{bmatrix}\cos( \pi\sqrt{1-t^2})& \frac{\sin( \pi\sqrt{1-t^2})}{\sqrt{1-t^2}}(t+1)\\
\frac{\sin( \pi\sqrt{1-t^2})}{\sqrt{1-t^2}}(t-1)&\cos( \pi\sqrt{1-t^2})\end{bmatrix}.
\notag
\end{align}%
Thus,
\begin{align}
\mu_{\mathrm L}(t\cdot \Phi)&=\log\Lexp (t\cdot\Phi)
\notag\\
&=\frac{\AC(-\cos (\pi\sqrt{1-t^2}))\sin (\pi\sqrt{1-t^2})  }{\sqrt{1-t^2}}
\begin{bmatrix}& -t-1\\-t+1&\end{bmatrix}
\notag\\
&=\pi\left(\frac1{\sqrt{1-t^2}}-1\right)\begin{bmatrix}& -t-1\\-t+1&\end{bmatrix}.
\notag
\end{align}
Consequently, if $t\in[0,1)$, then
\begin{align}
\|\mu_{\mathrm L}(t\cdot \Phi)\|_2
&=\pi\left(\frac1{\sqrt{1-t^2}}-1\right)(1+t)
\notag\\
&=\sqrt2\pi (1-t)^{-1/2}-2\pi-\frac{\sqrt2}4\pi(1-t)^{1/2}+O(1-t),
\notag
\end{align}
as $t\nearrow1$.
This indicates the rate of divergence of the Magnus expansion,
which is asymptotically $\sqrt2$ times the one in the previous example.

Regarding the individual terms of the Magnus expansion, we see that
\[\mu_{n,\mathrm L}(\Phi)=
\begin{cases}
0 & \text{if }n=1\\\mathrm G_{\lfloor n/2\rfloor}
\pi\tilde I&\text{if $n$ is even, $n\geq 2$}
\\
\mathrm G_{\lfloor n/2\rfloor}\pi(-\tilde K)&\text{if $n$ is odd, $n\geq 3$}
\end{cases}\]
where $\mathrm G_n$ is as in \eqref{eq:gfunc}.
Again, due to \eqref{eq:wallis},  we explicitly see that $\sum_{n=1}^\infty\mu_{n,\mathrm L}(\Phi)$ is divergent.
In particular, $\|\mu_{n,\mathrm L}(\Phi)\|_2=\pi\mathrm G_{\lfloor n/2\rfloor}$ for $n\geq2$; thus
\[\|\mu_{n,\mathrm L}(\Phi) \|_2=\sqrt{\frac{2\pi}{n}}(1+o(1)),\]
as $n\rightarrow\infty$.
\qedexer
\end{example}
~
\snewpage

\subsection{General remarks on the convergence of the Magnus expansion}~\\

As we know, the convergence radius of the Magnus expansion of $\phi$ depends
on the analytic continuability of the map $t\mapsto\log (\Rexp(t\cdot\phi))$ from $t\sim0$.

In the case of $2\times 2$ matrices, we have considered two different principles regarding this continuability:

(i) singularities can occur only if the $\Rexp(t\cdot\phi)$ have two equal eigenvalues
(which, in the traceless case can be only $1$ or $-1$);

(ii) singularities can occur only of after  the continuation of $\frac{\tr \Rexp(t\cdot\phi)}{2\sqrt{\det\Rexp(t\cdot\phi)}}$
hits $(-\infty,-1]$.

One might think that these observations could be combined into a single one, like that in the traceless case
singularities can occur at Jordan blocks $\bem-1&1\\&-1\eem$.
In this respect, the following example may be instructive:
\begin{example}
Consider $\Phi^2=\Phi\boldsymbol.\Phi$ (the concatenation of the Magnus critical case with itself).
Its Magnus expansion has convergence radius $1$.
In fact, $\log^{\ext}\Lexp(t\cdot\Phi^2)=2\log \Lexp(t\cdot\Phi)$,
it yields the very same convergence problem as the Magnus critical case.
The first singularities occur for $t=\pm1$, where $\Lexp (1\cdot\Phi^2)= \bem1&4\pi\\&1\eem$.
\qedexer
\end{example}
\begin{remark}
\plabel{rem:logradius2}
When we continue $\log^{\ext}\Lexp(t\cdot \Psi\boldsymbol.\Psi)$ for $t\in[0,1]$,
we find the first critical point at $t=\sqrt3/2$ with $\mathrm A$-index $1$, where analytic continuation is, of course, possible.
Then, we arrive to a critical point at $t=1$ with  $\mathrm A$-index $2$, where analytic continuation, is, of course, impossible.
\qedremark
\end{remark}

Moan, Niesen \cite{MN} conjectures criterion for the divergence of the Magnus expansion $\mu_{\mathrm L}(\phi)$
of $n\times n$ matrices in terms of $Y(t;\kappa)=\exp(\kappa\cdot\phi|_{(-\infty,t)})$;
 the reader is advised to it look up in \cite{MN} for its precise formulation.
But, it says that the divergence is caused by ``eigenvalue collisions'' under certain prescribed circumstances.
However, the precise formulation of such a statement
would require further care regarding what makes an eigenvalue collision  exactly ``lethal''.
\begin{example}
Consider
\[
\exp\begin{bmatrix}
&-\pi&-1&\\
\pi&&&-1\\
&&&-\pi\\
&&\pi&
\end{bmatrix}
=
\begin{bmatrix}
-1&&1&\\
&-1&&1\\
&&-1&\\
&&&-1
\end{bmatrix}.
\]

The exponential itself can be considered as a time-ordered exponential having a Magnus expansion with only one nontrivial term.
There is an ``eigenvalue collision'' for $t=1$ at the eigenvalue $\lambda=-1$ with nontrivial Jordan blocks,
 and by that time the eigenvalues encircling the origin.
Yet,  the Magnus expansion has convergence radius $+\infty$.
\qedexer
\end{example}

As Casas \cite{Ca} notes,
there are already some results in the literature regarding restrictive conditions for the general analytic continuation of the $\log$
of $n\times n$ complex matrix matrices (that is not necessarily in the context
of the Magnus expansion), see Yakubovich, Starzhinskii \cite{YS},
based on simple spectral principles like as in Lemma \ref{lem:rc}, but which can probably be developed further
in this more specific situation.

In general,  there is much yet to clarify about the development of the singularities in the analytical
continuation $\Rexp(t\phi)$, even in the matrix case.

\snewpage
\subsection{Unitary / quaternionic counterexamples to the Magnus expansion}
\begin{remark}\plabel{rem:physex}
The more detailed analysis of concrete systems in the unitary case was started by physicists.
As the Magnus expansion was applied in certain quantum mechanical investigations quite
early, this prompted some theoretical interest, see Pechukas, Light \cite{PL}, Wilcox \cite{Wi}
(leading to the work of  Mielnik, Pleba\'nski \cite{MP}).
Increasingly refined  analysis of some examples was taken up in
Fel'dman \cite{Fe}, Salzman \cite{Sa}, Maricq \cite{Mq}, Klarsfeld, Oteo \cite{KO}, etc.
(See also Blanes, Casas, Oteo, Ros \cite{BCOR}.)
Some aspects were scrutinized further in Moan, Niesen \cite{MN} and Casas \cite{Ca}.
(For physical motivations see Kuprov \cite{Kup}, although in applications the Magnus
expansion tends to used as theoretical background for related numerical ``microconvergence''.)
\qedremark
\end{remark}
\snewpage

Next, we show a variant of  Example
\ref{ex:critical}, such that the measure take values from $2\times 2$ skew-Hermitian matrices,
thus the time-ordered exponentials are $2\times2$ unitary matrices.
\begin{example}\plabel{ex:criticalquat}
(Magnus critical case, quaternionic version.)
 On the interval $[0,\pi]$, we consider the measure $\Phi^{\mathrm q}$, such that
\[\Phi^{\mathrm q}(\theta)=\mathrm i\cdot \Phi(\theta)=\mathrm i
\begin{bmatrix}
-\sin2\theta& \cos2\theta\\\cos2\theta&\sin2\theta
\end{bmatrix}
\,\mathrm d\theta|_{[0,\pi]}=
\exp(\theta I) K\exp(-\theta I)\,\mathrm d\theta|_{[0,\pi]}
.\]
Then,
\[\int\|\Phi^{\mathrm q}\|_2=\pi.\]

For $t\in\intD(0,1)$, therefore,
$\mu_{\mathrm L}(t\cdot \Phi^{\mathrm q})=\log\Lexp (t\cdot\Phi^{\mathrm q})$ holds.
Explicitly,
\begin{align}
\Lexp (t\cdot\Phi)&=\Lexp\left( t\mathrm i\begin{bmatrix}
-\sin2\theta& \cos2\theta\\\cos2\theta&\sin2\theta
\end{bmatrix}\mathrm d\theta|_{[0,\pi]}\right)
\notag\\
&=W(\mathrm i\pi t,\pi)
\notag\\
&=-\begin{bmatrix}\cos( \pi\sqrt{1+t^2})& \frac{\sin( \pi\sqrt{1+t^2})}{\sqrt{1+t^2}}(\mathrm it+1)\\
\frac{\sin( \pi\sqrt{1+t^2})}{\sqrt{1+t^2}}(\mathrm it-1)&\cos( \pi\sqrt{1+t^2})\end{bmatrix}.
\notag
\end{align}%
Thus,
\begin{align}
\mu_{\mathrm L}(t\cdot \Phi^{\mathrm q})&=\log\Lexp (t\cdot\Phi^{\mathrm q})
\notag\\
&=\frac{\AC(-\cos (\pi\sqrt{1+t^2}))\sin (\pi\sqrt{1+t^2})  }{\sqrt{1+t^2}}
\begin{bmatrix}& -\mathrm it-1\\-\mathrm it+1&\end{bmatrix}
\notag\\
&=\pi\left(\frac1{\sqrt{1+t^2}}-1\right)\begin{bmatrix}& -\mathrm it-1\\-\mathrm it+1&\end{bmatrix}.
\notag
\end{align}
This, of course,  continues radially for any $t\in\mathbb R$.

Regarding the individual terms of the Magnus expansion, we see that
\[\mu_{n,\mathrm L}(\Phi^{\mathrm q})=
\begin{cases}
0 & \text{if }n=1\\\mathrm G_{\lfloor n/2\rfloor}
\pi  I (-1)^{n/2}&\text{if $n$ is even, $n\geq 2$}
\\
\mathrm G_{\lfloor n/2\rfloor}\pi  K(-1)^{(n+1)/2}&\text{if $n$  is odd, $n\geq 3$}.
\end{cases}\]

The convergence radius is, obviously, $1$.
Due to the explicit formula, the Magnus expansion is Abel-summable.
In fact, due to its oscillatory nature it even converges but not absolutely.

In particular, multiplying $\Phi^{\mathrm q}$ with $1+\varepsilon$, we see that
convergence cannot be expected for cumulative norm $\pi(1+\varepsilon)$, not even in the unitary case.
\qedexer
\end{example}
(Here the  divergence was less pregnant than in the real case.)
\snewpage
\subsection{Cases of convergence analysis for the BCH expansion}
\plabel{ss:CABCH}~\\

The following proposition shows how Discussion \ref{disc:cont2} can be applied to establish
convergence radii in case of a relatively simple-looking example.

\begin{prop}
\plabel{pr:sharpunbal}
(a) In Example \ref{ex:vin}(a), the convergence radius of the Magnus (BCH) expansion is

(i) The unique $t_0$ such that $0<t_0<1$ and
\[-\frac1{t_0}\cot \frac\pi2{t_0}=\delta,\]
if $\delta<0$;

(ii) \[1,\] if $0<\delta\leq \coth\frac\pi2  $;

(iii) The unique $t_0$ such that $0<t_0<1$ and
\[\frac1{t_0}\coth \frac\pi2{t_0}=\delta,\]
if $\coth\frac\pi2<\delta $.

In case (ii), the Magnus (BCH) expansion is (completely) divergent.
The rate of divergence is indicated by

\[
\lim_{t\nearrow1}\mu_{\mathrm R}(t\cdot V^{[\delta]}_1\mathbf1\boldsymbol.V^{[\delta]}_2\mathbf1)\cdot\sqrt{1-t^2}
\equiv
\lim_{t\nearrow1}\log((\exp tV^{[\delta]}_1 )(\exp tV^{[\delta]}_2))\cdot\sqrt{1-t^2}
=2\sqrt{\delta\pi}\bem0&-1\\&0\eem.
\]

(b) In Example \ref{ex:vin}(b), the convergence radius of the Magnus (BCH) expansion is always ($\eta\neq 0$) smaller than $1$.

\begin{proof}
(a) Let $C_{12}^{V[\delta]}(t):=$
\begin{multline*}(\exp tV^{[\delta]}_1 )(\exp tV^{[\delta]}_2)=
\left[ \begin{matrix} \cos \pi t &-\sin \pi t \\\sin \pi t&\cos \pi t\end{matrix}\right]
\left[\begin{matrix} 1&2t\delta\\&1\end{matrix}\right]
=\bem \cos \pi t  &-\sin \pi t +2\delta t \cos \pi t\\\sin \pi t&\cos\pi t +2\delta t\sin \pi t\eem.
\end{multline*}

First, we remark that, considering the adiagonal terms in the matrix,

\begin{equation}\text{$C_{12}^{V[\delta]}(t)$ is a scalar matrix only for $t=0$.}\plabel{eq:vinscalar}\end{equation}

Taking logarithm, critical behaviour for
$\AC\left(\dfrac{\tr C_{12}^{V[\delta]}(t)}{2\sqrt{\dett C_{12}^{V[\delta]}(t)}^{\ext}}\right)\equiv
\AC\left(\dfrac{\tr C_{12}^{V[\delta]}(t)}{2}\right)$ occurs when $(t\neq0)$
\begin{equation}\cos\pi t+\delta t\sin\pi t=-1,\plabel{eq:vinprecrit}\end{equation}
or
\begin{equation}\cos\pi t+\delta t\sin\pi t=1.\plabel{eq:vinprecrit2}\end{equation}
For $t=\pm1$ equation \eqref{eq:vinprecrit} holds and $C_{12}^{V[\delta]}(t)$ is not a scalar matrix, thus the convergence radius cannot be greater than $1$.
(It is less, or we have explicit real singularities at $t=\pm1$.)
Thus we can restrict to the case $|t|<1$, where the critical behavior can be rephrased as
\begin{equation}-\frac1t\cot \frac\pi2t=\delta,\plabel{eq:vincrit}\end{equation}
or
\begin{equation}\frac1t\tan \frac\pi2t=\delta,\plabel{eq:vincrit2}\end{equation}
respectively.

Note that for $|t|<1$ the LHS of equation \eqref{eq:vincrit} is real if and only if $t$ is real or purely imaginary.
The same comment applies for  equation \eqref{eq:vincrit2}.
Thus it is sufficient to search for critical behaviour only along the real and purely imaginary axes.
If we apply conjugation by the matrix
\[G= \begin{bmatrix}1&\\&-\mathrm i\end{bmatrix},\]
then we find
\[ GV^{[\delta]}_1G^{-1}=\frac1{\mathrm i}\cdot \pi\begin{bmatrix}&1\\1&\end{bmatrix}
\qquad\text{and}\qquad
 GV^{[\delta]}_2G^{-1}=\frac1{\mathrm i}\cdot 2\delta\begin{bmatrix}0&-1\\&0\end{bmatrix}.\]
This shows that convergence along the imaginary axis has the same ``real'' qualitative features as convergence along the real axis.

Consequently, for a first non-continuable singularity it is not only necessary consider the solutions of
equation \eqref{eq:vinprecrit}/\eqref{eq:vincrit} (to go there or beyond),
but it is, due \eqref{eq:vinscalar}, is also sufficient.
(Thus, equations \eqref{eq:vinprecrit2}/\eqref{eq:vincrit2} are out of the play.)

In case (i), \eqref{eq:vincrit} has only real solutions with $t\in\intD(0,1)$, these are $t=\pm t_0$.
In case (ii),  \eqref{eq:vincrit} has  no solutions.
In case (iii), \eqref{eq:vincrit} has only purely imaginary solutions with $t\in\intD(0,1)$, these are $t=\pm t_0\mathrm i$.
This establishes the convergence radii.

In case (ii),
\begin{align*}
\mu_{\mathrm R}(t\cdot V^{[\delta]}_1\mathbf1\boldsymbol.V^{[\delta]}_2\mathbf1)
&=\log^{\ext}((\exp tV^{[\delta]}_1 )(\exp tV^{[\delta]}_2))\\
&=\AC\left(\cos\pi t+\delta t\sin\pi t\right)
\cdot\bem -\delta t\sin \pi t  &-\sin \pi t +2\delta t \cos \pi t\\\sin \pi t& \delta t\sin \pi t\eem.
\end{align*}
In fact, by Lemma \ref{lem:logreal}, $\log$ is proper for $t\in(-1,1)$.
(Remark: It would require more analysis, but this is also true for $t\in\intD(0,1)$.)
Taking the corresponding limit $t\nearrow1$ is a matter of elementary analysis.

(b) Here
\[\exp(t\tilde V^{[\eta]}_1)(\exp t\tilde V^{[\eta]}_2)=
\left[ \begin{matrix} \cos \pi t &-\sin \pi t \\\sin \pi t&\cos \pi t\end{matrix}\right]
\left[\begin{matrix} \mathrm e^{t\eta}&\\& \mathrm e^{-t\eta}\end{matrix}\right]
=
\bem \mathrm e^{t\eta}\cos \pi t &-\mathrm e^{-t\eta}\sin \pi t \\\mathrm e^{t\eta}\sin \pi t&\mathrm e^{-t\eta}\cos \pi t\eem.\]
Taking logarithm, the critical behaviour for $\AC$ is when $(t\neq0)$
\[\cos(\pi t)(\cosh \eta t)=-1.\]
We can assume $\eta>0$.
Then, by continuity, the equation above always has a root $0<t_0<1$.
On the other hand the matrix above is not a scalar matrix for $0<t<1$; thus  at critical point the $\log$
is not continuable.
\end{proof}

\end{prop}
\snewpage

\begin{prop} \plabel{pr:sharpbal}
(a) In Example \ref{ex:vinimp}.(a)
\[
A_1^{[\delta]}=\AC(\delta)\begin{bmatrix}\delta&-1\\1&-\delta\end{bmatrix},
\qquad\text{and}\qquad
A_2^{[\delta]}=\AC(\delta)\begin{bmatrix}-\delta&-1\\1&\delta\end{bmatrix}.
\]
As $\delta\searrow0$, or $\delta\nearrow0$, we have
\[\|A_1^{[\delta]}\|_2=\|A_2^{[\delta]}\|_2=\AC(\delta)(1+|\delta|)\searrow\frac\pi2.\]

For $\delta>0$,
the convergence radius of the Magnus expansion of $A_1^{[\delta]}\mathbf 1\boldsymbol. A_2^{[\delta]}\mathbf 1$ is exactly $1$
but the Magnus expansion is divergent (not even Abel summable). Here
\[\lim_{t\nearrow1}\,\,\, \log\left(\exp t A_1^{[\delta]}\exp tA_2^{[\delta]}\right) \cdot\sqrt{1-t^2}=
2\pi\sqrt{\frac{\delta}{\AC(\delta)}}\cdot\begin{bmatrix}0&-1\\0&0\end{bmatrix}\neq0\]
indicates the rate of divergence.

For $-1<\delta<0$,
the convergence radius of the Magnus expansion of $A_1^{[\delta]}\mathbf 1\boldsymbol. A_2^{[\delta]}\mathbf 1$ is
\[\frac{\AC(-\delta)}{\AC(\delta)}\equiv\frac{\AC(|\delta|)}{\AC(-|\delta|)}=\frac{\arccos(|\delta|)}{\pi- \arccos(|\delta|)}<1. \]

In fact, the cases $-1<\delta<0$ and $0<\delta<1$ are related to each other by
\begin{equation}
  A_1(\delta)\mathbf 1\boldsymbol. A_2(\delta)\mathbf 1=\frac{\AC(\delta)}{\AC(-\delta)}\cdot
 \tilde I\cdot A_1(-\delta)\mathbf 1\boldsymbol. A_2(-\delta)\mathbf 1 \cdot \tilde I^{-1};
\plabel{eq:dan1}
\end{equation}
i. e., through conjugating by $\tilde I=\bem&-1\\1&\eem$ and rescaling.

(b) In Example \ref{ex:vinimp}.(b)
\[
\tilde A_1^{[\eta]}=\frac\pi2\left[\begin{matrix}&-\mathrm e^{\eta} \\\mathrm e^{-\eta}&\end{matrix}\right],
\qquad\text{and}\qquad
\tilde A_2^{[\eta]}=\frac\pi2\left[\begin{matrix}&-\mathrm e^{-\eta} \\\mathrm e^{\eta}&\end{matrix}\right].
\]
As $|\eta|\searrow0$,
\[\|\tilde A_1^{[\eta]}\|_2=\|\tilde A_2^{[\eta]}\|_2=\frac\pi2\mathrm e^{|\eta|}\searrow\frac\pi2.\]
The convergence radius of in the Magnus expansion of $\tilde A_1^{[\eta]}\mathbf 1\boldsymbol.\tilde A_2^{[\eta]}\mathbf 1$ is
\[\frac2\pi\arccos\tanh|\eta|<1.\]

In fact,  case (a) $0<\delta<1$ and case (b) are related as follows:
Let us consider the matrix
\[\tilde L=\frac1{\sqrt{2}}\begin{bmatrix}1&-1\\1&1\end{bmatrix}.\]
Then we find
\begin{equation}
\tilde A_1^{[\eta]}\mathbf 1\boldsymbol.\tilde A_2^{[\eta]}\mathbf 1=
\underbrace{\frac1{\frac2\pi\arccos\tanh|\eta|}}_{>1}\cdot
\tilde L^{\sgn \eta}\cdot\mathbf 1A_1^{[\tanh|\eta|]}\mathbf 1\boldsymbol. A_2^{[\tanh|\eta|]}\mathbf 1\cdot \tilde L^{-\sgn\eta }.
\plabel{eq:dan2}
\end{equation}
I. e., relative to case (a), we apply $\delta=\tanh |\eta|$, conjugation by an orthogonal matrix, and simple rescaling.
\begin{commentx}
Hence, the essential difference between  Example \ref{ex:vinimp}(a) and (b) is only just that (a) picks up the very first not $\log$-able element in
the product exponential, while (b) picks up one rather later.
\end{commentx}

(c) In Example \ref{ex:vinimp}.(c), where $\delta\in\mathbb C\setminus(-\infty,-1]$ is allowed,
we still have
\[
A_1^{[\delta]}=\AC(\delta)\begin{bmatrix}\delta&-1\\1&-\delta\end{bmatrix},
\qquad\text{and}\qquad
A_2^{[\delta]}=\AC(\delta)\begin{bmatrix}-\delta&-1\\1&\delta\end{bmatrix}.
\]
As $\delta\rightarrow0$, we have
\[\|A_1^{[\delta]}\|_2=\|A_2^{[\delta]}\|_2\rightarrow\frac\pi2.\]

If $\Rea\delta=0$, such that $\delta=\mathrm i\varrho$ with $\varrho\in\mathbb R$, we have
\[
A_1^{[\mathrm i\varrho]}=\frac{\frac\pi2-\mathrm i\arsinh\varrho }{\sqrt{1+\varrho^2}}
\begin{bmatrix}\mathrm i\varrho&-1\\1&-\mathrm i\varrho\end{bmatrix},
\qquad\text{and}\qquad
A_2^{[\mathrm i\varrho]}=\frac{\frac\pi2-\mathrm i\arsinh\varrho }{\sqrt{1+\varrho^2}}
 \begin{bmatrix}-\mathrm i\varrho&-1\\1&\mathrm i\varrho\end{bmatrix};
\]
and, if $\varrho\neq0$, then the convergence radius of the Magnus expansion is $1$.

If  $\Rea\delta<0$, then the convergence radius of the Magnus expansion is strictly  less than $1$.

\proofremarkqed{In case (c),
the missing statement is that for  $\Rea\delta>0$,
 the   convergence radius of the Magnus expansion is $1$,
which is apparently true but
the actual computations seem to be tedious.
}
\begin{proof}

(a)
The logarithms and norms are straightforward;
only the convergence properties require nontrivial reasoning.

Let us assume $0<\delta<1$ (fix). Let $t\in\mathbb C$ (variable). Then
\[\exp t A_1^{[\delta]}=\cos\left(t\arccos \delta\right)
\begin{bmatrix}1&\\&1\end{bmatrix}+
\sin\left(t\arccos\delta\right)\frac1{\sqrt{1-\delta^2}}
\begin{bmatrix}\delta&-1\\1&-\delta\end{bmatrix},
\]
\[\exp tA_2^{[\delta]}=\cos\left(t\arccos\delta\right)
\begin{bmatrix}1&\\&1\end{bmatrix}+
\sin\left(t\arccos\delta\right)\frac1{\sqrt{1-\delta^2}}
\begin{bmatrix}-\delta&-1\\1&\delta\end{bmatrix}.
\]
Consequently,
\begin{multline*}
C_{12}^{[\delta]}(t):=\exp t A_1^{[\delta]}\exp tA_2^{[\delta]}=\\
=\cos^2\left(t\arccos\delta\right)\begin{bmatrix}1&\\&1\end{bmatrix}
-\frac{\sin^2\left(t\arccos\delta\right)}{1-\delta^2}\begin{bmatrix}1+\delta^2&2\delta
\\2\delta&1+\delta^2\end{bmatrix}
+\frac{\sin\left(2t\arccos\delta\right) }{\sqrt{1-\delta^2}}\begin{bmatrix}0&-1\\1&0\end{bmatrix}.
\end{multline*}
The possible obstacles to the analytic continuation of $\log$ are at
\[\frac{\tr C_{12}^{[\delta]}(t)}{2\sqrt{\det C_{12}^{[\delta]}(t)}^{\ext}}\equiv
\frac{\tr C_{12}^{[\delta]}(t)}{2}=
-\frac{1+\delta^2-2\cos^2(t\arccos \delta)}{1-\delta^2}=\pm1,\]
i. e. when
\[\cos(t\arccos \delta)=\pm1,\pm\delta.\]
We have to solve this for $t\neq0$.
For $t\in\intD(0,1)\setminus\{0\}$ this allows no solutions; thus the convergence radius of the Magnus expansion is at least $1$.
However, as we have seen, $C_{12}^{[\delta]}(1)$ does not allow a real logarithm, thus the Magnus expansion cannot be convergent at $t=1$.
For $t\in\intD(0,1)\setminus\{0\}$,
\begin{align*}
\log^{\ext} C_{12}^{[\delta]}(t)=&\AC\left( -\frac{1+\delta^2-2\cos^2(t\arccos \delta)}{1-\delta^2}\right)\cdot\\&
\left(-\frac{\sin^2\left(t\arccos\delta\right)}{1-\delta^2}\begin{bmatrix}0&2\delta
\\2\delta&0\end{bmatrix}
+\frac{\sin\left(2t\arccos\delta\right) }{\sqrt{1-\delta^2}}\begin{bmatrix}0&-1\\1&0\end{bmatrix}\right).
\end{align*}
In fact, by Lemma \ref{lem:logreal}, $\log$ is proper for $t\in(-1,1)$ (cf. next Remark).
Then, one can compute
\[\lim_{t\nearrow1} \left(\log  C_{12}^{[\delta]}(t)\right) \cdot\sqrt{1-t^2}=2\pi\sqrt{\frac{\delta\sqrt{1-\delta^2}}{\arccos\delta}}\cdot
\begin{bmatrix}0&-1\\0&0\end{bmatrix}\neq0.\]

Consider the case $\delta=1$. Here
\[\exp t A_1^{[1]}=
\begin{bmatrix}1&\\&1\end{bmatrix}+
t
\begin{bmatrix}1&-1\\1&-1\end{bmatrix},
\qquad
\text{and}
\qquad
\exp tA_2^{[1]}=
\begin{bmatrix}1&\\&1\end{bmatrix}+
t
\begin{bmatrix}-1&-1\\1&1\end{bmatrix}.
\]
Consequently,
\[
C_{12}^{[1]}(t):=\exp t A_1^{[1]}\exp tA_2^{[1]}
=\begin{bmatrix}1&\\&1\end{bmatrix}
-t^2\begin{bmatrix}2&2
\\2&2\end{bmatrix}
+2t\begin{bmatrix}0&-1\\1&0\end{bmatrix}.
\]
Similarly, searching for singularities, we have to solve
\[\frac{\tr C_{12}^{[1]}(t)}{2\sqrt{\det C_{12}^{[1]}(t)}^{\ext}}\equiv\frac{\tr C_{12}^{[1]}(t)}{2}=1-2t^2=\pm1.\]
for $t\neq0$.
Again,  we arrive to the critical value $t=1$; for $t\in\intD(0,1)\setminus\{0\}$,
\[\log^{\ext} C_{12}^{[1]}(t)=\AC\left( 1-2t^2\right)\cdot\left(
-2t^2\begin{bmatrix}0&1
\\1&0\end{bmatrix}
+2t\begin{bmatrix}0&-1\\1&0\end{bmatrix} \right);\]
$\log$ is proper for $t\in(-1,1)$; and
\[\lim_{t\nearrow1} \left(\log  C_{12}^{[1]}(t)\right) \cdot\sqrt{1-t^2}=2\pi
\begin{bmatrix}0&-1\\0&0\end{bmatrix}\neq0.\]

Assume that $\delta>1$.
Then
\[\exp t A_1^{[\delta]}=\cosh\left(t\arcosh \delta\right)
\begin{bmatrix}1&\\&1\end{bmatrix}+
\sinh\left(t\arcosh\delta\right)\frac1{\sqrt{\delta^2-1}}
\begin{bmatrix}\delta&-1\\1&-\delta\end{bmatrix},
\]
\[\exp tA_2(\delta)=\cosh\left(t\arcosh\delta\right)
\begin{bmatrix}1&\\&1\end{bmatrix}+
\sinh\left(t\arcosh\delta\right)\frac1{\sqrt{\delta^2-1}}
\begin{bmatrix}-\delta&-1\\1&\delta\end{bmatrix}.
\]
Consequently,
\begin{multline*}
C_{12}^{[\delta]}(t):=\exp t A_1^{[\delta]}\exp tA_2^{[\delta]}=\\
=\cosh^2\left(t\arcosh\delta\right)\begin{bmatrix}1&\\&1\end{bmatrix}
-\frac{\sinh^2\left(t\arcosh\delta\right)}{\delta^2-1}\begin{bmatrix}1+\delta^2&2\delta
\\2\delta&1+\delta^2\end{bmatrix}
+\frac{\sinh\left(2t\arcosh\delta\right) }{\sqrt{\delta^2-1}}\begin{bmatrix}0&-1\\1&0\end{bmatrix}.
\end{multline*}
The possible obstacles to the analytic continuation of $\log$ are at
\[\frac{\tr C_{12}^{[\delta]}(t)}{2\sqrt{\det C_{12}^{[\delta]}(t)}^{\ext}}
\equiv\frac{\tr C_{12}^{[\delta]}(t)}{2}=-\frac{1+\delta^2-2\cosh^2(t\arcosh \delta)}{1-\delta^2}=\pm1,\]
i. e., when
\[\cosh(t\arcosh \delta)=\pm1,\pm\delta.\]
We have to solve this for $t\neq0$.
For $t\in\intD(0,1)\setminus\{0\}$ this allows no solutions; but $t=1$ is already a point of divergence for the Magnus expansion.
For $t\in\intD(0,1)\setminus\{0\}$,
\begin{align*}
\log^{\ext} C_{12}^{[\delta]}(t)=&\AC\left( -\frac{1+\delta^2-2\cosh^2(t\arcosh \delta)}{1-\delta^2}\right)\cdot\\
&\cdot\left(-\frac{\sinh^2\left(t\arcosh\delta\right)}{\delta^2-1}\begin{bmatrix}0&2\delta
\\2\delta&0\end{bmatrix}
+\frac{\sinh\left(2t\arcosh\delta\right) }{\sqrt{\delta^2-1}}\begin{bmatrix}0&-1\\1&0\end{bmatrix}\right);
\end{align*}
it is easy to see is that $\log$ proper on $t\in(-1,1)$;
and
\[\lim_{t\nearrow1} \left(\log  C_{12}^{[\delta]}(t)\right) \cdot\sqrt{1-t^2}=2\pi\sqrt{\frac{\delta\sqrt{\delta^2-1}}{\arcosh\delta}}\cdot
\begin{bmatrix}0&-1\\0&0\end{bmatrix}\neq0.\]

Regarding the $-1<\delta<0$, the identity \eqref{eq:dan1} is easy to establish, then the rest follows from that.

(b) The identity \eqref{eq:dan2} can be checked in a straightforward manner,
from that the rest follows.

(c)
Only the statements regarding the convergence radii need much proof.
The formula \eqref{eq:dan1} also extends for $\delta\in\mathbb C\setminus\mathbb R$;
and then Lemma \ref{lem:arccos2}
implies that the convergence radius of the Magnus expansion is strictly less than $1$ for $0>\Rea \delta$.
Next, we deal with the case $\Rea \delta=0$. Let $\delta=\mathrm i\varrho$.
Then, similarly to our previous computations,  we find, for small $t$,
\begin{equation}
\frac{\tr C_{12}^{[\mathrm i\varrho]}(t)}{2\sqrt{\det C_{12}^{[\mathrm i\varrho]}(t)}^{\ext}}\equiv
\frac{\tr C_{12}^{[\mathrm i\varrho]}(t)}{2}=
-\frac{1-\varrho^2-2\cos^2\left(t\left(\frac\pi2-\mathrm i\arsinh \varrho\right) \right)}{1+\varrho^2}.
\plabel{eq:zon}
\end{equation}
Now, one can see that the lowest real value of $\cos^2\left(t\left(\frac\pi2-\mathrm i\arsinh \varrho\right)\right)$
for $t\in\Dbar(0,1)$
is
$\cos^2\left(\frac\pi2-\mathrm i\arsinh \varrho\right)$, which taken for some $t$ with $|t|=1$ (cf. Lemma \ref{lem:cos2}).
Thus smallest real value of \eqref{eq:zon} for $t\in\Dbar(0,1)$ is
\[
-\frac{1-\varrho^2-2\cos^2\left(\frac\pi2-\mathrm i\arsinh \varrho\right) }{1+\varrho^2}=
-\frac{1-\varrho^2-2\sin^2\left(-\mathrm i\arsinh \varrho\right) }{1+\varrho^2}
-\frac{1-\varrho^2+2\varrho^2 }{1+\varrho^2}
=-1,
\]
but it is taken only for some $t$ with $|t|=1$.
This implies that for $t\in\intD(0,1)$ the value of \eqref{eq:zon} is from $\mathbb\setminus(-\infty,0]$,
allowing the extension of the logarithm of the time-ordered exponential.
\end{proof}
\end{prop}

\begin{remark}\plabel{rem:logradiusexp}
Example \ref{ex:vinimp} (a) and (b) are cases when one
can (relatively easily) prove (cf. Remark \ref{rem:logradius}) that the $\log$-able radius of $\Rexp$ is equal to the convergence radius of
the Magnus expansion $\mu_{\mathrm R}$.

Indeed, it is sufficient to examine Example \ref{ex:vinimp}(a), $\delta>0$
(as the remaining cases are related by conjugation and rescaling).
For $0<\delta<1$
\[\frac{\tr C_{12}^{[\delta]}(t)}{2\sqrt{\det C_{12}^{[\delta]}(t)}^{\ext}}\equiv\frac{\tr C_{12}^{[\delta]}(t)}{2}
=-\frac{1+\delta^2-2\cos^2(t\arccos \delta)}{1-\delta^2}<-1,\]
means
\[\cos^2(t\arccos \delta)<-\delta^2.\]
We have to show its impossibility for $t\in\intD(0,1)$.
In general, $\cos z$ is real iff $z$ is real or $\Rea z\in\pi\mathbb Z$;
$\cos z$ is purely imaginary iff  $\Ima z\in\pi\left(\mathbb Z+\frac12\right)$.
(Cf. $\Rea \cos (x+y\mathrm i)=(\cos x)(\cosh y)$ and $\Ima \cos (x+y\mathrm i)=-(\sin x)(\sinh y)$ for $x,y\in\mathbb R$.)
As $\arccos \delta<\frac\pi2$, it sufficient to deal with case when $t$ is real,
in which case elementary calculus shows impossibility.
For $\delta=1$,
\[\frac{\tr C_{12}^{[1]}(t)}{2\sqrt{\det C_{12}^{[1]}(t)}^{\ext}}\equiv\frac{\tr C_{12}^{[1]}(t)}{2}=1-2t^2<-1\]
is clearly impossible for $t\in\intD(0,1)$.
For $\delta>1$,
\[\frac{\tr C_{12}^{[\delta]}(t)}{2\sqrt{\det C_{12}^{[\delta]}(t)}^{\ext}}
\equiv\frac{\tr C_{12}^{[\delta]}(t)}{2}
=-\frac{1+\delta^2-2\cosh^2(t\arcosh \delta)}{1-\delta^2}<-1,\]
means
\[\cosh^2(t\arcosh \delta)>\delta^2.\]
We have to show its impossibility  for $t\in\intD(0,1)$.
In general, $\cosh z$ is real iff $z$ is purely imaginary or $\Ima z\in\pi\mathbb Z$;
$\cosh z$ is purely imaginary iff  $\Rea z\in\pi\left(\mathbb Z+\frac12\right)$.
(Cf. $\Rea \cos (x+y\mathrm i)=(\cosh x)(\cos y)$ and $\Ima \cos (x+y\mathrm i)=(\sinh x)(\sin y)$ for $x,y\in\mathbb R$.)
In particular, $t$ should be real or purely imaginary.
If $t$ is purely imaginary, then  $\cosh^2(t\arcosh \delta)\leq1<\delta^2$,
while if $t$ is real, then elementary calculus shows impossibility.

An argument of this type was already used for purely imaginary $\delta$ in the proof of Proposition \ref{pr:sharpbal}(c).
\qedremark
\end{remark}

An uptake of the previous discussions is that the seemly uniform Examples \ref{ex:vin} and \ref{ex:vinimp}
are more fragmented than one naively expects.
\snewpage
\subsection{Skew-Hermitian balanced counterexamples in the BCH case}
\plabel{ss:SSBCH}~\\

We can modify the construction of Example \ref{ex:vinimp}(a) / Proposition \ref{pr:sharpbal}(a) as follows.
For $\delta\in\mathbb C\setminus(-\infty,0]$, let
\[
S_1^{[\delta]}=|\AC(\delta)|\begin{bmatrix}\delta&-1\\1&-\delta\end{bmatrix},
\qquad\text{and}\qquad
S_2^{[\delta]}=|\AC(\delta)|\begin{bmatrix}-\delta&-1\\1&\delta\end{bmatrix}.
\]
The measure $S_1^{[\delta]}\mathbf 1_{[0,1)}\boldsymbol.S_2^{[\delta]}\mathbf 1_{[1,2)}$ has  the same convergence radius as
$A_1^{[\delta]}\mathbf 1_{[0,1)}\boldsymbol.A_2^{[\delta]}\mathbf 1_{[1,2)}$ (as they differ only but a complex unit multiplier).
\begin{example}
\plabel{ex:vinimp12}
If $\delta\neq0$ but $\delta$ is purely imaginary, $\delta=\mathrm i\varrho$, then the
Magnus expansion of $S_1^{[\delta]}\mathbf 1_{[0,1)}\boldsymbol.S_2^{[\delta]}\mathbf 1_{[1,2)}$ has convergence radius $1$.
It is not absolutely convergent but Abel-summable.

Indeed, the Magnus expansion
of $S_1^{[\delta]}\mathbf 1_{[0,1)}\boldsymbol.S_2^{[\delta]}\mathbf 1_{[1,2)}$
differs from the Magnus expansion of
 $A_1^{[\delta]}\mathbf 1_{[0,1)}\boldsymbol.A_2^{[\delta]}\mathbf 1_{[1,2)}$
only by complex multipliers.
This implies the convergence radius and the lack of absolute convergence.
With some computation,
one can check that the extended $\log\Rexp(t\cdot S_1^{[\delta]}\mathbf 1_{[0,1)}\boldsymbol.S_2^{[\delta]}\mathbf 1_{[1,2)})$
has only four singularities on the unit circle, these are
$\pm\left(\frac{\AC(\delta)}{|\AC(\delta)|}\right)^{\pm1}
=\frac{\pm\frac\pi2\pm\mathrm i\arsinh\varrho}{\sqrt{(\frac\pi2)^2+(\arsinh\varrho)^2}}$.
At everywhere else the $\log$ can be continued over.
This also includes the value $1$, proving Abel-summability.

Here
$S_1^{[\delta]}\mathbf 1_{[0,1)}\boldsymbol.S_2^{[\delta]}\mathbf 1_{[1,2)}$ is skew-Hermitian, thus its time-ordered
exponential is unitary, and
\[
|S_1^{[\delta]}|= |S_2^{[\delta]}|= \sqrt{\left(\frac\pi2\right)^2+(\arsinh\varrho)^2}.
\]
which are $\sim\frac\pi2$ for small $\varrho$.

Again, the divergence is not so pronounced as in the original case, but one can amplify divergence by slightly
upscaling by a real multiplier.
\qedexer
\end{example}

\subsection{Arbitrary norm ratios in the BCH case}
\plabel{ss:ArbBCH}~\\

Here, an extension of the ``norm balanced'' Example \ref{ex:vinimp} to any other prescribed ``balance ratio'' is given.
(The case of arbitrary norm ratios for counterexamples to the BCH expansion is relatively complicated.
Our examples will be presented in analogy our earlier ones:
There is a ``parabolic'' counterexample which is relatively complicated, and a hyperbolic one, which is simpler;
 however, the parabolic one extends to the unitary / quaternionic case in a more straightforward manner.
Yet, the situation, in general, is not as nice as in the balanced case.
The parametrization we use in the parabolic case is not perfect but chosen for resemblance to the balanced case.)

\snewpage

\begin{exaprop}\plabel{ex:vinimp15}
(a) Let $\alpha\in(-\pi,\pi)$ be arbitrary.
Then, for $u>\cos\left(\pi\frac{\pi}{\pi+|\alpha|}\right)$ (the latter being a number from $[-1,0)$), we define
\[A_1^{\langle \alpha,u\rangle}=
\frac{\pi-\alpha}{2}\frac{\AC\left(u\right)}{\frac\pi2}\sqrt{\frac{1-u^2}{
\Cos\left(\tfrac{\alpha^2}{\pi^2}\AN( u)\right)^2-u^2
}}\left(\tilde I \Cos\left(\tfrac{\alpha^2}{\pi^2}\AN( u)\right)+u\tilde J \right)\]
and
\[A_2^{\langle \alpha,u\rangle}=
\frac{\pi+\alpha}{2}\frac{\AC\left(u\right)}{\frac\pi2}\sqrt{\frac{1-u^2}{
\Cos\left(\tfrac{\alpha^2}{\pi^2}\AN( u)\right)^2-u^2
}}\left(\tilde I \Cos\left(\tfrac{\alpha^2}{\pi^2}\AN( u)\right)-u\tilde J \right).\]
More precisely,
$\sqrt{\frac{1-u^2}{\Cos\left(\tfrac{\alpha^2}{\pi^2}\AN( u)\right)^2-u^2}}$
is resolved as $\frac{\pi}{\sqrt{\pi^2-\alpha^2}}$ for $u=1$; making
\[A_1^{\langle \alpha,1\rangle}=\sqrt{\frac{\pi-\alpha}{\pi+\alpha}}(\tilde I+\tilde J)
\qquad\text{and}\qquad
A_2^{\langle \alpha,1\rangle}= \sqrt{\frac{\pi+\alpha}{\pi-\alpha}}(\tilde I-\tilde J).\]
On the indicated domain, these define  smooth functions $A_i^{\langle \alpha,u\rangle}$ in $\alpha,u$.

(Remark: We will be primarily interested in this functions for $u\sim0$, for a fixed $\alpha$, only; where the finer considerations
 in the definition will be mostly unimportant.
In fact, for $u<1$, we have $\Cos\left(\frac{\alpha^2}{\pi^2}\AN( u)\right)\equiv\cos\left(\frac\alpha\pi\arccos u\right)$;
allowing all the computations to be done for $u\sim0$ in terms of $\arccos$.)

Then, one van see that $\|A_1^{\langle \alpha,u\rangle}\|_2:\|A_2^{\langle \alpha,u\rangle}\|_2=\frac{\pi-\alpha}2:\frac{\pi+\alpha}2$,
and as $u\rightarrow0$ we find that
$\left(\|A_1^{\langle \alpha,u\rangle}\|_2,\|A_2^{\langle \alpha,u\rangle}\|_2\right)
\rightarrow \left(\frac{\pi-\alpha}2,\frac{\pi+\alpha}2\right)$.

However, we claim that for $u\neq0$,
\[\Rexp\left(A_1^{\langle \alpha,u\rangle}\mathbf 1_{[0,1)}\boldsymbol.A_2^{\langle \alpha,u\rangle}\mathbf 1_{[1,2)}\right)
=\left(\exp A_1^{\langle \alpha,u\rangle}\right)\left(\exp A_2^{\langle \alpha,u\rangle}\right)\backsimeq
\bem-1&1\\&-1\eem.\]
(In particular it is not an exponential of any real $2\times2$ real matrix.)
This makes, in particular, the Magnus (BCH) expansion of
$A_1^{\langle \alpha,u\rangle}\mathbf 1_{[0,1)}\boldsymbol.A_2^{\langle \alpha,u\rangle}\mathbf 1_{[1,2)}$
completely divergent.
(This specializes to Example \ref{ex:vinimp}.(a) by $\alpha=0$ and $u=\delta$.)

(b) Suppose that $\alpha\in(-\pi,\pi)$ be arbitrary, and $\eta\in \mathbb R\setminus\{0\}$.
Let
\[\tilde A_1^{\langle\alpha,\eta\rangle}=\frac{\pi-\alpha}2\left[\begin{matrix}&-\mathrm e^{-\eta} \\\mathrm e^{\eta}&\end{matrix}\right],
\qquad\text{and}\qquad
\tilde A_2^{\langle\alpha,\eta\rangle}=\frac{\pi+\alpha}2\left[\begin{matrix}&-\mathrm e^{\eta} \\\mathrm e^{-\eta}&\end{matrix}\right].\]
Here,
\[\|\tilde A_1^{\langle\alpha,\eta\rangle}\|_2=\frac{\pi-\alpha}2\mathrm e^{|\eta|},
\qquad\text{and}\qquad
\|\tilde A_2^{\langle\alpha,\eta\rangle}\|_2=\frac{\pi+\alpha}2\mathrm e^{|\eta|}.\]

On the other hand, $(\exp A_1^{\langle\alpha,\eta\rangle})(\exp A_2^{\langle\alpha,\eta\rangle})$ is of negative strictly hyperbolic type.
In particular, it is not an exponential of a real $2\times2$ matrix.
In particular, the Magnus (BCH) expansion of
$A_1^{\langle\alpha,\eta\rangle}\mathbf 1_{[0,1)}\mathbf.A_2^{\langle\alpha,\eta\rangle}\mathbf 1_{[1,2)}$ is completely divergent.
In fact, its convergence radius is strictly smaller than $1$.
(This specializes to Example \ref{ex:vinimp}.(b) by $\alpha=0$.)

(c) Part (a) extends to $u\sim0$ (depending on $\alpha$), even if we allow $u$ to be complex.
\snewpage
\begin{proof}[Justification (Proof)]
(a) The function theoretic features, although not quite trivial can be established using standard analytical tools.
The spectral type of the exponential can be established as follows.
Using
\[\exp A_1^{\langle \alpha,u\rangle}=\Cos\left(\frac{(\pi-\alpha)^2}{\pi^2}\AN(u)\right)
+\Sin\left(\frac{(\pi-\alpha)^2}{\pi^2}\AN(u)\right) \exp A_1^{\langle \alpha,u\rangle}\]
and
\[\exp A_2^{\langle \alpha,u\rangle}=\Cos\left(\frac{(\pi+\alpha)^2}{\pi^2}\AN(u)\right)
+\Sin\left(\frac{(\pi+\alpha)^2}{\pi^2}\AN(u)\right) \exp A_2^{\langle \alpha,u\rangle}\]
we find
\begin{multline*}
(\exp A_1^{\langle \alpha,u\rangle})(\exp A_2^{\langle \alpha,u\rangle})=-\Id_2+2u\Cos\left(\tfrac{\alpha^2}{\pi^2}\AN( u)\right)\cdot
\\
\cdot\left(\sqrt{\frac{1-u^2}{\Cos\left(\tfrac{\alpha^2}{\pi^2}\AN( u)\right)^2-u^2}}\left(\tilde I-\tilde J\tfrac\alpha\pi\AC(u)
\Sin\left(\tfrac{\alpha^2}{\pi^2}\AN( u)\right) \right)- \tilde K\right).
\end{multline*}
Then one can check that $((\exp A_1^{\langle \alpha,u\rangle})(\exp A_2^{\langle \alpha,u\rangle})+\Id_2)^2=0$,
showing that the matrix  $(\exp A_1^{\langle \alpha,u\rangle})(\exp A_2^{\langle \alpha,u\rangle})$ has only eigenvalues $-1$.
(The computations are best to done for $u\sim0$ (depending on $\alpha$), then apply analytic extension in $u$.
For $u\sim0$, we have  $\Cos\left(\frac{\alpha^2}{\pi^2}\AN( u)\right)\equiv\cos\left(\frac\alpha\pi\arccos u\right)$
and $\frac\alpha\pi\AC(u)\Sin\left(\frac{\alpha^2}{\pi^2}\AN( u)\right)=
\frac{\sin\left(\frac\alpha\pi\arccos u\right)}{\sqrt{1-u^2}}$ .)
The nonvanishing of the coefficient of $\tilde K$ for $u\neq0$ shows that
$(\exp A_1^{\langle \alpha,u\rangle})(\exp A_2^{\langle \alpha,u\rangle})\neq-\Id_2$.

(b)
Only the establishing of the spectral type of the exponential require efforts.
Conjugating $\bem&1\\-1&\eem$ by $\bem\mathrm e^{\eta/2}&\\&\mathrm e^{-\eta/2}\eem^{\pm1}$, we see that
\[\exp \left(t\frac{\pi-\alpha}2\left[\begin{matrix}&-\mathrm e^{\eta} \\\mathrm e^{-\eta}&\end{matrix}\right]\right)=
\bem\cos\frac{\pi-\alpha}2t&-\mathrm e^{\eta} \sin\frac{\pi-\alpha}2t\\\mathrm e^{-\eta}\sin\frac{\pi-\alpha}2t&\cos\frac{\pi-\alpha}2t\eem,\]
and
\[\exp \left(t\frac{\pi+\alpha}2\left[\begin{matrix}&-\mathrm e^{-\eta} \\\mathrm e^{\eta}&\end{matrix}\right]\right)=
\bem\cos\frac{\pi+\alpha}2t&-\mathrm e^{-\eta} \sin\frac{\pi+\alpha}2t\\\mathrm e^{\eta}\sin\frac{\pi+\alpha}2t&\cos\frac{\pi+\alpha}2t\eem.\]
Hence, one can compute
$\tilde C_{12}=(\exp\tilde A_1^{\langle\alpha,\eta\rangle})\cdot(\exp\tilde A_2^{\langle\alpha,\eta\rangle})$,
whose characteristic equation (in $\lambda$) turns out to be
\[\lambda^2+2\left(1+(\sinh\eta)^2(1+\cos\alpha) \right)\lambda+1=0.\]
Now, $1+(\sinh\eta)^2(1+\cos\alpha) >1$ shows that $\tilde C_{12}$ has two distinct negative eigenvalues.
Consequently, $\tilde C_{12}$ is not an exponential.

If we replace $\tilde A_1^{\langle\alpha,\eta\rangle}$ and $\tilde A_2^{\langle\alpha,\eta\rangle}$ by
$\hat A_1^{\langle\alpha,\eta\rangle}=(1-\hat\varepsilon)A_1^{\langle\alpha,\eta\rangle}$
 and $\hat A_2^{\langle\alpha,\eta\rangle}=(1-\hat\varepsilon)A_2^{\langle\alpha,\eta\rangle}$
 (where $\hat\varepsilon$ is small),
then the characteristic polynomial of
$(\exp \hat A_1^{\langle\alpha,\eta\rangle})((\exp \hat A_2^{\langle\alpha,\eta\rangle}))$ is still of shape
$\lambda^2+2(1+\hat V(\alpha,\eta))\lambda+1$, where $\hat V(\alpha,\eta)>0$.
This implies that  the Magnus (BCH) expansion of
$\hat A_1^{\langle\alpha,\eta\rangle}\mathbf 1_{[0,1)}\mathbf.\hat A_2^{\langle\alpha,\eta\rangle}\mathbf 1_{[1,2)}$ is also divergent.
This shows that the convergence radius of the Magnus (BCH) expansion of
$\tilde  A_1^{\langle\alpha,\eta\rangle}\mathbf 1_{[0,1)}\mathbf. \tilde A_2^{\langle\alpha,\eta\rangle}\mathbf 1_{[1,2)}$
is strictly smaller than $1$.

(c) This is immediate.
\end{proof}
\end{exaprop}

\begin{example}\plabel{ex:vinquat}
Let us define (at least for $u\sim0$, depending on a fixed $\alpha\in(-\infty,\infty)$)
\begin{multline*}
\hat A_1^{\langle \alpha,u\rangle}=
\frac{\pi-\alpha}{2}
\left|\frac{\AC\left(u\right)}{\frac\pi2}
\sqrt{\frac{1-u^2}{\Cos\left(\tfrac{\alpha^2}{\pi^2}\AN( u)\right)^2-u^2}}
\Cos\left(\tfrac{\alpha^2}{\pi^2}\AN( u)\right)\right|
\cdot\\\cdot
\left(\tilde I+\frac{u}{\Cos\left(\tfrac{\alpha^2}{\pi^2}\AN( u)\right)} \tilde J \right)
\end{multline*}
and
\begin{multline*}
\hat A_2^{\langle \alpha,u\rangle}=
\frac{\pi+\alpha}{2}
\left|\frac{\AC\left(u\right)}{\frac\pi2}
\sqrt{\frac{1-u^2}{\Cos\left(\tfrac{\alpha^2}{\pi^2}\AN( u)\right)^2-u^2}}
\Cos\left(\tfrac{\alpha^2}{\pi^2}\AN( u)\right)\right|
\cdot\\\cdot
\left(\tilde I-\frac{u}{\Cos\left(\tfrac{\alpha^2}{\pi^2}\AN( u)\right)} \tilde J \right).
\end{multline*}
The difference between
$A_1^{\langle \alpha,u\rangle}\mathbf 1_{[0,1)}\boldsymbol.A_2^{\langle \alpha,u\rangle}\mathbf 1_{[1,2)}$
and
$\hat A_1^{\langle \alpha,u\rangle}\mathbf 1_{[0,1)}\boldsymbol.\hat A_2^{\langle \alpha,u\rangle}\mathbf 1_{[1,2)}$
is just multiplication by a complex unit vector.
That means the same cumulative norm and also that the
Magnus (BCH) expansion of
$\hat A_1^{\langle \alpha,u\rangle}\mathbf 1_{[0,1)}\boldsymbol.\hat A_2^{\langle \alpha,u\rangle}\mathbf 1_{[1,2)}$
is still divergent for $u\neq 0$ (although we cannot claim the things about the Jordan form and complete divergence anymore).

For a fixed $\alpha\in(-\pi,\pi)$, the function $f_\alpha:u\mapsto\frac{u}{\Cos\left(\tfrac{\alpha^2}{\pi^2}\AN( u)\right)}$
 takes $f_\alpha(0)=0$ with complex derivative $f'_\alpha(0)=\frac{1}{\cos\frac\alpha2}$.
Therefore, it can be inverted by the function $g_\alpha$ near $u\sim0$.
Here $g_\alpha(u)=u\cos\frac\alpha2+u^2\frac{\alpha\sin\alpha}{2\pi}+\ldots$.
This, for $\delta\sim0$, allows us to define
\[\hat A_i^{[\alpha,\delta]}:=\hat A_i^{\langle \alpha,g_\alpha(\delta)\rangle} ,\]
 (a slight change in parametrization).
This, for $\delta\sim0$, leads to matrices
\[\hat A_1^{[\alpha,\delta]}=\frac{\pi-\alpha}2S_\alpha(\delta)(\tilde I-\delta\tilde J)
\qquad\text{and}\qquad
\hat A_1^{[\alpha,\delta]}=\frac{\pi+\alpha}2S_\alpha(\delta)(\tilde I+\delta\tilde J)
\]
such that $S_\alpha(\delta)$ is real.
Altogether, one sees that
$\|A_1^{[\alpha,\delta]}\|_2:\|A_2^{[\alpha,\delta]}\|_2=\frac{\pi-\alpha}2:\frac{\pi+\alpha}2$,
and as $u\rightarrow0$ and
$\left(\|\hat A_1^{[\alpha,\delta]}\|_2,\|\hat A_2^{[\alpha,\delta]}\|_2\right)
\rightarrow \left(\frac{\pi-\alpha}2,\frac{\pi+\alpha}2\right)$.
At the same time, the Magnus (BCH) expansion of
$A_1^{[\alpha,\delta]} \mathbf 1_{[0,1)}\boldsymbol. A_2^{[\alpha,\delta]}\mathbf 1_{[1,2)}$
  divergent for $\delta\neq0$.
The additional feature here, however, is the following:
If $\delta$ is purely imaginary, then the matrices
$\hat A_i^{[\alpha,\delta]}$ are skew-Hermitian (purely imaginary quaternions)
making the time-ordered exponential special unitary (quaternionic with determinant $1$).
\qedexer
\end{example}
\snewpage
\begin{theorem}
If $n_1,n_2>0$ such that $n_1+n_2>\pi$, then we can find matrices
$\check{A_1}$ and $\check{A_2}$ such that

(i) $\|\check{A_1}\|_2=n_1$ and  $\|\check{A_2}\|_2=n_2$,

(ii) $\check{A_1}$ and $\check{A_2}$ are traceless real $2\times2$ matrices, or complex skew-Hermitian $2\times2$ matrices
(depending on request), such that

(iii)
the Magnus (BCH) expansion of
$ \check{A_1} \mathbf 1_{[0,1)}\boldsymbol. \check{A_2} \mathbf 1_{[1,2)}$ is divergent with convergent radius less than $1$.
\begin{proof}
Using the previous example, as $\delta\rightarrow0$, we can find matrices $\hat A_1$ and $\hat A_2$
such that $\|\hat{A_1}\|_2 : \|\hat{A_2}\|_2=n_1:n_2$ and $\|\hat{A_1}\|_2 + \|\hat{A_2}\|_2<n_1+n_2$.
The matrices are   traceless real or complex skew-Hermitian if the $\delta$ is real or purely imaginary in the process.
Then, linearly upscaling the matrices $\hat A_i$ by $\frac{n_1+n_2}{\|\hat{A_1}\| + \|\hat{A_2}\|}$ leads to the statement.
\end{proof}
\end{theorem}

\snewpage
\subsection{``Sharp'' counterexamples in the BCH case}
\plabel{ss:SharpBCH}~\\

The critical case with $\|A_1\|+\|A_2\|=\pi$, $\|A_1\|>0$, $\|A_2\|>0$ is trickier.
There are no counterexamples using finite matrices;
 Theorem \ref{th:better2} tells that in that case the convergence radius is greater than $1$.
However, there are ``sharp'' counterexamples if the dimension is allowed to be infinite.
\begin{example}\plabel{ex:vinimp2}
Let $\alpha\in(-\pi,\pi)$.
For $\delta\in\mathbb C$ we set
\begin{equation}
B_1^{[\alpha,\delta]}=\frac{\pi-\alpha}2\frac{1}{\Delta(\delta) }\bem \delta&-1\\1&-\delta\eem
\qquad\text{and}\qquad
 B_2^{[\alpha,\delta]}=\frac{\pi+\alpha}2\frac{1}{\Delta(\delta) }\bem -\delta&-1\\1&\delta\eem,
\plabel{eq:BB1}
\end{equation}
(cf.~\eqref{eq:uninorm}).
Here the cumulative radius of $ B_1^{[\alpha,\delta]}\mathbf 1_{[0,1)}\boldsymbol.  B_2^{[\alpha,\delta]}\mathbf 1_{[1,2)}$
 is $\frac{\pi-\alpha}2+\frac{\pi+\alpha}2=\pi$.
By Theorem \ref{th:better2}, the convergence radius of
$ B_1^{[\alpha,\delta]}\mathbf 1_{[0,1)}\boldsymbol.  B_2^{[\alpha,\delta]}\mathbf 1_{[1,2)}$
is strictly greater than $1$.
However, if $\delta$ is close to $0$, then, by Example (Proposition) \ref{ex:vinimp15} we know that the
convergence radius cannot be much greater $1$, because those counterexamples are just slightly upscaled versions of
this construction.

A variant of the construction is given by
\[
\tilde B_1^{[\alpha,\eta]}=\frac{\pi-\alpha}2 \mathrm{e}^{-|\Rea\eta|}
 \frac{|\cosh \eta|}{\cosh\eta}
\begin{bmatrix}&-\mathrm e^{-\eta}\\\mathrm e^{\eta}&\end{bmatrix}
\]
and
\[
\tilde B_2^{[\alpha,\eta]}=\frac{\pi+\alpha}2 \mathrm{e}^{-|\Rea\eta|}
\frac{|\cosh \eta|}{\cosh\eta}
\begin{bmatrix}&-\mathrm e^{\eta}\\\mathrm e^{-\eta}&\end{bmatrix},
\]
where $\eta\in\mathbb C\setminus \pi\mathrm i(\mathbb Z+\frac12)$.
Again, the cumulative radius of $\tilde B_1^{[\alpha,\eta]}\mathbf 1_{[0,1)}\boldsymbol. \tilde B_2^{[\alpha,\eta]}\mathbf 1_{[1,2)}$
 is $\frac{\pi-\alpha}2+\frac{\pi+\alpha}2=\pi$, etc.
The two constructions above are essentially equivalent,
\begin{equation}
\tilde B_1^{[\alpha,\eta]}\mathbf 1\boldsymbol.\tilde B_2^{[\alpha,\eta]}\mathbf 1=
\tilde L\cdot\mathbf 1B_1^{[\alpha,\tanh \eta ]}\mathbf 1\boldsymbol. B_2^{[\alpha,\tanh \eta ]}\mathbf 1\cdot\tilde L^{-1} .
\plabel{eq:dan3}
\end{equation}
For $\eta\sim0$ the correspondence $\delta=\tanh \eta$ is rather direct (apart from the conjugation by $\tilde L$);
real or purely imaginary $\eta$ corresponds to real or purely imaginary $\delta$, respectively.

Assume, for the sake of simplicity, that $\delta_*$ is sequence in $\mathbb C$ such that $\delta_i\neq0$ but $\delta_i\rightarrow0$
as $i\rightarrow+\infty$.
As  $\delta_i\rightarrow0$ we know that the lim sup of the convergence radius is at most $1$
(and, as we have discussed, the limit is actually $1$).
We can take the direct sums
\[  B_1^{[\alpha,\delta_*]}=\bigoplus_{n=1}^\infty B_1^{[\alpha,\delta_n]},
\qquad\text{and}\qquad
   B_2^{[\alpha,\delta_*]}=\bigoplus_{n=1}^\infty B_2^{[\alpha,\delta_n]}.\]
Here, $\|  B_1^{[\alpha,\delta_*]}\|_2=\frac{\pi-\alpha}2$ and $\|  B_2^{[\alpha,\delta_*]}\|_2=\frac{\pi+\alpha}2$.
Let $\phi=  B_1^{[\alpha,\delta_*]}\mathbf 1_{[0,1)}\boldsymbol.  B_2^{[\alpha,\delta_*]}\mathbf 1_{[1,2)}$.
Then the convergence radius of the Magnus (BCH) expansion of $\pi$ at most $1$.
Indeed, the various direct components (whose convergence radii close on $1$) prevent any greater value.
 On the other hand, according to Theorem \ref{th:MMNC}, the convergence radius is at least $1$, thus it is exactly $1$.
This statement, however, does not decide the convergence of the Magnus expansion yet.
For that we need an additional argument.
\qedexer
\end{example}

Assume that $\alpha\in(-\pi,\pi)$, and $\xi\in\mathbb C$ such that $|\xi|=1$.
Then we set
\begin{multline}
U_{\xi}^{[\alpha]}=\frac{\pi}{\sqrt{(\pi\Rea\xi)^2-4\xi^2(\cos\frac\alpha2)^2 }}\cdot
\left( \pi|\Rea\xi|\tilde I-2\xi\cos\frac\alpha2\left(\tilde J\sin\frac\alpha2 + \tilde K\cos\frac\alpha2\right)\right)
\\
\equiv\frac{\pi}{\sqrt{(\pi\Rea\xi)^2-\xi^2(\sin\alpha)^2-\xi^2(1+\cos\alpha)^2 }}\cdot
\left(\pi|\Rea\xi|\tilde I-(\xi\sin\alpha)\tilde J-\xi(1+\cos\alpha)\tilde K\ \right).\plabel{eq:xoro}
\end{multline}
(The expressions under the square root sign have real part $\geq 2(1+\cos\alpha)>0$,
 making the square root canonical, non-vanishing, and continuous).
We also define the set
\[\boldsymbol U^{[\alpha]}=\{U^{[\alpha]}_{\xi}\,:\,\xi\in\mathbb C\text{ and }|\xi|=1\}.\]
This is a compact set.
Now we can state
\begin{prop}\plabel{pr:vinimp2}
Assume that  $\alpha\in(-\pi,\pi)$.

Then the convergence radius of the Magnus expansion of
$  B_1^{[\alpha,\delta]}\mathbf 1_{[0,1)}\mathbf .  B_2^{[\alpha,\delta]}\mathbf1_{[1,2)}$ is greater than $1$.
For any $t\in\intD(0,1)$,
\[
\lim_{\delta\rightarrow0}\mu_{\mathrm R}(t\cdot   B_1^{[\alpha,\delta]}\mathbf 1_{[0,1)}\boldsymbol.
  B_2^{[\alpha,\delta]}\mathbf 1_{[1,2)} )\equiv
\lim_{\delta\rightarrow0} \log((\exp t\cdot   B_1^{[\alpha,\delta]})(\exp  t\cdot  B_2^{[\alpha,\delta]}))=t\cdot\pi\tilde  I.
\]
However (corresponding to case $t=1$), as $\delta\rightarrow 0$ ($\delta\neq0$) the values
\[\mu_{\mathrm R}(  B_1^{[\alpha,\delta]}\mathbf 1_{[0,1)}\boldsymbol.
  B_2^{[\alpha,\delta]}\mathbf 1_{[1,2)} )=\log((\exp  B_1^{[\alpha,\delta]})(\exp   B_2^{[\alpha,\delta]}))\]
limit to the set $\boldsymbol U^{[\alpha]}$.
The $\boldsymbol U^{[\alpha]}$ does not contain $\pi\tilde I$.
In fact, the distance $\boldsymbol U^{[\alpha]}$ and $\mathbb C\cdot \pi\tilde I$ is $2\cos\frac\alpha2>0$.
However, the elements $U^{[\alpha]}_{\xi}$ of $\boldsymbol U^{[\alpha]}$ satisfy $\exp U^{[\alpha]}_{\xi}=-\Id_2$.

\begin{commentx}
In particular,
\[\lim_{\delta\searrow0}\mu_{\mathrm R}(  B_1^{[\alpha,\delta]}\mathbf 1_{[0,1)}\boldsymbol.
  B_2^{[\alpha,\delta]}\mathbf 1_{[1,2)} )\equiv
\lim_{\delta\searrow0} \log((\exp  B_1^{[\alpha,\delta]})(\exp   B_2^{[\alpha,\delta]}))=  U_{1}^{[\alpha]},\]
and
\[\lim_{\delta\nearrow0}\mu_{\mathrm R}(  B_1^{[\alpha,\delta]}\mathbf 1_{[0,1)}\boldsymbol.
  B_2^{[\alpha,\delta]}\mathbf 1_{[1,2)} )\equiv
\lim_{\delta\nearrow0} \log((\exp  B_1^{[\alpha,\delta]})(\exp   B_2^{[\alpha,\delta]}))=  U_{-1}^{[\alpha]},\]
and
\[\lim_{\delta\searrow0}\mu_{\mathrm R}(  B_1^{[\alpha,\mathrm i\delta]}\mathbf 1_{[0,1)}\boldsymbol.
  B_2^{[\alpha,\mathrm i\delta]}\mathbf 1_{[1,2)} )\equiv
\lim_{\delta\searrow0} \log((\exp  B_1^{[\alpha,\mathrm i\delta]})(\exp   B_2^{[\alpha,\mathrm i\delta]}))
=  U_{\mathrm i}^{[\alpha]},\]
and
\[\lim_{\delta\nearrow0}\mu_{\mathrm R}(  B_1^{[\alpha,\mathrm i\delta]}\mathbf 1_{[0,1)}\boldsymbol.
  B_2^{[\alpha,\mathrm i\delta]}\mathbf 1_{[1,2)} )\equiv
\lim_{\delta\nearrow0} \log((\exp  B_1^{[\alpha,\mathrm i\delta]})(\exp   B_2^{[\alpha,\mathrm i\delta]}))
=  U_{-\mathrm i}^{[\alpha]},\]
where
\begin{align*}
  U_{\pm1}^{[\alpha]}&=\frac\pi{\sqrt{\pi^2-2-2\cos\alpha}}\cdot
\left(\pi \tilde I\pm\left(-\tilde  J\sin\alpha - \tilde K(1+\cos\alpha)\right)  \right).
\\
&=\frac\pi{\sqrt{\pi^2-\left(2\cos\frac\alpha2\right)^2}}\cdot
\left(\pi \tilde I\pm
\left(2\cos\frac\alpha2\right)\cdot\left(-\tilde J\sin\frac\alpha2 - \tilde K\cos\frac\alpha2\right)
 \right),
\end{align*}
and
\[  U_{\pm\mathrm i}^{[\alpha]}=\pm\mathrm i\pi\left(- \tilde J\sin\frac\alpha2 - \tilde K\cos\frac\alpha2\right).\]
We note that
\[\|  U_{\pm1}^{[\alpha]}\|_2=\pi\sqrt{\frac{\pi+2\cos\frac\alpha2}{\pi-2\cos\frac\alpha2}}>\pi\]
and
\[\|  U_{\pm\mathrm i}^{[\alpha]}\|_2=\pi.\]
\end{commentx}

\begin{proof}
We have already seen the statement about the convergence radius.
For $\delta\neq0$, the complex pencil generated by  $B_1^{[\alpha,\delta]}$ and  $B_2^{[\alpha,\delta]}$ is the same as the one
 generated by $\tilde I$ and $\tilde J$, excluding any common eigenvector for $B_1^{[\alpha,\delta]}$ and  $B_2^{[\alpha,\delta]}$.
This and the cumulative norm, and Theorem \ref{th:better1} implies that
 the logarithm of the time-ordered exponential can be taken.
The limit for $t\in\intD(0,1)$ follows from the continuity of $\log$ at $\exp\left(t\pi\tilde I\right)$.
This argument, of course, does not work for $t=1$.
For $\delta\neq0$, (and $\delta\neq\pm1$ but those values can be resolved), we find
\begin{multline}
\log((\exp  B_1^{[\alpha,\delta]})(\exp   B_2^{[\alpha,\delta]}))=
\AC\left(\frac{\Cos\left(\frac{\pi^2(1-\delta^2)  }{\Delta(\delta)^2}\right)
-\delta^2\Cos\left(\frac{\alpha^2(1-\delta^2)  }{\Delta(\delta)^2}\right)}{1-\delta^2}\right) |\delta| \cdot
\\
\cdot\frac1{|\delta|}\Biggl(\frac\pi{\Delta(\delta)}\Sin\left(\frac{\pi^2(1-\delta^2)  }{\Delta(\delta)^2}\right)\tilde I
-\frac{\alpha\delta}{\Delta(\delta)}\Sin\left(\frac{\alpha^2(1-\delta^2)  }{\Delta(\delta)^2}\right)\tilde J
\\+\delta\left(\frac{\Cos\left(\frac{\pi^2(1-\delta^2)  }{\Delta(\delta)^2}\right)
-\Cos\left(\frac{\alpha^2(1-\delta^2)  }{\Delta(\delta)^2}\right)}{1-\delta^2}\right)\tilde K
\Biggr).
\plabel{eq:aoso}
\end{multline}
Let us write $\delta=\xi|\delta|$.
Then \eqref{eq:aoso} is as \eqref{eq:xoro} but $+o(|\delta|)$ for $\delta\sim0$.
(One can use separation of cases for $\Rea\xi\geq0$ and $\Rea\xi\leq0$.)
This proves the limiting statement.
It is easy to see that $\boldsymbol U^{[\alpha]}$ does not contain real multiples of $\pi\tilde I$.
Regarding the distance, after conjugating by $\exp(\frac\alpha4 \tilde I)$ and multiplying by $\xi^{-1}$, we
have to minimize $\|\lambda\tilde I-2\cos\frac\alpha2\tilde K\|_2$.
Having antidiagonal matrices, this is easy to investigate, it is minimal for $\lambda=0$.

The exponential statement is straightforward.
\end{proof}
\end{prop}

\snewpage
\begin{theorem}\plabel{th:vinimp2}
In Example \ref{ex:vinimp2}, the Magnus expansion of
$B_1^{[\alpha,\delta_*]}\mathbf 1_{[0,1)}\boldsymbol.  B_2^{[\alpha,\delta_*]}\mathbf 1_{[1,2)}$ is (completely) divergent;
yielding a counterexample to the convergence of the Magnus (BCH) expansion with
$\|  B_1^{[\alpha,\delta_*]}\|_2=\frac{\pi-\alpha}2$ and $\|  B_2^{[\alpha,\delta_*]}\|_2=\frac{\pi+\alpha}2$ where $\alpha\in(-\pi,\pi)$.

If the $\delta_i$ are all purely imaginary, then the $B_i^{[\alpha,\delta_*]}$ are skew-Hermitian (and the exponentials are unitary).
In that case, the realifications are skew-symmetric and the exponential are orthogonal.

\begin{proof}
Let $D_\alpha$ be the $\|\cdot\|_2$-distance of the compact sets $\boldsymbol U^{[\alpha]}$ and $[0,1] \cdot\pi\tilde I$
(which is at least $2\cos\frac\alpha2$).
We will prove

(X) ``For any $0<t<1$ and any $\varepsilon>0$ there is a $k\geq1$ and $t<\tilde t<1$ such that
\[\left\|
\Rexp\left(\tilde t\cdot B_1^{[\alpha,\delta_k]}\mathbf 1_{[0,1)}\boldsymbol.  B_2^{[\alpha,\delta_k]}\mathbf 1_{[1,2)}\right)
-
\Rexp\left(t\cdot B_1^{[\alpha,\delta_k]}\mathbf 1_{[0,1)}\boldsymbol.  B_2^{[\alpha,\delta_k]}\mathbf 1_{[1,2)}\right)
\right\|_2>D_\alpha-\varepsilon
\]
holds.''

Statement (X) obviously precludes the radial convergence of the Magnus expansion of
$B_1^{[\alpha,\delta_*]}\mathbf 1_{[0,1)}\boldsymbol.  B_2^{[\alpha,\delta_*]}\mathbf 1_{[1,2)}$,
thus it is ideal for our purposes. We can prove it as follows:

If $k$ is sufficiently large, then
if $|\delta_k|$ is sufficiently small, and due to continuity of
$\log$ at $\exp(t\cdot\pi\tilde I)$, we have
\[\left\|
\log\Rexp\left(t\cdot B_1^{[\alpha,\delta_k]}\mathbf 1_{[0,1)}\boldsymbol.  B_2^{[\alpha,\delta_k]}\mathbf 1_{[1,2)}\right)
-t\cdot\pi\tilde I
\right\|_2<\dfrac\varepsilon3.\]
However, also, if $\delta_k$ is sufficiently small, then
\[\left\|
\log\Rexp\left(  B_1^{[\alpha,\delta_k]}\mathbf 1_{[0,1)}\boldsymbol.  B_2^{[\alpha,\delta_k]}\mathbf 1_{[1,2)}\right)
-U^{[\alpha,\delta_k]}
\right\|_2<\dfrac\varepsilon3,\]
where $U^{[\alpha,\delta_k]}$ is an appropriate element of $\boldsymbol U^{[\alpha]}$.
(We can use $U^{[\alpha,\delta_k]}=U^{[\alpha]}_{\delta_k/|\delta_k|}$.)
Finally, we can choose  $t<\tilde t<1$ such that
\[\left\|
\log\Rexp\left(\tilde t  B_1^{[\alpha,\delta_k]}\mathbf 1_{[0,1)}\boldsymbol.  B_2^{[\alpha,\delta_k]}\mathbf 1_{[1,2)}\right)
-
\log\Rexp\left(  B_1^{[\alpha,\delta_k]}\mathbf 1_{[0,1)}\boldsymbol.  B_2^{[\alpha,\delta_k]}\mathbf 1_{[1,2)}\right)
\right\|_2<\dfrac\varepsilon3\]
Thus for a small $\delta_k$ chosen, and an also an appropriate  $t<\tilde t<1$, we find
\begin{multline*}\Biggl\|
\left(\Rexp\left(\tilde t\cdot B_1^{[\alpha,\delta_k]}\mathbf 1_{[0,1)}\boldsymbol.  B_2^{[\alpha,\delta_k]}\mathbf 1_{[1,2)}\right)
-
\Rexp\left(t\cdot B_1^{[\alpha,\delta_k]}\mathbf 1_{[0,1)}\boldsymbol.  B_2^{[\alpha,\delta_k]}\mathbf 1_{[1,2)}\right)
\right)
-\\-
\left( U^{[\alpha,\delta_k]}-t\cdot\pi\tilde I\right)
\Biggr\|_2<\varepsilon.
\end{multline*}
This latter inequality implies the inequality of statement (X).

The comment about the unitary / orthogonal cases is immediate.
\end{proof}
\end{theorem}
Note that we have divergence only in norm topology but not in strong topology.
This still leaves a window for further investigations.

\snewpage\section{The conformal range (continued)}
\plabel{sec:ConformalRangeTwo}
The conformal range we have introduced is a particular aspect (in fact, a projection of) the so-called Davis--Wielandt
shell, cf. Davis \cite{D1} (1968), Davis \cite{D2} (1970), Wielandt \cite{Wie} (1953).
A more proper notation for it would be $\DW^{\mathbb R}_{\mathrm{PH}}$, the 2-dimensional Davis--Wielandt shell in the (asymptotically closed)
Poincar\'e half-plane model.
As this name is too long, we retain the name `conformal range'.
Here we include some more advanced aspects of the conformal range, hence the  Bolyai--Lobachevski\u{\i} hyperbolic geometry will also be applied.
(See Berger \cite{Ber} for a standard account of hyperbolic geometry.)
One familiar with \cite{D1}, \cite{D2} will find the subsequent discussion very easy.
In fact, for comparison, we review some basic properties of the Davis--Wielandt shell
in Appendices \ref{app:hyprev} and \ref{app:DWrev}.

\begin{lemma}\plabel{lem:conf} (Conformal invariance.)
Suppose that $g(x)=\frac{ax+b}{cx+d}$ is a real fractional linear function, $ad-bc\neq 0$.
Assume $A\in\mathcal B(\mathfrak H)$ and that $cA+d\Id$ is invertible.

(a) If $\mathbf x\in\mathfrak H\setminus 0$ and $\mathbf y= (cA+d\Id)^{-1}\mathbf x$, then
\[g(A)\mathbf x:\mathbf x= g(A\mathbf y:\mathbf y  )^{\textrm{conjugated if } ad-bc<0}.\]

(b) Consequently,
\[\CR(g(A))=g(\CR(A))^{\textrm{conjugated if } ad-bc<0}.\]
\[\CRext(g(A))=g(\CRext(A)).\]
\begin{proof} (a) The elementary rules
\begin{align}
\alpha\mathbf y:\mathbf x&=\alpha \cdot (\mathbf y:\mathbf x)^{\textrm{conjugated if } \alpha<0}
&&(\alpha\in \mathbb R),\notag\\
(\mathbf y+\beta\mathbf x):\mathbf x&=\mathbf y:\mathbf x+\beta
&&(\beta\in \mathbb R),\notag\\
\gamma\mathbf y:\gamma\mathbf x&=\mathbf y:\mathbf x
&&(\gamma\in \mathbb R\setminus\{0\}),\notag\\
\mathbf y:\mathbf x &= \overline{(\mathbf x:\mathbf y)}^{-1}
&&(\mathbf y\neq 0)\notag
\end{align}
are easy to check. If $g$ is linear ($c=0$), then the statement follows from from the first three rules.
If $g$ is not linear ($c\neq 0$), then $g(x)=\frac ac-\frac{ad-bc}{c^2}\left(x+\frac dc\right)^{-1}$, and
\begin{align}
g(A)\mathbf x:\mathbf x&=\frac ac-\frac{ad-bc}{c^2}
\left( \mathbf x:\left(A+\frac dc\Id\right)^{-1}\mathbf x
\right)^{-1,\textrm{ conjugated if } ad-bc<0} \notag\\
&=\frac ac-\frac{ad-bc}{c^2}
\left( \left(A+\frac dc\Id\right)\mathbf y: \mathbf y\right)^{-1,\textrm{ conjugated if } ad-bc<0}\notag\\
&=\left(\frac ac-\frac{ad-bc}{c^2}  \left( \left(A\mathbf y: \mathbf y\right)
+\frac dc\right)^{-1}\right)^{\textrm{ conjugated if } ad-bc<0}\notag\\
&=g \left(A\mathbf y: \mathbf y\right)^{\textrm{ conjugated if } ad-bc<0}.\notag
\end{align}
(b) This follows from the previous part and the conjugational symmetry of  $\CRext(A)$.
\end{proof}
\end{lemma}
\snewpage
Let $z_1,z_2\in\mathbb C$ such that $\Ima z_1,\Ima z_2\geq0$.
We say that the $h$-segment $[z_1,z_2]_h$ is the circular or straight segment connecting $z_1$ and $z_2$,
whose circle or line is perpendicular to the real axis,
and lies in the upper half plane $\upper=\{z\in\mathbb C\,:\,\Ima z\geq 0\}$.
\begin{lemma}\plabel{lem:conv} ($h$-Convexity.)
Suppose that $A\in\mathcal B(\mathfrak H)$, and $\dim_{\mathbb R}\mathfrak H\neq 2$.
Then $\CR(A)$ is $h$-convex, i.~e.~$z_1,z_2\in\CR(A)$
implies $[z_1,z_2]_h\subset \CR(A)$.
\begin{proof} We can suppose that $z_1\neq z_2$.
Applying linear conformal transformations to $A$, we can assume that $\Rea z_1=\Rea z_2=0$ (lineal case) or $|z_1|=|z_2|=1$ (circular case).
Assume that $ A\mathbf x_1:\mathbf x_1=z_1$, $ A\mathbf x_2:\mathbf x_2=z_2$.
Extend the span of $\{\mathbf x_1,\mathbf x_2\}$ to a $3$-dimensional space $V\subset\mathfrak H$.
Consider the quadratic form defined by
\[q(\mathbf x)=\begin{cases} \langle A\mathbf x,\mathbf x\rangle_{\mathrm{real}}&\text{(lineal case),}\\
\langle A\mathbf x,A\mathbf x\rangle_{\mathrm{real}}-\langle \mathbf x,\mathbf x\rangle_{\mathrm{real}}&\text{(circular case)}.\end{cases}\]
The nullset $V_q$ of $q$ on $V$ is either $V$, a plane, or a double cone (cf. $\mathbf x_1,\mathbf x_2\in V_q$).
In any case,
\[\{(A\mathbf x:\mathbf x)\,:\, \mathbf x\in V_q\setminus\{0\}\}\]
is a connected set (cf. $(A\mathbf x:\mathbf x)=(A(-\mathbf x):(-\mathbf x))$), which is contained in
\[L_h=\begin{cases} \{z\in\overline{\mathbb C}^+\,:\,\Rea z=0\}&\text{(lineal case),}\\
\{z\in\overline{\mathbb C}^+\,:|z|=1\}&\text{(circular case)}.\end{cases}\]
The connectedness implies $[z_1,z_2]_h\subset L_h$.
\end{proof}
\end{lemma}
(If $A$ is a linear relation, then dimension of the graph of $A$ is relevant, and the vertical segments $[z,\infty]_h$
also come to play.)

\snewpage
The previous lemma suggests that $\CR(A)\subset \upper$ is best to be interpreted as the
asymptotically closed Poincar\'e half-plane, which is a standard model for the hyperbolic plane.
In that regard, it is useful to consider the diagram
\[\xymatrix{{\mathbb C^\star}\ar[rr]^{\cq\qquad}&&(\upper)^\star\equiv\overline{H^2}_{\mathrm{PH}}\ar[rr]_{\simeq}^{
\tfrac{\mathrm{CKB}}{\mathrm{PH}}}&&\Dbar(0,1) }\equiv\overline{H^2}_{\mathrm{CKB}} .\]
Here the first map from the Riemann sphere to the asymptotically closed Poincar\'e half-plane is just factorization up to
conjugation, and the second map is the passage to the Cayley--Klein--Beltrami plane model, given by
\[\tfrac{\mathrm{CKB}}{\mathrm{PH}}: \qquad   (u_1,u_2)\mapsto\left(\frac{2u_1}{u_1^2+u_2^2+1} ,\frac{u_1^2+u_2^2-1}{u_1^2+u_2^2+1}\right),\qquad \infty\mapsto(0,1).\]
The CKB model is advantageous because there, $h$-segments correspond to ordinary segments, showing that $h$-convexity is
qualitatively not different from ordinary convexity.

Another variant is the  parabolic version of the Cayley--Klein--Beltrami model with
\[\xymatrix{{\mathbb C^\star}\ar[rr]^{\cq\qquad}&&(\upper)^\star\equiv\overline{H^2}_{\mathrm{PH}}\ar[rr]_{\simeq\qquad}^{
\tfrac{\mathrm{CKB(P)}}{\mathrm{PH}}\qquad}&&\{(x,z)\,:\,z\geq x^2\}\cup\{\infty\} }\equiv\overline{H^2}_{\mathrm{CKB(P)}},\]
given by
\[\tfrac{\mathrm{CKB(P)}}{\mathrm{PH}}: \qquad   (u_1,u_2)\mapsto\left(u_1,u_1^2+u_2^2\right),\qquad \infty\mapsto\infty.\]

We will not use it computations, but it is useful for certain visualizations to consider mapping to the Poincar\'e disk model,
\[\xymatrix{{\mathbb C^\star}\ar[rr]^{\cq\qquad}&&(\upper)^\star\equiv\overline{H^2}_{\mathrm{PH}}\ar[rr]_{\simeq}^{
\tfrac{\mathrm{P}}{\mathrm{PH}}}&&\Dbar(0,1) }\equiv\overline{H^2}_{\mathrm{P}},\]
given by
\[\tfrac{\mathrm{P}}{\mathrm{PH}}: \qquad   (u_1,u_2)\mapsto\left(\frac{2u_1}{u_1^2+(u_2+1)^2} ,\frac{u_1^2+u_2^2-1}{u_1^2+(u_2+1)^2}\right),\qquad \infty\mapsto(0,1).\]

As we will consider basically only bounded operator, the point $\infty$ will appear in the conformal range
(using $\upper$ instead of $(\upper)^\star$ will be sufficient), but, from time to time, it appears in geometric considerations.

If we deal with complex numbers (that id the Poincar\'e half-plane model), then we keep the terminology `conformal range',
but otherwise  we prefer to think about the range as the real Davis--Wielandt shell.
That is
$\CR(A)\equiv\DW_{\mathrm{PH}}^{\mathbb R}(A)$,
$\DW_{\mathrm{CKB}}^{\mathbb R}(A)\equiv \tfrac{\mathrm{CKB}}{\mathrm{PH}}(\CR(A) ) $,
$\DW_{\mathrm{CKB(P)}}^{\mathbb R}(A)\equiv\tfrac{\mathrm{CKB(P)}}{\mathrm{PH}}(\CR(A) )$,
$\DW_{\mathrm{P}}^{\mathbb R}(A)\equiv\tfrac{\mathrm{P}}{\mathrm{PH}}(\CR(A) )$.
It is a detour here, but the following lemma may help to get a feel with respect to various models.
For the purposes of the statements coordinates $x_{\mathrm{CKB}}$ etc. will be considered as functions
on the models (but, in reality, in the abstract conformal range).
\snewpage
\begin{lemma}[Extremal values in ranges]
\plabel{lem:extran}

(o) In the $\mathrm{PH}$ model, let $\mathrm n_{\mathrm{PH}}$ denote the distance from the origin.
Then
\[\sup\mathrm n_{\mathrm{PH}}( \CR(A) ) =\|A\|_2,\]
\[\inf\mathrm n_{\mathrm{PH}}( \CR(A) ) =\|A\|_2^-.\]

(a) In the $\mathrm{PH}$/$\mathrm{CKB(P)}$  models,
\[\sup x_{\mathrm{PH}}( \CR(A) )
=\sup x_{\mathrm{CKB(P)}}(\DW_{\mathrm{CKB(P)}}^{\mathbb R}(A))
=\max\spec\left(\frac {A+A^*}2\right),\]
\[\inf x_{\mathrm{PH}}( \CR(A) )
=\sup x_{\mathrm{CKB(P)}}(\DW_{\mathrm{CKB(P)}}^{\mathbb R}(A))
=\min\spec\left(\frac {A-A^*}2\right) .\]
(These are the real supremum and real infimum.)

(b) In the $\mathrm{CKB(P)}$  model,
\[\sup y_{\mathrm{CKB(P)}}(\DW_{\mathrm{CKB(P)}}^{\mathbb R}(A))
=\left(\|A\|_2\right)^2,\]
\[\inf y_{\mathrm{CKB(P)}}(\DW_{\mathrm{CKB(P)}}^{\mathbb R}(A))
=\left(\|A\|_2^-\right)^2.\]

(c) In the $\mathrm{CKB}$  model,
\[\sup y_{\mathrm{CKB}}(\DW_{\mathrm{CKB}}^{\mathbb R}(A))
=\frac{\left(\|A\|_2\right)^2-1}{\left(\|A\|_2\right)^2+1},\]
\[\inf y_{\mathrm{CKB}}(\DW_{\mathrm{CKB}}^{\mathbb R}(A))
=\frac{\left(\|A\|_2^-\right)^2-1}{\left(\|A\|_2^-\right)^2+1}.\]

(d) In the $\mathrm{CKB}$  model,
\[\sup x_{\mathrm{CKB}}(\DW_{\mathrm{CKB}}^{\mathbb R}(A))
=\frac{1-\biggl(\left\|\dfrac{\Id-A}{\Id+A}\right\|_2^-\biggr)^2}{1+\biggl(\left\|\dfrac{\Id-A}{\Id+A}\right\|_2^-\biggr)^2}
=\frac{\left(\left\|\dfrac{\Id+A}{\Id-A}\right\|_2\right)^2-1}{\left(\left\|\dfrac{\Id+A}{\Id-A}\right\|_2\right)^2+1}
\]
(the result is $1$ if $1\in\spec(A)$),
\[\inf x_{\mathrm{CKB}}(\DW_{\mathrm{CKB}}^{\mathbb R}(A))
=\frac{1-\left(\left\|\dfrac{\Id-A}{\Id+A}\right\|_2\right)^2}{1+\left(\left\|\dfrac{\Id-A}{\Id+A}\right\|_2\right)^2}
=\frac{\biggl(\left\|\dfrac{\Id+A}{\Id-A}\right\|_2^-\biggr)^2-1}{\biggl(\left\|\dfrac{\Id+A}{\Id-A}\right\|_2^-\biggr)^2+1}
\]
(the result is $-1$ if $-1\in\spec(A)$).
Altogether, the information here is obtained from the Cayley transform of $A$.
\begin{proof}
(o) This is just \eqref{eq:wrsup}--\eqref{eq:wrsupmin} repackaged.

(a) For $\mathbf x\in\mathfrak H$, $\Rea\langle A\mathbf x,\mathbf x\rangle=\left\langle\frac {A+A^*}2\mathbf x,\mathbf x\right\rangle$.

(b) $(x_{\mathrm{PH}})^2+(y_{\mathrm{PH}})^2=\lambda^2$ translates to $y_{\mathrm{CKB(P)}}=\lambda^2$.

(c)  $(x_{\mathrm{PH}})^2+(y_{\mathrm{PH}})^2=\lambda^2$ translates to $y_{\mathrm{CKB}}=\frac{\lambda^2-1}{\lambda^2+1}$.

(d) The (involutive) Cayley transform $\lambda\mapsto \frac{1-\lambda}{1+\lambda}$ in the $\mathrm{PH}$ model induces
reflexion to the line $x_{\mathrm{CKB}}+y_{\mathrm{CKB}}=0 $ in the $\mathrm{CKB}$ model.
Thus it is induced from case (c).
\end{proof}
\end{lemma}

\begin{lemma}\plabel{lem:join}  (a)
Suppose that $A_1\in\mathcal B(\mathfrak H_1)$, $A_2\in\mathcal B(\mathfrak H_2)$.
Let us consider the direct sum  $A_1\oplus A_2\in\mathcal B(\mathfrak H_1\oplus\mathfrak H_2)$.
Then
\[\CR(A_1\oplus A_2)
=\bigcup\{[z_1,z_2]_h \,:\,z_1 \in\CR(A_1),z_2\in \CR(A_2)\}.\]

(b) If $\dim_{\mathbb R}\mathfrak H\neq 2$, then complexification of $A$ does not change $\CR(A)$. Else,
\[\CR(A^{\mathbb C})
=\bigcup\{[z_1,z_2]_h \,:\,z_1 ,z_2\in \CR(A)\} ,\]
and this is already $h$-convex.
\begin{proof} (a) Suppose that $\mathbf x_1\in\mathfrak H_1$,  $\mathbf x_2\in\mathfrak H_2$, $A\mathbf x_i:\mathbf x_i=z_i$.
Let $b$ be a real number such that $\Rea z_1+b=\Rea z_2+b=0$ or $|z_1+b|=|z_2+b|$.
Then  $(A+b\Id)\mathbf x_i:\mathbf x_i=z_i+b$.
Now, it is simple geometry that $(A+b\Id)(\sqrt{1-t^2}\mathbf x_1+ t\mathbf x_2):(\sqrt{1-t^2}\mathbf x_1+ t\mathbf x_2)$
runs along $[z_1+b,z_2+b]_h$ for   $t\in[0,1]$.
This implies that  $A(\sqrt{1-t^2}\mathbf x_1+ t\mathbf x_2):(\sqrt{1-t^2}\mathbf x_1+ t\mathbf x_2)$
runs along $[z_1,z_2]_h$.

(b) $\CR(A^{\mathbb C})=\CR(A\oplus A)$; and if $\CR(A)$ is already $h$-convex, then it does not change.
\end{proof}
\end{lemma}

In fact, the only non-convex case, $\dim\mathfrak H=2$ is easy to overview:
\begin{lemma}\plabel{lem:2key}
Consider the real matrix
\begin{equation}A=\begin{bmatrix}a&b\\c&d\end{bmatrix}.\plabel{eq:Areal}\end{equation}

(a) For $A$ acting on $\mathbb R^2$,
\[
\CRext(A^{\mathbb{R}})=
\partial \Dbar\left(\tfrac{a+d}2+\tfrac{c-b}2\mathrm i,
\sqrt{\left(\tfrac{a-d}2\right)^2+\left(\tfrac{b+c}2\right)^2} \right)
\cup\partial \Dbar\left(\tfrac{a+d}2-\tfrac{c-b}2\mathrm i,
\sqrt{\left(\tfrac{a-d}2\right)^2+\left(\tfrac{b+c}2\right)^2} \right).
\]

(b) For $A$ acting on $\mathbb C^2$,
\begin{align}\notag
\CRext(A^{\mathbb C})=
& \Dbar\left(\tfrac{a+d}2+\tfrac{c-b}2\mathrm i,
\sqrt{\left(\tfrac{a-d}2\right)^2+\left(\tfrac{b+c}2\right)^2} \right)
\setminus
\intD\left(\tfrac{a+d}2-\tfrac{c-b}2\mathrm i,
\sqrt{\left(\tfrac{a-d}2\right)^2+\left(\tfrac{b+c}2\right)^2} \right)
\\
&\notag\cup \Dbar\left(\tfrac{a+d}2-\tfrac{c-b}2\mathrm i,
\sqrt{\left(\tfrac{a-d}2\right)^2+\left(\tfrac{b+c}2\right)^2} \right)
\setminus
\intD\left(\tfrac{a+d}2+\tfrac{c-b}2\mathrm i,
\sqrt{\left(\tfrac{a-d}2\right)^2+\left(\tfrac{b+c}2\right)^2} \right)
.
\end{align}
This is $\CRext(A^{\mathbb{R}})$ but with the components of $\mathbb C\setminus \CRext(A^{\mathbb{R}})$ disjoint from
$\mathbb R$ filled in.
\begin{proof}
(a) $\mathbb R^2$ can be identified $\mathbb C$.
One can check that for $|w|=1$,
\[\frac{Aw}w=\left(\frac{a+d}2+\frac{c-b}2\mathrm i\right )
+\frac1{w^2}\left(\frac{a-d}2+\frac{b+c}2\mathrm i\right ).\]
As an image, this yields a circle, the boundary of chiral disk.
The statement is an immediate consequence of this picture:
$\CR(A)$ is obtained, by conjugation-factorization, $\CRext(A)$ is obtained by conjugational doubling.
(In particular, the information in $\CR(A)$ is the same as in the principal disk; and the boundary of the chiral disk acts as an oriented conformal range.)

(b) This is a consequence of $\CRext(A^{\mathbb{C}})\cap\mathbb C^+=\CRext(A^{\mathbb{R}}\oplus A^{\mathbb{R}})\cap\mathbb C^+$.
\end{proof}
\end{lemma}
Thus,  for $\dim\mathfrak H=2$,  in terms of hyperbolic geometry,
$\CR(A^\mathbb C)$ yields points and disks around them; lines and distance bands around them, asymptotic points and corresponding horodisks.

\begin{proof}[Alternative proof to Lemma \ref{lem:normcompute}.]
$\CRext(A^{\mathbb R})$ is made of circles.
The farthest distance from the origin gives the norm;
and the closest distance from the origin gives the, say, co-norm.
These distances, however, can immediately be read off from
the center and the radius of the chiral disk (or the principal disk).
\end{proof}

Taking arbitrary Hilbert spaces $\mathfrak H$ now, the main point is that we see that complexification
does not lead to essential changes; thus general questions can be treated in the complex (or complexified) setting.
(Previously, this led to some awkwardness in Sections \ref{sec:ConformalRange} and \ref{sec:MagnusHilbert}.)

In the finite dimensional case
\begin{equation}\qquad\qquad\qquad\qquad\qquad\qquad\|A\|_2^-=\left(\|A^{-1}\|_2\right)^{-1}\qquad\qquad\qquad\qquad
(\dim\mathfrak H<\infty)\plabel{eq:finkey}\end{equation}
holds with  $\|A\|_2^-=0$ being equivalent to non-invertibility.
In particular, not only $\|A\|_2=\|A^*\|_2$ but also
\begin{equation}\qquad\qquad\qquad\qquad\qquad \qquad\|A\|_2^-=\|A^*\|_2^- \qquad\qquad\qquad\qquad\qquad
(\dim\mathfrak H<\infty)\plabel{eq:finadjkey}\end{equation}
holds in general. As consequence,
\begin{lemma}
\plabel{lem:CRfinadj}\plabel{lem:finAdj}
 Suppose that $\dim\mathfrak H<\infty$. Then
\[\CR(A)=\CR(A^*).\]
\begin{proof}
It is sufficient to prove $\CR(A^*)\subset\CR(A)$.
If $\dim_{\mathbb R} H\neq 2$; then we can use convexity:
Assume that $\CR(A^*)\nsubseteq\CR(A)$.
Then it is easy to see in CKB(P), that there is a point $P_{\mathrm{CKB(P)}}\in \DW_{\mathrm{CKB(P)}}(A^*)$
which is strictly separated by a line from the compact convex set $\DW_{\mathrm{CKB(P)}}(A)$.
Even more, we can assume that the line is not vertical.
The line corresponding in the PH model is then of shape $|z-\lambda|=r$ with some $\lambda\in\mathbb R$, $r>0$.
This mean depending on the position of the point and set that either
($\|A-\lambda\Id\|_2<r$ and $\|A^*-\lambda\Id\|_2>r$) or ($\|A-\lambda\Id\|_2^->r$ and $\|A^*-\lambda\Id\|_2^-<r$ holds).
This is a contradiction to the equality of norms and co-norms.
If $\dim_{\mathbb R} \mathfrak H= 2$ the we have to pass to $h$-convex closure.
There we see that $h$-convex closure contains the same information as the original conformal range.
Otherwise, the argument is the same. (Alternatively, the explicit shape can be examined.)
\end{proof}
\end{lemma}
\snewpage

Note that in the case of complex Hilbert spaces most of the results above
can be derived from the properties of the Davis--Wielandt shell; and in general,
the real case is not much different anyway.
We will not use the properties of Davis--Wielandt shell in our discussion,
but as exception we quote the following three statements
as corollaries of the information presented in Appendix \ref{app:DWrev}.

\begin{cor}
\plabel{lem:ellipdisk}
Suppose that $A$ acts on a 2-dimensional complex Hilbert space.
Then the conformal range of $A$ is the $h$-convex closure of a possibly degenerate $h$-ellipse
on the asymptotically closed hyperbolic plane (minus the $\infty$ point).

That is, in the asymptotically closed $\mathrm{CKB(P)}$ it is a possibly degenerate elliptical disk.
(The same applies in the $\mathrm{CKB}$ model but there the point $\infty_{\mathrm{CKB}}=(0,1)$
is to be avoided.)
\proofremark{
In the $\mathrm{CKB}$ / $\mathrm{CKB(P)}$ models, an $h$-ellipse (allowing an $h$-circle) is simply an ellipse in the interior
of the model. All other ``degenerate'' situations are limits of those.
}
\begin{proof}
$\DW_{*}^{\mathbb R}(A)$ is an orthogonal $h$-projection of $\DW_{*}(A)$.
(This is best to be visualized in the projective models, where this orthogonal $h$-projection
is represented by an ordinary orthogonal $h$-projection.)
In this case a possible degenerate $h$-tube (may be a $h$-line, a $h$-horosphere or asymptotic point) is projected.
In the projective models this can easily be seen to lead to the result indicated.
\end{proof}
\end{cor}
(The statement above was presented due to the elegance of the proof, but it can established in several other ways.)
\begin{cor}
\plabel{lem:CRvertex1}
Suppose $\dim\mathfrak H<\infty$.
Assume that $(x_0,z_0)$ is an angular boundary point of $\DW_{\mathrm{CKB(P)}}(A) $.
Let $\lambda=x_0+\mathrm i \sqrt{z_0-(x_0)^2}$.
Alternatively put: Assume that $\lambda$ is a $h$-angular boundary point of $\CR(A)$.

Then the generalized $\{\lambda,\bar \lambda\}$-eigenspace $V$ of $A$ is non-trivial,
$V\oplus V^\bot_{\mathfrak H}$ forms an $A$-invariant decomposition of $\mathfrak H$, and
$A|_V$ is normal.

In particular, $\DW_{\mathrm{CKB(P)}}(A) $ may have only  finite many angular vertices;
at any other boundary points the support lines are unique.
\begin{proof}
In the complex case is an immediate consequence of Corollary \ref{cor:Wielandt3}.
The real case follows by complexification.
\end{proof}
\end{cor}
(Again, other linear algebraic arguments are either more advanced or they would reproduce the argument of
used for the Davis--Wielandt shell.)
\begin{cor}
\plabel{lem:CRvertex2}
Assume that $(x_0,z_0)$, with $z_0>(x_0)^2$,  is an angular boundary point of $\overline{\DW_{\mathrm{CKB(P)}}(A)} $.
Let $\lambda=x_0+\mathrm i \sqrt{z_0-(x_0)^2}$.
Alternatively put: Assume that $\Ima \lambda>0$, and
$\lambda$ is an $h$-angular boundary point of $\overline{\CR(A)}$.

(i) If $(x_0,z_0)\in{\DW_{\mathrm{CKB(P)}}(A)}$. i. e.  $\lambda\in\CR(A)$, then
$\lambda$ or $\bar\lambda$ is an element of $\spec_p(A)\cup\spec_r(A)$ (cf. Discussion \ref{disc:decSpectrum}).

(ii) In any case, $\lambda$ or $\bar\lambda$ is an element of $\spec(A)$.

\begin{proof}
In the complex case is an immediate consequence of Lemma \ref{lem:Wielandt}, Lemma \ref{lem:Wielandt2}, Theorem \ref{thm:davisline}.
The real case follows by complexification.
\end{proof}
\end{cor}
\snewpage

The proof of Lemma \ref{lem:CRfinadj} leads to
\begin{disc}
\plabel{disc:DualView}
(The dual viewpoint -- PH model.)
As $\CR(A^{\mathbb C})$ is a $h$-convex set (the complexification is only for $\dim_{\mathbb R}\mathfrak H\neq 2$),
one can characterize its closure by duality, i. e. from its supporting $h$-lines.
Then the dual characterization yields
\begin{multline}
\overline{\CRext(A^{\mathbb C})} =\\=\{ z\in\mathbb C \,:\, |z-\lambda|\leq \|A-\lambda\Id\|_2,
|z-\lambda|\geq \|A-\lambda\Id\|_2^-\text{ for all }\lambda\in\mathbb R\}.
\plabel{eq:spectrange00}
\end{multline}
In fact, in this characterization above we have not even used the lineal $h$-lines,
but this is not a problem, the circular $h$-lines are sufficient, indeed.

Another way to express this is that the information contained in $\overline{\CRext(A^{\mathbb C})}$
is equivalent to the knowledge of $\|A-\lambda\Id\|_2$ and $\|A-\lambda\Id\|_2^-$ for all $\mathbb R$.

As for one more viewpoint:
The conditions in \eqref{eq:spectrange00} can be expressed as the complement of the conditions
\[|z-\lambda|^2\Id+\varepsilon\Id\leq (A^*-\lambda\Id)(A-\lambda\Id),\]
\[|z-\lambda|^2\Id-\varepsilon\Id\geq (A^*-\lambda\Id)(A-\lambda\Id),\]
for $\varepsilon>0$, $\lambda\in\mathbb R$. I. e. $\overline{\CRext(A^{\mathbb{C}})}$ contains the $z$ where none of the conditions above hold.
These conditions are linear inequalities for $\Id, A+A^*$, and $A^*A$.
(In the case of lineal relations the inequalities were for only $\Id$ and $ A+A^*$.)
\end{disc}
\begin{disc} (The dual viewpoint -- CKB(P) model.)
\plabel{disc:DualView2}
For this reason, it may be more convenient to use a projective model, in particular the CKB(P) model.
So, let  us consider
\[\DW_{\mathrm{CKB(P)}}^{\mathbb R}(A)=\left\{\left(\frac{\langle A\mathbf x,\mathbf x\rangle_{\mathrm{real}}}{\|\mathbf x\|_2},\frac{\|A\mathbf x\|_2}{\|\mathbf x\|_2}\right)\,:\,0\neq\mathbf x\in\mathfrak H\right\}.\]
Then (for $u,v,w\in\mathbb R$,)
\[u x_{\mathrm{CKB(P)}}+v y_{\mathrm{CKB(P)}} +w\geq0\qquad\text{for all}\qquad
( x_{\mathrm{CKB(P)}},  y_{\mathrm{CKB(P)}})\in \DW_{\mathrm{CKB(P)}}^{\mathbb R}(A)\]
is equivalent to
\[u \langle A\mathbf x,\mathbf x\rangle_{\mathrm{real}} +v \|A\mathbf x\|_2  +w\|A\mathbf x\|_2\geq0\qquad\text{for all}\qquad
0\neq\mathbf x\in\mathfrak H,\]
i. e.
\[u\frac{A+A^*}{2}+vA^*A+w\Id\geq0.\]

Let
$\Dual_{\mathrm{CKB(P)}}^{\mathbb R}(A)$
be the set of all supporting lines $(u:v:w)'$ (using projective line coordinates now).
Then we find that
\begin{equation}
\overline{\conv \DW_{\mathrm{CKB(P)}}^{\mathbb R}(A)}=\mathcal X( \Dual_{\mathrm{CKB(P)}}^{\mathbb R}(A)  )
\plabel{eq:DWRrec}
\end{equation}
when $\mathcal X$ is the operation given by
\begin{multline*}
\mathcal X(L):=\bigcap_{\phi\in[0,2\pi)}\{( x_{\mathrm{CKB(P)}},  y_{\mathrm{CKB(P)}})\in\mathbb R^2\,:\,\\
(\cos\phi)x_{\mathrm{CKB(P)}}+(\sin\phi)  y_{\mathrm{CKB(P)}}+\inf\{s\,:\, (\cos\phi : \sin\phi :s)'\in L\}\geq0.
\}\end{multline*}
\snewpage

Note, however, that the variant of  \eqref{eq:DWRrec},
\begin{equation}
\overline{\conv \DW_{\mathrm{CKB(P)}}^{\mathbb R}(A)}=\mathcal X( \Dual_{\mathrm{CKB(P)}}^{\mathbb R,\mathrm{var}}(A)  )
\plabel{eq:DWRrecvar}
\end{equation}
 remains valid
even if $\Dual_{\mathrm{CKB(P)}}^{\mathbb R}(A)$ extended by an arbitrary set of lines intersecting $ \DW_{\mathrm{CKB(P)}}^{\mathbb R}(A)$,
and after that only  the lines parallel to a dense set of directions selected to be in
$ \Dual_{\mathrm{CKB(P)}}^{\mathbb R,\mathrm{var}}(A)$.
(In the first step $\inf$'s are not affected, in the second step some $\inf$'s are set to $-\infty$ but still sufficiently many directions are left.)
In particular, the operation $\mathcal X$ can be defined more even more flexibly, but it will be sufficient for our purposes.

The connection to picture with (co-)norms is straightforward.
If $N$ denotes the square of the norm or co-norm, $\lambda\in\mathbb R$, then
\[|(x+\mathrm iy)-\lambda|^2\lesseqgtr N(A-\lambda\Id)\]
is equivalent to
\[-2x_{\mathrm{CKB(P)}}+y_{\mathrm{CKB(P)}}+\lambda^2-N(A-\lambda\Id)\lesseqgtr0.\]
We find that
\[\Dual_{\mathrm{CKB(P)}}^{\mathbb R,\mathrm{ncn}}(A)=\{ (-2\lambda: 1 : \lambda^2-N(A-\lambda\Id) )'  : \lambda\in\mathbb R, \text{ choice for }N\}\]
is almost the same as $ \Dual_{\mathrm{CKB(P)}}^{\mathbb R}(A)$ but the vertical lines
(i. e. supporting lines of the real supremum and infimum) are omitted.
In particular it fits to \eqref{eq:DWRrecvar}.

In the setting of CKB(P), it is particulary nice that the slopes of the supporting lines are directly connected to the
poles $\lambda$.
\end{disc}
\begin{disc} ($C^*$-algebras.)
In particular, $\overline{\CRext(A^{\mathbb{C}}})$ is characterized in purely $C^*$-algebraic terms.
As the valid linear inequalities for $\Id, A+A^*, A^*A$ form a convex set, it would be possible
to develop a theory for
\[\ccCRext(A):=\overline{\CRext(A^{\mathbb{C}})},\]
or
\[\ccCR(A):=\overline{\CR(A^{\mathbb{C}})},\]
etc., in purely $C^*$-algebraic terms.
It is, however, more economical to use the Gelfand--Naimark representation theorem, which
is $\leq$-compatible; showing that the closed $h$-convex version of the conformal range
can be defined for abstract $C^*$-algebras invariantly.
\end{disc}
\snewpage
\begin{disc} (The dual viewpoint -- the finite dimensional case, PH.)
In the finite dimensional case we spell out
\begin{lemma}\plabel{lem:finchar}
 If $\dim \mathfrak H<\infty$, then
\begin{multline}
{\CRext(A^{\mathbb C})}=\{ z\in\mathbb C \,:\, |z-\lambda|\leq \|A-\lambda\Id\|_2\text{ for all }\lambda\in\mathbb R;\text{ and }\\
|z-\lambda|\geq \|(A-\lambda\Id)^{-1}\|_2^{-1}\text{ for all }\lambda\in\mathbb R\setminus \spec(A) )\}
\plabel{eq:spectrange35}
\end{multline}
\begin{proof}
As $\CRext(A)$ is (compact) closed, the dual viewpoint and  \eqref{eq:finkey} immediately implies the characterization.
\end{proof}
\end{lemma}
(We can also say that the spectrum is finite, thus there are no point and residual lacunas.)
\end{disc}
\begin{disc} (The dual viewpoint -- the finite dimensional case, CKB(P).)
\plabel{disc:DualView3}

Assume that $(u:s:w)'$ is a supporting line such that
$ux+sy+w\geq 0$ for all $(x,y)\in \DW_{\mathrm{CKB}}^{\mathbb R} (A)$ but with equality somewhere.
this is to say that \[u\Rea\langle A\mathbf x,\mathbf x\rangle+s\|A\mathbf x\|_2 +w\|\mathbf x\|_2\equiv
u\mathbf x^*\frac{A+A^*}2\mathbf x +s\mathbf x^* A^* A\mathbf x+w\mathbf x^*\mathbf x\geq 0\] for all $\mathbf x\neq0$
but with at least with one $\mathbf x_0\neq0$ having equality. At such an $\mathbf x_0$ the value $0$
is an extremal value, thus, by differentiation,
\begin{equation}
u\frac{A+A^*}2\mathbf x_0 +s A^* A\mathbf x_0+w\mathbf x_0=0.
\plabel{eq:van1}
\end{equation}
That is
\begin{equation}
0\neq\mathbf x_0\in \ker\left( u\frac{A+A^*}{2}+sA^*A+w\Id\right)
\plabel{eq:van2}
\end{equation} holds.
Consequently
\begin{equation}
\det\left(u\frac{A+A^*}{2}+sA^*A+w\Id\right)=0.\plabel{eq:Kip}
\end{equation}
(This is in analogue of
\[\det\left(u\frac{A+A^*}{2}+v\frac{A-A^*}{2\mathrm i}+ w\Id\right)=0\]
from Kippenhahn \cite{Ki1} / \cite{Ki2}.)
Conversely, if \eqref{eq:Kip} holds ($(u,s,w)\neq(0,0,0)$), then with some $\mathbf x_0\neq 0$ \eqref{eq:van2} / \eqref{eq:van1} hold.
According to this, $(u:s:w)'$ may not be a supporting line but passes  through a point of $\DW_{\mathrm{CKB(P)}}^{\mathbb R}(A)$
(namely, the image of $\mathbf x_0$).

Hence
\[\Dual_{\mathrm{CKB(P)}}^{\mathbb R,\mathrm{alg}}(A)=\left\{(u: s: w)'\,:\, \det\left(u\frac{A+A^*}{2}+sA^*A+w\Id\right)=0  \right\}\]
 fits to \eqref{eq:DWRrecvar}, again.

If we obtain supporting lines from norms and co-norms, then we can use the fact that
$N(A-\lambda\Id)$ is a maximal or minimal eigenvalue of  $(A-\lambda\Id)^*(A-\lambda\Id)$
(all eigenvalues are real, in fact, nonnegative).
Thus we can define
\[\Dual_{\mathrm{CKB(P)}}^{\mathbb R,\mathrm{ncn},\mathrm{alg}}(A)=
\left\{(-2\lambda: 1: \lambda^2-\nu)'\,:\, \det\left(\nu\Id- (A^*-\lambda\Id)(A-\lambda\Id)\right)=0  \right\}.\]
Again, we may have added some intermediate lines between norm and co-norm lines, and we have none of vertical lines,
but it still fits into  \eqref{eq:DWRrecvar}.

In particular, we find that the real homogeneous (in $u,s,w$) polynomial
\[K^{\CR}_A(u,s,w)=\det\left(u\frac{A+A^*}{2}+sA^*A+w\Id\right)\]
(naturally normalized by $K_A(0,0,1)=1$)  and / or the real polynomial
\[F^{\CR}_A(\lambda,\nu)=\det\left(\nu\Id- (A^*-\lambda\Id)(A-\lambda\Id)\right)\]
determine the conformal range.
(Equivalence is by
$F^{\CR}_A(\lambda,\nu)=K^{\CR}_A(2\lambda,-1,\nu-\lambda^2)$,
and, formally, $K^{\CR}_A\left(-\dfrac us,-1,-\dfrac ws\right)=F^{\CR}_A\left( -\dfrac{u}{2s},-\dfrac ws+\dfrac{u^2}{4s^2}\right) $.
For real arguments, having determinants of self-adjoint operators, the value are real.
Consequently, the polynomials are real.)
Practical aspects of obtaining the conformal range from the polynomials this will be considered later.

We can note, however, that angular vertices in $\DW_{\mathrm{CKB}}^{\mathbb R} (A)$
will result dual segments (in line space) in $\Dual_{\mathrm{CKB(P)}}^{\mathbb R}(A)$,
thus, by its linear algebraic nature, dual lines in $\Dual_{\mathrm{CKB(P)}}^{\mathbb R,\mathrm{alg}}(A)$.
That is a homogeneous component $ux_0+sy_0+w$.

It is natural to call the polynomial $K^{\CR}_A(u,s,w)$ as the algebraic conformal range,
and the polynomial
\[K_A^{\mathrm W}(u,v,w)=\det\left(u\frac{A+A^*}{2}+v\frac{A-A^*}{2\mathrm i}+ w\Id\right)\]
as the algebraic numerical range.
In the finite dimensional case the algebraic ranges are finer invariants than the geometrical ones.
In contrast to
\[K^{\mathrm W}_A(u,v,w)=K^{\mathrm W}_{A^*}(u,-v,w);\]
due to the identity $\det(\nu\Id-ST)=\det(\nu\Id-TS)$,
we have
\[F^{\CR}_A(\lambda,\nu)=F^{\CR}_{A^*}(\lambda,\nu),\]
or expressed otherwise,
\[K^{\CR}_A(u,s,w)=K^{\CR}_{A^*}(u,s,w).\]
This can be considered as a refinement of Lemma \ref{lem:CRfinadj}.
Unsurprisingly, the conformal invariance properties also extend to $K^{\CR}$ in appropriate form.

In the case of the algebraic numerical range, according to Kippenhahn's observation,
it is easy to see that the eigenvalues of $A$ are given by the curve theoretic foci
(``as obtained from the equation written in line coordinates'')
\[K^{\mathrm W}_A(1,\mathrm i,-\lambda)\equiv \det(A -\lambda\Id)=0.\]
(Curve theoretic foci are an old idea of Pl\"ucker, about which do not have to be concerned,
as the computation is on a purely algebraic level here.)
In the case of the algebraic conformal range,
$\lambda$ or $\bar\lambda$ is an eigenvalue of $A$ if and only if
\[K^{\CR}_A(2\lambda,-1,-\lambda^2)\equiv F^{\CR}_A(\lambda,0)\equiv \det\left(-(A^* -\lambda\Id)(A -\lambda\Id)\right)=0.\]
We may call these as the ``hyperbolic foci'' (regarding the given setting).
\end{disc}
\snewpage
\begin{disc}[The decomposition of the spectrum in general]
\plabel{disc:decSpectrum}
It is well-known that the spectrum of the linear operator $A$ can be decomposed as
\[\spec(A)=\spec_p(A)\cup\spec_c(A)\cup\spec_r(A),\]
i. e. point, continuous, and residual spectrum.
Following Davis \cite{D1} we define,

$\lambda\in\spec_p(A)$ iff  $\ker A-\lambda\Id\neq 0$;

$\lambda\in\spec_c(A)$ iff  $A-\lambda\Id$ restricted from $(\ker A-\lambda\Id)^\bot_{\mathfrak H}$ to
$\overline{\mathrm{im}\,  A-\lambda\Id}$ is not invertible;

$\lambda\in\spec_r(A)$ iff  $\overline{\mathrm{im}\, A-\lambda\Id}\neq \mathfrak H$.
This is not compatible to the vast majority of the literature, where the spectral decomposition considered is
\[\spec(A)=\boldsymbol\sigma_p(A)\,\dot\cup\,\boldsymbol\sigma_c(A)\,\dot\cup\,\boldsymbol\sigma_r(A)\]
with $\boldsymbol\sigma_p(A)=\spec_p(A)$, $\boldsymbol\sigma_c(A)=\spec_c(A)\setminus \spec_p(A)\setminus \spec_r(A)$, $\boldsymbol\sigma_r(A)=\spec_r(A)\setminus \spec_p(A)$.
The name compression spectrum is used for our residual spectrum, $\boldsymbol\sigma_{cp}(A)=\spec_r(A)$.

However, the notation according to Davis has the advantage that $\spec_p(A^*)=\spec_r(A)$, $\spec_c(A^*)=\spec_c(A)$, $\spec_r(A^*)=\spec_p(A)$;
furthermore, the approximate point spectrum is given conveniently as $\boldsymbol\sigma_{ap}(A)=\spec_p(A)\cup \spec_c(A)$.
\end{disc}
\begin{disc}[The spectrum]
\plabel{disc:Spectrum}
Now, Lemma \ref{lem:spectrange} is but an immediate consequence of the relations,
\[\spec_p(A)\subset\CRext(A^\mathbb C) \]
\[\spec_c(A)\subset\overline{\CRext(A^\mathbb C)} \cap \overline{\CRext((A^\mathbb C)^*)}\]
\[\spec_r(A)\subset\CRext((A^\mathbb C)^*) \]
which, in turn, are rather immediate from the definitions.
(Note that $\CRext$ is conjugation-invariant.)
Furthermore,
\[\spec_p(A)\cap\mathbb R=\CRext(A^\mathbb C) \cap\mathbb R.\]
and
\begin{equation}
\spec(A)\cap\mathbb R=\left(\overline{\CRext(A^\mathbb C)} \cup \overline{\CRext((A^\mathbb C)^*)}\right)\cap \mathbb R.\plabel{eq:dane}
\end{equation}
and
\[\spec_r(A)\cap\mathbb R=\CRext((A^\mathbb C)^*) \cap\mathbb R.\]
(These are sufficient to check for $0$ as a spectrum point.)
\end{disc}

The generalization of Lemma \ref{lem:spectrealtemp} is
\begin{lemma}\plabel{lem:spectrange15}
Let $\smc(\overline{\CRext(A)})$ denote the simply connected closure of $\overline{\CRext(A)}$, i.~e.~the
complement of the infinite component of $\mathbb C\setminus \overline{\CRext(A)}$.
Then
\begin{equation} \spec(A)\subset\smc(\overline{\CRext(A)}).\plabel{eq:spectrange15} \end{equation}
\begin{proof}
Complexification does not change $\smc(\overline{\CRext(A)})$, thus we can assume that the complex case.
Indeed, indirectly, suppose that $C$ is a polygonal chain from $\infty$ to $\xi$ in the
complement $\mathbb C\setminus \overline{\CRext(A)}$.
It can be assumed that $\xi$ is the first and (last) element of $C$ such that $A-\xi\Id$ is not invertible.
According to \eqref{eq:wrav}, the inverse $(A-\lambda\Id)^{-1}$ is bounded by $\dist(C,\CRext(A) )^{-1}$ for $\lambda\in C\setminus\{\xi\}$.
Hence, its derivative $(A-\lambda\Id)^{-2}A$ is bounded by $\dist((-\infty,0],\CRext(A) )^{-2}\|A\|_2$  for $\lambda\in C\setminus\{\xi\}$.
This, however, implies that the inverse extends to $A-\xi\Id$; which is a contradiction.
\end{proof}
\end{lemma}
An alternative line of argument in this direction, which yields a bit more is as follows:
\begin{lemma}\plabel{lem:smcReal}
\[\smc(\overline{\CRext(A)} )\cap \mathbb R=\overline{\CRext(A)} \cap \mathbb R=\conv(\spec(A)\cap\mathbb R )\]
($\conv$ means ordinary convex hull); and
\begin{align}
\smc(\overline{\CRext(A)} )=\{ z\in\mathbb C \,:&\, |z-\lambda|\leq \|A-\lambda\Id\|_2\text{ and }\plabel{eq:smcchar}\\
&|(\nu-z)^{-1}-\lambda|\leq \|(\nu\Id-A)^{-1}-\lambda\Id\|_2\notag\\
&\text{ for all }\lambda\in\mathbb R\text{ and }\nu\in\mathbb R\setminus \conv(\spec(A)\cap\mathbb R )\}.\notag
\end{align}
\begin{proof}
The first formula can derived directly as in Lemma \ref{lem:spectrealtemp}.

Regarding the second:
Let us define
\[\ccCRext_\infty(A)=\{ z\in\mathbb C \,:\, |z-\lambda|\leq \|A-\lambda\Id\|_2\text{ for all }\lambda\in\mathbb R\}.\]
Now, $\ccCR_\infty(A)$ is not the same as $\ccCR(A)$ because it also contains every point ``below''.
From hyperbolic point of view, however,  $\ccCR_\infty(A)$ is the $h$-convex view of  $\ccCR(A)$ from $\infty$.
(This particularly transparent in the Cayley--Klein--Beltrami model.)
However, any point $\nu$ of $\mathbb R\setminus \smc(\overline{\CRext(A)} )$
can be moved to infinity of $\infty$
using the fractional linear transformation $x\mapsto\frac1{x-\lambda}$.
Thus $\ccCRext_\infty$ applied to $(A-\nu\Id)^{-1}$, and transformed back, leads to the same information as in
$\ccCRext_\nu(A)$, which is  the $h$-convex view of  $\ccCRext(A)$ from $\nu$.
There,
\[\ccCRext_\nu(A)=\{ z\in\mathbb C \,:\,|(\nu-z)^{-1}-\lambda|\leq \|(\nu\Id-A)^{-1}-\lambda\Id\|_2\text{ for all }\lambda\in\mathbb R\}.\]
Intersecting these $h$-convex views (again, the viewpoint of the Cayley--Klein--Beltrami model may be useful),
we obtain $\smc(\ccCRext(A))\equiv \smc(\overline{\CRext(A)} )$.
\end{proof}
\end{lemma}

\begin{lemma}\plabel{lem:smcadj}
\[\smc(\overline{\CRext(A)} )=\smc(\overline{\CRext(A^*)}\]
\begin{proof}
Note $\conv(\spec(A)\cap\mathbb R)= \conv(\spec(A^*)\cap\mathbb R)$ because of $\eqref{eq:dane}$.
Then the statement follows from \eqref{eq:smcchar}, where  the same norms are yielded for $A$ and $A^*$.
(In fact, in terms of the previous proof, $\ccCRext_\nu(A)=\ccCRext_\nu(A^*)$ as long as $\nu\in (\mathbb R\cup\{\infty\})\setminus\spec(A)$.)
\end{proof}
\end{lemma}
\snewpage

This situation of Lemma \ref{lem:smcadj} can be analysed further:
\begin{disc}[Lacunas]
\plabel{disc:lacunas}
Consider now $\overline{\CR((A^\mathbb C)^*)}$.
Then $\upper\setminus \overline{\CR((A^\mathbb C)^*)}$
has countable many connected components, and due to convexity, in bijection to open segments on the asymptotic boundary.
Thus
\[\upper\setminus\overline{\CR(A^\mathbb C)}=L_\infty\,\dot\cup\,\bigcup_I^{\boldsymbol.} L_{I},\]
we have a disjoint union of lacunas; one infinite and some other ones corresponding to some
finite open intervals $I\subset\mathbb R$.
The situation is pretty much the same if consider $\mathbb C\setminus \overline{\CRext(A^\mathbb C)}$,
but the lacunas are conjugationally symmetric.
\end{disc}
\begin{lemma}
\plabel{lem:lacunas}
Assume that
\[\mathbb C\setminus\overline{\CRext(A^\mathbb C)}=L^{\ext}_\infty\,\dot\cup\,\bigcup_I^{\boldsymbol.} L^{\ext}_{I}\]
is a decomposition to lacunas.
Then, for each $L^{\ext}_I$,
\begin{equation}
L^{\ext}_I\cap(\spec_p(A)\cup \spec_c(A))=\emptyset
\plabel{eq:lac1}
\end{equation}
holds.
Also, for each $L^{\ext}_I$,
\begin{equation}
L^{\ext}_I\cap  \spec_r(A)=\emptyset
\qquad\text{or}\qquad
L^{\ext}_I\subset  \spec_r(A)
\plabel{eq:lac2}
\end{equation}
holds (this is a dichotomy, depending on $I$; in case of $L^{\ext}_\infty$ obviously with the first case holds).
\begin{proof}
\eqref{eq:lac1} is immediate, we have to prove only \eqref{eq:lac2}.
If $z\in L^{\ext}_I\setminus \spec_r(A)$, then $z\in\mathbb C\setminus\spec(A)$.
Then any other point $w$ in the lacuna $L^{\ext}_I$ can be connected to $z$ by  polygonal chain still in the lacuna.
Then, as in the proof of Lemma \ref{lem:spectrange15}, we also obtain $w\in\mathbb C\setminus\spec(A)$, in particular, $w\notin\spec_r(A)$
\end{proof}
\end{lemma}
It is reasonable to call the lacuna $L^{\ext}_I$ a residual lacuna of $A$ if $L^{\ext}_I\subset  \spec_r(A)$ holds,
and we can call other lacunas $L^{\ext}_I$ as non-spectral lacunas of $A$.

We can also consider the same picture for $A^*$ with lacunas $\tilde L^{\ext}_J$.
According to Lemma \ref{lem:smcadj}, $L^{\ext}_\infty=\tilde L^{\ext}_\infty$; the infinite lacunas are the same.
It is reasonable to call a residual lacuna $\tilde L^{\ext}_J$ of $A$ as a point lacuna of $A$.
\snewpage

\begin{prop}
\plabel{pr:lacunas}
(a)
The closed $h$-convex set
\[\overline{\CRext(A^\mathbb C)} \cap \overline{\CRext((A^\mathbb C)^*)}\]
is non-empty.

(b) Consider the lacunar decomposition
\[\mathbb C\setminus( \overline{\CRext(A^\mathbb C)} \cap \overline{\CRext((A^\mathbb C)^*)} ) =\hat L^{\ext}_\infty\,\dot\cup\,\bigcup_K^{\boldsymbol.} \hat L^{\ext}_{K}\]
Then, for each $\hat L^{\ext}_{K}$ one of the following three cases holds:
\[\text{ $\hat L^{\ext}_{K}$ is a non-spectral lacuna of $A$ (non-spectral lacuna of $A^*$) , $\subset \mathbb C\setminus(\spec(A) )$}\]
or
\[\text{ $\hat L^{\ext}_{K}$ is a residual lacuna of $A$ (point lacuna of $A^*$) , $\subset \spec_r(A)\setminus(\spec_p(A) \cup \spec_c(A) )$}\]
or
\[\text{ $\hat L^{\ext}_{K}$ is a point lacuna of $A$ (residual lacuna of $A^*$) , $\subset \spec_p(A)\setminus(\spec_r(A) \cup \spec_c(A) )$.}\]

(The set of such lacunas is the union of the sets  of all lacunas of $A$ and $A^*$.)
$\hat L^{\ext}_\infty$ is equal to the infinite lacunas of $A$ and $A^*$.
\begin{proof}
(a)
Otherwise, $\overline{\CR(A^\mathbb C)}$ and $\overline{\CR((A^\mathbb C)^*)}$ could be separated by a $h$-line which is actually a
semicircle, which quickly leads to contradiction in terms of norms.
Another way to argue is based on the observation $\partial  \smc(\overline{\CR(A^\mathbb C)}) \subset \overline{\CR(A^\mathbb C)}$
and Lemma \ref{lem:smcadj}.

(b)
Point lacunas are subsets of  $\spec_p(A)\setminus \spec_r(A)$, and
residual lacunas are subsets of  $\spec_r(A)\setminus \spec_p(A)$.
Thus lacunas of $A$ and $A^*$ are pairwise disjoint with the exception of the
possibility that a non-spectral lacuna $L^{\ext}_I$ of $A$ and a non-spectral lacuna $\hat L^{\ext}_J$ of $A^*$ intersects.
In that case, however,   $\partial L^{\ext}_I\cap \partial\hat L^{\ext}_J\subset \overline{\CRext(A^\mathbb C)}$,
and due to convexity of  $\overline{\CRext(A^\mathbb C)}$, we conclude that the supporting real intervals $I$ and $J$ intersect.
If $\lambda\in I\cap J$, then it is non-spectral, and one can apply the fractional linear transformation $x\mapsto\frac1{x-\lambda}$
and Lemma \ref{lem:smcadj} to see that $L^{\ext}_I=\hat L^{\ext}_J$.
On the other hand, by the same argument, any non-spectral $\lambda\in\mathbb R$ belongs to non-spectral lacuna of $A$ and also of $A^*$,
this establishes the trichotomy.
\end{proof}
\end{prop}

\begin{cor}
\plabel{cor:lacunas}
\begin{multline*}
\left(\CRext(A^\mathbb C) \setminus \overline{\CRext((A^\mathbb C)^*)}\right)
\cup
\left(\CRext((A^\mathbb C)^*) \setminus \overline{\CRext(A^\mathbb C)}\right)
\subset
\spec(A)\subset  \\
\left(\CRext(A^\mathbb C) \setminus \overline{\CRext((A^\mathbb C)^*)}\right)
\cup
\left(\overline{\CRext(A^\mathbb C)} \cap \overline{\CRext((A^\mathbb C)^*)}\right)
\cup
\left(\CRext((A^\mathbb C)^*) \setminus \overline{\CRext(A^\mathbb C)}\right)
\\\subset\smc(\overline{\CRext(A)} )= \smc(\overline{\CRext(A^*)} ).  \plabel{eq:spectrange2}
\end{multline*}
\begin{proof}
This is immediate from the previous theorem.
\end{proof}
\end{cor}
\snewpage
\begin{example}\plabel{exa:CR}
(a) If $\mathfrak H=\ell^2(\mathbb Z;\mathbb C)$, and let $A$ be the unilateral shift $A\mathbf e_n=\mathbf e_{n+1}$.
Then
\[\CRext(A)=\partial\Dbar(0,1)\setminus \{-1,1\},\]
\[\CRext(A^*)=\partial\Dbar(0,1)\setminus \{-1,1\},\]
\[\spec(A)=\partial\Dbar(0,1),\]
\[\smc(\overline{\CRext(A)})=\Dbar(0,1).\]
In this latter case, beyond the infinite lacuna, $A$ has a single non-spectral lacuna.

(b)
Let $\mathfrak H=\ell^2(\mathbb N;\mathbb C)$, and let $A$ be the unilateral shift $A\mathbf e_n=\mathbf e_{n+1}$.
Then
\[\CRext(A)=\partial\Dbar(0,1)\setminus \{-1,1\},\]
\[\CRext(A^*)=\Dbar(0,1)\setminus \{-1,1\},\]
\[\spec(A)=\Dbar(0,1),\]
\[\smc(\overline{\CRext(A)})= \Dbar(0,1).\]
In this latter case, beyond the infinite lacuna, $A$ has a single residual lacuna
(which is a point lacuna of $A^*$).\qedexer
\end{example}

It is easy to see
that any $h$-convex bounded closed set of $\upper$ can be obtained as $\overline{\CR(A^\mathbb C) }$ of a bounded linear operator $A$ on a separable Hilbert space;
and any $h$-convex bounded set of $\upper$ can be obtained as $\CR(A^\mathbb C)$ of a bounded linear operator $A$ on a possibly inseparable Hilbert space.
In the finite case, $\CR(A^\mathbb C)$ has a more specific character:
\snewpage

As we have seen, in theory, we can determine the closure of the conformal range using norms and co-norms.
Now the standard differential geometric construction is as follows.

\begin{theorem}
\plabel{th:envel}
(The standard $h$-horo-translation induced tracing.)

Let $N(\cdot)$ denote the square of the norm or the co-norm.
Then $\partial\CR(A)$ is the enveloping curve of the semicircles
\[(x-\lambda)^2+y^2=N(A+\lambda\Id),\qquad\qquad y\geq0.\]
This curve can be computed as
\[\lambda\mapsto E^A(\lambda)=\left(\lambda-\frac12 \frac{\mathrm d N(A-\lambda\Id)}{\mathrm d\lambda}\right)+
\mathrm i  \sqrt{ N(A-\lambda\Id)- \left(\frac12\frac{\mathrm d N(A-\lambda\Id)}{\mathrm d\lambda}\right)^2}.\]
The norm produces  the upper part, the co-norm produces the lower part.
(The joins correspond to $\lambda=\pm\infty$.)
The expression is defined almost everywhere (restricted to $N(A-\lambda\Id)\neq0$), but discontinuities can occur, which
should be bridged by $h$-segments.

In the CKB(P) model we have to consider the enveloping curves of the lines
\[-2\lambda x +y+\lambda^2-N(A-\lambda\Id)=0.\]
(Note that the steepness of the lines is $2\lambda$.)
The enveloping curves can be written as
\[\lambda\mapsto E_{\mathrm{CKB(P)}}^A(\lambda)=\left(\lambda-\frac12 \frac{\mathrm d N(A-\lambda\Id)}{\mathrm d\lambda},\lambda^2-\lambda\frac{\mathrm d N(A-\lambda\Id)}{\mathrm d\lambda}+ N(A-\lambda\Id)\right).\]

\begin{proof}
In the  CKB(P) model, the standard (but not smooth) differential geometric arguments apply.
Then, the result can be transcribed to PH.
\end{proof}
\end{theorem}
It is easy to see that discontinuities in the enveloping curves correspond to $1$-di\-men\-sio\-nal faces of the conformal range (in $h$-sense),
and plateaus correspond to angular vertices or asymptotic points of the closed conformal range (in $h$-sense).
The `values $\lambda=\pm\infty$' would correspond to the real supremum and infimum of $A$
(that is the vertical supporting lines in $\upper$); those might be singular or not.
Nevertheless, it is sufficient to know the value of the enveloping upper curves on a dense subset once restricted to $N(A-\lambda\Id)\neq0$.
The enveloping curves go counterclockwise, i. e. to the left on the upper part, and to the right on the lower part.

The norm branch is the upper boundary curve and corresponds to the $h$-convex view from $\infty$, which is part of the boundary of the infinite lacuna.
\snewpage

\begin{example}\plabel{ex:envel1}
Generally, in the real case, $A= \tilde a\Id +\tilde b\tilde I+\tilde c\tilde J+\tilde d\tilde K$
yields
\[N(A-\lambda\Id_2)=(\tilde a-\lambda)^2+\tilde b^2+\tilde c^2+\tilde d^2\pm2\sqrt{(\tilde a-\lambda)^2+\tilde b^2}\sqrt{\tilde c^2+\tilde d^2}\]
(sign $+$ for the  norm branch, sign $-$ for the co-norm branch).
The enveloping curve in the $\mathrm{CKB(P)}$ model is given by
\[ E_{\mathrm{CKB(P)}}^A(\lambda)=
\left(a\pm\frac{(\tilde a-\lambda)\sqrt{\tilde c^2+\tilde d^2}}{\sqrt{(\tilde a-\lambda)^2+\tilde b^2}},
\tilde a^2+\tilde b^2+\tilde c^2+\tilde d^2
\pm\frac{2(\tilde a^2-\lambda \tilde a+\tilde b^2)\sqrt{\tilde c^2+\tilde d^2}}{\sqrt{(\tilde a-\lambda)^2+\tilde b^2}}\right).\]
Transcribed to the $\mathrm{PH}$ model (that is essentially complex numbers), it yields
\[E^A(\lambda)=\left(\tilde a\pm\sqrt{\tilde c^2+\tilde d^2}\frac{(\tilde a-\lambda)}{
\sqrt{(\tilde a-\lambda)^2+\tilde b^2}}
,
\left|\tilde b+\sqrt{\tilde c^2+\tilde d^2}\frac{\tilde b}{
\sqrt{(\tilde a-\lambda)^2+\tilde b^2}}
\right|
\right).\]

But, e. g., for $A=\begin{bmatrix}1&\\&-1\end{bmatrix}$, i. e.
 for $\tilde a=\tilde b=\tilde d=0$, $\tilde c=1$, this degenerates to $\lambda\mapsto(-\sgn \lambda,0)$ in the norm case,
 and  to $\lambda\mapsto(\sgn \lambda,0)$ in the co-norm case; and almost the whole
conformal range (which is the $h$-segment connecting $-1$ to $1$) comes from a discontinuity.
The constant plateaus in the enveloping curve, however, here, correspond to some asymptotic points
(which are also angular boundary points here) in the closed conformal range.
\qedexer
\end{example}

Variants of the enveloping construction above are possible:
\snewpage
\begin{theorem}
\plabel{th:envelrot}
(An $h$-rotation induced tracing.)

In the CKB model, the $\partial\DW_{\mathrm{CKB}}^{\mathbb R} (A)$ is the enveloping curve of the lines
\begin{equation}
-(\sin\omega)x_{\mathrm{CKB}}+(\cos \omega)y_{\mathrm{CKB}}
-\frac{\left(\left\|\dfrac{(\cos\frac\omega2)A-(\sin\frac\omega2)\Id}{(\sin\frac\omega2)A+(\cos\frac\omega2)\Id}\right\|_2\right)^2-1
}{\left(\left\|\dfrac{(\cos\frac\omega2)A-(\sin\frac\omega2)\Id}{(\sin\frac\omega2)A+(\cos\frac\omega2)\Id}\right\|_2\right)^2+1}=0.
\plabel{eq:ctran1}
\end{equation}

It is given analytically, almost everywhere, by
\begin{equation}
\omega \mapsto \widehat E^A_{\mathrm{CKB}}(\omega)=
\left(
\frac{\sin\omega -(\sin\omega)\widehat N^2-2(\cos\omega)\dfrac{\mathrm d\widehat N}{\mathrm d\omega}}{(1+\widehat N)^2}
,
\frac{-\cos\omega +(\cos\omega)\widehat N^2-2(\sin\omega)\dfrac{\mathrm d\widehat N}{\mathrm d\omega}}{(1+\widehat N)^2}
\right),
\plabel{eq:ctran2}
\end{equation}
where `$\widehat N$' is an abbreviation for $\left(\left\|\dfrac{(\cos\frac\omega2)A-(\sin\frac\omega2)\Id}{(\sin\frac\omega2)A+(\cos\frac\omega2)\Id}\right\|_2\right)^2$.
\begin{proof}
\eqref{eq:ctran1} simply specifies the tangent lines with outward pointing normal vector $(-\sin\omega,\cos\omega)$;
while \eqref{eq:ctran2} is the analytical solution.
\end{proof}
\end{theorem}
The latter is scheme is, however, less practical than the first one.
\begin{example}\plabel{ex:horoenv}
Let us consider the case of $A=\begin{bmatrix}0&1\\&0\end{bmatrix}$.
As we have seen, $\partial\CR(A)$ is the $(0,\frac12)$ centered circle with radius $\frac12$.
In this case,
\[N(A-\lambda\Id_2)=\lambda^2+\frac12\pm\frac{\sqrt{1+4\lambda^2}}2;\]
the standard enveloping construction yields
\[ E_{\mathrm{CKB(P)}}^A(\lambda)=
\left(\mp\frac{\lambda}{\sqrt{1+4\lambda^2}},\frac12\pm \frac12\frac{1}{\sqrt{1+4\lambda^2}}\right),\]
and
\[ E^A(\lambda)=\left(\mp\frac{\lambda}{\sqrt{1+4\lambda^2}},\frac12\pm\frac12 \frac{1}{\sqrt{1+4\lambda^2}}\right).\]

On the other hand,
\[\left(\left\|\dfrac{(\cos\frac\omega2)A-(\sin\frac\omega2)\Id}{(\sin\frac\omega2)A+(\cos\frac\omega2)\Id}\right\|_2\right)^2
=\frac{2-\cos\omega+\sqrt{1+(\sin\omega)^2}}{2+\cos\omega-\sqrt{1+(\sin\omega)^2}};\]
and the rotational construction yields
\[\widehat E_{\mathrm{CKB(P)}}^A(\omega)=\left(-\frac{\sin\omega}{\sqrt{1+(\sin\omega)^2}},
-\frac12+\frac12\frac{\cos\omega}{\sqrt{1+(\sin\omega)^2}} \right).\]
Transcribed to the $\mathrm{PH}$ model, it yields
\[\widehat E^A(\omega)=\left(\frac{-\sin\omega}{3\sqrt{1+(\sin\omega)^2}-\cos\omega },
\frac{\sqrt{1+(\sin\omega)^2}+\cos\omega }{3\sqrt{1+(\sin\omega)^2}-\cos\omega }\right).\]
At first sight it is not obvious that this traces out a circle.
\qedexer
\end{example}
\snewpage
\begin{commentx}
\begin{example}
(a) One (relatively) simple case is, the complex matrix
\[A=\begin{bmatrix}
\cos\alpha+\mathrm i\sin\alpha&2t\\&-\cos\alpha+\mathrm i\sin\alpha
\end{bmatrix}\]
$(\alpha\in[0,\pi/2] , \,t\geq0)$.
In this case
\[\left(\left\|\dfrac{(\cos\frac\omega2)A-(\sin\frac\omega2)\Id}{(\sin\frac\omega2)A+(\cos\frac\omega2)\Id}\right\|_2\right)^2
=\frac{\sqrt{1+t^2}+\sqrt{(\cos\alpha)^2(\sin\omega)^2+t^2}}{\sqrt{1+t^2}-\sqrt{(\cos\alpha)^2(\sin\omega)^2+t^2}} ;\]
and
\[\widehat E^A(\omega)=\left(\frac{-((\cos\alpha)^2+t^2)\sin\omega}{\sqrt{1+t^2} \sqrt{(\cos\alpha)^2(\sin\omega)^2+t^2}},
\frac{t^2\cos\omega}{\sqrt{1+t^2} \sqrt{(\cos\alpha)^2(\sin\omega)^2+t^2}}
\right).\]
This implies that $\partial\DW_{\mathrm{CKB}}^{\mathbb R} (A)$ is a possibly degenerate ellipse in the CKB model, whose axes are
\begin{equation*} \left[-\frac{\sqrt{(\cos\alpha)^2+t^2 }}{\sqrt{1+t^2}} ,\frac{\sqrt{(\cos\alpha)^2+t^2 }}{\sqrt{1+t^2}} \right]\times\{0\}\qquad\text{and}\qquad  \{0\}\times\left[-\frac{t}{\sqrt{1+t^2}},\frac{t}{\sqrt{1+t^2}}\right];
\end{equation*}
Again, this not entirely trivial to notice, although done so, it is easy to prove.
Using this line of argument, one can prove that the conformal range of
$2\times2$ complex matrices yields possibly degenerate $h$-ellipses.
In Part IIA several proofs will presented, and, in particular, this family of matrices will be reexamined.

However, it is easy to see that typically (if $\alpha\in(0,\pi/2)$ and $t>0$ hold) the ellipses in this example will not already yield ordinary
quadrics in the Poincar\'e half plane model, i. e. on the complex plane.

(b) A generalization of  Example \ref{ex:horoenv} applies to
\[A=\begin{bmatrix}0&\cos\beta\\0&\mathrm i\sin\beta\end{bmatrix}\]
with $\beta\in[0,\pi]$.
Then
\[\left(\left\|\dfrac{(\cos\frac\omega2)A-(\sin\frac\omega2)\Id}{(\sin\frac\omega2)A+(\cos\frac\omega2)\Id}\right\|_2\right)^2
=\frac{2-\cos\omega+\sqrt{2(\sin\omega)^2(\cos\beta)^2+(\cos\omega)^2}}{2+\cos\omega-\sqrt{2(\sin\omega)^2(\cos\beta)^2+(\cos\omega)^2}};\]
and the rotational construction yields
\[\widehat E_{\mathrm{CKB(P)}}^A(\omega)=\left(-\frac{(\cos\beta)^2\sin\omega}{\sqrt{2(\sin\omega)^2(\cos\beta)^2+(\cos\omega)^2}},
-\frac12+\frac12\frac{\cos\omega}{\sqrt{2(\sin\omega)^2(\cos\beta)^2+(\cos\omega)^2}} \right).\]
This implies that $\partial\DW_{\mathrm{CKB}}^{\mathbb R} (A)$ is a possibly degenerate ellipse in the CKB model, whose axes are
\[ \left[-\frac{\sqrt2}{2}\cos\beta,\frac{\sqrt2}{2}\cos\beta \right]\times\left\{-\frac12\right\}\qquad\text{and}\qquad  \{0\}\times[-1,0].\]
Again, for $\beta\in(0,\pi/2)$ these will not be quadrics in the PH model.

Together with the null matrix $0_2$ the matrices of this will cover every complex $2\times2$ matrix up to real M\"obius transformations.
\qedexer
\end{example}
\end{commentx}
More generally, Theorem \ref{th:envel} makes the finite dimensional case computable algebraically (at least, as much as such a computation is possible at all).
\begin{example}\plabel{ex:envel2}
Consider the matrix
\[A=\bem1&1&0\\&&-1\\&&-1\eem.\]
Then $N(A-\lambda I)$ is the maximum or minimum of $\spec((A-\lambda\Id)^*(A-\lambda\Id))$.
Thus, it satisfies the characteristic equation of $(A-\lambda\Id)^*(A-\lambda\Id)$ , yielding
\begin{equation}
N^3+(-3\,{\lambda}^{2}-4)N^2+(3\,{\lambda}^{4}+2\,{\lambda}^{2}+4)N+(-{\lambda}^{6}+2\,{\lambda}^{4}-{\lambda}^{2})=0
\plabel{eq:granchar}
\end{equation}
with $N=N(A-\lambda I)$.
As this equation is of order $3$ (the problem is real  $3$ dimensional), there are $3$ branches $N$ (corresponding to the $3$ nonnegative characteristic values).
For us, these are $N^+$ (the norm branch), $N^0$ (the middle branch) , $N^-$ (the co-norm branch).
For us only $N^+$ and $N^-$ are geometric significance, but, for certain computations we not have to specify a branch necessarily.
The branches are completely well-defined, but where they meet \eqref{eq:granchar} has double roots.
The equation for such $\lambda$ is the discriminant equation of \eqref{eq:granchar}.
It is actually
\begin{equation}
864\,{\lambda}^{6}-99\,{\lambda}^{4}+288\,{\lambda}^{2}=0.
\plabel{eq:grandisk}
\end{equation}
In the present case this is very convenient, as the only real discriminant root is $\lambda=0$.
Thus the branches of $N$ are analytic on parts of $\mathbb R$ divided only by $-\infty,0, \infty$.
In fact, in order to avoid any other complications, we also put $\lambda=-1,1,0$ (the real spectrum of $A$, where the conformal range meets the real line)
into the discriminant set.
Thus, ultimately, we consider $-\infty,-1,0,1,\infty$ as the extended discriminant set.

Apart from the locus of the discriminant, taking the derivative of \eqref{eq:granchar},
\begin{multline}\frac{\mathrm dN}{\mathrm d\lambda}\cdot\left((3)N^2+(-6\,{\lambda}^{2}-8
)N+(3\,{\lambda}^{4}+2\,{\lambda}^{2}+4)\right)+\\+\left((-6\,\lambda
)N^2+(12\,{\lambda}^{3}+4\,\lambda)N+(-6\,{\lambda}^{5}+8\,{\lambda}^{3}-2\,\lambda)\right)=0.
\end{multline}
Now it might be case that the coefficient of $\frac{\mathrm dN}{\mathrm d\lambda}$ is identically $0$ on some pieces (they would correspond to locally multiple branches),
but this is not the case, thus we simply obtain
\[\frac{\mathrm dN}{\mathrm d\lambda}=2\,{\frac {\lambda\, \left( 3\,{\lambda}^{4}-6\,{\lambda}^{2}N+3\,{N}^
{2}-4\,{\lambda}^{2}-2\,N+1 \right) }{3\,{\lambda}^{4}-6\,{\lambda}^{2
}N+3\,{N}^{2}+2\,{\lambda}^{2}-8\,N+4}}.
\]
Hence, the corresponding enveloping curve in the CKB(P) model is
\begin{multline}(X(\lambda),Y(\lambda))
=\Biggl(\,{\frac { -3\lambda\, \left( -2\,{\lambda}^{2}+2\,N-1 \right) }{3\,{
\lambda}^{4}-6\,{\lambda}^{2}N+3\,{N}^{2}+2\,{\lambda}^{2}-8\,N+4}},\\{
\frac {-3\,{\lambda}^{6}+9\,{\lambda}^{4}N-9\,{\lambda}^{2}{N}^{2}+10
\,{\lambda}^{4}+3\,{N}^{3}-2\,{\lambda}^{2}N-8\,{N}^{2}+2\,{\lambda}^{
2}+4\,N}{3\,{\lambda}^{4}-6\,{\lambda}^{2}N+3\,{N}^{2}+2\,{\lambda}^{2
}-8\,N+4}}
\Biggr).\end{multline}
We can plot this as in Figure \ref{fig:figB01}(a)(b).

\begin{figure}[H]
   \begin{subfigure}[b]{3.5in}
    \includegraphics[width=3.5in]{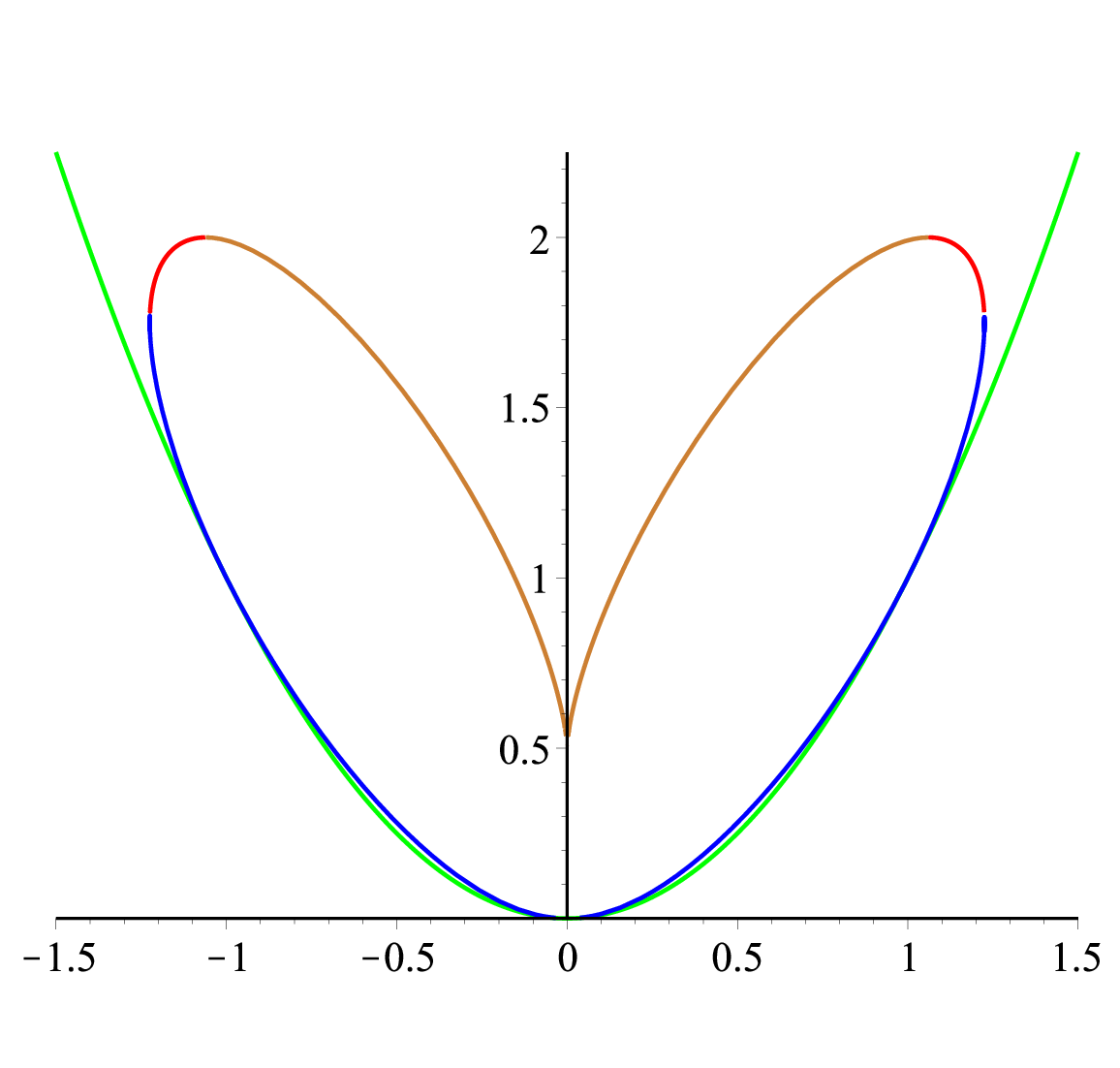}
    \caption*{Fig. \ref{fig:figB01}(a) enveloping construction, CKB(P) model}
  \end{subfigure}
  \begin{subfigure}[b]{2.1875in}
    \includegraphics[width=2.1875in]{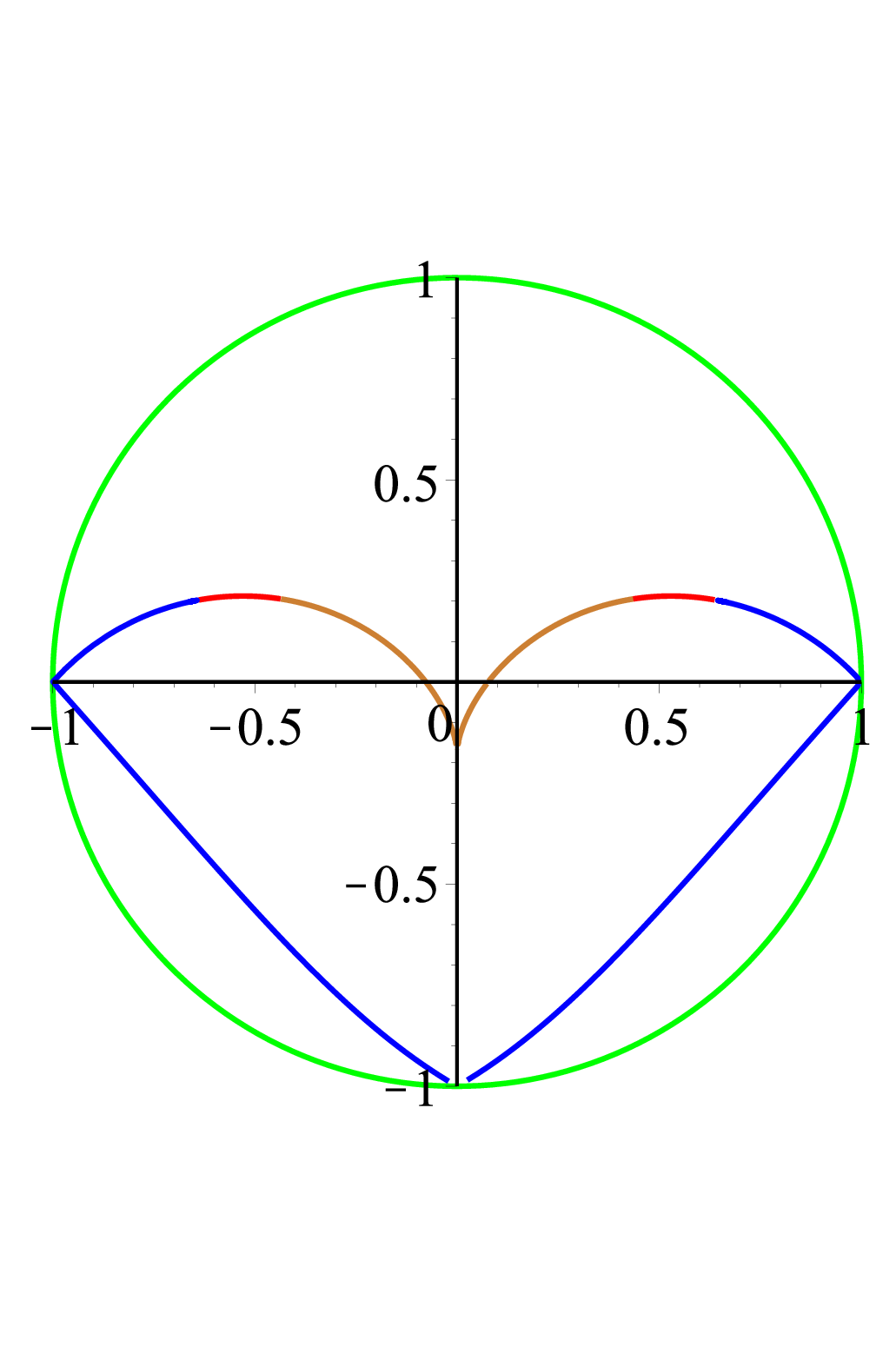}
    \caption*{\ref{fig:figB01}(b) (Poincar\'e disk model)}
  \end{subfigure}
\phantomcaption
\plabel{fig:figB01}
\end{figure}~\\[-0.5cm]
Here the upper (norm) branch is given by red, the lower (co-norm) branch is given by blue, and the middle (artifact) branch is by brown;
green is the asymptotic boundary of the CKB(P) mode.
The co-norm branch meets the asymptotic boundary  at points corresponding to $\lambda=-1,0,1$.
Note, however, that the enveloping construction is not even piecewise continuous, as $\frac{\mathrm dN}{\mathrm d\lambda}$
might have singularities at the points of the discriminant locus.
Indeed, this is the case, for the upper (norm) branch at $\lambda=0$,
where the discontinuity must be bridged by a segment.
\\[-0.5cm]
\begin{figure}[H]
   \ContinuedFloat
   \begin{subfigure}[b]{3.5in}
    \includegraphics[width=3.5in]{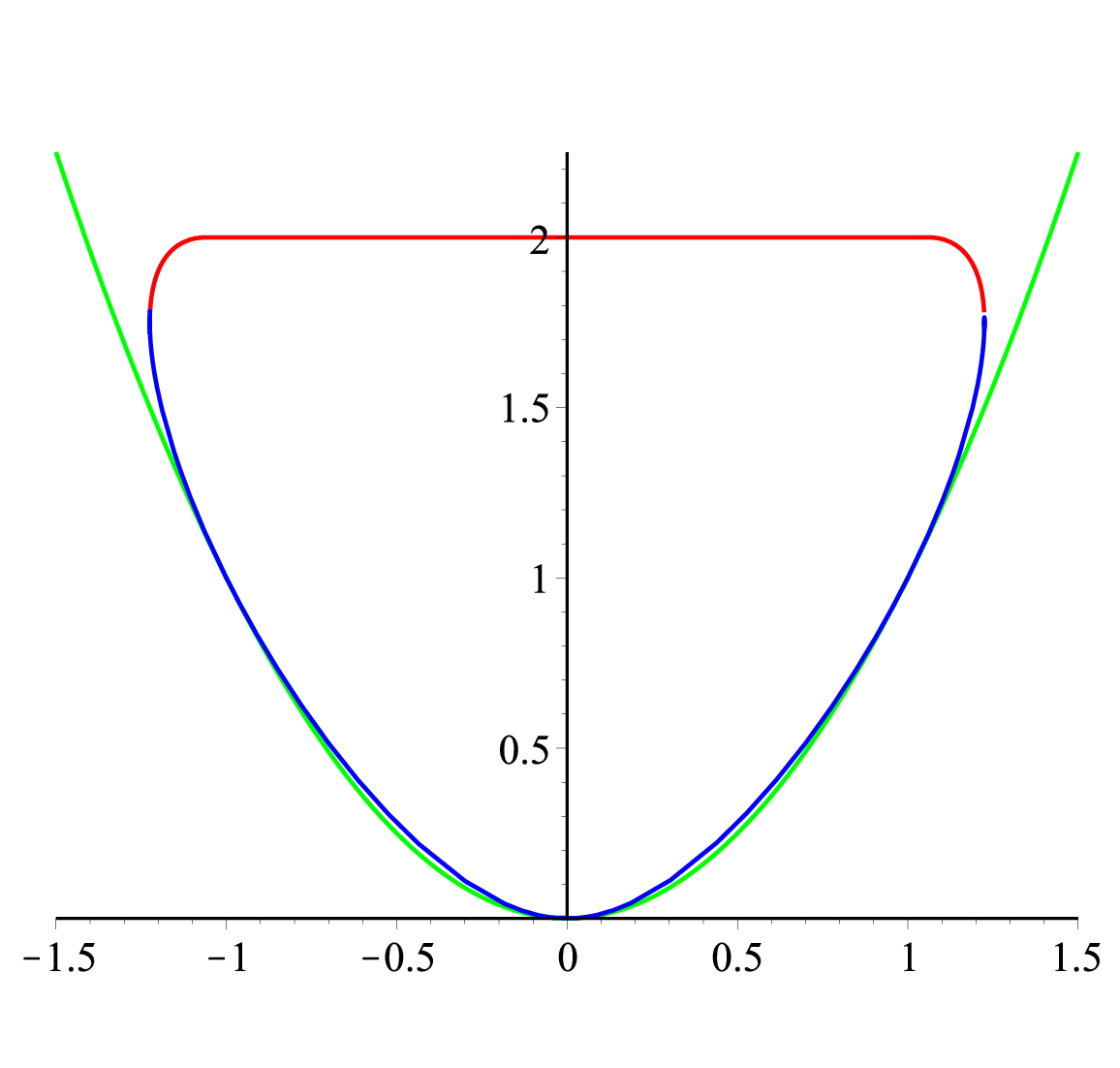}
    \caption*{Fig. \ref{fig:figB01}(c) full conformal range, CKB(P) model}
  \end{subfigure}
  \begin{subfigure}[b]{2.1875in}
    \includegraphics[width=2.1875in]{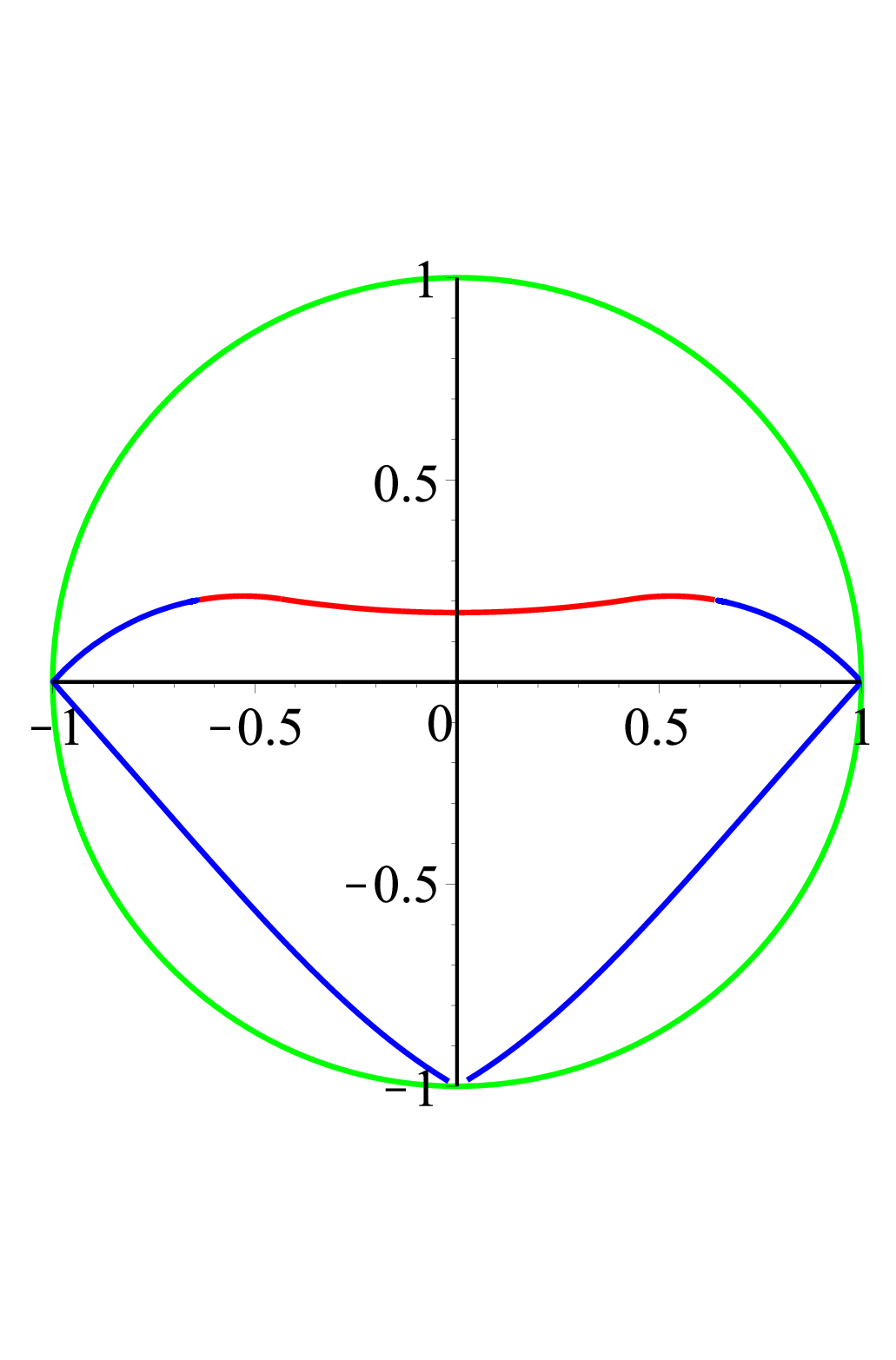}
    \caption*{\ref{fig:figB01}(d) (Poincar\'e disk model)}
  \end{subfigure}
\phantomcaption
\end{figure}

At that point we have to compute norm (or co-norm, in other cases), and the corresponding $h$-line
will give the line of the segment.
Here this is the $h$-line $y=2$.
The endpoints occur as limits from non-discriminant points of branches, at this point these are at $x=\pm\frac{3\sqrt2}{4}$.
Thus, ultimately, the $h$-segment between $(-\frac{3\sqrt2}{4},2)$ and $(\frac{3\sqrt2}{4},2)$ should be added to the upper (norm) branch.
Similar considerations for $\lambda=\pm\infty$ are not necessary here, as the norm and co-norm branches
meet at $\left(\pm \frac{\sqrt6}{2},\frac74\right)$.
After eliminating the artificial branches but bridging the discontinuities, we obtain
the proper enveloping construction as shown in Figure \ref{fig:figB01}(c)(d).

Due to \eqref{eq:granchar}, however, $N$ and $\lambda$ are algebraically not independent.
Indeed, after a bit of computation, we find that $X(\lambda)$ and $Y(\lambda)$ satisfy the algebraic
equation
\[p_{\mathrm K}(x,y)\equiv128\,{x}^{4}-11\,{x}^{2}{y}^{2}+24\,{y}^{4}-124\,{x}^{2}y-36\,{y}^{3}+
4\,{x}^{2}+18\,{y}^{2}-3\,y=0.
\]
However, this is not the whole story yet, as we have to add the boundary segment corresponding to the discontinuity at the discriminant locus $\lambda=0$.
The resulting algebraic equation is
\[(y-2)(128\,{x}^{4}-11\,{x}^{2}{y}^{2}+24\,{y}^{4}-124\,{x}^{2}y-36\,{y}^{3}+
4\,{x}^{2}+18\,{y}^{2}-3\,y)=0;\]
the boundary of the conformal range is on that algebraic curve.
It is clear, however, that the boundary is not the whole curve because only part of $y-2=0$ is included in the picture, and the
artificial branches must also be eliminated from consideration.
See Figure \ref{fig:figB01}(e)(f). (The Poincar\'e picture is self-inverted.)
\\[-0.5cm]
\begin{figure}[H]
   \ContinuedFloat
   \begin{subfigure}[b]{3.5in}
    \includegraphics[width=3.5in]{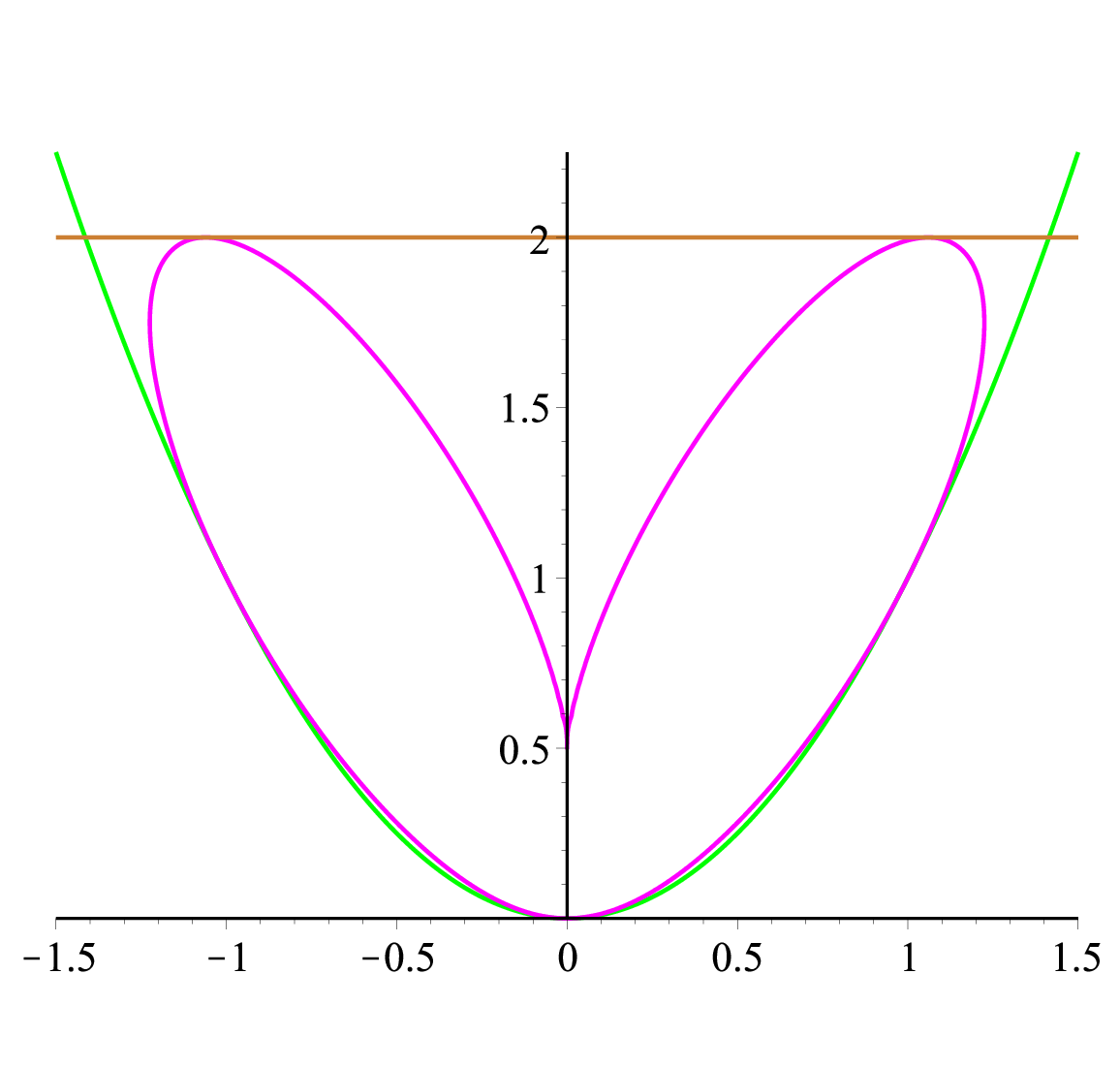}
    \caption*{Fig. \ref{fig:figB01}(e)  algebraic data, CKB(P) model}
  \end{subfigure}
  \begin{subfigure}[b]{2.1875in}
    \includegraphics[width=2.625in]{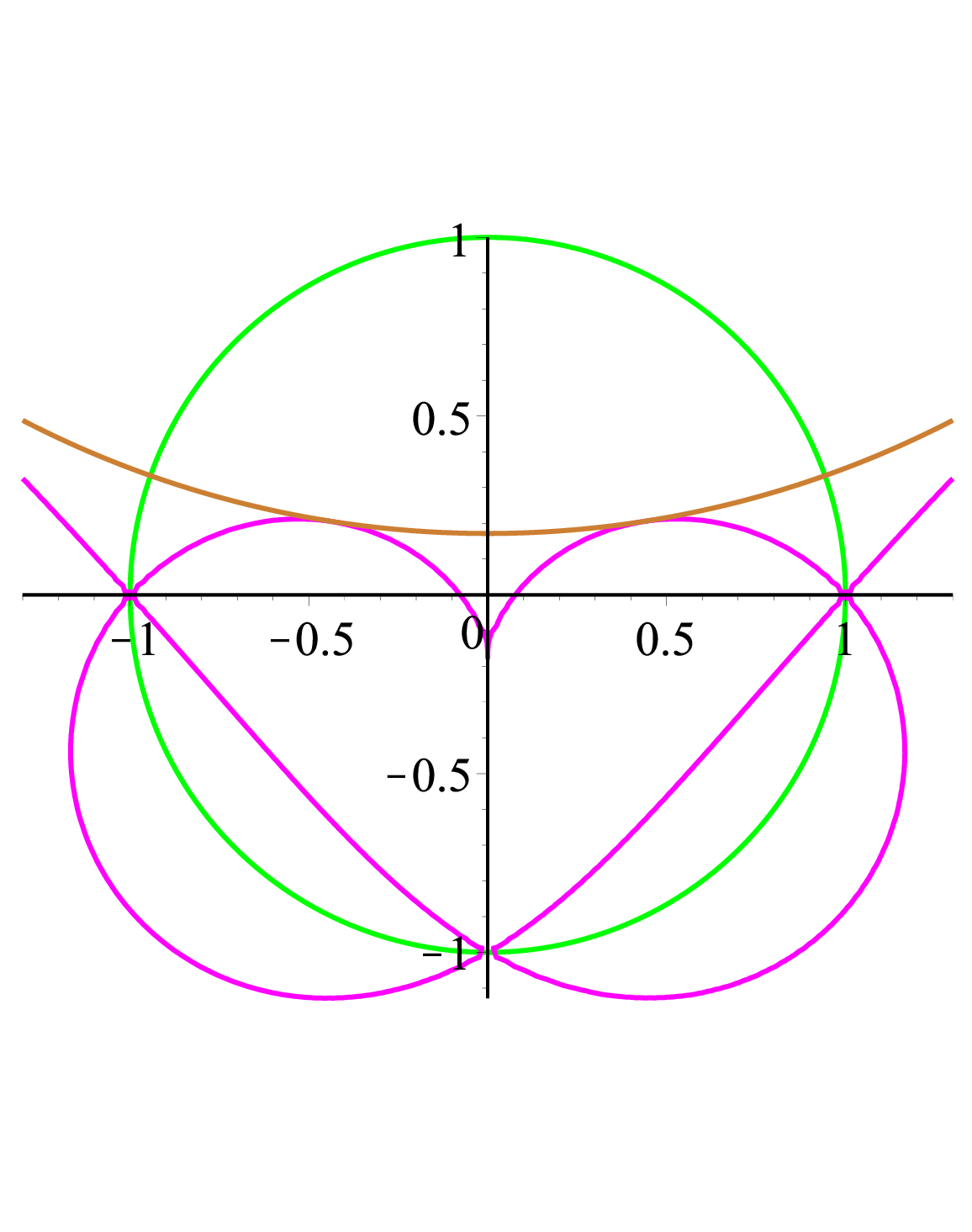}
    \caption*{\qquad\quad\ref{fig:figB01}(f)  (Poincar\'e disk model)}
  \end{subfigure}
\phantomcaption
\end{figure}

Purely algebraically, we have solved the system
\[-2\lambda x+y+\lambda^2-N=0\]
\[-2x+2\lambda-N'=0\]
\[F(\lambda ,N)=0\]
\[\partial_1F(\lambda ,N)+N'\partial_2F(\lambda,N)=0\]
with $F(\lambda,N)$ as \eqref{eq:granchar}.
(Thus in the present case it happens to be the $F_A^{\CR}(\lambda,N)$, which is already reduced.)
This leads to the discriminant of $F(\lambda,-2\lambda x+y+\lambda^2)$ in $\lambda$.
This amounts having not only to the smooth enveloping construction but also to the
the bitangents (that is to discontinuities in $N'$).
In the present case
\begin{multline}
F(\lambda,-2\lambda x+y+\lambda^2)=\left( 8\,{x}^{3}-12\,x \right) {\lambda}^{3}+\\+ \left( -12\,{x}^{2}y+
16\,{x}^{2}+6\,y-3 \right) {\lambda}^{2}+ \left( 6\,x{y}^{2}-16\,xy+8
\,x \right) \lambda-{y}^{3}+4\,{y}^{2}-4\,y;
\end{multline}
and the discriminant made equal to $0$ is
\[36\, \left( y-2 \right) ^{2} \left( 128\,{x}^{4}-11\,{x}^{2}{y}^{2}+24
\,{y}^{4}-124\,{x}^{2}y-36\,{y}^{3}+4\,{x}^{2}+18\,{y}^{2}-3\,y
 \right)=0.\]
The linear factors can be omitted as convex closure takes care of them.
\qedexer
\end{example}
\begin{theorem}
\plabel{th:qKippenhahn}
If $\dim\mathfrak H<\infty$, then $\partial\CR(A)$ is the union of finitely
many algebraic arcs,
data obtained from the polynomials
$F_A^{\CR}(\lambda,\nu)$ and / or $K^{\CR}_A(u,s,w)$.
\begin{proof}
The process we have seen in the previous example works in general.
Multiple branches may need to be reduced out in order to obtain $\frac{\mathrm dN}{\mathrm d\lambda}$.
Special considerations about $\lambda=\pm\infty$ are in fact unnecessary, as the operator can be
subjected to conformal (fractional linear) transformations making the situation in the real suprema and infima smooth.)
\end{proof}
\end{theorem}
\begin{disc}
\plabel{rem:qKippenhahn}
Let us summarize the computation of the conformal range in the $\mathrm{CKB(P)}$ model.

(o) Consider $F^{\CR}_A(\lambda,\nu)$.
We can observe that any root in $\nu$ is real; and any locally polynomial (in $\lambda$) dependence of a
root $\nu_i$ of $F^{\CR}_A(\lambda,\nu)$ can occur only in form $\nu_i=\lambda^2-2x_i\lambda+y_i$.
Indeed, local polynomial dependence extends to a global dependence and the norm / co-norm bounds for the eigenvalues allow only $\lambda^2$
asymptotics.

(a) The analytic process of finding the conformal range can be summarized  as follows:

First we reduce out multiplicities in $F^{\CR}_A$, in order to get single branches for $\nu$.
Then we eliminate the $\lambda$ which give multiple roots in $\nu$.
(Take discriminant in $\nu$, solve for real $\lambda$.)
This eliminates finitely many $\lambda$, on the complement $\nu$ has simple, disjoint branches.
(The eliminated $\lambda$ lead to ``algebraic faces''.)
For polynomial branches the enveloping construction gives constant curves (vertices), or,
otherwise, smooth curves.
If an artificial branch would get out of the conformal range, then at one of its ``farthest'' point,
the corresponding tangent would avoid the conformal range, but this would be in contradiction
to that the corresponding eigenvalue is not extremal.
Now the closure of the convex hull of the enveloping curves must be taken.

Eliminated $\lambda$ are not problematic because the end points of the corresponding geometric faces
will be
approached by non-artificial smooth simple branches with continuous tangents
or vertices.
But taking closure is necessary, as the end points are typically just approached;  and this also applies to ``artificial geometric faces''
corresponding $\lambda=+\infty$, which is missing as a slope.

(b) We can proceed algebraically as follows: (This is when we replace actual differentiation by
implicit differentiation, that is essentially by algebraic trickery.)

By (o), $F^{\CR}_A(\lambda,\nu)$ decomposes as
\[F^{\CR}_A(\lambda,\nu)=\prod_{i=0}^s(\nu-(\lambda^2-2x_i\lambda+y_i))\prod_{j=0}^r f_j(\lambda,\nu),\]
where the $f_j$ are real irreducible such that all their roots are of locally nonpolynomial dependence from $\lambda$.

In terms of homogeneous polynomials, this corresponds to
\[K^{\CR}_A(u,s,w)=\prod_{i=0}^s(x_iu+y_is+w)\prod_{j=0}^r k_j(u,s,w),\]
where the normalization $k_j(0,0,1)=1$ can be assumed.

The linear factors lead to vertices $(x_j,y_j)$.
(All angular corners of the conformal range lead to algebraic vertices, cf. also
$F^{\CR}_A(\lambda,\lambda^2-2x\lambda+y)=K^{\CR}_A(2\lambda,-1,-2\lambda x +y)$.)
The non-linear factors lead to non-normal curves (in $x,y$) which occur from the discriminants of
$f_j(\lambda,\lambda^2-2x\lambda+y )$ in $\lambda$
(i. e.  $k_j( 2\lambda,-1,-2\lambda x +y)$ in $\lambda$) but the linear factors should be removed.
Then the conformal range is the convex hull of the vertices and the non-normal curves.
(Actual decomposition to irreducibles is not necessary, we can proceed with the reduced product.)
\end{disc}

Thus in that, and other properties, the conformal range is a proper analogue
of the numerical range and also some other types of range, cf. Horn, Johnson \cite{HJ}, and
Lins, Spitkovsky, Zhong \cite{LSZ}.

\begin{example}\plabel{ex:envel3}
Consider the matrix
\[\tilde A=\bem 1&\sqrt2&\\&-1&\\&&0\eem.\]
Compared to matrix $A$ of Example \ref{ex:envel2}, we see that their
algebraic numerical ranges are the same:
\[K^{\mathrm W}_A(u,v,w)=K^{\mathrm W}_{\tilde A}(u,v,w)=w\left(w^2-\frac32u^2-\frac12v^2\right).\]
Geometrically, regarding the numerical range, the first component gives the point $(0,0)$, while the
second component yields the ellipse
$\dfrac23x^2+\dfrac21y^2-1=0$.
Taking convex closure, the numerical range is given by
\[\dfrac{x^2}{\sqrt{\dfrac32}^2}+\dfrac{y^2}{\sqrt{\dfrac12}^2}-1\leq0\]
(on the complex plane, for $x+\mathrm iy$).

Regarding the conformal range, however, the algebraic conformal range is given by
\[K^{\CR}_{\tilde A}(u,s,w)=w\left(w^2-\frac32u^2+s^2+4s\right).\]
The first component gives the point $(0,0)$, and the second one gives the ellipse with equation
\[\frac23(x_{\mathrm{CKB(P)}})^2+\frac13(y_{\mathrm{CKB(P)}}-2)^2=1.\]
Thus $\DW_{\mathrm{CKB(P)}}^{\mathbb R}(\tilde A)$ is the convex closure of the union of ellipse
and $\{(0,0)\}$.
(In Figure \ref{fig:figB011}, the corresponding algebraic data is indicated.)
\begin{figure}[H]
   \begin{subfigure}[b]{3.5in}
    \includegraphics[width=3.5in]{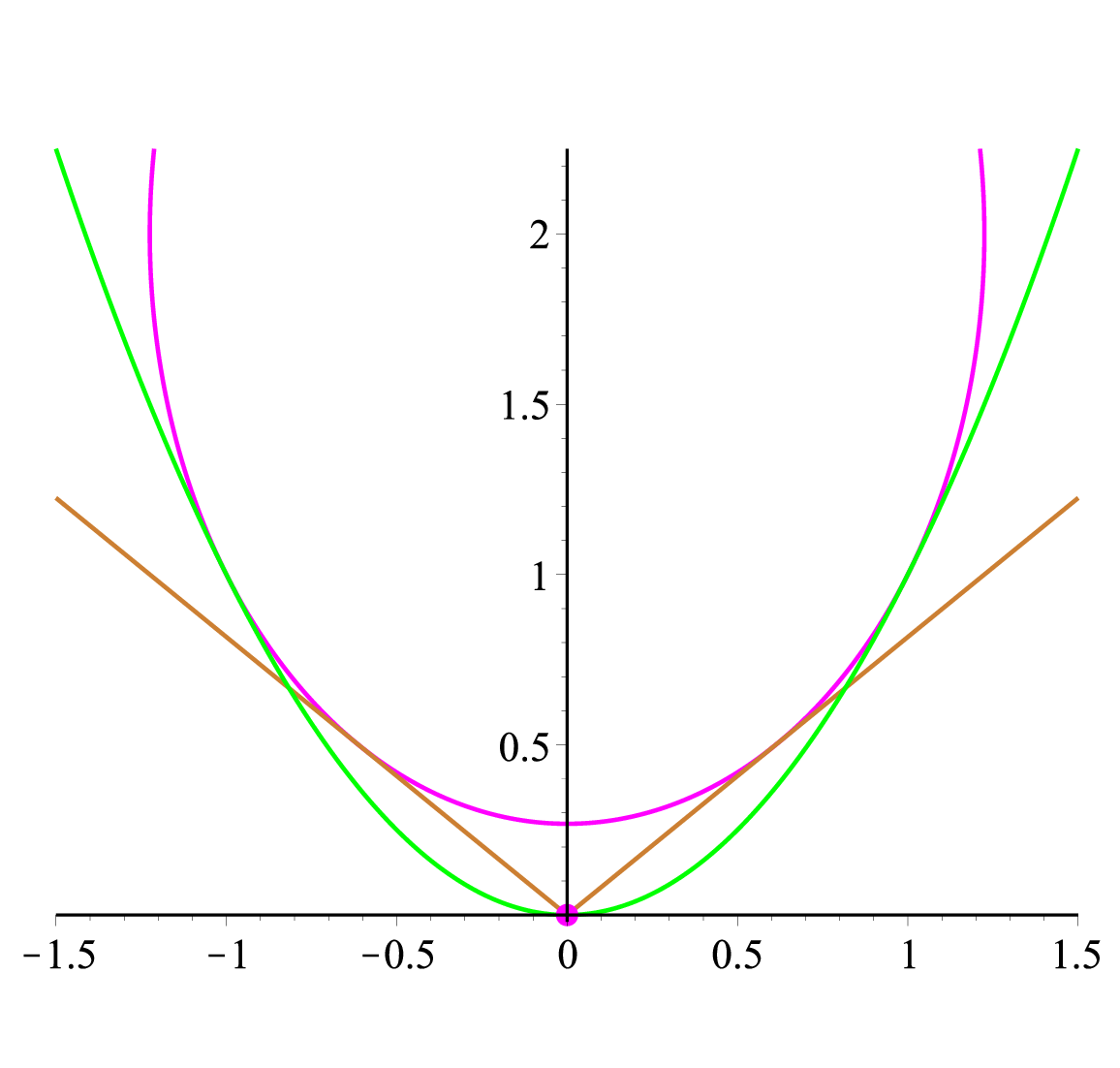}
    \caption*{Fig. \ref{fig:figB011}(a)  algebraic data, CKB(P) model}
  \end{subfigure}
  \begin{subfigure}[b]{2.1875in}
    \includegraphics[width=2.625in]{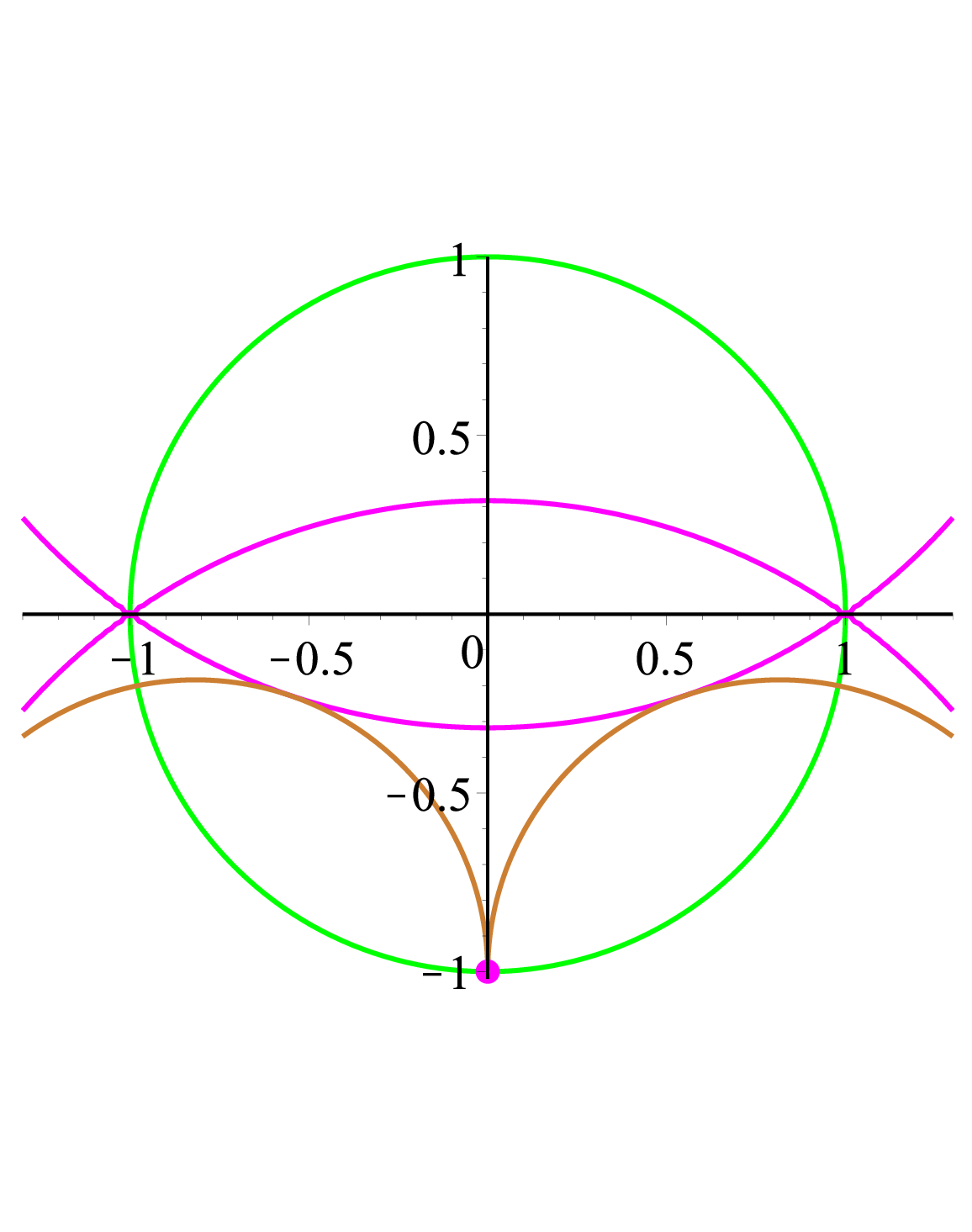}
    \caption*{\qquad\quad\ref{fig:figB011}(b)  (Poincar\'e disk model)}
  \end{subfigure}
\phantomcaption
\plabel{fig:figB011}
\end{figure}

This shows that, beyond the real suprema and infima, and  focal properties, there are little
obvious relations between the numerical range and the conformal range.

Thus, even if the algebraic numerical ranges are the same, the conformal ranges may be different.
For this, the following simpler argument also suffices:
One can check that
\[\|A\|_2=\sqrt2\qquad<\quad \|\tilde A\|_2=\frac{\sqrt6+\sqrt2}2.\]
As the norm can be read off from the conformal range in general, the conformal ranges of $A$ and $\tilde A$ must be different.
\qedexer
\end{example}
\begin{example}
\plabel{ex:envel4}
Conversely, it may happen that the conformal ranges are the same but the numerical ranges differ.
The trivial choice is the pair of matrices $\bem\mathrm i\eem$ and $\bem-\mathrm i\eem$.
However, the phenomenon can also occur for real matrices:
Consider
\[A_1=\bem&-1&1&\\1&&&1\\&&&-1\\&&1&\eem
\qquad\text{and}\qquad
A_2=\bem&-1&1&\\1&&&1\\&&&1\\&&-1&\eem.
\]
Although these are real matrices, we may notice the unitary equivalences
\begin{equation}
A_1\simeq\bem\mathrm i&1&&\\&\mathrm i&&\\&&-\mathrm i&1\\&&&-\mathrm i\eem
\qquad\text{and}\qquad
A_2\simeq\bem\mathrm i&1&&\\&-\mathrm i&&\\&&\mathrm i&1\\&&&-\mathrm i\eem.
\plabel{eq:A12deco}
\end{equation}
The algebraic conformal ranges are the same for $A_1$ and $A_2$:
\[K^{\CR}_{A_1}(u,s,w)=K^{\CR}_{A_2}(u,s,w)=\left(w^2-\frac14u^2+s^2+3sw\right)^2.\]
This yields two copies of an ellipse for the conformal range.
Taking convex closure, this yields
\[4(x_{\mathrm{CKB(P)}})^2+\frac45\left(y_{\mathrm{CKB(P)}}-\frac32 \right)^2\leq1\]
for the conformal range.

The numeral ranges can simply be recovered from the decomposition
\eqref{eq:A12deco} and the elliptic range theorem.
For $A_2$ the numerical range will be the elliptical disk
\[\mathrm W(A_2)=\left\{x+y\mathrm i \,:\, 4x^2+\frac45 y^2\leq 1\right\}.\]
On the other hand,
\[\mathrm W(A_1)=\conv\left( \Dbar\left(\mathrm i,\frac12\right) \cup \Dbar\left(-\mathrm i,\frac12\right)\right),\]
it is ``stadium shaped''.
$\mathrm W(A_2)$ is strictly contained in $\mathrm W(A_1)$.
\qedexer
\end{example}
\begin{remark}
For $\dim_{\mathbb C}\mathfrak H\leq 2$, however, there is a strong connection between the
numerical range and the conformal range of linear operators.
Then the numerical range is equivalent to the Davis-- Wielandt
shell, and conformal contains only slightly less in\-for\-ma\-tion; cf. \cite{L2}.
The various ranges of $3\times3$ matrices still seem to be quite manageable,
cf. also Keeler, Rodman, Spitkovsky \cite{KRS}, but already much more difficult.
\qedremark
\end{remark}
\snewpage
\section{Conformal range and time-ordered exponentials (continued)}
\plabel{sec:MagnusHilbertTwo}

For the sake of visualization, let us consider the boundary of $\exp \Dbar(0,p)\cap \overline{\mathbb C}^+$ for $p>0$.
Consider the curve
\[\gamma_p(t)=\exp(p(\cos t+\mathrm i\sin t))=\mathrm e^{p\cos(t)}\cos(p\sin t)+\mathrm i \mathrm e^{p\cos(t)}\sin(p\sin t) ).\]
If $0\leq p\leq \pi$, the boundary is given by this curve, $t\in\mathcal I_p=[0,\pi]$.
For $p>\pi$, it is also given by this curve, but for $t\in\mathcal I_p=\left[0,\arcsin\frac\pi p\right]\cup \left[\pi-\arcsin\frac\pi p,\pi\right]$.
In the following Figure \ref{fig:figB02}, we depict this for $p=\frac\pi4,\boldsymbol{\frac\pi2},\frac{3\pi}4,\boldsymbol\pi,\frac{5\pi}4$.\\
\begin{figure}[H]
   \begin{subfigure}[b]{2.5in}
    \includegraphics[width=2.5in]{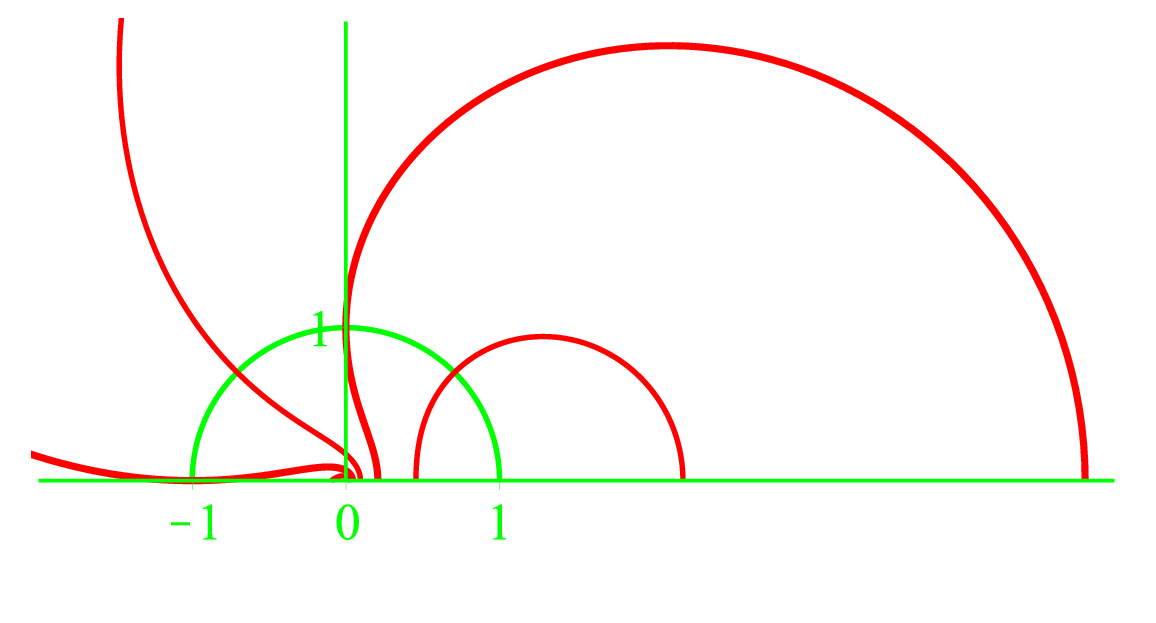}
    \caption*{Fig. \ref{fig:figB02}(a)  Poincar\'e half-plane, large scale}
  \end{subfigure}
  \begin{subfigure}[b]{2.5in}
    \includegraphics[width=2.5in]{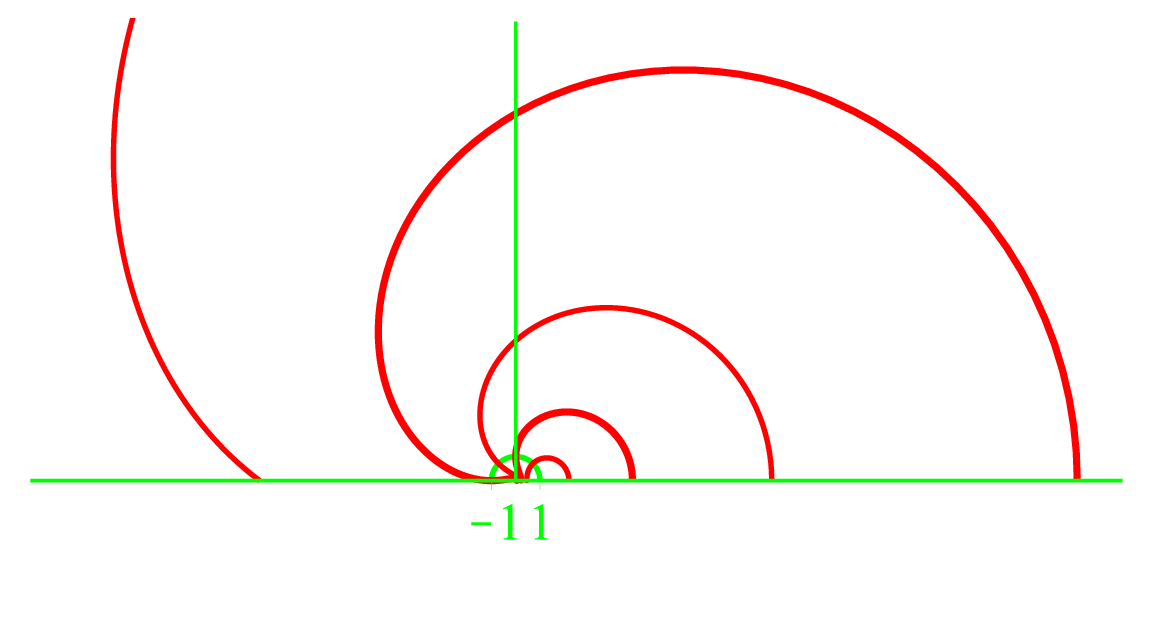}
    \caption*{\ref{fig:figB02}(b)  Poincar\'e half-plane, small scale}
  \end{subfigure}
\phantomcaption
\plabel{fig:figB02}
\end{figure}
We also show the corresponding picture in the Cayley--Klein--Beltrami model, where $h$-convexity is apparent;
and we also include the view in the Poincar\'e disk model (which shows certain details better).
\begin{figure}[H]
   \ContinuedFloat
   \begin{subfigure}[b]{2.5in}
    \includegraphics[width=2.5in]{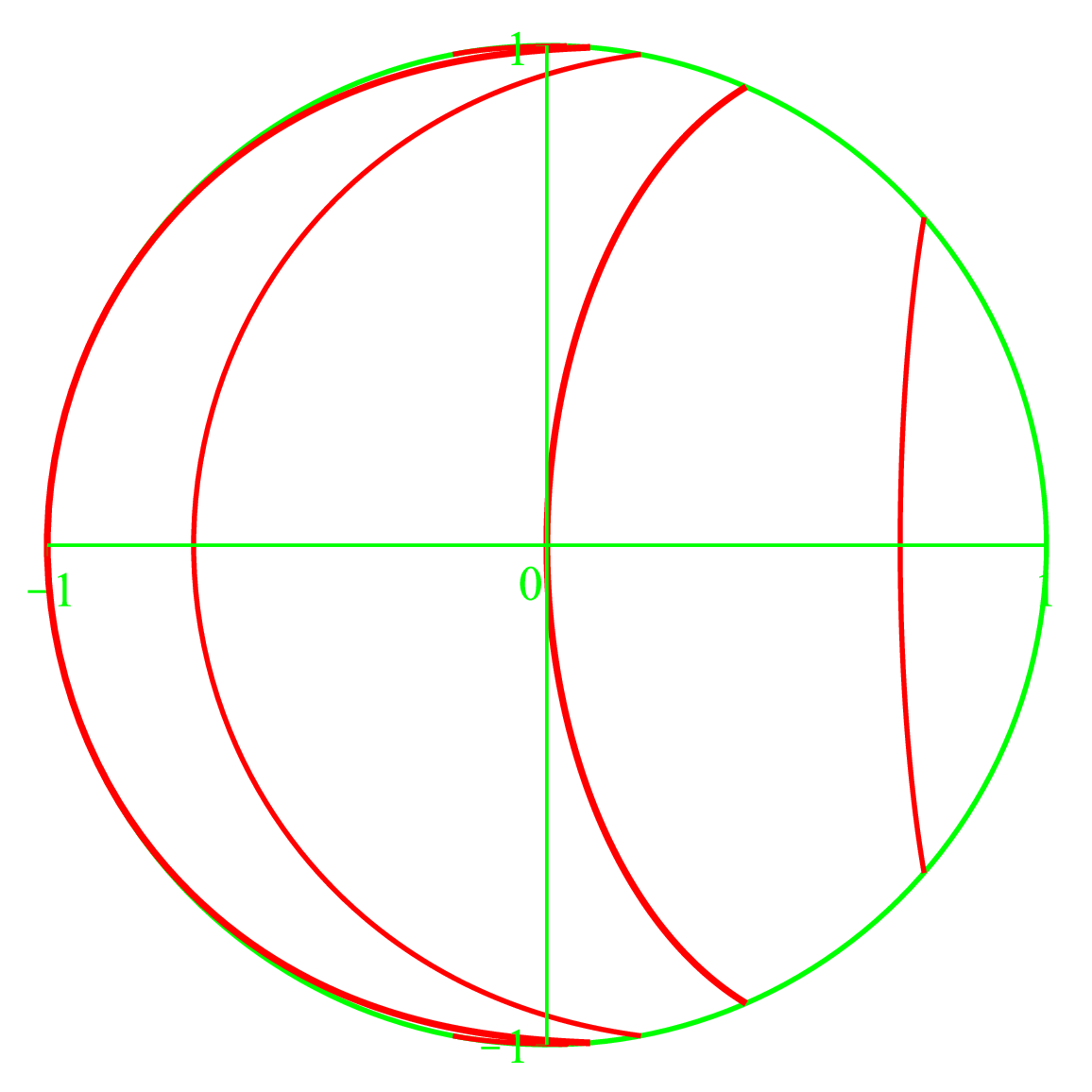}
    \caption*{\ref{fig:figB02}(c)   Cayley--Klein--Beltrami model}
  \end{subfigure}
  \begin{subfigure}[b]{2.5in}
    \includegraphics[width=2.5in]{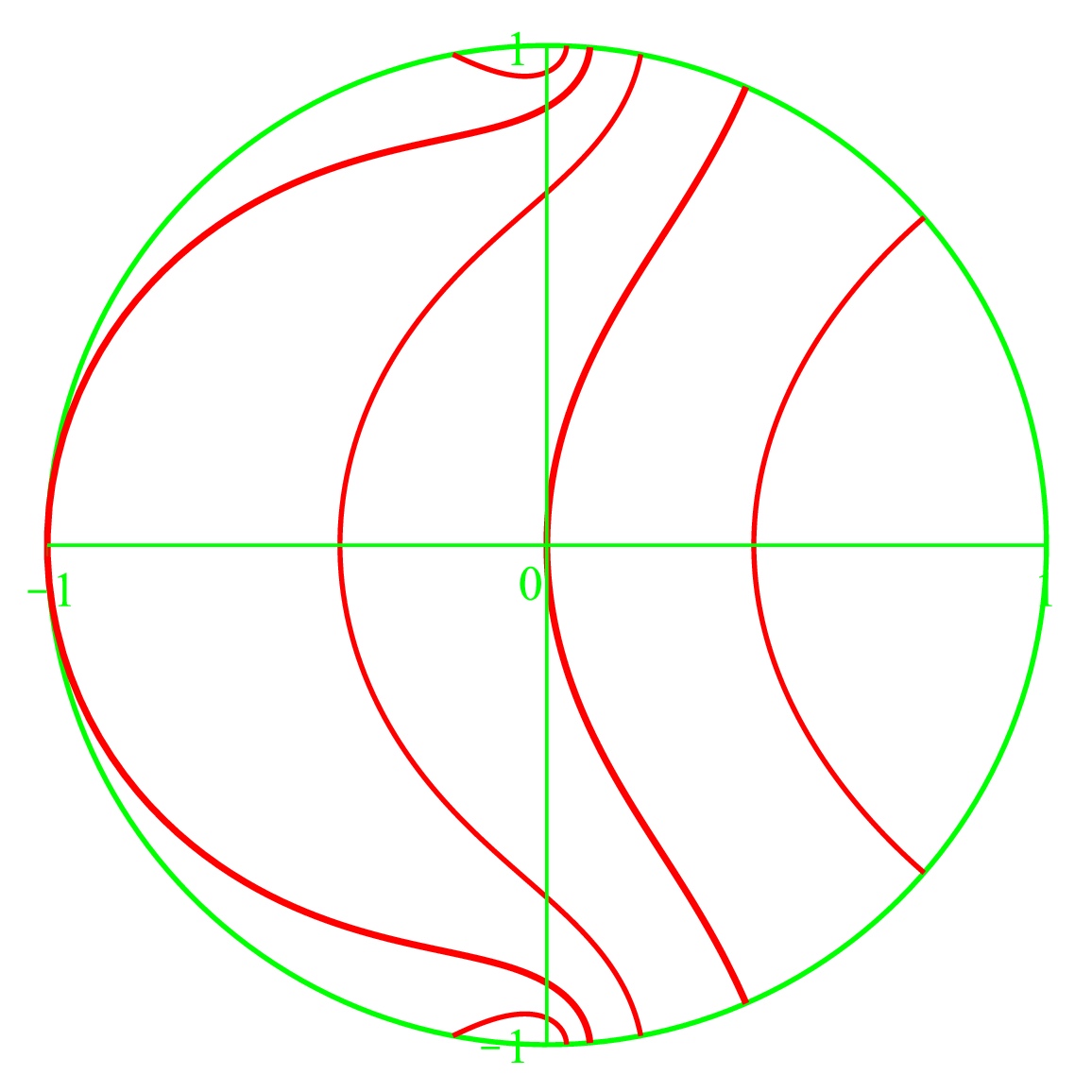}
    \caption*{\ref{fig:figB02}(d)   Poincar\'e disk model}
  \end{subfigure}
  \phantomcaption
\end{figure}

\snewpage
For a subset, whenever
\[a\hat x+b\hat y+c\geq0\]
holds in the CKB model, it means that
\[a \,2 x+b(x^2+y^2-1)+c(x^2+y^2+1)\geq0\]
holds in the Poincar\'e half-plane model.
With respect to (the half-planes containing) $\CR(A)$, it means that
\[a (A+A^*)+b(A^*A-\Id)+c(A^*A+\Id)\geq0,\]
or, equivalently,
\[(b+c) A^*A+ a (A+A^*)+(c-b)\Id\geq0\]
should hold.
This can be applied to the half-planes induced by the tangent lines of $\exp \Dbar(0,p)\cap \overline{\mathbb C}^+$ for $p>0$.
After some computation, this yields that  $\CRext(A)\subset \exp\Dbar(0,p)$ is equivalent to the collection of operator inequalities
\[\mathcal E_{p,t}(A):\,\,
-{\frac {{{\rm e}^{-p\cos t }}\sin \left( t+p\sin
 \left( t \right)  \right) }{\sin \left( t \right) }}
A^*A+(A^*+A)-
{\frac {{{\rm e}^{p\cos t }}\sin \left( t-p\sin
 \left( t \right)  \right) }{\sin \left( t \right) }}
\Id\geq0
\]
for $t\in\mathcal I_p$. (For $t=0$, this is to be  understood as
\[-(1+p)\mathrm e^{-p}A^*A+(A+A^*)-(1-p)\mathrm e^{p}\Id\geq0;\]
and for $t=\pi$, this is
\[-(1-p)\mathrm e^{p}A^*A+(A+A^*)-(1+p)\mathrm e^{-p}\Id\geq0.)\]

\snewpage
It is natural to ask whether Theorem \ref{th:CRrange} holds on the level of the operator inequalities.
The answer is affirmative:

\begin{theorem}\plabel{th:CRrangeInv}
Suppose that $t\in[0,\pi]$; $0\leq p\leq q\leq \frac\pi{\sin t}$.
If $\phi$ is $\mathcal B(\mathfrak H)$-valued,
\[B=\Lexp(\phi),\qquad\textstyle{\int \|\phi\|_2}\leq q-p,\]
then
\[\mathcal E_{p,t}(A)\Rightarrow\mathcal E_{q,t}(BA).\]
\begin{proof}
For a fixed $t$, take all tangent $h$-lines $L_{p,t}$ at $\gamma_p(t)$, for $p\in\left(0,\frac\pi{\sin t}\right)$.
Then, according to  Theorem \ref{th:cs}, the question is whether $\log(L_{p,t})$ and $\log(L_{q,t})$ are sufficiently
great (i. e. at least $q-p$) distance from each other.

For the sake of visualization, for some fixed $t$, we draw some ($p$-equidistant) tangent lines in the
CKB model, the Poincar\'e half-plane model (i. e.  on the upper complex halfplane), and their logarithms.
For $t=\frac\pi2$, they yield (Figure \ref{fig:figB03}):
\begin{figure}[H]
   \begin{subfigure}[b]{2in}
    \includegraphics[width=2in]{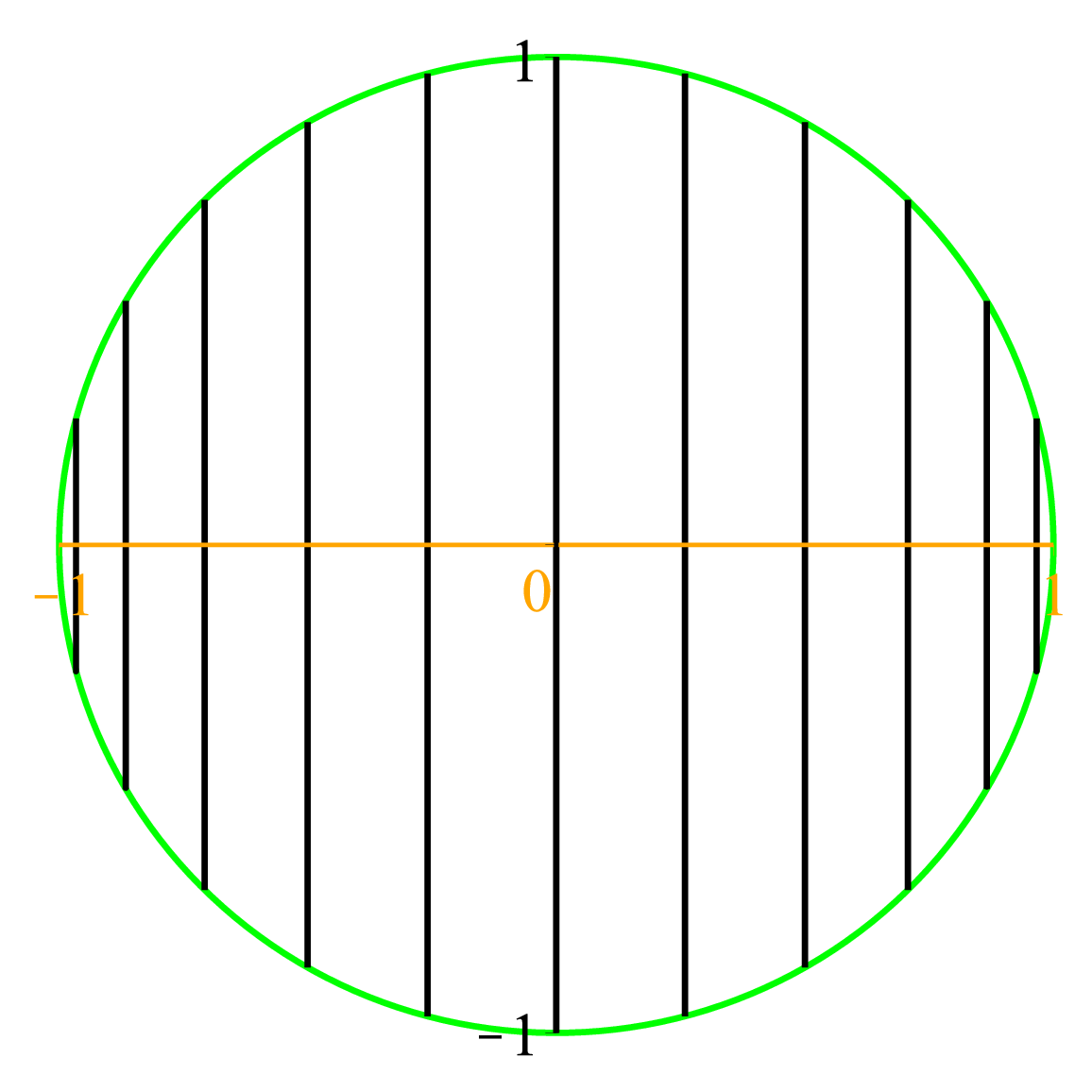}
    \caption*{Fig. \ref{fig:figB03}(a)   CKB model}
  \end{subfigure}
  \begin{subfigure}[b]{3in}
    \includegraphics[width=3in]{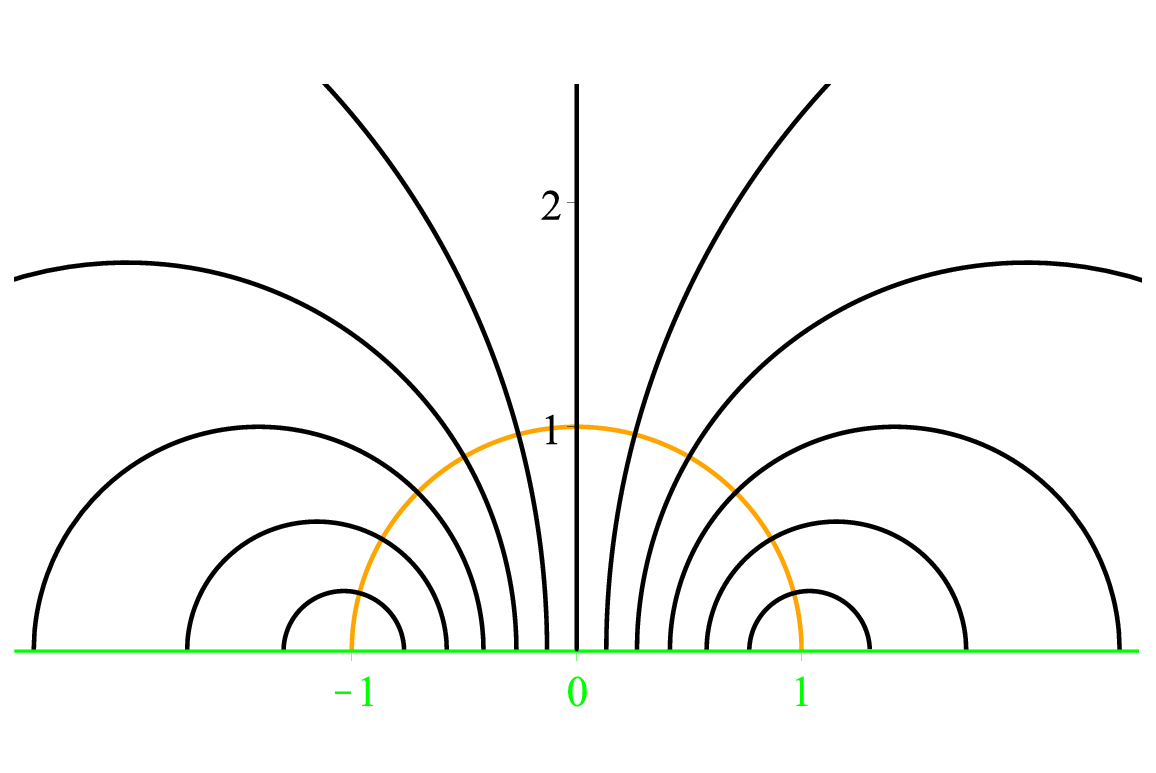}
    \caption*{\ref{fig:figB03}(b)  complex / Poincar\'e half-plane}
  \end{subfigure}
\phantomcaption
\plabel{fig:figB03}
\end{figure}
\begin{figure}[H]
  \ContinuedFloat
  \begin{subfigure}[b]{2in}
    \includegraphics[width=2in]{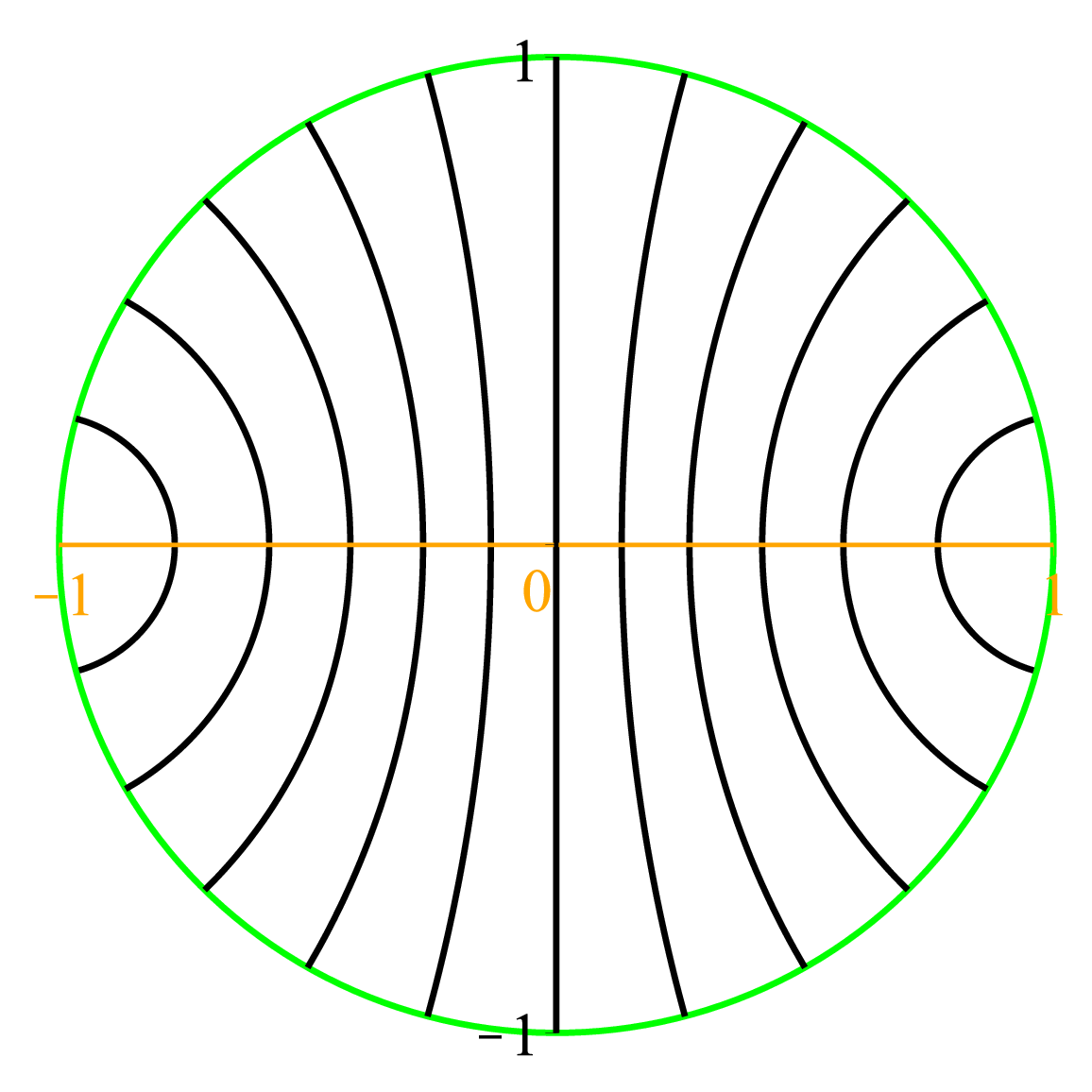}
    \caption*{\ref{fig:figB03}(c)   Poincar\'e disk model}
  \end{subfigure}
  \begin{subfigure}[b]{3in}
    \includegraphics[width=3in]{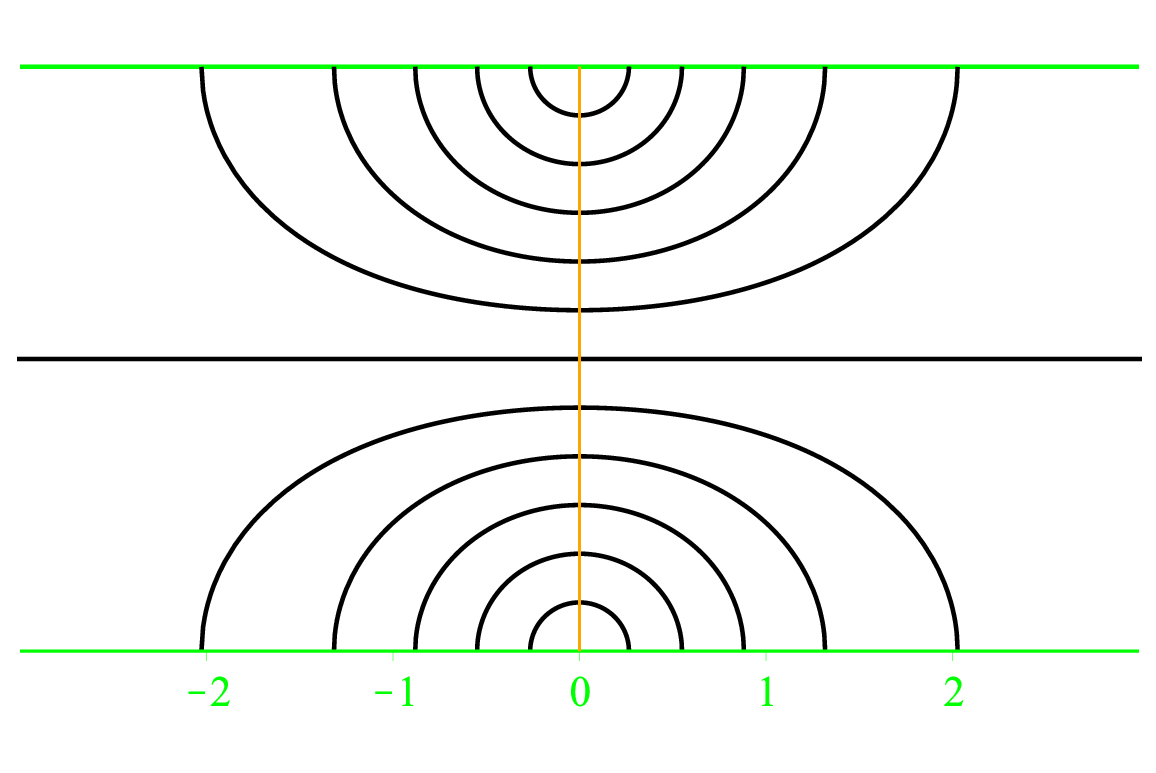}
    \caption*{\ref{fig:figB03}(d) log complex / logPHP}
  \end{subfigure}
\phantomcaption
\end{figure}

For $t=\frac\pi4$, they yield (Figure \ref{fig:figB04}):
\begin{figure}[H]
   \begin{subfigure}[b]{2in}
    \includegraphics[width=2in]{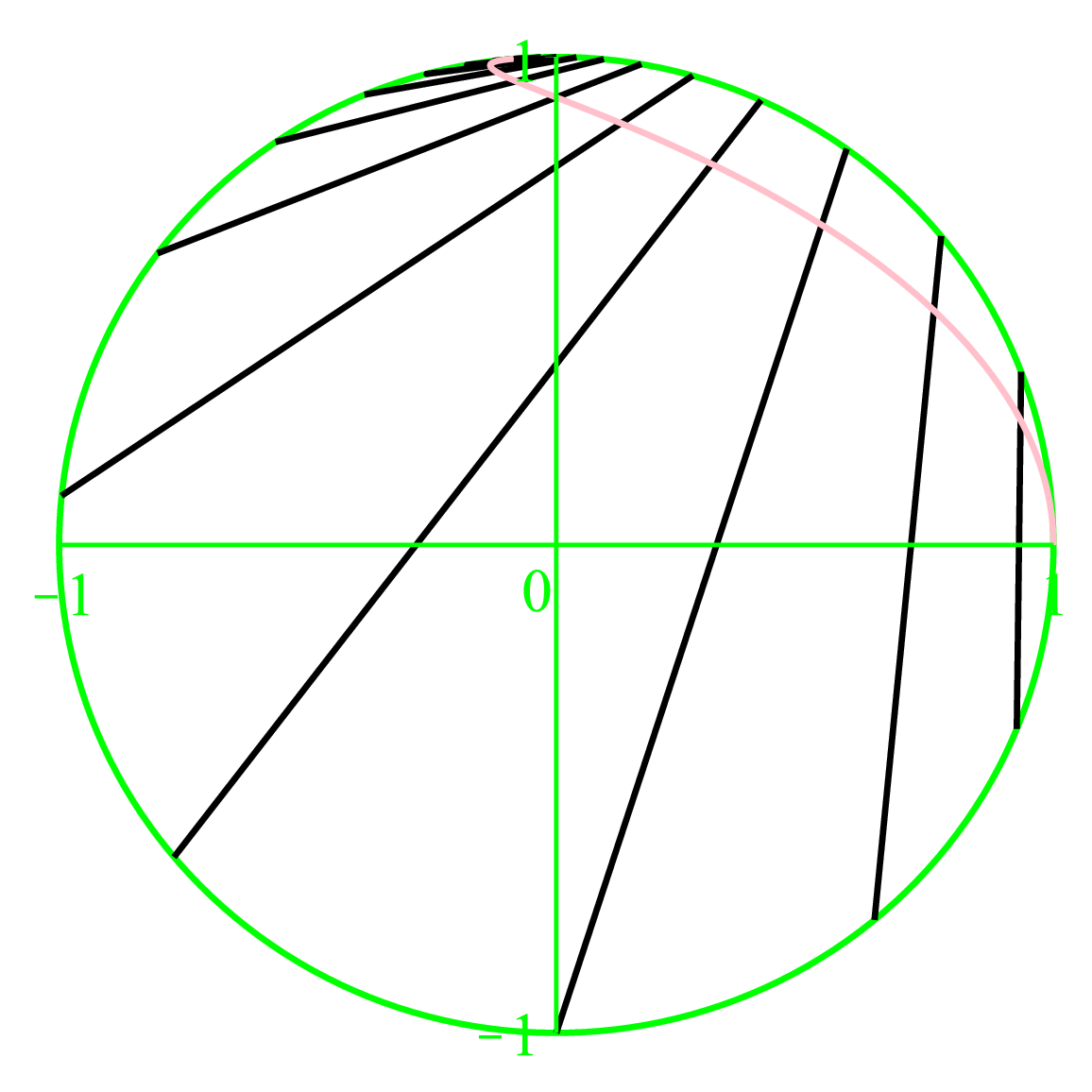}
    \caption*{Fig. \ref{fig:figB04}(a)   CKB model}
  \end{subfigure}
  \begin{subfigure}[b]{3in}
    \includegraphics[width=3in]{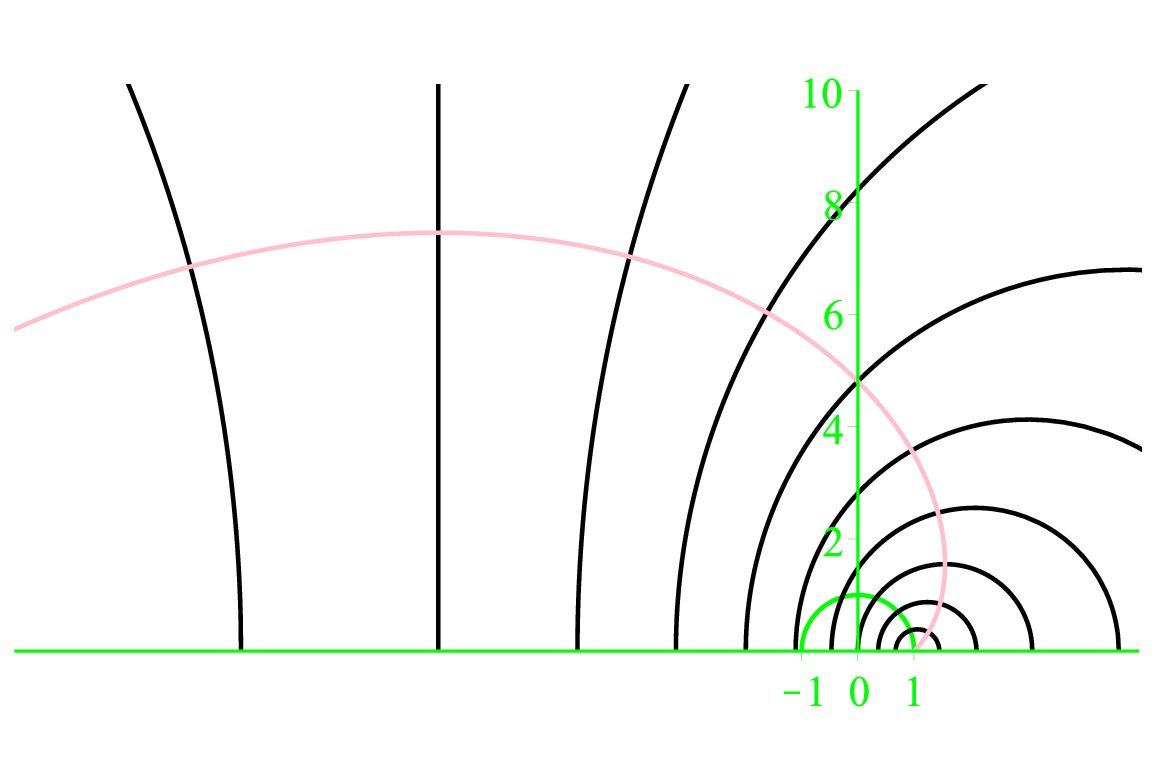}
    \caption*{\ref{fig:figB04}(b)  complex / Poincar\'e half-plane}
  \end{subfigure}
\phantomcaption
\plabel{fig:figB04}
\end{figure}
\begin{figure}[H]
  \ContinuedFloat
  \begin{subfigure}[b]{2in}
    \includegraphics[width=2in]{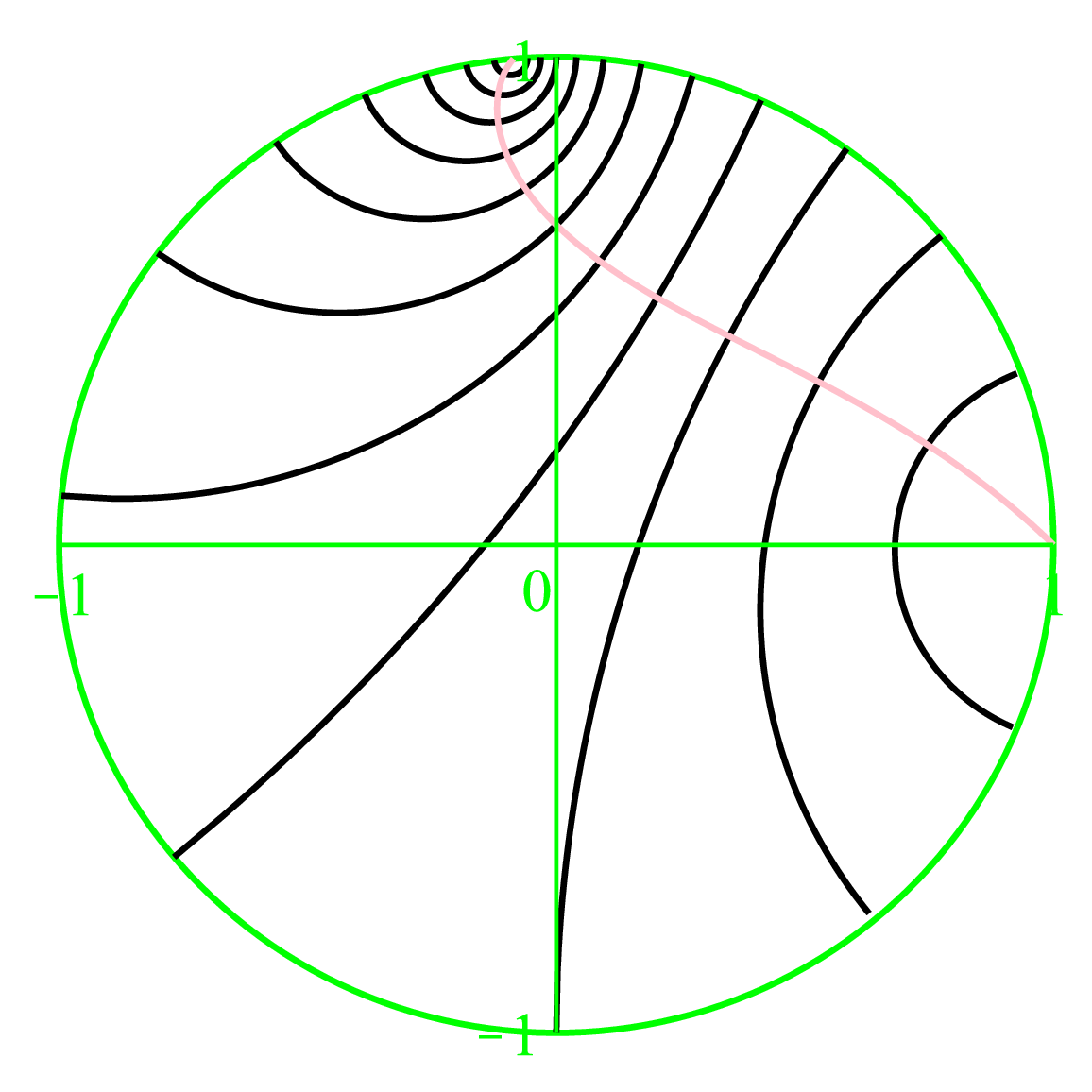}
    \caption*{\ref{fig:figB04}(c)   Poincar\'e disk model}
  \end{subfigure}
  \begin{subfigure}[b]{3in}
    \includegraphics[width=3in]{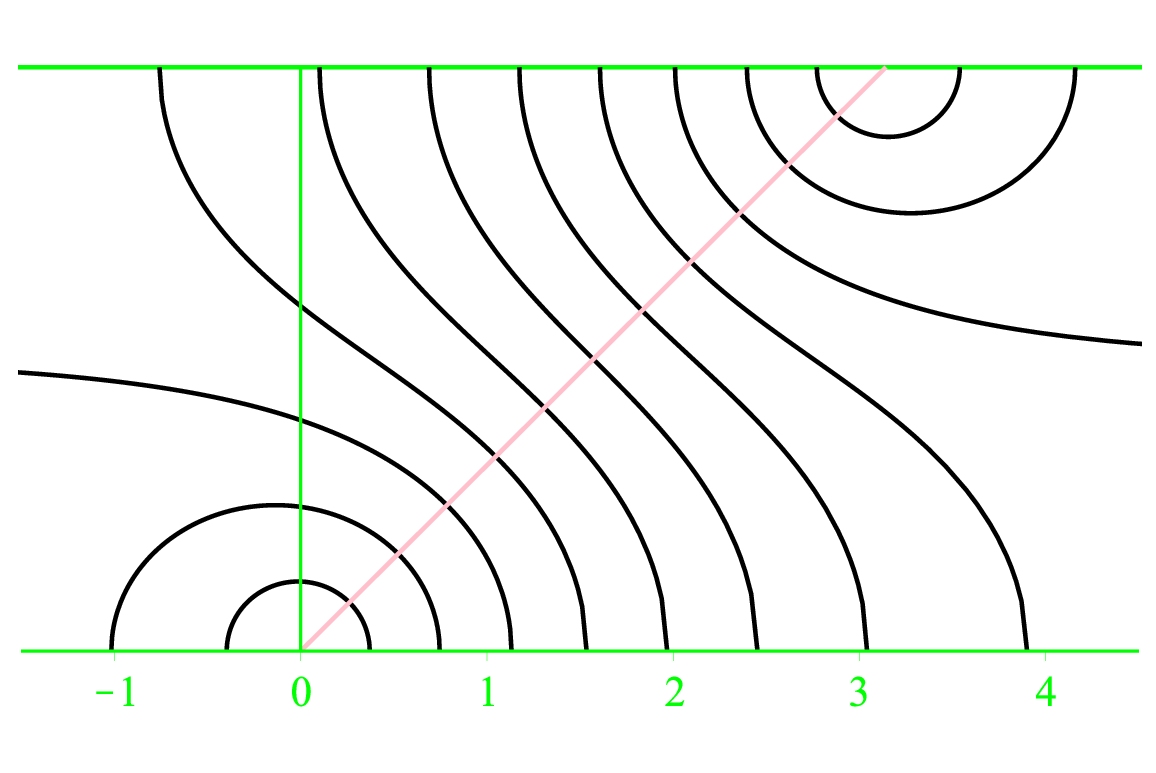}
    \caption*{\ref{fig:figB04}(d) log complex / logPHP}
  \end{subfigure}
\phantomcaption
\end{figure}

Suppose that the tangent lines are parametrized consistently by $\tau\in[-1,1]$ such that
$L_{p,t}(\tau)\equiv L_t(p,\tau)$ yields a regular mapping with respect to $(p,\tau)\in \left(0,\frac\pi{\sin t}\right)\times [-1,1]$.
Let $F_t(p,\tau)=\log L_t(p,\tau)$ (extended to the boundary in a natural manner).
Then, from simple plane-topological considerations, it is sufficient to prove that
\[\proj^\bot_{\frac{\partial F_t}{\partial \tau}}\frac{\partial F_t}{\partial p}\geq 1\]
holds in general. Equivalently,
\[\frac{\left| \dfrac{\partial F_t}{\partial \tau}\wedge \dfrac{\partial F_t}{\partial p} \right|}
{\left| \dfrac{\partial F_t}{\partial \tau}\right|}
\geq1,\qquad\text{ or, yet, equivalently }\qquad\frac{\left| \dfrac{\partial F_t}{\partial \tau}\wedge \dfrac{\partial F_t}{\partial p} \right|^2}
{\left| \dfrac{\partial F_t}{\partial \tau}\right|^2}
\geq1.\]
Taking into account that we deal with logarithmic derivatives, this is
\[\frac{\left| \dfrac{\partial L_t}{\partial \tau}\wedge \dfrac{\partial L_t}{\partial p} \right|^2}{\left| L_t\right|^2\left|\dfrac{\partial L_t}{\partial \tau}\right|^2}\geq1.\]
This is sufficient to prove in the interior. (In fact, the boundary behaviour is particularly nice:
the Poincar\'e half-plane model and the conformality of $\log$ show that the curves are perpendicular to the boundary, except the two critical points.)

In order to express the formulas with a greater efficiency, let use the abbreviations
\begin{equation}P=\mathrm e^{p\cos t}, \qquad C=\tan \frac{p\sin t}2,\qquad T=\tan \frac t2.
\plabel{eq:PCT}\end{equation}

Tangent line are the easiest to parametrize in the CKB model.
(For later reference, we include some comments regarding the orientation of the curves.)
 Then
\[\tfrac{\mathrm{CKB}}{\mathrm{PH}}\circ\gamma(t)=\left(\frac{\cos(p\sin t)}{\cosh(p\cos t)},\tanh (p\cos t) \right).\]
(This is a downward oriented curve.)
The corresponding tangent $h$-line at $\tfrac{\mathrm{CKB}}{\mathrm{PH}}\circ\gamma(t)$ is given by
\begin{align*}
\tilde H_{p,t}(\tau)=&\,\frac{1-\tau}2\cdot\left(\frac{-2P(C^2T^2-1)}{P^2(CT+1)^2+(CT-1)^2},\frac{P^2(CT+1)^2-(CT-1)^2}{P^2(CT+1)^2+(CT-1)^2}   \right)
\\&+\frac{1+\tau}2\cdot\left(\frac{-2P(C^2-T^2)}{P^2(C-T)^2+(C+T)^2},\frac{P^2(C-T)^2-(C+T)^2}{P^2(C-T)^2+(C+T)^2}
\right),
\end{align*}
where $\tau\in[-1,1]$. (This is also downward oriented.) The tangency occurs at
\[\tau=-\frac{P^2(CT+1)^2-P^2(C-T)^2+(CT-1)^2-(C+T)^2  }{(P^2+1)(C^2+1)(T^2+1)}.\]

We could proceed using this $\tilde H_{p,t}(\tau)$ perfectly well, however, the formulas would be quite long.
Thus, for the moment we abandon the cases $t=0$ and $t=\pi$, and reparametrize to
\begin{align*}
H_{p,t}(\tau)=&\biggl(\frac{-2P\left(\frac{1-\tau}2(C^2T^2-1)+\frac{1+\tau}2(C^2-T^2) \right)  }{
\frac{1-\tau}2(  P^2(CT+1)^2+(CT-1)^2)+\frac{1+\tau}2(  P^2(C-T)^2+(C+T)^2)},\\&\,
\frac{ \frac{1-\tau}2(  P^2(CT+1)^2-(CT-1)^2)+\frac{1+\tau}2(  P^2(C-T)^2-(C+T)^2)}{
\frac{1-\tau}2(  P^2(CT+1)^2+(CT-1)^2)+\frac{1+\tau}2(  P^2(C-T)^2+(C+T)^2)}\biggr),
\end{align*}
where $\tau\in[-1,1]$. This is tangent at $\tau=0$.
(The endpoints are the same as before, in particular, this is still downward oriented.
Even more relevantly, one can say that $\tfrac{\mathrm{CKB}}{\mathrm{PH}}(\exp\Dbar(0,p) \cap \upper )$
lies on left of the oriented tangent line $H_{p,t}(\tau)$.)

Hence $L_{p,t}(\tau)=\left(\tfrac{\mathrm{CKB}}{\mathrm{PH}}\right)^{-1}(H_{p,t}(\tau))$ yields the desired
parametrization in the PH model. In concrete terms,
\[L_{p,t}(\tau)=\left(\frac{-P\left(\frac{1-\tau}2(  C^2T^2-1))+\frac{1+\tau}2(C^2-T^2)\right)}{\frac{1-\tau}2  (CT-1)^2+\frac{1+\tau}2(C+T)^2},
\frac{\sqrt{1-\tau^2}PC(1+T^2)}{\frac{1-\tau}2  (CT-1)^2+\frac{1+\tau}2(C+T)^2}\right).\]
(As $\tfrac{\mathrm{CKB}}{\mathrm{PH}}$ is orientation preserving itself, it is still true that
$\exp\Dbar(0,p) \cap \upper$ lies on left of the oriented tangent $h$-line $L_{p,t}(\tau)$.)

Then direct computation yields
\[\frac{\left| \dfrac{\partial  L_t}{\partial \tau}\wedge \dfrac{\partial  L_t}{\partial p} \right|^2}{\left|  L_t\right|^2\left|\dfrac{\partial  L_t}{\partial \tau}\right|^2}
-1=\frac{4\tau^2C^2T^2
}{\left(\frac{1-\tau}2  (CT-1)^2+\frac{1+\tau}2(C+T)^2\right)\left(\frac{1-\tau}2  (CT+1)^2+\frac{1+\tau}2(C-T)^2\right)},\]
which is indeed nonnegative.
Thus the same kind of behaviour  applies to $\tilde L_{p,t}(\tau)=\left(\tfrac{\mathrm{CKB}}{\mathrm{PH}}\right)^{-1}(\tilde H_{p,t}(\tau))$, which extends to $t=0,\pi$ by
continuity. (These latter are the norm and co-norm cases anyway, which are also easy to see directly.)
\end{proof}
\end{theorem}
\begin{commentx}
\begin{remark}\plabel{rem:CRrangeInv}
In the previous proof,
\[\frac{\left| \dfrac{\partial \tilde L_t}{\partial \tau}\wedge \dfrac{\partial \tilde L_t}{\partial p} \right|^2}{\left| \tilde L_t\right|^2\left|\dfrac{\partial \tilde L_t}{\partial \tau}\right|^2}
-1=\]
\[\frac1{ \frac{1-\tau}2(CT-1)^2( P^2(C-T)^2+(C+T)^2) + \frac{1+\tau}2(T+C)^2(P^2(CT+1)^2+(CT-1)^2)  }{}\]
\[\cdot\frac1{\frac{1-\tau}2(CT+1)^2( P^2 (C-T)^2+(C+T)^2) + \frac{1+\tau}2(C-T)^2(P^2(CT+1)^2+(CT-1)^2)}\]
\[\cdot4C^2T^2
\cdot\left(\frac{1-\tau}2 (P^2(C-T)^2+(C+T)^2)-\frac{1+\tau}2( P^2(CT+1)^2+(CT-1)^2 )     \right)^2.\]
\qedremark
\end{remark}
\end{commentx}
\snewpage
\section{Growth estimates for the Magnus expansion}
\plabel{sec:MagnusGrowth}
As a technique, we will estimate $\log A$ from $\CRext(A)$.
This is certainly doable, but the computations are not always simple.

The expression of $\log$ can be rewritten ($s=\frac{\lambda}{\lambda-1}$) as
\begin{equation}
\log A=\int^{0}_{s=-\infty}\frac{A-\Id}{(1-s)(A-s\Id)}\, \mathrm ds.\plabel{eq:logdef3}
\end{equation}
(This refers more directly to the resolvents.)

The resolvent term can be estimated already in the traditional angular localization \eqref{eq:mapcrude},
thus, as warm-up, we start with that.
\begin{theorem}\plabel{th:angestimate}
Assume $\CRext(A)\subset \exp \{z\in\mathbb C\,:\,|\Rea z|\leq p,|\Ima z|\leq p \}$, $\pi/2< p<\pi$.

Then
\[\hat f:\left[\frac{\mathrm e^{p}}{\cos p},\frac{\mathrm e^{-p}}{\cos p}\right]\rightarrow
\left[\frac{-1+\mathrm e^{p}\cos p}{1-\mathrm e^{-p}\cos p},\frac{-1+\mathrm e^{-p}\cos p}{1-\mathrm e^{p}\cos p}\right]\]
\[c\mapsto\hat f(c)=\frac{c(1-c\cos^2 p)}{1-c}\]
is monotone increasing diffeomorphism. We claim,

for $c\in \left[\frac{\mathrm e^{p}}{\cos p},\frac{\mathrm e^{-p}}{\cos p}\right]$,
\[\left\|\frac{A-\Id}{A-\hat f(c)\Id}\right\|_2\leq\frac{c-1}{c\sin p};\]
for $s\in\left[\frac{-1+\mathrm e^{-p}\cos p}{1-\mathrm e^{p}\cos p}, 0\right]$,
\[\left\|\frac{A-\Id}{A-s\Id}\right\|_2\leq  \frac{\sqrt{((\mathrm e^{-p}\cos p)-1)^2+( \mathrm e^{-p}\sin p)^2  } }{\sqrt{((\mathrm e^{-p}\cos p)-s)^2+( \mathrm e^{-p}\sin p)^2  }  };\]
for $s\in\left(\infty,\frac{-1+\mathrm e^{p}\cos p}{1-\mathrm e^{-p}\cos p}\right]$,
\[\left\|\frac{A-\Id}{A-s\Id}\right\|_2\leq  \frac{\sqrt{((\mathrm e^{p}\cos p)-1)^2+( \mathrm e^{p}\sin p)^2  } }{\sqrt{((\mathrm e^{p}\cos p)-s)^2+( \mathrm e^{p}\sin p)^2  }  }.\]
\snewpage
\begin{proof}
Firstly, let us estimate  $\|(A-c\Id)^{-1}\|_2$ for $c\in(\infty,0]$.
(Let us note that, by Lemma \ref{lem:spectrange15}, the inverse exists.)
This involves computing the inverse of the co-norm of $A-c\Id$.
We have to take the smallest $c$ centered circle which still intersects $R=\exp \{z\in\mathbb C\,:\,|\Rea z|\leq p,|\Ima z|\leq p \}$.

It is easy to that for $c\in\left[\frac{\mathrm e^{p}}{\cos p},\frac{\mathrm e^{-p}}{\cos p}\right]$ this ``co-norm circle'' is tangent to the radial boundary of  $R$,
and it is of radius $-c\sin p$; consequently,
\[\|(A-c\Id)^{-1}\|_2\leq\frac1{-c\sin p}.\]
For $c\in\left[\frac{\mathrm e^{-p}}{\cos p},0\right]$, ``co-norm circle'' goes though the corner $\left(\mathrm e^{-p}\cos p, \mathrm e^{-p}\sin p \right)$,
thus it of radius $\sqrt{((\mathrm e^{-p}\cos p)-c)^2+( \mathrm e^{-p}\sin p)^2  }$;
consequently,
\[\|(A-c\Id)^{-1}\|_2\leq\frac1{\sqrt{((\mathrm e^{-p}\cos p)-c)^2+( \mathrm e^{-p}\sin p)^2  }  }.\]
For $c\in\left(\infty,\frac{\mathrm e^{p}}{\cos p},0\right]$, ``co-norm circle'' goes though the corner $\left(\mathrm e^{p}\cos p, \mathrm e^{p}\sin p \right)$,
thus it of radius $\sqrt{((\mathrm e^{p}\cos p)-c)^2+( \mathrm e^{p}\sin p)^2  }$;
consequently,
\[\|(A-c\Id)^{-1}\|_2\leq\frac1{\sqrt{((\mathrm e^{p}\cos p)-c)^2+( \mathrm e^{p}\sin p)^2  }  }.\]

However we want to estimate
$\left\|\frac{A-\Id}{A-s\Id}\right\|_2$.
Assume that $|\mathbf x|_2=1$.
According to the Lemma \ref{lem:conf},
\[\left|\frac{A-\Id}{A-\lambda\Id}\mathbf x\right|_2\in \left\{\frac{|\omega-1|}{|\omega-s|} :\omega\in\CRext(A)\right\}. \]

So, we can estimate $\left|\frac{A-\Id}{A-s\Id}\right|_2$ as follows:
Take the Apollonian circles relative to $s$ and $1$, and take the one closest to $s$ but which still intersects $R$.
Then the characteristic ratio of this Apollonian circle provides an upper estimate.

If $c\in\left[\frac{\mathrm e^{p}}{\cos p},\frac{\mathrm e^{-p}}{\cos p}\right]$, then the conjugate (ie. inverse) of $1$ with respect to given circle is
\[s=\hat f(c)=\frac{c(1-c\cos^2 p)}{1-c}.\]
(It is easy to see that $\hat f$ is monotone increasing with the given range.)
In this case the Apollonian ratio is $\frac{c-1}{c\sin p}$.
For $s\in\left[\frac{-1+\mathrm e^{-p}\cos p}{1-\mathrm e^{p}\cos p}, 0\right]$, the Apollonian circle goes though the corner $\left(\mathrm e^{-p}\cos p, \mathrm e^{-p}\sin p \right)$;
for $s\in\left(\infty,\frac{-1+\mathrm e^{p}\cos p}{1-\mathrm e^{-p}\cos p}\right]$,  the Apollonian circle goes though the corner $\left(\mathrm e^{p}\cos p, \mathrm e^{p}\sin p \right)$.
This is not entirely trivial as Apollonian circles with center outside of $(-\infty,0]$ enter into the picture, but it can be checked.
(It is easier to establish for small  $|s|$, then use inversion symmetry to the unit circle).
\end{proof}
\end{theorem}

\snewpage
\begin{theorem}
\plabel{th:angloge}
Assume $\CRext(A)\subset \exp \{z\in\mathbb C\,:\,|\Rea z|\leq p,|\Ima z|\leq p \}$, $\pi/2< p<\pi$.
\begin{align*}
\|\log A\|_2\leq&\int_{s=-\infty}^{\frac{-1+\mathrm e^{p}\cos p}{1-\mathrm e^{-p}\cos p}}\frac1{1-s}\cdot
\frac{\sqrt{((\mathrm e^{p}\cos p)-1)^2+( \mathrm e^{p}\sin p)^2  } }{\sqrt{((\mathrm e^{p}\cos p)-s)^2+( \mathrm e^{p}\sin p)^2  }  }\,\mathrm ds
\\&+\int_{c=\frac{\mathrm e^{p}}{\cos p} }^{\frac{\mathrm e^{-p}}{\cos p} }\frac1{1-\hat f(c)}\cdot\frac{c-1}{c\sin p}\,\mathrm d\hat f(c)
\\&+\int_{s=\frac{-1+\mathrm e^{-p}\cos p}{1-\mathrm e^{p}\cos p} }^{0}
\frac{\sqrt{((\mathrm e^{-p}\cos p)-1)^2+( \mathrm e^{-p}\sin p)^2  } }{\sqrt{((\mathrm e^{-p}\cos p)-s)^2+( \mathrm e^{-p}\sin p)^2  }  }\,\mathrm ds
\\=&\left [
{\rm arcsinh} \left({\frac {s{{\rm e}^{-p}}-\cos \left( p \right) s-
\cos \left( p \right) +{{\rm e}^{p}}}{ \left( 1-s \right) \sin \left(
p \right) }}\right)
\right]_{s=-\infty}^{\frac{-1+\mathrm e^{p}\cos p}{1-\mathrm e^{-p}\cos p}}
\\&+\left[2\,{\rm arctanh} \left({\frac {\sin \left( p \right) }{c\, \left( \cos
 \left( p \right)  \right) ^{2}-1}}\right)-{\frac {\ln  \left( -c
 \right) }{\sin \left( p \right) }}
\right]_{c=\frac{\mathrm e^{p}}{\cos p} }^{\frac{\mathrm e^{-p}}{\cos p} }
\\&+\left[{\rm arcsinh} \left({\frac {s{{\rm e}^{p}}-\cos \left( p \right) s-
\cos \left( p \right) +{{\rm e}^{-p}}}{ \left( 1-s \right) \sin
 \left( p \right) }}\right)
\right]_{s=\frac{-1+\mathrm e^{-p}\cos p}{1-\mathrm e^{p}\cos p} }^{0}
\\=&2\,{\frac {p}{\sin \left( p \right) }}+2\,{\rm arctanh} \left({\frac {
\sin \left( p \right)  \left( {{\rm e}^{2\,p}}-1 \right) }{-{{\rm e}^{
2\,p}}+2\,\cos \left( p \right) {{\rm e}^{p}}-1}}\right)
\\&2\cdot\left({\rm arcsinh} \left({\frac {-\cos \left( p \right) +{{\rm e}^{-p}}}{
\sin \left( p \right) }}\right)+{\rm arcsinh} \left({\frac {\sin
 \left( p \right)  \left( {{\rm e}^{-2\,p}}-1 \right) }{{{\rm e}^{-2\,
p}}\cos \left( p \right) +\cos \left( p \right) -2\,{{\rm e}^{-p}}}}
\right)\right)
\\=&\frac{2\pi}{\pi-p}-2\log\frac1{\pi-p}+O(1)\qquad\text{as}\quad p\nearrow\pi.
\end{align*}
\begin{proof}
This is the direct application of the estimates in Theorem \ref{th:angestimate} to \eqref{eq:logdef3}.
\end{proof}
\end{theorem}
There are many reasons why Theorem \ref{th:angloge} cannot work as a really good estimate for the Magnus expansion,
but an obvious one is that it uses a suboptimal localization of the conformal range.
\snewpage
\begin{theorem}\plabel{th:monoestimate}
Assume $\CRext(A)\subset \exp \Dbar(0,p)$, $0< p<\pi$.
Then the strictly increasing diffeomorphism
\[f_p:[0,\pi]\rightarrow\left[-\frac{1-\mathrm e^{p}(1-p)}{1-\mathrm e^{-p}(1+p)},
-\frac{1-\mathrm e^{-p}(1+p)}{1-\mathrm e^{p}(1-p)}\right] \Subset[-\infty,0]\]
given by
\[f_p(t)=-\frac{\sin t-\mathrm e^{p\cos t}\sin(t-p\sin t)}{\sin t-\mathrm e^{-p\cos t}\sin(t+p\sin t)}\]
has the property that for any $t\in(0,\pi)$,
\[\left\|\frac{A-\Id}{A-f_p(t)\Id}\right\|_2\leq\frac{|\gamma_p(t)-1|}{|\gamma_p(t)-f_p(t)|}=
\frac{\sin t-\mathrm e^{-p\cos t}\sin(t+p\sin t)}{\sin(p\sin t)}.\]

For $s\in(-\infty,f_p(0)] \cup (\mathrm e^p,+\infty)$,
\[\left\|\frac{A-\Id}{A-s\Id}\right\|_2\leq\frac{|\gamma_p(0)-1|}{|\gamma_p(0)-s|}=
\frac{|\mathrm e^p-1|}{|\mathrm e^p-s|}.\]

For $s\in[f_p(\pi)\mathrm e^{-p})$,
\[\left\|\frac{A-\Id}{A-s\Id}\right\|_2\leq\frac{|\gamma_p(\pi)-1|}{|\gamma_p(\pi)-s|}=
\frac{|\mathrm e^{-p}-1|}{|\mathrm e^{-p}-s|}.\]
\begin{proof}[Remark]
The existence of $\frac{A-\Id}{A-s\Id}$ for $s\in\mathbb R\setminus [\mathrm e^{-p},\mathrm e^{p} ] $ follows from Lemma \ref{lem:spectrange15},
although applied in the setting of Theorem \ref{th:CRrange} also follows from \eqref{eq:spectrange}.
\renewcommand{\qedsymbol}{}
\end{proof}
\begin{proof}
Again, estimating $\left\|\frac{A-\Id}{A-s\Id}\right\|_2$ involves taking Apollonian circles with respect to the the point $s$ and $1$,
taking the one closest to $s$ but still intersecting $\exp \Dbar(0,p)$, and computing the Apollonian ratio.

This leads to considering circles (and lines) which are tangent to the curve
\[\gamma_p(t)=\mathrm e^{p\cos(t)}\cos(p\sin t)+\mathrm i \mathrm e^{p\cos(t)}\sin(p\sin t) )\]
($t\in[0,\pi]$),  and their center is on  the real axis (or in the infinity).
If $t\in(0,\pi)$, then the normal line at $\gamma_p(t)$ intersects the real axis
at \[C_p(t)=\frac{\mathrm e^{p\cos t}\sin t }{\sin(t+p\sin t)},\]
the center of the circle. This leads to radius
\[r_p(t)= \frac{\mathrm e^{p\cos t}\sin (p\sin t) }{\sin(t+p\sin t)}.\]
(The sign counts the touching orientation to $\gamma_p$.)
Taking the inverse of $1$, relative to the circle above, leads to the Apollonian pole
\[f_p(t)=-\frac{\sin t-\mathrm e^{p\cos t}\sin(t-p\sin t)}{\sin t-\mathrm e^{-p\cos t}\sin(t+p\sin t)}\]
conjugate to $1$.
The functions $C_p$ and $r_p$ are singular, but $f_p$ is not.
This can be seen from
\[ \sin t-\mathrm e^{p\cos t}\sin(t-p\sin t)=\int_{q=0}^p  \mathrm e^{q\cos t}\sin (q\sin t)\,\mathrm dq>0,\]
\[\sin t-\mathrm e^{-p\cos t}\sin(t+p\sin t)=\int_{q=0}^p  \mathrm e^{-q\cos t}\sin (q\sin t)\,\mathrm dq>0.\]
In fact, $f_p$ is strictly increasing.
Indeed,
\[f'_p(t)=\frac{(\sin(p\sin t)-(p \sin t) \cos (p \sin t))(\mathrm e^{p\cos t} +\mathrm e^{-p\cos t}-2\cos (p \sin t))}{(\sin t-\mathrm e^{-p\cos t}\sin(t+p\sin t))^2}=\]
\[=\frac{(\sin t)^2\left( \int_{q=0}^p q\sin (q\sin t)\,\mathrm dq \right)
\cdot 2(\cosh (p\cos t)-\cos (p \sin t))
}{\left( \int_{q=0}^p \mathrm e^{-q\cos t}\sin (q\sin t)\,\mathrm dq \right)^2}>0.\]
It is easy to see that the range of $f_p$ on $(0,\pi)$ is
\[(f_p(0+),f_p(\pi-))=\left(-\frac{1-\mathrm e^{p}(1-p)}{1-\mathrm e^{-p}(1+p)},
-\frac{1-\mathrm e^{-p}(1+p)}{1-\mathrm e^{p}(1-p)}\right).\]
The characteristic ratio belonging to the relevant Apollonian circle is
\[\chi_p(t)=\frac{|\gamma_p(t)-1|}{|\gamma_p(t)-f_p(t)|}=
\frac{\sin t-\mathrm e^{-p\cos t}\sin(t+p\sin t)}{\sin(p\sin t)}.\]
The values $t=0$ and $t=\pi$ exceptional, because tangent circles there always have their
centers on the real axis.

Let $s\in\mathbb R\setminus [\gamma_p(\pi),\gamma_p(0)]$. Consider the Apollonian circles between $s$ and $1$,
and consider the one closest to $s$ but still touching $\gamma_p$.
From geometrical considerations (the injectivity of $f_p$) we can devise that the
closest touching circle touches at
\begin{align}\gamma_p(0)&\quad\text{if}\quad  s\in(-\infty,f_p(0+)]\cup (\gamma(0),+\infty),\notag\\
\gamma_p(t)&\quad\text{if}\quad  s=f_p(t)\in(f_p(0+),f_p(\pi-)),\notag\\
\gamma_p(\pi)&\quad\text{if}\quad  s\in[f_p(\pi-),\gamma_p(\pi)].\notag
\end{align}
These observations together yield the statement of the theorem.
\end{proof}
\end{theorem}

Through \eqref{eq:logdef3}, the theorem above can be used to estimate $\log A$:
\snewpage
\begin{theorem}\plabel{th:Lestimate}
If $\CRext(A)\subset \exp \Dbar(0,p)$, $0< p<\pi$, then
\begin{equation}\|(\log A)\|_2\leq H(p),\plabel{eq:Hestimate}\end{equation}
where
\begin{equation}
H(p)=p-2\log\left(2\cosh \frac p2-\frac2p\sinh\frac p2\right)
+\int_{t=0}^\pi HH(p,t)\,\mathrm dt\plabel{eq:hexp}
\end{equation}
with
\begin{equation}HH(p,t)=\frac{(\sin(p\sin t)-(p \sin t) \cos (p \sin t))
(\mathrm e^{p\cos t} +\mathrm e^{-p\cos t}-2\cos (p \sin t))}{
(\sin(p\sin t))(2\sin t+\mathrm e^{p\cos t}\sin(-t+p\sin t)-\mathrm e^{-p\cos t}\sin(t+p\sin t))}.\plabel{eq:hexp2}
\end{equation}
$H(p)$ and $HH(p,t)$ are positive and finite for $0<p<\pi$.

The statement (trivially) extends to $p=0$ with $H(p)=0$.
\begin{proof}[Remark] $HH(p,t)$  can be rewritten as
\begin{equation}
HH(p,t)=
\frac{(p^2\sin t)\left(\frac1{p^3\sin t} \int_{q=0}^p q\sin (q\sin t)\,\mathrm dq \right)
\cdot \left(\frac1{p^2}(\cosh (p\cos t)-\cos (p \sin t))\right)}{\left(\frac{\sin(p\sin t)}{p\sin t}\right)
\left( \frac1{p^2\sin t}\int_{q=0}^p \cosh(q\cos t)\sin (q\sin t)\,\mathrm dq \right)}.
\end{equation}
From the power series expansion, it is easy to see that the expressions in the big parentheses are actually entire functions of
$p$ and $t$. Moreover, one can see that these entire functions are positive for $(p,t)\in[0,\pi)\times[0,\pi]$.
In fact, what prevents the smooth extension to $(p,t)\in[0,\pi]\times[0,\pi]$ is only the singularity in $\frac{\sin(p\sin t)}{p\sin t}$.
\renewcommand{\qedsymbol}{}
\end{proof}
\begin{proof} Continuing the proof of the previous theorem,
this provides the estimate
\begin{align}|(\log A)\mathbf x|_2\leq&
\int_{s=-\infty}^{f_p(0+)} \frac{|\gamma_p(0)-1|}{(1-s)|\gamma_p(0)-s|} \,\mathrm ds+
\int_{t=0}^\pi\frac{\chi_p(t)}{1-f_p(t)}\,\mathrm df_p(t)+\notag\\
&+\int_{s=f_p(\pi-)}^0 \frac{|\gamma_p(\pi)-1|}{(1-s)|\gamma_p(\pi)-s|} \,\mathrm ds.\notag
 \end{align}
The first and third integrals expands as
\[\int_{s=-\infty}^{f_p(0+)}{\frac {{{\mathrm e}^{p}}-1}{
\left( 1-s \right)  \left( {{\mathrm e}^{p}}-s \right) }}\,\mathrm ds
=\left[\log  \left( {\frac {{{\rm e}^{p}}-s}{1-s}} \right)\right]_{s=-\infty}^{f_p(0+)}
=\log\frac{p}{p-1+\mathrm e^{-p}(p+1)},
\]
\[\int_{s=f_p(\pi-)}^0 {\frac { 1-{{\rm e}^{-p}}}{ \left( 1-s \right)  \left( {{\rm e}^{-p}}-
s \right) }} \,\mathrm ds=\left[\log  \left( {\frac {1-s}{{{\rm e}^{-p}}-s}} \right) \right]_{s=f_p(\pi-)}^0
=\log\frac{p}{p-1+\mathrm e^{-p}(p+1)}.\]
Note that
\[\log\frac{p}{p-1+\mathrm e^{-p}(p+1)}=
\frac p2-\log\frac{\mathrm e^{\frac p2}(p-1)+\mathrm e^{-\frac p2}(p+1)}{p}.\]
The integrand in the second integral expands as indicated in \eqref{eq:hexp}.
\end{proof}
\end{theorem}

\begin{remark}
\plabel{rem:Hestimate}
The estimate \eqref{eq:Hestimate} is not sharp.
A simple indication for that is as follows:
For example, in the proof, applied to an individual vector $\mathbf x$ with $\|\mathbf x\|_2=1$,
we estimated $|A^{-1}\mathbf x|_2$ by $\mathrm e^{p}$, which belongs to
$A(A^{-1}\mathbf x): A^{-1}\mathbf x=\mathrm e^{-p}$,  i. e. $A^{-1}\mathbf x=\mathrm e^{p}\mathbf x$.
But then $|(\log A)\mathbf x|_2=|-p\mathbf x|_2=p<H(p)$ would hold.
In general, there is a penalty or gain (depending on the viewpoint) for approaching the real axis
in $\CRext(A)$, for which we have not accounted.
\qedremark
\end{remark}

\begin{theorem}\plabel{th:Hasympt}
(a) As $p\searrow0$,
\begin{equation}H(p)=p+\frac14{p}^{2}+{\frac {23}{864}}{p}^{4}+\mathrm O(p^6).\plabel{eq:san1}\end{equation}

(b) As $p\nearrow\pi$
\[H(p)=\frac{2\pi^2}{\sqrt{\pi^2-p^2}}+H_\pi+o(1)=\frac{\sqrt{2}\pi^{3/2}}{\sqrt{\pi-p}}+H_\pi+o(1)=
p\sqrt{\frac{\pi+p}{\pi-p}}+H_\pi+o(1),\]
where
\begin{equation}
H_\pi=\pi-2\log\left(2\cosh \frac \pi2-\frac2\pi\sinh\frac \pi2\right)
+\int_{t=0}^\pi \left(HH(\pi,t)-\frac{2}{\cos^2 t} \right)\,\mathrm dt\notag
\end{equation}
(and the integrand is actually a smooth function of $t$). Numerically, $H_\pi=-2.513\ldots$
\proofremark{As $p\nearrow\pi$,
\[p\sqrt{\frac{\pi+p}{\pi-p}}
\equiv\pi\sqrt{\frac{\pi+p}{\pi-p}}
\equiv\frac{\sqrt{2}\pi^{3/2}}{\sqrt{\pi-p}}
\equiv\frac{2\pi^2}{\sqrt{\pi^2-p^2}}\modu O(\sqrt{\pi-p})
.\]
}
\begin{proof} Consider \eqref{eq:hexp}.
One finds
\begin{equation}p-2\log\left(2\cosh \frac p2-\frac2p\sinh\frac p2\right)=p-{\frac {5}{12}}{p}^{2}+{\frac {49}{1440}}{p}^{4}+O \left( {p}^{6}\right).\plabel{eq:san2}\end{equation}
Regarding $HH(p,t)$, one  can see that
\begin{align}
\frac1{p^3\sin t} \int_{q=0}^p q\sin (q\sin t)\,\mathrm dq&=\frac13-\frac{\sin^{2} t  }{30}\,{p}^{2}+O(p^4),\notag\\
\frac{\cosh (p\cos t)-\cos (p \sin t)}{p^2}&=\frac12+ \frac{ \cos^{2} t  -  \sin^{2} t  }{24} \,{p}^{2}+O(p^4),\notag\\
\frac{\sin(p\sin t)}{p\sin t}&=1-\frac{\sin^{2}t}6\,{p}^{2}+O(p^4),\notag\\
\frac1{p^2\sin t}\int_{q=0}^p \cosh(q\cos t)\sin (q\sin t)\,\mathrm dq&=\frac12+{\frac {3\,   \cos^{2} t - \sin^{2} t }{24 }}\,{p}^{2}+O(p^4).\notag
\end{align}
Consequently,
\[HH(p,t)=\frac{\sin t}3p^2+{\frac {(2\sin^{2} t-5 \cos^{2} t )\sin t }{90}}{p}^{4}+O(p^6).\]
Integrating this for $t\in[0,\pi]$, it gives
\begin{equation}\int_{t=0}^\pi HH(p,t)\,\mathrm dt=\frac23\,{p}^{2}-{\frac {1}{135}}{p}^{4}+O(p^6).\plabel{eq:san3}\end{equation}
Adding \eqref{eq:san2} and \eqref{eq:san3} yields \eqref{eq:san1}.

(b) Notice that $\frac{\sin x}{(\pi^2-x^2)}$ is analytic function, which is positive
on $x\in[-\pi,\pi]$. Consequently,
$\frac{\sin(p\sin t)}{(\pi^2- p^2\sin^2 t)p\sin t }$ is an entire function of $p,t$ such that
it is positive for $(p,t)\in[0,\pi]\times[0,\pi]$. Hence
\[HH(p,t)=\frac{1}{\pi^2- p^2\sin^2 t } \widetilde{HH}(p,t),\]
where $\widetilde{HH}(p,t)$ is smooth
on  $(p,t)\in[0,\pi]\times[0,\pi]$.
Due to symmetry for $t\leftrightarrow\pi-t$, $\widetilde{HH}(p,t)- \widetilde{HH}(p,\pi/2)$ not only vanishes at
$t=\pi/2$ but $\cos^2 t$ can be factored out.
Thus
\[\widehat{HH}(p,t)=\frac{\widetilde{HH}(p,t)- \widetilde{HH}(p,\pi/2)}{\pi^2\cos^2 t}\]
can also be considered as a smooth function on  $(p,t)\in[0,\pi]\times[0,\pi]$.
Now we have
\[HH(p,t)=\frac{1}{\pi^2- p^2\sin^2 t } \widetilde{HH}(p,\pi/2)+
\frac{\pi^2\cos^2 t}{\pi^2- p^2\sin^2 t }\widehat{HH}(p,t).\]
For a fixed $p$ the first summand integrates to
\begin{multline}\notag
\int_{t=0}^\pi \frac{1}{\pi^2- p^2\sin^2 t } \widetilde{HH}(p,\pi/2)\,\mathrm dt=
\frac{\widetilde{HH}(p,\pi/2)}{\sqrt{\pi^2-p^2}}=\sqrt{\pi^2-p^2}\,\frac{\sin p-p\cos p}{\sin p}=\\
=\frac{2\pi^2}{\sqrt{\pi^2-p^2}}+o(1)=\frac{\sqrt{2}\pi^{3/2}}{\sqrt{\pi-p}}+H_\pi+o(1)=
p\sqrt{\frac{\pi+p}{\pi-p}}+H_\pi+o(1).
\end{multline}

The function $\frac{\pi^2\cos^2 t}{\pi^2- p^2\sin^2 t }=\frac{\pi^2- \pi^2\sin^2 t}{\pi^2- p^2\sin^2 t } $
is uniformly bounded by $0$ and $1$, and, in fact
\[\lim_{ p \nearrow\pi}\frac{\cos^2 t}{\pi^2- p^2\sin^2 t }=1 \qquad\text{ for  }t\in [0,\pi]\setminus\left\{\frac{\pi}2\right\} \]
pointwise.
Thus, by Lebesgue's dominated convergence theorem, the integral of the second summand is
\[\int_{t=0}^\pi \widehat{HH}(\pi,t)\,\mathrm dt+o(1).\]
Notice that $\widehat{HH}(\pi,t)$ is a smooth function.
Taking limit with $p\nearrow\pi$ we find that
\[\widehat{HH}(\pi,t)=HH(\pi,t)-\frac{2}{\cos^2 t}.\]
The numerical evaluation of $H_\pi$ can be realized by various methods.
\end{proof}
\end{theorem}
\begin{remark}
In Example \ref{ex:critical}, using the notation $p=\pi t$, we find
\[\int\|p/\pi\cdot \Phi\|_2=p\]
and
\begin{equation}
\|\mu_{\mathrm L}(p/\pi\cdot \Phi)\|_2=\sqrt2\pi^{3/2} (\pi-p)^{-1/2}-2\pi-\frac{\sqrt2}4\pi^{1/2} (\pi-p)^{1/2}+O(\pi-p),
\plabel{eq:lower1}
\end{equation}
as $p\nearrow\pi$.
Thus, despite Remark \ref{rem:Hestimate},
we find that the estimate in Theorem \ref{th:Lestimate} is already asymptotically sharp as $p\nearrow\pi$;
which is quite good, after all.
(Later, we will see slightly better examples and  estimates.)
\qedremark
\end{remark}

\snewpage
Before we provide an improved estimate, it is illuminating to  rephrase the statement of
Theorem \ref{th:monoestimate} with respect to
$\log$ rewritten by $\nu=2\lambda-1$,
\begin{equation}
\log A=\int_{\nu=-1}^1 \frac{A-\Id}{(1+\nu)\Id +(1-\nu)A }\,d\nu.\plabel{eq:logdef2}
\end{equation}
The theorem says, in the same context, that the strictly decreasing diffeomorphism
\[g_p=\frac{f_p+1}{f_p-1}:[0,\pi]\rightarrow
\left[-\frac{p\cosh p-\sinh p}{1+p\sinh p-\cosh p},\frac{p\cosh p-\sinh p}{1+p\sinh p-\cosh p}\right]\Subset[-1,1]\]
\[g_p(t)=\frac{\mathrm e^{p\cos t}\sin(-t+p\sin t)+\mathrm e^{-p\cos t}\sin(t+p\sin t)
}{2\sin t+\mathrm e^{p\cos t}\sin(-t+p\sin t)-\mathrm e^{-p\cos t}\sin(t+p\sin t)}\]
has the property that for $t\in(0,\pi)$,
\begin{multline}\left\|\frac{A-\Id}{(1+g_p(t))\Id+(1-g_p(t))A}\right\|_2\leq \psi_p(t)\equiv\\
\equiv\frac{2\sin t+\mathrm e^{p\cos t}\sin(-t+p\sin t)-\mathrm e^{-p\cos t}\sin(t+p\sin t)}{2\sin(p\sin t)};\plabel{eq:resmono}\end{multline}
etc.
Using the previous terminology, $f_p(\pi-t)=\frac1{f_p(t)}$ was apparent; here not only $g_p(\pi-t)=-g_p(t)$ holds but $\psi_p(\pi-t)=\psi_p(t)$.
This reflects the symmetry $A\leftrightarrow A^{-1}$. Considering the other cases, too, Theorem \ref{th:monoestimate}
provides estimates of $\left\| \frac{A-\Id}{(1+\nu)\Id+(1-\nu)A}\right\|$ for $\nu\in\left(-\coth p,\coth p\right)\Supset[-1,1]$.

\begin{theorem}\plabel{th:biestimate}
Let $\CRext(A)\subset \exp \Dbar(0,p)$, $0< p<\pi$.
Then the function
\[h_p:[0,\pi]\rightarrow\left(-\coth p,\coth p\right) \Supset[-1,1]\]
given by
\[h_p(t)={\frac {{{\rm e}^{p\cos  t }}\sin \left( p\sin  t
   \right) +{{\rm e}^{p\cos  t  }}\cos  t
 +\cos \left( t+p\sin  t   \right) }{{{\rm e}^{p
\cos  t  }}\sin \left( p\sin  t  \right) +{
{\rm e}^{p\cos  t  }}\cos  t  -\cos \left( t
+p\sin  t  \right) }}
\]
is a strictly decreasing diffeomorphism.
Let us also consider
\[\varpi_p(t)=-\frac14\,{\frac {{{\rm e}^{-p\cos t }} \left(
\left( {{\rm e}^{p
\cos t  }}
 \sin \left( p\sin  t
 \right)\right)^2-\left( {{\rm e}^{p
\cos  t  }}\cos  t    -\cos \left( t+p\sin  t
  \right)  \right) ^{2}
 \right) }{\sin \left( p\sin  t   \right) \cos \left( t+p
\sin  t   \right) }}
\]
and
\[\phi_p(t)=\frac14\,{\frac {{{\rm e}^{-p\cos t  }} \left( {{\rm e}^{p
\cos t  }}
 \sin \left( p\sin  t
 \right) +{{\rm e}^{p
\cos  t  }}\cos  t    -\cos \left( t+p\sin  t
  \right)  \right) ^{2}}{\sin \left( p\sin  t
 \right)  \left( \sin \left( t+p\sin  t   \right) +1
 \right) }};
\]
these are also smooth functions on $(-1,1)$, the latter one is positive.
They satisfy the symmetry properties
\[h_p(t)=-h_p(\pi-t),\]
\[\varpi_p(t)=-\varpi_p(\pi-t),\]
\[\phi_p(t)=\phi_p(\pi-t). \]

Now, we claim, for any $t\in(0,\pi)$,
\[\left\|\frac{A-\Id}{(1+h_p(t))\Id+(1-h_p(t))A}-\varpi_p(t)\right\|_2\leq \phi_p(t).\]
\begin{proof}
Let us use he notation \eqref{eq:PCT} from the proof of Theorem \ref{th:CRrangeInv}.
Note that  in the give range $P>0$, $C>0$, $T>0$ but $(P-1)(T-1)\leq0$, thus
\[CPT+CP+CT-PT+C+P+T-1=C \left( T+1 \right)  \left( P+1 \right) - \left( T-1 \right)  \left(
P-1 \right)>0.\]
This shows that the formulas
\[h_p(t)={\frac {CPT+CP-CT-PT-C+P-T+1}{CPT+CP+CT-PT+C+P+T-1}},\]
\[\varpi_p(t)=-\frac18\,{\frac { \left( PC-PT+C+T \right) ^{2}-(CPT+CT+P-1)^2}{PC \left( {T}^{
2}+1 \right) }},\]
\[\phi_p(t)=
\frac18\,{\frac { \left( CPT+CP+CT-PT+C+P+T-1 \right) ^{2}}{P \left( {T}^{
2}+1 \right) C}}\]
do provide smooth functions; and $\phi_p(t)>0$.
Regarding the monotonicity properties,
\[ h_p'(t)=-\frac{2P( C^2(T+1)^2p+2C(T^2+1)+(T-1)^2p) }{ (CPT+CP+CT-PT+C+P+T-1)^2 }<0;\]
the limits can be checked easily.

As for the main statement, it is sufficient to prove that
$\CRext\left(\frac{A-\Id}{(1+h_p(t))\Id+(1-h_p(t))A}\right)\subset \Dbar((\varpi(t),0),\phi_p(t))$.
We can show this by demonstrating that $\partial\Dbar((\varpi(t),0),\phi_p(t))\cap \overline{\mathbb C}^+$ is an
actual $h$-tangent line to the image of $\exp(\Dbar(0,p))$ under the conformal map
$U_{p,t}:u\mapsto\frac{u-1}{(1+h_p(t))1+(1-h_p(t))u}$,
touching from the correct side.
In the terminology of the proof of Theorem \ref{th:CRrangeInv}, it is sufficient to show that the curve
\[S_{p,t}(\tau):=U_{p,t}(L_{p,t}(\tau))=\frac{L_{p,t}(\tau)-1}{(1+h_p(t))1+(1-h_p(t))L_{p,t}(\tau)}\]
satisfies
\begin{equation}S_{p,t}(0)=\varpi_p(t)+\mathrm i\phi_p(t)\plabel{eq:subkey}\end{equation}
and
\begin{equation}\frac{\partial S_{p,t}}{\partial \tau}(0)<0\plabel{eq:topkey}\end{equation}
(including the statement that the latter one is real).
Indeed; so $S_{p,t}(\tau)$ must yield a semicircle, oriented from right to left, with $S_{p,t}(0)$ as its top point;
then the M\"obius map $U_{p,t}$ is orientation-preserving as $\det\begin{bmatrix}1&-1\\1-h_p(t)&1+h_p(t)  \end{bmatrix}=2>0$
shows that;
thus $U_{p,t}(\exp(\Dbar(0,p))\cap{\mathbb C}^+)$ still lies on the left of the oriented $S_{p,t}(\tau)$;
in case of \eqref{eq:topkey}, this is the bounded region.
As an example, Figure 7 depicts the semicircle $S_{p,t}(\tau)$, ($\tau\in[-1,1]$) and
$U_{p,t}(\gamma_p(\theta))$ $(\theta\in[0,\pi])$, the image of the boundary curve for $p=\pi$ and $t=\pi/3$:

\begin{figure}[H]
   \begin{subfigure}[b]{4in}
    \includegraphics[width=4.0in]{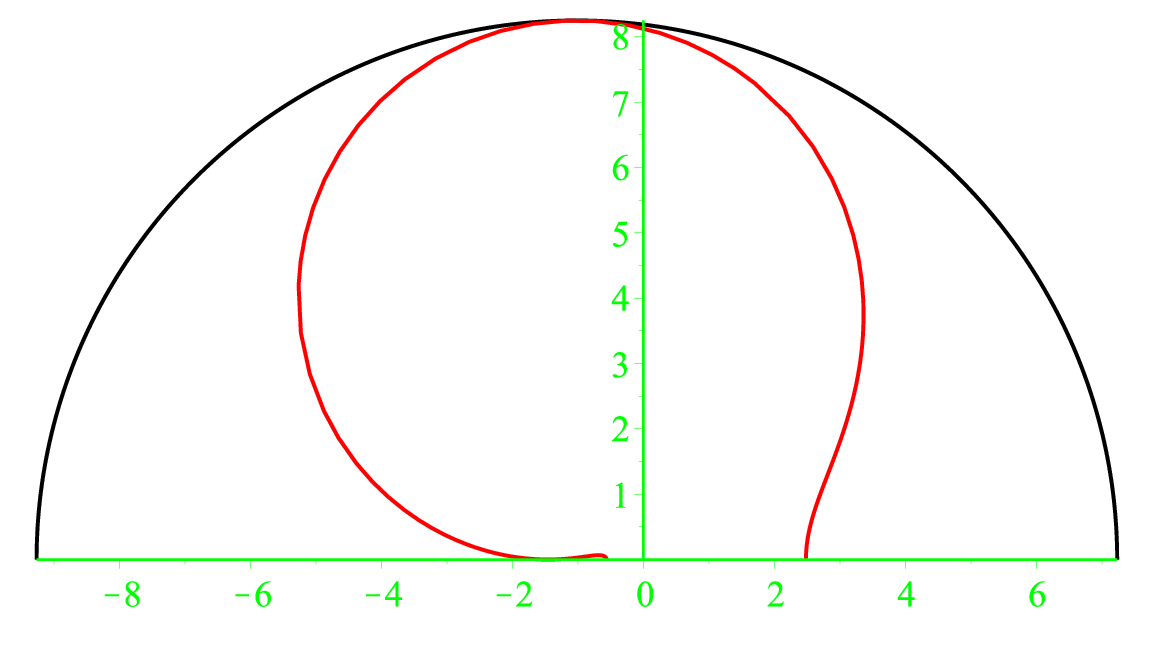}
    \caption*{Fig. \ref{fig:figB05}.  boundary curve with a tangent line, displaced in PHP}
   \end{subfigure}
\phantomcaption
\plabel{fig:figB05}
\end{figure}

Now, \eqref{eq:subkey} and  \eqref{eq:topkey} follows from
\[ S_{p,t}(\tau)=\varpi_p(t)+(-\tau+\mathrm i\sqrt{1-\tau^2})\phi_{p}(t),\]
which is easy to check.
\end{proof}
\end{theorem}
\begin{remark}\plabel{rem:monoestimate}
At first sight, the estimate of Theorem \ref{th:monoestimate}
is less artificial than the estimate of Theorem \ref{th:biestimate}.
In fact, the opposite view is also reasonable, as in  the former case we estimate a simple
but sort of arbitrary resolvent term, while in the latter case, the use
of conformal range allows to separate the ``central'' and the ``real scalar'' parts of the
resolvent terms naturally.
\qedremark
\end{remark}
\begin{remark}
\plabel{rem:biestimate}
A consequence of the previous statement that for $t\in(0,\pi/2]$ the estimate
\begin{equation}
\left\|\frac{A-\Id}{(1+h_p(t))\Id+(1-h_p(t))A} + \frac{A-\Id}{(1-h_p(t))\Id+(1+h_p(t))A} \right\|_2\leq 2\phi_p(t)
\plabel{eq:respair}
\end{equation}
holds.
Let us assume now that for a unit vector $\mathbf x$ in the Hilbert space the equality
\[\mathcal J_{p,t}(\mathbf x):\quad
\left|\frac{A-\Id}{(1+h_p(t))\Id+(1-h_p(t))A}\mathbf x + \frac{A-\Id}{(1-h_p(t))\Id+(1+h_p(t))A}\mathbf x \right|_2= 2\phi_p(t)\]
holds for some $t\in (0,\pi)$. According to the previous discussion, and strict convexity, this implies that there is a unit vector $\mathbf y$,
perpendicular to $\mathbf x$, such that
\begin{align}
\frac{A-\Id}{(1+h_p(t))\Id+(1-h_p(t))A}\mathbf x=\varpi_p(t)\mathbf x+\phi_p(t)\mathbf y,\plabel{eq:ram1}\\
\frac{A-\Id}{(1-h_p(t))\Id+(1+h_p(t))A}\mathbf x=-\varpi_p(t)\mathbf x+\phi_p(t)\mathbf y\plabel{eq:ram2}
\end{align}
hold.
From this, it is easy to conclude that
\begin{align}
A\mathbf x=&\alpha_{11}(p,t) \mathbf x+\alpha_{21}(p,t) \mathbf y,\plabel{eq:lam1}\\
A\mathbf y=&\alpha_{12}(p,t) \mathbf x+\alpha_{22}(p,t) \mathbf y;\plabel{eq:lam2}
\end{align}
where
\[\alpha_{11}(p,t)={\frac {{{\rm e}^{p\cos t }}\cos \left( -t+p\sin  t   \right) +{{\rm e}^{-p\cos t }}
\cos \left( t+p\sin  t   \right) }{2\cos t }} ,\]
\[\alpha_{21}(p,t)= {\frac {{{\rm e}^{p\cos t }}\sin \left( -t+p\sin  t   \right) +{{\rm e}^{-p\cos t }}\sin
 \left( t+p\sin  t   \right) +{{\rm e}^{p\cos  t  }}-{{\rm e}^{-p\cos t }}}{2\cos  t  }},\]
\[\alpha_{12}(p,t)={\frac {-{{\rm e}^{p\cos t }}\sin \left( -t+p\sin t   \right) -{{\rm e}^{-p\cos t }}\sin
 \left( t+p\sin \left( t \right)  \right) +{{\rm e}^{p\cos t  }}-{{\rm e}^{-p\cos t }}}{2\cos  t  }} ,\]
\[\alpha_{22}(p,t)={\frac {{{\rm e}^{p\cos t }}\cos \left( -t+p\sin t   \right) +{{\rm e}^{-p\cos t }}
\cos \left( t+p\sin  t   \right) }{2\cos t }} .\]

Then, one can see that the contribution $\CRext(A|_{\mathbb C\mathbf x+\mathbb C\mathbf y})$ comes from the
principal disk with center
\[\Omega_{p}(t)=\frac{1}{2\cos t}\left(\mathrm e^{p\cos t  -\mathrm i t+\mathrm i p\sin t}+\mathrm e^{-p\cos t  +\mathrm i t+\mathrm i p\sin t}  \right)\]
and radius
\[\omega_p(t)=\frac{\sinh(p\cos t)}{\cos t}.\]
Later, in Part IIA, we will see that these are ``maximal disks'' (of hyperbolic type).
Let us restrict our attention to $0<p<\pi$.
We see that $\mathcal J_{p,t}(\mathbf x)$ can be true for a pair $(p,t)$,
but more than it is impossible because its sharpness already implies strong arithmetic requirements.
In fact, we see later that a (maximal) disk cannot be tangent to the boundary at more than 2 places (with multiplicities), and even that
can hold only in restricted circumstances.

Although it is somewhat complicated, it can be quantified that even just nearing $\mathcal J_{p,t}(\mathbf x)$ must impose restrictions for
other  $\mathcal J_{p,t'}(\mathbf x)$'s.
Indeed,  the norm of $(U_{p,t}(A)-\varpi_{p(t)})\mathbf x+(U_{p,-t}(A)-\varpi_{p(-t)})\mathbf x$ can approach
$2\phi_p(t)$ only if the norm of $(U_{p,t}(A)-\varpi_{p(t)})\mathbf x-(U_{p,-t}(A)-\varpi_{p(-t)})\mathbf x$
is small (cf. $|\mathbf u-\mathbf v|^2=2|\mathbf u|^2+2|\mathbf v|^2-|\mathbf u+\mathbf v|^2$).
Moreover, taking the strict convexity of the boundary into account  also forces the value of $(U_{p,t}(A)-\varpi_{p(t)})\mathbf x:\mathbf x$
up to a small error term.
Ultimately, this implies that (\ref{eq:ram1}--\ref{eq:ram2}) hold but with small error terms on the right (with an appropriate $\mathbf y$).
Consequently, (\ref{eq:lam1}--\ref{eq:lam2}) also hold with small error terms.
Thus, approaching a particular $\mathcal J_{p,t}(\mathbf x)$  ($p$ fixed, although the date depends continuously on it)  will impose explicit restrictions to other ones
(as everything can be expressed a linear combination of $\mathbf x$ and $\mathbf y$ plus an error term).
\qedremark
\end{remark}
\snewpage

For $p>0$, let us define $0<\ell(p)<\frac\pi2$ as the solution of the equation
\[\ell(p)+p\sin \ell(p)=\frac\pi2.\]
Then $\ell:(0,\infty)\rightarrow (0,\pi/2)$ is a decreasing diffeomorphism.
Using a few Newton iterations, it is easy to see that
\[\ell(\pi)=0.386519539\ldots \]

Related to the previous theorem, it is easy to see that for $0<p<\pi$,
\[\ell(p)=(h_p)^{-1}(1),\qquad \pi-\ell(p)=(h_p)^{-1}(-1).\]
\begin{theorem}\plabel{th:biestimate2}
Let $\CRext(A)\subset \exp \Dbar(0,p)$, $0< p<\pi$.
Then
\begin{equation}
\|\log A\|_2\leq J(p)=\int_{t=\ell(p)}^{\pi-\ell(p)}
\underbrace{{\frac {p+\sin \left( p\sin  t  \right)-\cos \left( p\sin  t  \right) p\sin t}{2\sin
 \left( p\sin  t   \right) }}}_{JJ(p,t)}\,\mathrm dt.
 \plabel{eq:biestimate}
\end{equation}
\begin{proof}
Applying \eqref{eq:logdef2}, we find
\[\log A=\int_{t=\ell(p)}^{\pi-\ell(p)}\frac{A-\Id}{(1+h_p(t))\Id+(1-h_p(t))A}\,\left|h'_p(t)\right|\mathrm dt.\]
Applying the estimate of Theorem \ref{th:biestimate}, and using the antisymmetry properties of $h_p(t)$ and $\omega_p(t)$,
we find
\[\|\log A\|_2\leq\int_{t=\ell(p)}^{\pi-\ell(p)}\phi_p(t)\,\left|h'_p(t)\right|\mathrm dt.\]
Now, $\phi_p(t)\,\left|h'_p(t)\right|$ simplifies as indicated in  \eqref{eq:biestimate}.
\end{proof}
\end{theorem}
\begin{remark}\plabel{rem:biestimate2}
Again,  \eqref{eq:biestimate} is not sharp.
The previous estimate can be checked through  $(\log A)\mathbf x$ evaluated and
estimated by using $\mathcal J_{p,t}(\mathbf x)$.
Remark \ref{rem:biestimate} shows, however, that collective equality (in almost everywhere sense) cannot hold under the integral sign.
In fact, based on this, one could quantify stronger estimates but that would require
considerably more computation.
\qedremark
\end{remark}
\snewpage
\begin{theorem}\plabel{th:Jasympt}
(a) As $p\searrow0$,
\begin{equation}J(p)=p+\frac16\,{p}^{3}-{\frac {1}{72}}{p}^{5}+  \frac{211}{15120} p^7 +O({p}^{9}).\plabel{eq:Jsan1}\end{equation}

(b) As $p\nearrow\pi$,
\begin{equation}
J(p)=
\pi\sqrt{\frac{\pi+p}{\pi-p}}+J_\pi+O({\pi-p}),
\plabel{eq:Jsan2}
\end{equation}
where
\begin{equation}
J_\pi=-4\tan\ell(\pi)
+\int_{t=\ell(\pi)}^{\pi-\ell(\pi)} \left(JJ(\pi,t)-\frac{2}{\cos^2 t} \right)\,\mathrm dt\notag
\end{equation}
(and the integrand is actually a smooth function of $t$). Numerically, $J_\pi=-3.0222\ldots$

(c) As a crude estimate, for $0<p<\pi$,
\begin{equation}J(p)< p\sqrt{\frac{\pi+p}{\pi-p}}.\plabel{eq:crude}\end{equation}

\begin{proof}(a) As $p\searrow0$, it is easy to see that
\[\ell(p)=\frac\pi 2-p+\frac12\,{p}^{3}-{\frac {13}{24}}\,{p}^{5}+\frac{541}{720}\,p^7+O({p}^{9})\]
and
\begin{align*}
JJ_p(t)&=\frac1{2\sin t}+\sum_{k=1}^\infty\frac{2^{2k-1}\sin^{2k}t+  (2^{2k-1}-1)\sin^{2k-1}t }{(2k)!}(-1)^{k+1}B_{2k} p^{2k}\\
&=\frac1{2\sin t}+ \frac{(2\sin t+1)\sin t}{12}p^2
+\frac{(8\sin t+7)\sin^3 t}{720}p^4\\
&\quad+ \frac{(32\sin t+31)\sin^5 t}{30240}p^6
+\frac{(128\sin t+127)\sin^7 t}{1209600}p^8+O(p^{10}).
\end{align*}
 Then \eqref{eq:Jsan1} follows from evaluating  \eqref{eq:biestimate}, essentially formally in $p$.
(Here the possible singularity in the integrand is  from $\frac{1}{\sin t}$,
but the integration domain $t\in(\ell(p),\pi-\ell(p))$ is well kept away from $0$ and $\pi$, due to $\ell(p)\sim\pi/2$. )

(b)  The singular behaviour (in the limit) comes from $p\sin t \sim \pi$. Fix $p_0\in(0,\pi)$.
Let $D_0=\{(p,t)\,:\, p_0\leq p\leq \pi,\quad   \ell(p)\leq t\leq \pi-\ell(p)\}$.
Then $p\sin t$ takes positive values on $D_0$, in fact, the compact interval of $[\pi/2-\ell(p_0),\pi]$.
Notice that $\frac{\sin x}{(\pi^2-x^2)}$ is an analytic function, which is positive on $x\in(0,\pi]$.
Following this idea, computing the appropriate residues, we find that
\[JJ(p,t)=p\left( \frac{\pi}{\pi^2-p^2\sin^2 t}+J_1(p\sin t) \right)+\left( \frac{\pi^2}{\pi^2-p^2\sin^2 t}+J_2(p\sin t) \right)\]
where $J_1$ and $J_2$ are analytic in a neighborhood of  $[\pi/2-\ell(p_0),\pi]$.

Now we have
\[JJ(p,t)=\frac{p\pi+\pi^2}{\pi^2- p^2\sin^2 t }+\underbrace{pJ_1(p\sin t)+J_2(p\sin t)}_{\widehat{JJ}(p,t)}.\]

Then
\[
J(p)=\int_{t=0}^\pi \frac{p\pi+\pi^2}{\pi^2- p^2\sin^2 t }\,\mathrm dt
-\left(\int_{t=0}^{\ell(p)}
+\int_{t=\pi-\ell(p)}^\pi\right) \frac{p\pi+\pi^2}{\pi^2- p^2\sin^2 t }\,\mathrm dt+\int_{t=\ell(p)}^{\pi-\ell(p)} \widehat{JJ}(p,t)\,\mathrm dt.\]

As $p\nearrow\pi$, the first summand of $J(p)$ yields
\[
\frac{p\pi+\pi^2}{\sqrt{\pi^2-p^2}}=
\pi\sqrt{\frac{\pi+p}{\pi-p}}.
\]
As $p\nearrow\pi$, the middle summand of $J(p)$ yields
\[
 -\left(\int_{t=0}^{\ell(\pi)}
+\int_{t=\pi-\ell(\pi)}^\pi\right) \frac{2}{\cos^2t}\,\mathrm dt
+O(\pi-p)
=-4\tan \ell(\pi)+O(\pi-p). \]
As $p\nearrow\pi$, the third summand of $J(p)$ yields
\[\int_{t=\ell(\pi)}^{\pi-\ell(\pi)} \widehat{JJ}(\pi,t)\,\mathrm dt+O({\pi-p})=
\int_{t=\ell(\pi)}^{\pi-\ell(\pi)}\left( {JJ}(\pi,t)-\frac{2}{\cos^2 t}\right)\,\mathrm dt+O({\pi-p}).\]
Adding up these terms yields the corresponding statement.
The numerical evaluation of $J_\pi$ can be realized by various methods.

(c) As $p\searrow0$,
\[p\sqrt{\frac{\pi+p}{\pi-p}}=p+\frac1\pi p^2+O(p^3);\]
thus, considering  \eqref{eq:Jsan1}, we see that  \eqref{eq:crude} holds for $p\sim0$.
As $p\nearrow\pi$, considering  \eqref{eq:Jsan2} and $J_\pi<0$, we see that \eqref{eq:crude} holds for $p\sim\pi$.
Thus it is sufficient to check \eqref{eq:crude} on a properly specified compact subinterval of $(0,\pi)$.
This is, however, basically a numerical problem, so we omit the details.
\end{proof}
\end{theorem}
\begin{remark}
\plabel{rem:Jasympt}
As the resolvent pair estimate \eqref{eq:respair} is better than the individual resolvent estimate \eqref{eq:resmono},
Theorem \ref{th:biestimate2}, predictably, leads to a better estimate than Theorem \ref{th:Lestimate}.
This shows up in comparing Theorem \ref{th:Hasympt} and Theorem \ref{th:Jasympt}.
This leads to a modest gain as $p\nearrow \pi$ (as $J_\pi<H_\pi$).
What is more conspicuous is the improvement in the estimate for short term behaviour $p\searrow0$.
The lack of the quadratic term in \eqref{eq:Jsan1}  is a bit striking.
In fact, as we will see in Part IIA, the short-term estimate is sharp up to
$\frac1{270}p^7+o(p^7)$.
\qedremark
\end{remark}
\begin{remark}
\plabel{rem:astrem}
We can define $K^{\CR}(p)=\sup\{\|\log A\|_2\,:\, \CRext(A)\subset\exp\Dbar(0,p)\}$
for $p\in[0,\pi)$.
Then there is a smallest constant $K^{\CR}_\pi$ such that
\[ K^{\CR}(p)\leq \frac{\sqrt{2}\pi^{3/2}}{\sqrt{\pi-p}}+K^{\CR}_\pi+o(1)\]
holds as $p\nearrow\pi$.
Similarly, we can define $K^{\mathrm M}(p)=\sup\{\|\log(\Lexp(\phi))\|_2\,:\,\smallint\|\phi\|_2\leq p\}$,
and the corresponding best constant  $K^{\mathrm M}_\pi$.
We can summarize our findings as follows.

According to the previous discussions,
\[-2\pi\leq K^{\mathrm M}_\pi\leq K^{\CR}_\pi<J_\pi=-3.0222\ldots\quad.\]

Furthermore, in terms of some explicit expressions, (cf. Example \ref{ex:critical})
\[\pi\sqrt{\frac{{\pi+p}}{{\pi-p}}}-p-\pi<K^{\mathrm M}(p)\leq K^{\CR}(p)<  p\sqrt{\frac{\pi+p}{\pi-p}} =\pi\sqrt{\frac{\pi+p}{\pi-p}}-\sqrt{\pi^2-p^2}.\]

\qedremark
\end{remark}
\snewpage
As a corollary, we obtain
\begin{theorem}
\plabel{th:Magnusest}
 If $\phi$ is $\mathcal B(\mathfrak H)$-valued, and $\int\|\phi\|_2<\pi$, then the following hold:

(a) Regarding the norm of the Magnus expansion,
\[\|\mu_{\mathrm R}(\phi)\|_2
\equiv\left\|\sum_{k=1}^\infty\int_{t_1\leq\ldots\leq t_k\in I}\phi(t_1)\cdot\ldots\cdot\phi(t_k)\right\|_2
\leq J\left(\int\|\phi\|_2\right).\]

(b) Regarding the $k$th term of the Magnus expansion,
\[\|\mu_{k,\mathrm R}(\phi)\|_2\equiv
\left\|\int_{t_1\leq\ldots\leq t_k\in I}\mu_k(\phi(t_1),\ldots,\phi(t_k))\right\|_2\leq
\pi^{-k+1}2\sqrt{\mathrm e k}\cdot\left(\int\|\phi\|_2\right)^k.\]

\begin{proof}
(a) This follows from Theorems \ref{th:CRrange} and \ref{th:biestimate2}.
(b) $\int\|\phi\|_2>0$ can be assumed.
Consider the operator valued  function $\eta$ given by
\[\eta(z)=\log\Rexp\left(\frac{z}{\smallint\|\phi\|_2}\phi\right).\]
This is analytic in $\Dbar(0,\pi)$, moreover,
\[\|\eta(z)\|_2\leq J(|z|).\]
Applying the generalized Cauchy theorem with
$\partial \Dbar\left(0,\pi-\frac1{2k}\pi\right)$,
we estimate the $k$th power series coefficient $\eta_k$ of $\eta$ at $z=0$, by
\begin{align}
\|\eta_k\|_2&\leq \left(\pi-\frac1{2k}\pi\right)^{-k}J\left(\pi-\frac1{2k}\right)
\notag\\
&\leq \left(\pi-\frac1{2k}\pi\right)^{-k}
\left(\pi-\frac1{2k}\pi\right)\sqrt{\frac{2\pi-\frac1{2k}\pi}{\frac1{2k}\pi}}
\notag\\
&=\pi^{-k+1}\left(1-\frac1{2k}\right)^{-k+1}2\sqrt{k-\frac14}
\notag\\
&\leq\pi^{-k+1}2\sqrt{\mathrm e k}.
\notag
\end{align}
On the other hand,
\[\eta_k=\left(\int\|\phi\|_2\right)^{-k}\mu_{k,\mathrm R}(\phi).\]
This proves the statement.
\end{proof}
\end{theorem}

\begin{remark}
\plabel{rem:MagnusOff}
Compared to \eqref{eq:critast} in Example \ref{ex:critical}, the estimate above is essentially off
by a linear factor $\sqrt{2\mathrm e\mathrm \pi}k$.
This is not very good, but considering the crudeness of the method of estimate, it is fair. \qedremark
\end{remark}

If we can want to extend our results for $C^*$-algebras, then it requires no particular efforts.
According to the Gelfand--Naimark representation theorem, $C^*$-algebras have isometric
$*$-representations on Hilbert spaces, thus our estimates carry over directly.
The only concern occurs only if one wants to use the terminology of conformal range (or
Davis--Wielandt range) in the the $C^*$-algebraic setting.
In this case the original graph-induced version of  $\CRext(A)$ cannot be used.
However, the closed convex closure $\ccCRext(A)$ can be characterised using $A^*,\leq$, etc.,
thus carries over; it is independent of the isometric $*$-representation.

\snewpage

Note, however, the immediate corollary (but which uses Theorem \ref{th:Magnusest} only in a qualitative way)
\begin{theorem}
\plabel{th:Magnusstrong}
Suppose that $\phi$ is $\mathcal B(\mathfrak H)$-valued, and $\int\|\phi\|_2<\pi$.
Assume that the projections $P_\lambda$ (in a sequence or net) strongly converge to the identity $\Id$.
Then the Magnus expansions $\mu_{\mathrm L}(P_\lambda\cdot\phi\cdot P_\lambda)$
strongly converge to the  Magnus expansion $\mu_{\mathrm L}(\phi)$.
\begin{proof}
According to Theorem \ref{th:Magnusest}, the Magnus expansions are uniformly absolute convergent.
Meanwhile, the individual terms converge strongly.
This implies global convergence strongly.
\end{proof}
\end{theorem}
The theorem above applies, in particular, when the $P_\lambda$ form the natural net
of finite dimensional projections.
So we see that while there are several qualitative differences between the
finite dimensional and infinite dimensional cases; there are no
quantitative differences
between the general finite dimensional case and the infinite dimensional case
in a set of problems like absolute norm estimates.

Also note that what we did with $\log$ here is not truly specific to that function,
as one can use resolvent expansions to other analytic functions too, cf. Mielnik,  Pleba\'{n}ski \cite{MP} or Part I.
It is true, however, that the most natural class of functions, the power functions, are closely related to $\log$.
\snewpage
\appendix
\section{The hyperbolic space}
\plabel{app:hyprev}

The Davis--Wielandt shell is a construction introduced by Wielandt \cite{Wie} and Davis \cite{D1}, \cite{D2}.
Its input is a (multivalued) linear operator on a complex linear Hilbert space
and its output is a subset of the points of a canonical model of the asymptotically closed hyperbolic $3$-space.

\begin{point}
\plabel{po:3Dmodels}
The notable canonical models of the hyperbolic $3$-space which we use
  are represented by subsets $\mathbb R^3$
as the set of points with appropriate $h$-lines and $h$-planes distinguished.
These are the:

$\bullet$ The Cayley--Klein--Beltrami model $H^3_{\mathrm{CKB}}$.
Its points $(x_{\mathrm{CKB}}, y_{\mathrm{CKB}},z_{\mathrm{CKB}})\in H^3_{\mathrm{CKB}}$ satisfy
\begin{equation}
(x_{\mathrm{CKB}})^2+(y_{\mathrm{CKB}})^2+(z_{\mathrm{CKB}})^2<1,
\plabel{eq:A1}\end{equation}
and $h$-lines and $h$-planes are  the non-empty traces of ordinary lines and planes.

$\bullet$ The Poincar\'e model $H^3_{\mathrm{P}}$.
Its points $(x_{\mathrm{P}}, y_{\mathrm{P}},z_{\mathrm{P}})\in H^3_{\mathrm{P}}$ satisfy
\begin{equation}
(x_{\mathrm{P}})^2+(y_{\mathrm{P}})^2+(z_{\mathrm{P}})^2<1,
\plabel{eq:A2}\end{equation}
again, but $h$-lines and $h$-planes are  the traces of lines or circles and planes or spheres perpendicular to the unit sphere.

$\bullet$ The parabolic version of the Cayley--Klein--Beltrami model $H^3_{\mathrm{CKB(P)}}$.
Its points $\left(x_{\mathrm{CKB(P)}}, y_{\mathrm{CKB(P)}},z_{\mathrm{CKB(P)}}\right)\in H^3_{\mathrm{CKB(P)}}$ satisfy
\begin{equation}
(x_{\mathrm{CKB(P)}})^2+(y_{\mathrm{CKB(P)}})^2<z_{\mathrm{CKB(P)}},
\plabel{eq:A3}\end{equation}
and $h$-lines and $h$-planes are  the non-empty traces of ordinary lines and planes.

$\bullet$ The Poincar\'e half-space model $H^3_{\mathrm{PH}}$.
Its points $(x_{\mathrm{PH}}, y_{\mathrm{PH}},z_{\mathrm{PH}})\in H^3_{\mathrm{PH}}$ satisfy
\begin{equation}
0<z_{\mathrm{PH}},
\plabel{eq:A4}\end{equation}
and $h$-lines and $h$-planes are  the traces of lines or circles and planes or spheres perpendicular to the $xy$-plane.

In the case of hyperbolic space, the incidence structure already determines the natural metric structure (in contrast to the Euclidean space where
it determines only an affine linear structure). Thus the specifications above are sufficient.
\snewpage
There is a rather canonical correspondence between these models such that in terms of coordinates
\begin{align*}
(x_{\mathrm{CKB}},y_{\mathrm{CKB}},z_{\mathrm{CKB}})
&=\frac{2\cdot(x_{\mathrm{P}},y_{\mathrm{P}},z_{\mathrm{P}})} {1+(x_{\mathrm{P}})^2+(y_{\mathrm{P}})^2+(z_{\mathrm{P}})^2 }\\
&=\frac{(2x_{\mathrm{CKB(P)}},2y_{\mathrm{CKB(P)}},z_{\mathrm{CKB(P)}}-1)}{z_{\mathrm{CKB(P)}}+1}           \\
&=\frac{(2x_{\mathrm{PH}},2y_{\mathrm{PH}},(x_{\mathrm{PH}})^2+(y_{\mathrm{PH}})^2+(z_{\mathrm{PH}})^2-1) }{(x_{\mathrm{PH}})^2+(y_{\mathrm{PH}})^2+(z_{\mathrm{PH}})^2+1 }      \\
&   \\
(x_{\mathrm{PH}},y_{\mathrm{PH}},z_{\mathrm{PH}})
&=\frac{(x_{\mathrm{CKB}},y_{\mathrm{CKB}},\sqrt{1-(x_{\mathrm{CKB}})^2-(y_{\mathrm{CKB}})^2-(z_{\mathrm{CKB}})^2 })}{1-z_{\mathrm{CKB}}}           \\
&=\frac{(2x_{\mathrm{P}},2y_{\mathrm{P}},1-(x_{\mathrm{P}})^2-(y_{\mathrm{P}})^2-(z_{\mathrm{P}})^2) }{(x_{\mathrm{P}})^2+(y_{\mathrm{P}})^2+(z_{\mathrm{P}}-1)^2 }           \\
&=\left(x_{\mathrm{CKB(P)}},y_{\mathrm{CKB(P)}},\sqrt{z_{\mathrm{CKB(P)}}-(x_{\mathrm{CKB(P)}})^2 -(y_{\mathrm{CKB(P)}})^2}   \right)  \\
&   \\
(x_{\mathrm{CKB(P)}},y_{\mathrm{CKB(P)}},z_{\mathrm{CKB(P)}})
&=\frac{(x_{\mathrm{CKB}},y_{\mathrm{CKB}},1+z_{\mathrm{CKB}})}{1-z_{\mathrm{CKB}}}           \\
&=\frac{(2x_{\mathrm{P}},2y_{\mathrm{P}},(x_{\mathrm{P}})^2+(y_{\mathrm{P}})^2+(z_{\mathrm{P}}+1)^2   ) }{(x_{\mathrm{P}})^2+(y_{\mathrm{P}})^2+(z_{\mathrm{P}}-1)^2 }           \\
&=\left( x_{\mathrm{PH}},y_{\mathrm{PH}},(x_{\mathrm{PH}})^2+(y_{\mathrm{PH}})^2+(z_{\mathrm{PH}})^2 \right)  \\
&   \\
(x_{\mathrm{P}},y_{\mathrm{P}},z_{\mathrm{P}})
&=\frac{(x_{\mathrm{CKB}},y_{\mathrm{CKB}},z_{\mathrm{CKB}})}{1+\sqrt{1-(x_{\mathrm{CKB}})^2-(y_{\mathrm{CKB}})^2-(z_{\mathrm{CKB}})^2 }} \\
&=\frac{(2x_{\mathrm{CKB(P)}},2y_{\mathrm{CKB(P)}},z_{\mathrm{CKB(P)}}-1)}{1+z_{\mathrm{CKB(P)}}+2\sqrt{z_{\mathrm{CKB(P)}}- (x_{\mathrm{CKB(P)}})^2- (y_{\mathrm{CKB(P)}})^2 }}  \\
& =\frac{(2x_{\mathrm{PH}},2y_{\mathrm{PH}},(x_{\mathrm{PH}})^2+(y_{\mathrm{PH}})^2+(z_{\mathrm{PH}})^2-1) }{(x_{\mathrm{PH}})^2+(y_{\mathrm{PH}})^2+(z_{\mathrm{PH}}+1)^2 }      \\
\end{align*}

If, in general, $(x_{\mathrm A},y_{\mathrm A},z_{\mathrm A})=f(x_{\mathrm B},y_{\mathrm B},z_{\mathrm B})$
then $f=\frac{\mathrm A}{\mathrm B}$ is the function of the canonical correspondence between the models.
A good alternative viewpoint is that there is a single abstract instance $H^3$ of the hyperbolic $3$-space,
with coordinate functions $x_{\mathrm{CKB}},$ etc., and whose point $(x_{\mathrm{CKB}},y_{\mathrm{CKB}},y_{\mathrm{CKB}})^{\mathrm{CKB}}$
is the point which corresponds to $(x_{\mathrm{CKB}},y_{\mathrm{CKB}},y_{\mathrm{CKB}})$ in $H^3_{\mathrm{CKB}}$, etc.
\end{point}

\snewpage
\begin{point}
\plabel{po:3Dbound}
The asymptotic closure  $\overline{H}^3$ of the hyperbolic space  $H^3$ is its augmentation by its asymptotic points.
In terms of models it can be represented as follows:
The set of points of $\overline{H}^3_{\mathrm{CKB}}$ and $\overline{H}^3_{\mathrm{P}}$ is just the
closure of the original set of points in ordinary sense.
Similar comment applies to $\overline{H}^3_\mathrm{CKB(P)}$ and $\overline{H}^3_{\mathrm{PH}}$ except the additional
points $\infty_\mathrm{CKB(P)}$ and $\infty_{\mathrm{PH}}$ should also be taken, respectively.
In terms of correspondences, the transformations above are still valid, except we have the additional correspondences
$(0,0,1)_{\mathrm{CKB}}\leftrightarrow(0,0,1)_{\mathrm{P}}\leftrightarrow \infty_\mathrm{CKB(P)}\leftrightarrow\infty_{\mathrm{PH}} $.
The $h$-lines and $h$-planes in the asymptotically closed cases are just the natural closures of the original ones.

Here $\overline{H}^3_{\mathrm{PH}}$ can be imagined as the
the factor of the quaternionic Riemann sphere $\mathbb H^\star$ up to conjugation by complex numbers.
In this case the factorization operation can be represented as
\[\cq^3:\qquad q_0+q_1\mathrm i+q_2\mathrm j+q_3\mathrm k\mapsto \left(q_0,q_1,\sqrt{(q_2)^2+(q_3)^2}\right)\in \overline{H}^3_{\mathrm{PH}}\]
\[\infty_{\mathbb H}\mapsto\infty_{\mathrm{PH}}\,.\]
Another viewpoint is that $(x_{\mathrm{PH}},y_{\mathrm{PH}},z_{\mathrm{PH}})\in \overline{H}^3_{\mathrm{PH}}$ simply identifies  to the quaternion $x_{\mathrm{PH}}+y_{\mathrm{PH}}\mathrm i+z_{\mathrm{PH}}\mathrm j$.
In any case, the Riemann sphere $\mathbb C^\star$ parametrizes the set of asymptotic points in $\overline{H}^3_{\mathrm{PH}}$ (i. e. the asymptotic boundary). This yields the correspondence
\[\iota_{\mathrm{PH}}:\quad \lambda\in\mathbb C\mapsto(\Rea\lambda,\Ima\lambda,0)\in\overline{H}^3_{\mathrm{PH}}
,\quad \infty\mapsto\infty_{\mathrm{PH}}\in\overline{H}^3_{\mathrm{PH}}.\]
In terms of the other models, the mappings are
\[\iota_{\mathrm{CKB(P)}}:\quad\lambda\in\mathbb C\mapsto(\Rea\lambda,\Ima\lambda,|\lambda|^2)\in\overline{H}^3_{\mathrm{CKB(P)}}
,\quad \infty\mapsto\infty_{\mathrm{CKB(P)}}\in\overline{H}^3_{\mathrm{CKB(P)}};\]
and
\[\iota_{\mathrm{CKB}}:\quad\lambda\in\mathbb C\mapsto
\left(\frac{2\Rea\lambda}{|\lambda|^2+1},\frac{2\Ima\lambda}{|\lambda|^2+1},\frac{|\lambda|^2-1}{|\lambda|^2+1}\right)\in\overline{H}^3_{\mathrm{CKB}}
,\quad \infty\mapsto(0,0,1)\in\overline{H}^3_{\mathrm{CKB}};\]
and similarly,
\[\iota_{\mathrm{P}}:\quad\lambda\in\mathbb C\mapsto
\left(\frac{2\Rea\lambda}{|\lambda|^2+1},\frac{2\Ima\lambda}{|\lambda|^2+1},\frac{|\lambda|^2-1}{|\lambda|^2+1}\right)\in\overline{H}^3_{\mathrm{P}}
,\quad \infty\mapsto(0,0,1)\in\overline{H}^3_{\mathrm{P}}.\]

The collineations of $H^3$ (which are the same as isometries) are in bijective correspondence to the collineations of $\overline{H}^3$.
In the case of projective models ${\mathrm{CKB}}$ and ${\mathrm{CKB(P)}}$,
those appear as the projective transformations leaving the base set invariant.
In the case of conformal models ${\mathrm{PH}}$ and ${\mathrm{P}}$,
those appear as the conformal transformations leaving the base set invariant.
\end{point}
\snewpage
\begin{point}
\plabel{po:3Dconformal}
If $a,b,c,d\in\mathbb C$ and $ad-bc\neq0$, then we consider the M\"obius transformation $f$ on the Riemann sphere given by
$\lambda\mapsto f(\lambda)=\frac{a\lambda+b}{c\lambda +d}$.
By the transformation $(\mathbf x,\mathbf y)\mapsto (c\mathbf y+d\mathbf x ,a\mathbf y+b\mathbf x)$, this acts on any linear relation $A$ taking it to $f(A)$.
This is compatible with the customary case when $A$ is a linear transformation.

Now, the transformation $f$ also acts on the (models of) the hyperbolic 3-space.
The simplest case to motivate this $\overline{H}^3_{\mathrm{PH}}$.
On a quaternion $q$, the transformation $f$ acts by
\begin{equation}
f(q)=\begin{cases}
\sqrt{\dfrac ad}\cdot q \cdot \sqrt{\dfrac ad}+\dfrac bd&\text{ if $c=0$,}\\
&\\
\dfrac ac-\dfrac{\sqrt{ad-bc}}{c}\cdot\left(q+\dfrac dc\right)^{-1}\cdot\dfrac{\sqrt{ad-bc}}{c}&
\text{ if $c\neq0$}.
\end{cases}
\plabel{eq:quatmobius}
\end{equation}
where the square roots are chosen consistently.
Note that
\[(aq+b)(cq+d)^{-1}=\sqrt{ad-bc}\cdot f(q)\cdot (\sqrt{ad-bc})^{-1}\qquad\phantom{\,.} \]
\[(qc+d)^{-1}(qa+b)=(\sqrt{ad-bc})^{-1}\cdot f(q)\cdot \sqrt{ad-bc}\, ; \qquad \]
thus the action extends some naive formulas from the case $ad-bc=0$.

This is compatible to conjugation by complex numbers, and one can check that this descends (in other viewpoint: restricts) to an associative action $f_{\mathrm{PH}}$ on $\overline{H}^3_{\mathrm{PH}}$,
and, in this way, there is a bijection between the orientation-preserving conformal transformations of $\mathbb C^\star$ and the
orientation-preserving collineation / conformal group on $\overline{H}^3_{\mathrm{PH}}$.
This is not surprising as it is easy to see that collineations transformations are determined by their actions on the asymptotic boundary.
(Orientation-reversing collineations can be obtained in combination with, say,  the map $q\mapsto -\bar q $. )
\snewpage

In this way, we obtain the corresponding orientation-preserving collineations in other models.
One can transcribe this action in terms of projective transformations:
If \[(\tilde x_{\mathrm{PH}},\tilde y_{\mathrm{PH}},\tilde z_{\mathrm{PH}})=f_{\mathrm{PH}}\left((x_{\mathrm{PH}},y_{\mathrm{PH}},z_{\mathrm{PH}})\right)\]
($z_{\mathrm{PH}}\neq0$), then
\[\begin{bmatrix}\tilde x_{\mathrm{PH}}/\tilde z_{\mathrm{PH}}\\\tilde y_{\mathrm{PH}}/\tilde z_{\mathrm{PH}}\\
\frac{(\tilde x_{\mathrm{PH}})^2+(\tilde y_{\mathrm{PH}})^2+ (\tilde z_{\mathrm{PH}})^2-1 }2/\tilde z_{\mathrm{PH}}\\
\frac{(\tilde x_{\mathrm{PH}})^2+(\tilde y_{\mathrm{PH}})^2+ (\tilde z_{\mathrm{PH}})^2+1 }2/\tilde z_{\mathrm{PH}}\end{bmatrix}
\quad =\quad R_{\mathrm{CKB}}(f)
\begin{bmatrix} x_{\mathrm{PH}}/z_{\mathrm{PH}}\\ y_{\mathrm{PH}}/z_{\mathrm{PH}}\\
\frac{( x_{\mathrm{PH}})^2+( y_{\mathrm{PH}})^2+ ( z_{\mathrm{PH}})^2-1 }2/z_{\mathrm{PH}}\\
\frac{( x_{\mathrm{PH}})^2+( y_{\mathrm{PH}})^2+ ( z_{\mathrm{PH}})^2+1 }2/z_{\mathrm{PH}}\end{bmatrix}
\]
where
\[R_{\mathrm{CKB}}(f)=\frac1{|ad-bc|}\cdot\frac12
\begin{bmatrix}
&1&1&\\&-\mathrm i&\mathrm i&\\1&&&-1\\1&&&1
\end{bmatrix}
\begin{bmatrix}
a\bar a&a\bar b&b\bar a&b\bar b\\a\bar c&a\bar d&b\bar c&b\bar d\\c\bar a&c\bar b&d\bar a&d\bar b\\c\bar c&c\bar d&d\bar c&d\bar d
\end{bmatrix}
\begin{bmatrix}
&&1&1\\1&\mathrm i&&\\1&-\mathrm i&&\\&&-1&1
\end{bmatrix}
\]
\[=\frac1{|ad-bc|}
\begin{bmatrix}
\re(\bar ad+\bar cb)&\ima(\bar ad+\bar cb)&\re(\bar ca-\bar db)&\re(\bar ca+\bar db)
\\
\ima(\bar da-\bar bc)&\re(\bar da-\bar bc)&\ima(\bar ca-\bar db)&\ima(\bar ca+\bar db)
\\
\re(\bar ab-\bar cd)&\ima(\bar ab-\bar cd)&\frac{|a|^2-|b|^2-|c|^2+|d|^2}{2}&\frac{|a|^2+|b|^2-|c|^2-|d|^2}{2}
\\
\re(\bar ab+\bar cd)&\ima(\bar ab+\bar cd)&\frac{|a|^2-|b|^2+|c|^2-|d|^2}{2}&\frac{|a|^2+|b|^2+|c|^2+|d|^2}{2}
\end{bmatrix}.\]

Thus, if \[(\tilde x_{\mathrm{CKB}},\tilde y_{\mathrm{CKB}},\tilde z_{\mathrm{CKB}})=f_{\mathrm{CKB}}\left((x_{\mathrm{CKB}},y_{\mathrm{CKB}},z_{\mathrm{CKB}})\right),\] then
\[\begin{bmatrix}\tilde x_{\mathrm{CKB}}\\\tilde y_{\mathrm{CKB}}\\\tilde z_{\mathrm{CKB}}\\1\end{bmatrix}
\quad \sim\quad R_{\mathrm{CKB}}(f)
\begin{bmatrix}x_{\mathrm{CKB}}\\ y_{\mathrm{CKB}}\\ z_{\mathrm{CKB}}\\1\end{bmatrix},
\]
where $\sim$ means proportional, yielding a projective representation (the relation extending to the asymptotic boundary).
The matrices $R_{\mathrm{CKB}}(f)$ can be utilized to describe the slightly more complicated
(in the conform case: quadratic rational) transformations in the other models.
(As $R_{\mathrm{CKB}}(f)$, we obtain all matrices from $\SO^{\uparrow}(2,1)$.
Using more general matrices $R\in\mathrm O^{\uparrow}(2,1)$,  we can account for all collineations.)
\end{point}

\snewpage
\begin{point}
Now, we consider the relationship of (the analytical descriptions of) the hyperbolic plane and 3-space.
The plane can be embedded into the space. In turn, the hyperbolic 3-space can be projected to the embedded hyperbolic plane.
In the $\mathrm{CKB}$ and $\mathrm{CKB(P)}$ models this means discarding second coordinate.
In concrete terms, the projection maps in the relevant models are
\begin{align*}
\pi^{[2]}_{\mathrm{CKB}}:&\quad(x_{\mathrm{CKB}},y_{\mathrm{CKB}},z_{\mathrm{CKB}})\mapsto(x_{\mathrm{CKB}},0,z_{\mathrm{CKB}}),\\
\pi^{[2]}_{\mathrm{P}}:&\quad(x_{\mathrm P},y_{\mathrm P},z_{\mathrm P})\mapsto\\
&\qquad\frac{(2x_{\mathrm P},0,2z_{\mathrm P})}{1+(x_{\mathrm P})^2+(y_{\mathrm P})^2+(z_{\mathrm P})^2+\sqrt{
(1-(x_{\mathrm P)}^2-(y_{\mathrm P})^2-(z_{\mathrm P})^2)^2+(2y_{\mathrm P})^2  }},\\
\pi^{[2]}_{\mathrm{PH}}:&\quad(x_{\mathrm{PH}},y_{\mathrm{PH}},z_{\mathrm{PH}})\mapsto(x_{\mathrm{PH}},0,\sqrt{(y_{\mathrm{PH}})^2+(z_{\mathrm{PH}})^2}),\\
\pi^{[2]}_{\mathrm{CKB(P)}}:&\quad(x_{\mathrm{CKB(P)}},y_{\mathrm{CKB(P)}},z_{\mathrm{CKB(P)}})\mapsto(x_{\mathrm{CKB(P)}},0,z_{\mathrm{CKB(P)}}),
\end{align*}
respectively.
The images of the projection maps can be identified with the corresponding canonical models for $\overline{H}^2_*$.

In the formulas above, it is reflected that, especially with respect to the PH model,
in relation of the $3$ and $2$ dimensional models, is better to use the coordinates $(x_*, z_*)$ for the plane model.
Working purely in terms of planar models, the notation `$z_*$' can be changed to `$y_*$', but this is only a minor inconvenience.

As we have seen before, \textit{all} isometries in the PH model can be represented by real fractional linear
transformations $f:\lambda\mapsto \frac{ax+b}{cx+d}$ where $a,b,c,d\in\mathbb R$, $ad-bc\neq0$.
Their effect for $w=x_{\mathrm{PH}}+\mathrm iz_{\mathrm{PH}}$ is given by
\[ f_{\mathrm PH}^{[2]}: w\mapsto \left(\frac{aw+b}{cw+d}\right)^{\text{ conjugated if } ad-bc<0}.\]
This action is inherited directly from \eqref{eq:quatmobius}, except in the result
`$\mathrm j$' is transcribed to `$\mathrm i$'.

Projective representations of the isometries can be also be obtained by restriction.
For example, in terms of the CKB model,
\[R_{\mathrm{CKB}}^{[2]}(f)=R_{\mathrm{CKB}}(f)|_{\{1,3,4\}\times \{1,3,4\}}\]
can be taken. (Remember, this applies in the case when $a,b,c,d\in\mathbb R$.)
This yields
\begin{equation}
R_{\mathrm{CKB}}^{[2]}(f)=\frac1{|ad-bc|}
\begin{bmatrix}
 ad+ cb& ca- db& ca+ db
\\
 ab- cd&\frac{a^2-b^2-c^2+d^2}{2}&\frac{a^2+b^2-c^2-d^2}{2}
\\
 ab+ cd&\frac{a^2-b^2+c^2-d^2}{2}&\frac{a^2+b^2+c^2+d^2}{2}
\end{bmatrix}.
\plabel{eq:projrepCKB}
\end{equation}
As the only nontrivial omitted term was $R_{\mathrm{CKB}}(f)_{22}=\frac{ad-bc}{|ad-bc|}$,
it is easy to see that $\det R_{\mathrm{CKB}}^{[2]}(f)=\frac{ad-bc}{|ad-bc|} $.
In fact, we obtain all elements of $\mathrm O^\uparrow(2,1)$ in this way.
(At some places, it is customary to omit the absolute values in \eqref{eq:projrepCKB}, restoring the property $\det=1$ but possibly spoiling ortochronality.)
Otherwise, similar comments apply as in the higher dimensional case.
\end{point}
\snewpage
\section{The Davis--Wielandt shell}
\plabel{app:DWrev}
\begin{point}
\plabel{po:DWgraph}
We assume that $\mathfrak H$ is a complex Hilbert space with product $\langle\cdot,\cdot\rangle$
(linear in the first variable, conjugate-linear in the second).
Then a linear operator $A$ on $\mathfrak H$ can be identified by its graph $\{(\mathbf x,\mathbf y)\in\mathfrak H\times\mathfrak H\,:\,A\mathbf x=\mathbf y \}$.
More generally, a multivalued linear operator $A$ on $\mathfrak H$ is just a linear subspace of $\mathfrak H\times\mathfrak H$.
The norm $|A|=\sup\{|\mathbf y|/|\mathbf x|\,:\,(\mathbf x,\mathbf y)\in A\}$ and inverse
$A^{-1}=\{(\mathbf y,\mathbf x)\,:\,(\mathbf x,\mathbf y)\in A\}$ can be defined as usual, and thus the co-norm $\lfloor A\rfloor=\left(|A|^{-1}\right)^{-1}$, too.
\end{point}

\begin{point}
\plabel{po:DWmap}
Now, if $A$ is a linear operator on $\mathfrak H$, then to every pair $(\mathbf x,\mathbf y)\neq(0,0)$, $\mathbf y=A\mathbf x$ we can associate a point

\[\DW_{\mathrm{CKB}}((\mathbf x,\mathbf y))=
\left(\frac{2\Rea\langle\mathbf  y,\mathbf x\rangle}{|\mathbf y|_2^2+|\mathbf x|_2^2} ,\frac{2\Ima\langle \mathbf y,\mathbf x\rangle}{|\mathbf y|_2^2+|\mathbf x|_2^2} , \frac{|\mathbf y|_2^2-|\mathbf x|_2^2}{|\mathbf y|_2^2+|\mathbf x|_2^2}\right)\in
 \overline{H}^3_{\mathrm{CKB}}\]
or
\[\DW_{\mathrm{P}}((\mathbf x,\mathbf y))=
\frac{\left(2\Rea\langle \mathbf y,\mathbf x\rangle ,2\Ima\langle \mathbf y,\mathbf x\rangle , |\mathbf y|_2^2-|\mathbf x|_2^2\right)}{|\mathbf y|_2^2+|\mathbf x|_2^2+2\sqrt{|\mathbf x|_2^2|\mathbf y|_2^2-|\langle \mathbf y,\mathbf x\rangle|^2   }}\in
 \overline{H}^3_{\mathrm{P}}\]
or
\[\DW_{\mathrm{CKB(P)}}((\mathbf x,\mathbf y))=
\left(\frac{\Rea\langle \mathbf y,\mathbf x\rangle}{|\mathbf x|_2^2} ,\frac{\Ima\langle \mathbf y,\mathbf x\rangle}{|\mathbf x|_2^2} , \frac{|\mathbf y|_2^2}{|\mathbf x|_2^2}\right)\in
 \overline{H}^3_{\mathrm{CKB(P)}}\]
or
\[\DW_{\mathrm{PH}}((\mathbf x,\mathbf y))=
\left(\frac{\Rea\langle \mathbf y,\mathbf x\rangle}{|\mathbf x|_2^2} ,\frac{\Ima\langle \mathbf y,\mathbf x\rangle}{|\mathbf x|_2^2} , \frac{\sqrt{|\mathbf x|_2^2|\mathbf y|_2^2-|\langle \mathbf y,\mathbf x\rangle|^2}}{|\mathbf x|_2^2}\right)\in
 \overline{H}^3_{\mathrm{PH}}\]
depending on the model we use (the points are corresponding according to our conventions). The set of points obtained in this way is the Davis--Wielandt shell $\DW(A)$.
In practice, it is $\DW_{\mathrm{CKB}}(A)$ or $\DW_{\mathrm{P}}(A)$ or $\DW_{\mathrm{CKB(P)}}(A)$ or $\DW_{\mathrm{PH}}(A)$.
The same can be said if $A$ is just a linear relation. In the case of pairs $(\mathbf x,\mathbf y)$, $\mathbf x=0$, $\mathbf y\neq0$ we associate the infinite points
$(0,0,1)_{\mathrm{CKB}}\leftrightarrow(0,0,1)_{\mathrm{P}}\leftrightarrow \infty_\mathrm{CKB(P)}\leftrightarrow\infty_{\mathrm{PH}} $.

The definition by Wielandt is $\DW_{\mathrm{CKB(P)}}(A)$, and this is the one generally referred by the linear algebraic community, cf. \cite{HJ}.
The definition by Davis is $\DW_{\mathrm{CKB}}(A)$, which is more geometric.
This comes handy in developing some deeper properties of Davis--Wielandt shell.
These projective models are also nice because in those cases $h$-convexity coincides with ordinary convexity.
$\DW_{\mathrm{PH}}(A)$ of the Poincar\'e half-space model is, however, probably the most advantageous for the purposes of analytical computations.
In contrast,  $\DW_{\mathrm{P}}(A)$ using the other conformal model, is rather ugly to compute with,
however, it is most advantageous for visualization
as it is more proportional than the other ones.
\end{point}
\snewpage

The fundamental properties of the Davis--Wielandt  were established by Wielandt \cite{Wie} and Davis \cite{D1}, \cite{D2}.
It must be noted, however, that  Wielandt \cite{Wie} is more like a report of results, in essence it refers back
Wielandt \cite{Wie0}, where the background is explained.
The deeper analysis is due to Davis.
We summarize some of the fundamental results as follows:

\begin{theorem}[Wielandt, Davis]
\plabel{th:DWbasic}
Suppose that $A$ is a linear relation. Then:

(a) $\iota_*(\lambda)\in\DW_*(A)$ if and only if $ \lambda$ is an eigenvalue of $A$.

(b) If $\dim A<+\infty$, then $\DW_*(A)$ is compact.

(c) If $f$ is a complex M\"obius transformation, then $\DW_{*}(f(A))=f_*(\DW_{*}(A))$.

Suppose that $A_1$ and $A_2$ are linear relations, $\dim A_1,\dim A_2\geq1$. Then:

(d)  $\DW_*(A_1\oplus A_2)$ is the union of $h$-segments connecting the points of $\DW_*(A_1)$ to the points of $\DW_*(A_2)$.
\end{theorem}

If $\dim A=0$, then $\DW_*(A)$ is empty.
If $\dim A=1$, then $\DW_*(A)$ contains only a point. Otherwise
\begin{theorem}[Davis]
\plabel{th:Davis}
 Suppose that $A$ is a linear relation. Then:

(a) If $\dim A=2$, then  $\DW_{\mathrm{CKB}}(A)$ is a possibly degenerate ellipsoid in $\overline{H}^3_{\mathrm{CKB}}$.
Thus, in that sense, it can be said that $\DW_{*}(A)$ is a possibly degenerate $h$-ellipsoid.

(b) If $\dim A\geq3$, then $\DW_*(A)$ is convex.

Suppose that $A_1$ and $A_2$ are linear relations, $\dim A_1,\dim A_2\geq1$. Then:

(c) $\DW_*(A_1\oplus A_2)$ is the $h$-convex hull of  $\DW_*(A_1) \cup \DW_*(A_2)$.
\begin{proof}[Sketch of proofs]
Regarding \ref{th:Davis}(a):
Take an orthonormal basis $(\mathbf x_1,\mathbf y_1)$, $(\mathbf x_2,\mathbf y_2)$ in the graph
with respect to the  restriction of the natural product n $\mathfrak H\oplus\mathfrak \mathfrak H.$
(Thus $|\mathbf x_1|_2^2+|\mathbf y_1|^2=|\mathbf x_1|_2^2+|\mathbf y_1|^2=1$ and
$\langle\mathbf x_1,\mathbf x_2\rangle+\langle\mathbf y_1,\mathbf y_2\rangle=0$)
Consider the image of the (graph) unit sphere
$\{z_1(\mathbf x_1,\mathbf y_1)+ z_2(\mathbf x_2,\mathbf y_2)\:,\:z_1,z_2\in\mathbb C, |z_1|^2+|z_2|^2=1\}$.
Then the image through $\DW_{\mathrm{CKB}}$
will be a linear combination of $\Rea (z_1z_2),$ $ \Ima (z_1z_2)$, $|z_1|^2-|z_2|^2$ and
$|z_1|^2+|z_2|^2=1$.
As the $( 2\Rea (z_1\bar z_2), 2\Ima (z_1\bar z_2), |z_1|^2-|z_2|^2 )$  ($=\iota_{\mathrm{CKB}}(z_1/z_2)$) form a sphere,
we see that the image is an affine linear image of a sphere.
(Taking orthonormal basis and CKB is not essential but in that way projective geometry can be avoided.)

Regarding \ref{th:Davis}(b): In the presence an extra dimension, the image of
$\{z_1t(\mathbf x_1,\mathbf y_1)+ z_2t(\mathbf x_2,\mathbf y_2)+\sqrt{1-t^2}(\mathbf x_3,\mathbf y_3)\:,\:z_1,z_2\in\mathbb C, |z_1|^2+|z_2|^2=1, t\in[0,1]\}$ can be considered.
We see that the original ellipsoid can be contracted in the shell.
As $\mathbb S^2$ is not contractible (a topological argument!), this means, that in the contraction, the whole interior of the ellipsoid must
be taken as image.
Thus, in the shell, the ellipsoid is filled in.

The rest, conformal invariance, etc., is easy.
\end{proof}
\end{theorem}

Although linear relations were useful in the proof of  Theorem \ref{th:Davis},
in what follows $A$ will always be linear operator (on a complex Hilbert space).

If $\dim_{\mathbb C}\mathfrak H=2$, then a nice geometrical picture emerges:

\begin{theorem}[Wielandt, Davis]
\plabel{thm:DWconc}
Suppose that $A$ is a linear operator on a $2$-dimensional Hilbert space.
[It is sufficient to assume that $A$ is a linear relation on $\mathfrak H$ such that
$\dim A=\dim\mathfrak H=2$ but $A$ is not an extension from a $1$-dimensional Hilbert space.]

We have the following possibilities:

(i) $A$ has a double eigenvalue $\lambda$, and $A$ is normal (thus $A=\lambda \Id$).

Then $\DW_*(A)$ contains only the point $\iota_*(\lambda)$.

(ii)  $A$ has two different eigenvalues $\lambda_1\neq\lambda_2$, and $A$ is normal.

Then $\DW_*(A)$ is the asymptotically closed $h$-line connecting $\iota_*(\lambda_1)$ and $\iota_*(\lambda_2)$.

(iii) $A$ has a double eigenvalue $\lambda$, and $A$ is not normal.

Then $\DW_*(A)$ is an asymptotically closed $h$-horosphere with  asymptotical point $\iota_*(\lambda)$.
In the $\mathrm{CKB}$ model this is an ellipsoid, whose equation is linearly
generated by the quadratic equation of the unit sphere
and the equation of the double plane tangent to unit sphere at $\iota_*(\lambda)$.

(iv) $A$ has two different eigenvalues $\lambda_1\neq\lambda_2$, and $A$ is not normal.

Then $\DW_*(A)$ is the an asymptotically closed $h$-tube around the $h$-line connecting $\iota_*(\lambda_1)$ and $\iota_*(\lambda_2)$.
In the $\mathrm{CKB}$ model this is an ellipsoid, whose equation is linearly
generated by the quadratic equation of the unit sphere
and the quadratic equation of the union of planes tangent to unit sphere at $\iota_*(\lambda_1)$ and $\iota_*(\lambda_2)$ .
\begin{proof}
The cases (a) and (b) are trivial.

(c) By a conformal transformation and unity equivalence one can assume that
$
A=L_t\equiv\begin{bmatrix}1&2t\\&-1\end{bmatrix}
$
We can recognize that the M\"obius transformations $f:x\mapsto\frac{ax+b}{bx+a}$ leave $L_t$ invariant,
and these correspond to the translation-rotation group around the line with asymptotic  points $\iota_*(\pm1)$.
Thus the non-asymptotic part  of the shell is a union of possible degenerate tubes around that line, but we already know
that it must be one tube exactly.

(d) Here the representative $A=S_0\equiv\begin{bmatrix}0&1\\&0\end{bmatrix}$ can be taken.
We can recognize that the M\"obius transformations $f:x\mapsto\frac{x}{cx+1}$ leave $S_0$ invariant,
and these correspond to the horo-translation group around asymptotic  point $\iota_*(0)$.
Thus the non-asymptotic part of the shell is a union of horospheres asymtotic at $\iota_*(0)$,
but we already know that it must only one horosphere exactly.

The generation of the tubes and horospheres from a pencil is a standard hyperbolic geometry.

(Here is the critical part of the argument in the style of Davis:

(c) We can take $A$ to be the linear relation
generated by $(\mathbf e_1,0)$ and $(0,t\mathbf e_1-\mathbf e_2)$ as a representative.
[The involutive Cayley transform of $L_t$.]
Then (for $z_1\neq0$),
\[\DW_{\mathrm {PH}}( z_1(\mathbf e_1,0) +z_2 (0,t\mathbf e_1-\mathbf e_2))=
\left( \frac{2t \Rea(\bar z_1 z_2)}{|z_1|^2},\frac{2t \Ima(\bar z_1 z_2)}{|z_1|^2},
\frac{|z_2|}{|z_1|} \right).
\]
While this may look like to yield the cone $x^2+y^2-(4t)^2z^2=0$, but it is $h$-tube with asymptotical points $\iota_{\mathrm{PH}}(0)$
and $\iota_{\mathrm{PH}}(\infty)$.

(d) We can take $A$ to be the linear relation
generated by $(0,\mathbf e_1)$ and $(\mathbf e_1,\mathbf e_2)$ as a representative.
[The inverse of $S_0$.]
Then (for $z_2\neq0$),
\[\DW_{\mathrm {PH}}( z_1(0,\mathbf e_1)  +z_2 (\mathbf e_1,\mathbf e_2) )=
\left( \frac{2 \Rea(\bar z_1 z_2)}{|z_2|^2},\frac{2 \Ima(\bar z_1 z_2)}{|z_2|^2},
1 \right).
\]
This may look like to yield the plane  $z-1=0$, but it is $h$-horoshere with asymptotical point $\iota_{\mathrm{PH}}(\infty)$.

Remark: One can also use $\iota_{\mathrm{CKB}}(z_2/z_1)=( 2\Rea (\bar z_1 z_2), 2\Ima (\bar z_1 z_2), |z_2|^2-|z_1|^2 )=
( 2\Rea (z_1\bar z_2), -2\Ima (z_1\bar z_2), -(|z_1|^2-|z_2|^2) )$  to produce the unit sphere.)
\end{proof}
\end{theorem}
Lins, Spitkovsky, Zhong \cite{LSZ} gives the quadratic equation for shell; cf. also \cite{L2}.

\begin{commentx}
\begin{theorem}(Lins, Spitkovsky, Zhong \cite{LSZ}.)
Let $A$ be a $2\times2$ complex matrix.
Then  $\mathrm{DW}_{\mathrm{CKB(P)}}(A)$ is given by the quadratic equation
\[\begin{bmatrix}x_{\mathrm{CKB(P)}}\\ y_{\mathrm{CKB(P)}}\\ z_{\mathrm{CKB(P)}}\\1\end{bmatrix}^\top
\begin{bmatrix}A_{11}&A_{12}&A_{13}&A_{14}\\A_{12}&A_{22}&A_{23}&A_{24}\\A_{13}&A_{23}&A_{33}&A_{34}\\A_{14}&A_{24}&A_{34}&A_{44}\end{bmatrix}
\begin{bmatrix}x_{\mathrm{CKB(p)}}\\ y_{\mathrm{CKB(P)}}\\ z_{\mathrm{CKB(P)}}\\1\end{bmatrix}=0\]
where
\[A_{11}=\tr(A^*A)+2\Rea(\det A)\]
\[A_{22}=\tr(A^*A)-2\Rea(\det A)\]
\[A_{33}=1\]
\[A_{12}=2\Ima(\det A)\]
\[A_{13}=-\Rea(\tr A)\]
\[A_{23}=-\Ima(\tr A)\]
\[A_{14}=-\Rea(\tr A^*\det A)\]
\[A_{24}=-\Ima(\tr A^*\det A)\]
\[A_{34}=\frac12\left(|\tr A|^2-\tr(A^*A)\right)\]
\[A_{44}=|\det A|^2;\]
and solution set should be restricted to asymptotically closed CKB model
(which restriction is relevant iff $A$ is normal).
\qed
\end{theorem}
\end{commentx}
\snewpage
\begin{point}
Next we consider what happens if we discard the imaginary part of the scalar product in Davis-Wielandt shell.
In the $\mathrm{CKB}$ and $\mathrm{CKB(P)}$ models this means discarding second coordinate.
This can be considered as the projection $\pi^{[2]}_*$ of the hyperbolic $3$-space to a hyperbolic $2$-space.
The maps are compatible to the definition of the Davis-Wielandt shell.

The compositions $\pi^{[2]}_*\circ\DW_*$ can de identified as the real Davis-Wielandt shell $\DW^{\mathbb R}_{*}$
(except that the latter one can also be defined when there is no underlying complex structure).
The theory of $\DW^{\mathbb R}_{*}$ is quite similar but less informative and simpler than complex case.
We will not say much more because we have already met $\DW^{\mathbb R}_{\mathrm{PH}}$ as $\CR$.

If $A$, in particular, is a linear operator on a $2$-dimensional Hilbert space, then we can consider the same operator
acting on the same but $4$-dimensional real Hilbert space.
As we already know the complex case, we immediately see that in the CKB  and CKB(P) models the real Davis--Wielandt shell is a
projection of a possibly degenerate tube of horosphere, thus a
possibly degenerate elliptical disk in the model space, yielding Lemma \ref{lem:ellipdisk}.
\end{point}
\begin{point}
Obviously, the vertical projection (discarding the third coordinate) in CKB(P) model gives the numerical range.
(This is one primary reason to use the $\mathrm{CKB(P)}$ model.)
In terms of hyperbolic geometry, this is the central projection $\pi^{\infty}_*$ from the
the asymptotic point $\infty_*$ to the all the other asymptotic points (and composed with the inverse of $\iota_*$).
\end{point}

\begin{theorem}[Davis]
\plabel{th:DWspect}
$\iota_*(\spec_p(A)\cup \spec_c(A))=\iota_*(\mathbb C)\cap \overline{\DW(A)}$.

\begin{proof}
Due to translation invariance by complex numbers (in the Poincar\'e half space model)
it is sufficient to test for $\lambda =0$.
Then $0\in\spec_p(A)\cup \spec_c(A)$ is equivalent to the co-norm being $0$.
\end{proof}
\end{theorem}
\begin{theorem}[Wielandt]
\plabel{th:Wielandt}
Suppose that $A$ is a linear operator on a finite dimensional Hilbert space.
Then $A$ is normal if and only if $\DW_*(A)$ is an (asymptotical) $h$-polytope
(spanned by the $\iota_*$ of the spectrum).
\begin{proof}
We will consider the $\mathrm{CKB}$  model.
In the normal case this is certainly the situation, as it can be thought as a direct sum of $1$-dimensional
operators. If the operator is not normal, then
it has $2$-dimensional invariant subspace where it acts as a non-normal operator
(due to the a non-trivial Jordan block, or by the non-orthogonality of the pure eigenspaces).
Restricted to this subspace, it yields case (iii) or (iv) of the Theorem \ref{thm:DWconc}.
But then $\DW_{\mathrm{CKB}}(A)$ is ``rounded''  at the corresponding asymptotic points, so it cannot yield an $h$-polytope.
\end{proof}
\end{theorem}
All our canonical models of the hyperbolic space come a(n equicalent) canonical reflection, which
sign change in the second coordinate.
With respect to this reflection, using the standard separation argument, one can prove
\begin{theorem}[Li, Poon, Sze,\cite{LPS}]
\plabel{th:LPS}
If $\dim_{\mathbb C} \mathfrak H<\infty$, then
$\DW_*(A)=\DW_*(A^*)^{\mathrm{refl}}$.
\qed
\end{theorem}
Other, still quite elementary arguments are as follows.
\begin{lemma}
\plabel{lem:Wielandt}
Assume that $\iota_{\mathrm{CKB(P)}}(\lambda)\in\DW_{\mathrm{CKB(P)}}(A)$, and $\DW_{\mathrm{CKB(P)}}(A)$ allows more than one
supporting plane at $\iota_{\mathrm{CKB(P)}}(\lambda)$.

Suppose that $\iota_{\mathrm{CKB(P)}}(\lambda)=\DW_{\mathrm{CKB(P)}}((\mathbf x,A\mathbf x))$.
Then $\mathbb C\mathbf x\oplus \mathbf x^{\bot}_{\mathfrak H} $ is an invariant orthogonal decomposition of $\mathfrak H$,
with $A\mathbf x=\lambda\mathbf x$.

Alternatively put: The $\lambda$-eigenspace $V_{\lambda}^A$ of $A$ yields an $A$-invariant orthogonal decomposition
$V_{\lambda}^A \oplus (V_{\lambda}^A)^{\bot}_{\mathfrak H}$ of $\mathfrak H$, and $A$ restricted to  $V_{\lambda}^A$ is normal.
\begin{proof}
Passing to an appropriate linear transform $\mathrm e^{\mathrm i\theta}(A-\lambda)$, if necessary, we can assume that
$\lambda=0$, and $z=\varepsilon x$ ($\varepsilon\neq0$) is another supporting plane beside $z=0$.
Then the derivative of the first coordinate in  $\DW_{\mathrm{CKB(P)}}((\mathbf x,A\mathbf x))$ in $\mathbf x$ must be $0$, thus
we obtain $\frac{A+A^*}2\mathbf x=0$.
Together with the fact $A\mathbf x=0$, this implies $A^*\mathbf x=0$.
Then, by the standard linear arguments, the invariance of the decomposition
$\mathbb C\mathbf x\oplus \mathbf x^{\bot}_{\mathfrak H} $
follows.
\end{proof}
\end{lemma}
\begin{proof}[Alternative proof to Theorem \ref{th:Wielandt}]
For the main argument:
Assume $\DW_{\mathrm{CKB(P)}}((\mathbf x,A\mathbf x))$ allows a second supporting plane at
every point of $\iota_{\mathrm{CKB(P)}}(\spec(A))$.
Then Lemma \ref{lem:Wielandt} can be applied to splitting $\mathfrak H$ away inductively.
\end{proof}

\begin{lemma}
\plabel{lem:Wielandt2}
(a) Assume that $P$ is an interior point of the
CKB(P) model (i. e. $P\notin \iota_{\mathrm{CKB(P)}}(\mathbb C)$), and  $P\in\DW_{\mathrm{CKB(P)}}(A)$, and
$\DW_{\mathrm{CKB(P)}}(A)$ allows at least two
supporting planes $\pi_1,\pi_2$ at $P$.

Then $\pi_1\cap\pi_2$ intersects $\iota_{\mathrm{CKB(P)}}(\mathbb C)$ at two points, $ \iota_{\mathrm{CKB(P)}}(\lambda_1)$ and $ \iota_{\mathrm{CKB(P)}}(\lambda_2)$.
(In particular, $\pi_1\cap\pi_2$ is not vertical.)
Furthermore,  $\bar \lambda_1$ and $\bar \lambda_2$ are eigenvalues of $A^*$ associated to some eigenvectors.

(b) Assume additionally, that $\overline{\DW_{\mathrm{CKB(P)}}(A)}=\overline{\DW_{\mathrm{CKB(P)}}(A^*)}^{\mathrm{refl}}$
(this holds automatically in the finite dimensional case).

Then the segment connecting $ \iota_{\mathrm{CKB(P)}}(\lambda_1)$ and $ \iota_{\mathrm{CKB(P)}}(\lambda_2)$ (and containing $P$)
is contained in $\DW_{\mathrm{CKB(P)}}(A)$.

Suppose that $\iota_{\mathrm{CKB(P)}}(\lambda)=\DW_{\mathrm{CKB(P)}}((\mathbf x,A\mathbf x))$.
Then there exist two nonzero, orthogonal vectors $\mathbf x_1,\mathbf x_2\in\mathfrak H$ (uniquely), such that
$\mathbb C\mathbf x_1\oplus \mathbb C\mathbf x_2\oplus\{\mathbf x_1,\mathbf x_2\}^{\bot}_{\mathfrak H} $ is an invariant orthogonal decomposition of $\mathfrak H$,
with $A\mathbf x_i=\lambda_i\mathbf x_i$; and $\mathbf x=\mathbf x_1\oplus\mathbf x_2$.

Alternatively put: The $\{\lambda_1,\lambda_2\}$-eigenspace $V_{\{\lambda_1,\lambda_2\}}^A$ of $A$ yields an $A$-invariant
orthogonal decomposition
$V_{\{\lambda_1,\lambda_2\}}^A \oplus (V_{\{\lambda_1,\lambda_2\}}^A)^{\bot}_{\mathfrak H}$ of $\mathfrak H$, and $A$ restricted to  $V_{\{\lambda_1,\lambda_2\}}^A$ is normal, both eigenvalues are present.
\begin{proof}
(a) First, we assume that $\pi_1\cap\pi_2$ intersects $\iota_{\mathrm{CKB(P)}}(\mathbb C)$ at two points (i. e. it is not vertical).
Passing to an appropriate linear transform of $A$, if necessary, we can assume that
$(z_{\mathrm{CKB(P)}}-1) =m_1(x_{\mathrm{CKB(P)}}-0)$
and
$(z_{\mathrm{CKB(P)}}-1) =m_2(x_{\mathrm{CKB(P)}}-0)$
are the two supporting planes, $m_1\neq m_2$.
Then $\lambda_1=\mathrm i$ and $\lambda_2=-\mathrm i$ can be chosen.
Then the derivatives of the first and third coordinates in  $\DW_{\mathrm{CKB(P)}}((\mathbf x,A\mathbf x))$ in $\mathbf x$ must be $0$, thus
we obtain $\frac{A+A^*}2\mathbf x=0$ and  $A^*A\mathbf x=\mathbf x$.
This implies that the linear span of $\mathrm{Span}\{\mathbf x,A\mathbf x\}=\mathrm{Span}\{\mathbf x,A^*\mathbf x\}$ is
2-dimensional vector space invariant for $A^*$ with $(A^*)^2+\Id$ vanishing on it.
Note that $\mathbf x$ cannot be an eigenvector of $A^*$.
Indeed, if, for example, $A^*\mathbf x=\mathrm i\mathbf  x$ holds, then $A\mathbf x=-\mathrm i\mathbf  x$ holds, and
$P=\iota_{\mathrm{CKB(P)}}(-\mathrm i)$ is a contradiction to our assumptions.
Thus there are $\mathrm{Span}\{\mathbf x,A\mathbf x\}=\mathrm{Span}\{\mathbf x,A^*\mathbf x\}$
has two (unique) nonzero $A^*$-eigenvectors $\mathbf x_1$ and $\mathbf x_2$
such that $A^*\mathbf x_i=\bar \lambda_i$ and $\mathbf x_1+\mathbf x_2=\mathbf x$.
(The choice of normalization somewhat confuses, say, $\lambda_1$ and $\bar\lambda_2$, but after transformed back
the eigenvalues are as they should be.)

Consider now the situation when $\pi_1\cap\pi_2$ is vertical. Then an appropriate fractional
linear transform of $A$, very close to the original, has the property that its adjoint has a very large eigenvalue.
This contradicts to the boundedness of $A$.
(An alternative, and more generalizable argument is the following:
We can assume that the line is $x_{\mathrm{CKB(P)}}=y_{\mathrm{CKB(P)}}=0$.
Then, by differentiation, $\frac{A+A^*}{2}\mathbf x=0$ and $\frac{A-A^*}{2\mathrm i}\mathbf x=0$ holds,
contradicting to $\mathbf x^*A^*A\mathbf x>0$.)

(b) Let us continue the first half of the argument for (a).
Due to the extra condition, Lemma \ref{lem:Wielandt} can be applied to $A^*$ with respect to $\mathbf x_1$
and $\mathbf x_2$.
This proves that the decomposition is orthogonal, and in particular, $A\mathbf x_i=\lambda_i\mathbf x_i$ holds.
\end{proof}
\end{lemma}

\begin{cor}
\plabel{cor:Wielandt3}
Assume that $A$ acts on a finite dimensional Hilbert space $\mathfrak A$. Then

(i) $\DW_{\mathrm{CKB(P)}}(A)$ has at most finitely many angular vertices
(at those points with infinitely many support planes whose intersection is single point)
corresponding to some $\iota_{\mathrm{CKB(P)}}(\lambda)$.

(ii)
In the interior of the $\mathrm{CKB(P)}$ model, $\DW_{\mathrm{CKB(P)}}(A)$ may have finitely many angular edges
(at whose points with infinitely many support planes whose intersection is a line)
connecting certain angular vertices.

(iii) Apart from, at any other boundary point $P$ of $\DW_{\mathrm{CKB(P)}}(A)$ (either on the boundary, or in the interior)
 there is only one supporting plane.

(iv) The generalized eigenspace $V_{\mathrm{vertex}}^A$ allows the $V_{\mathrm{vertex}}^A \oplus (V_{\mathrm{vertex}}^A)^{\bot}_{\mathfrak H}$ of $\mathfrak H$, and $A$ restricted to  $V_{\mathrm{vertex}}^A$ is normal, all eigenvalues corresponding the vertices to are present.
\begin{proof}
It is an immediate consequence of
Lemma \ref{lem:Wielandt}
and
Lemma \ref{lem:Wielandt2}.
\end{proof}
\end{cor}

\begin{theorem}
\plabel{thm:davisline}
Assume that $P$ is an interior point of the
CKB(P) model (i. e. $P\notin \iota_{\mathrm{CKB(P)}}(\mathbb C)$), and  $p\in\overline{\DW_{\mathrm{CKB(P)}}(A)}$, and
$\overline{\DW_{\mathrm{CKB(P)}}(A)}$ allows at least two
supporting planes $\pi_1,\pi_2$ at $p$.

Then $\pi_1\cap\pi_2$ intersects $\iota_{\mathrm{CKB(P)}}(\mathbb C)$ at two points, $ \iota_{\mathrm{CKB(P)}}(\lambda_1)$
(in particular, $\pi_1\cap\pi_2$ is not vertical),
and $ \iota_{\mathrm{CKB(P)}}(\lambda_2)$; and $\lambda_1,\lambda_2\in\spec(A)$.
\begin{proof}
First we prove that  $\pi_1\cap\pi_2$ cannot be vertical.
Assume that it is so.
By taking an appropriate linear transform of $A$, we can assume that for shell
\[|y_{\mathrm{CKB(P)}}|\leq k_0 x_{\mathrm{CKB(P)}}\]
with it $k_0\geq 0$, and $(0,0)$ is in the vertical projection of the closure of the shell.
This makes
\[E:=\frac{A+A^*}{2}\geq0,
\qquad\text{and}\qquad
F:=\frac{A-A^*}{2\mathrm i}\text{ self-adjoint}.\]
Assume that $P=(0,0,p)$, $p>0$.
Then, there is sequence $\mathbf x_i$ such that $\langle A \mathbf x_i,\mathbf x_i\rangle\rightarrow0$
but $\langle A \mathbf x_i,A\mathbf x_i\rangle\rightarrow p$ while $\mathbf x_i=1$.
In particular,
\[\langle E \mathbf x_i,\mathbf x_i\rangle\rightarrow0
\qquad\text{ and}\qquad
\langle F \mathbf x_i,\mathbf x_i\rangle\rightarrow0.\]
Using the Cauchy--Schwarz inequality for $\langle E\cdot,\cdot\rangle$, we have the general inequality
\[|\langle E\mathbf x_i, H\mathbf x_i\rangle|^2\leq
\langle E\mathbf x_i, \mathbf x_i\rangle\langle EH\mathbf x_i, H\mathbf x_i\rangle\leq
\langle E\mathbf x_i, \mathbf x_i\rangle \|E\|_2 \|H\|_2^2 |\mathbf x_i|_2^2.
 \]
Consequently,
\[\langle E\mathbf x_i, E\mathbf x_i\rangle\rightarrow0
\qquad\text{and}\qquad
\langle E\mathbf x_i, F\mathbf x_i\rangle\rightarrow0.\]
Considering  $\langle A \mathbf x_i,A\mathbf x_i\rangle=
\langle (E+\mathrm iF) \mathbf x_i,(E+\mathrm iF)\mathbf x_i\rangle$, it is not hard to see that
\begin{equation}\langle F\mathbf x_i, F\mathbf x_i\rangle\rightarrow p.\plabel{eq:xFF}\end{equation}

Let $k>k_0$ arbitrary. Consider the polynomial (in $t$).
\begin{equation}f_i(t):=
\langle F(\mathbf x_i+tF\mathbf x_i),(\mathbf x_i+tF\mathbf x_i)\rangle
-k\langle E(\mathbf x_i+tF\mathbf x_i),(\mathbf x_i+tF\mathbf x_i)\rangle.
\plabel{eq:daviscorner}
\end{equation}
Then the $t^2$ coefficient is bounded, the $t$ coefficient limits to $2p$, the $t^0$ coefficient limits to $0$.
Thus for any sufficiently large $i$ and we can chose a value $t_i$ such that
\begin{equation}t_i\rightarrow0\plabel{eq:davisroot}\end{equation}
 and
\[f_{i}(t_i)=0.\]
Assume, for a moment, that it happen infinitely many times that $|\mathbf x_{i}+t_iF\mathbf x_{i}|_2=0$.
Then for this subsequence, $(t_i)^2\rightarrow p$ by \eqref{eq:xFF};
 which is in  contradiction to \eqref{eq:davisroot}.

Thus, in general $\mathbf x_{i}+t_iF\mathbf x_{i}\neq0$.
But then it induces an element in  $\DW_{\mathrm{CKB(P)}}(A)$
such that $y_{\mathrm{CKB(P)}}= k x_{\mathrm{CKB(P)}}$.
As $k>k_0$, this can happen only if $x_{\mathrm{CKB(P)}}=y_{\mathrm{CKB(P)}}=0$.
Consequently for large $i$, this yields a  point of the shell in the interior of the vertical line.
Then Lemma \ref{lem:Wielandt2} shows that this is impossible.

As for the general case, $\lambda_i\notin\spec(A)$ would yields a contradiction to $(A-\lambda_i \Id)^{-1}$.
\end{proof}
\end{theorem}
One may wonder about possible strengthenings of the theorem above.
Davis \cite{D2} claims a stronger result but seems to prove only Lemma \ref{lem:Wielandt2}(a).

The dual viewpoint is put forward by Li, Poon, Sze \cite{LPS}
who, in particular, emphasize the upper boundary (that is convex view from infinity, i. e. the norm branch of the boundary).
$C^*$-algebraic view is advocated in Aramba\v{s}i\'c, Beri\'c, Raji\'c \cite{ABR}.

In general, all constructions exposed regarding the dual view in the case
of conformal range, apply more generally to the Davis--Wieland shell.
With respect to the finite dimensional case, this, leads in particular to
the real homogeneous polynomal
\[K_A^{\DW}(u,v,s,w)=\det\left(u\frac{A+A^*}{2}+v\frac{A-A^*}{2\mathrm i}+sA^*A+w\Id \right)=0, \]
and / or
\[F_A^{\DW}(\lambda_1,\lambda_2,\nu)=\det\left(\nu\Id-(A^*-(\lambda_1-\mathrm i\lambda_2))(A-(\lambda_1+\mathrm i\lambda_2)) \right)=0, \]
which determine
$\Dual_{\mathrm{CKB(P)}}^{\mathrm{alg}}(A)$ and $\Dual_{\mathrm{CKB(P)}}^{\mathrm{ncn},\mathrm{alg}}(A)$ respectively, and
thus the Davis--Wielandt shell. Here
\[F_A^{\DW}(\lambda_1,\lambda_2,\nu)=F_{A^*}^{\DW}(\lambda_1,-\lambda_2,\nu)\]
and
\[K_A^{\DW}(u,v,s,w)=K_{A^*}^{\DW}(u,-v,s,w)\]
play the role of Theorem \ref{th:LPS}.

More about the  spectral properties of the Davis--Wielandt shell
can be found in Davis \cite{D1}, \cite{D2}, and Li, Poon, Sze \cite{LPS}.
The lacunar properties of the Davis--Wielandt shell are similar to the ones of the conformal range,
but they appear to be not particularly studied.
The Davis--Wielandt shell is generalized in Davis \cite{D3}.
The Davis--Wielandt shell is related to other types of ranges in
Lins, Spitkovsky, Zhong \cite{LSZ}.
\snewpage

We emphasize that the enveloping construction can also be carried out with respect to the Davis--Wielandt shell,
also making finite dimensional case (more or less) computable.
\begin{example}
\plabel{ex:DWenvel}
We will consider the case
\[A=\bem1&1&\\&&-1\\&&-1\eem.\]
This case is somewhat special as $A$ is real, thus the shell will be symmetric in the second coordinate,
but, otherwise,  it can demonstrate the general principles.

In the CKB(P) model, we have to consider the enveloping surface of the planes
\[-2\lambda_1 x_{\mathrm{CKB(P)}} -2\lambda_2 y_{\mathrm{CKB(P)}}+z_{\mathrm{CKB(P)}}+\lambda_1^2+\lambda_2^2-N(A-(\lambda_1+\mathrm i\lambda_2)\Id)=0,\]
where $N$ is the square of the norm or the co-norm.

Using the abbreviation  $N\equiv N(A-(\lambda_1+\mathrm i\lambda_2)\Id)$,  the enveloping surface is  parametrized as
\[\lambda_1+\mathrm i\lambda_2\mapsto E_{\mathrm{CKB(P)}}(\lambda_1+\mathrm i\lambda_2)=
\left(\lambda_1-\frac12 \frac{\mathrm d N}{\mathrm d\lambda_1},
\lambda_2-\frac12 \frac{\mathrm d N}{\mathrm d\lambda_2},
\lambda_1^2+\lambda_2^2-\lambda_1\frac{\mathrm d N}{\mathrm d\lambda_1}-\lambda_2\frac{\mathrm d N}{\mathrm d\lambda_2}+ N\right).\]

Again, instead of just the norm and co-norm, we can apply to this to all branches of the solutions $N$ of
\[\underbrace{\det\left(N\Id- (A-(\lambda_1+\mathrm i\lambda_2)\Id)^*(A-(\lambda_1+\mathrm i\lambda_2)\Id)   \right)}
_{F_A^{\DW}(\lambda_1,\lambda_2,N)}=0.\]

Now, $\frac{\mathrm d N}{\mathrm d\lambda_1}$ and  $\frac{\mathrm d N}{\mathrm d\lambda_2}$
can be (generically) be expressed as rational functions of $\lambda_1,\lambda_2,N$.
This yields
\[E_{\mathrm{CKB(P)}}(\lambda_1+\mathrm i\lambda_2)=\left(\frac{X_1}{X_4},\frac{X_2}{X_4},\frac{X_3}{X_4}\right)\]
where
\begin{align*}
X_1\equiv&-3\,\lambda_{{1}} \left( -2\,{\lambda_{{1}}}^{2}-2\,{\lambda_{{2}}}^{2}+2\,N-1 \right),
\\
X_2\equiv&-\lambda_{{2}} \left( -2\,{\lambda_{{1}}}^{2}-2\,{\lambda_{{2}}}^{2}+2\,N-3 \right),
\\
X_3\equiv&-3\,{\lambda_{{1}}}^{6}-9\,{\lambda_{{1}}}^{4}{\lambda_{{2}}}^{2}-9\,{
\lambda_{{1}}}^{2}{\lambda_{{2}}}^{4}-3\,{\lambda_{{2}}}^{6}+9\,N{
\lambda_{{1}}}^{4}+18\,N{\lambda_{{1}}}^{2}{\lambda_{{2}}}^{2}+9\,N{
\lambda_{{2}}}^{4}\\&-9\,{N}^{2}{\lambda_{{1}}}^{2}-9\,{N}^{2}{\lambda_{{
2}}}^{2}+10\,{\lambda_{{1}}}^{4}+8\,{\lambda_{{1}}}^{2}{\lambda_{{2}}}
^{2}-2\,{\lambda_{{2}}}^{4}\\&+3\,{N}^{3}-2\,N{\lambda_{{1}}}^{2}+10\,N{
\lambda_{{2}}}^{2}-8\,{N}^{2}+2\,{\lambda_{{1}}}^{2}+2\,{\lambda_{{2}}
}^{2}+4\,N,
\\
X_4\equiv&3\,{\lambda_{{1}}}^{4}+6\,{\lambda_{{1}}}^{2}{\lambda_{{2}}}^{2}+3\,{
\lambda_{{2}}}^{4}-6\,N{\lambda_{{1}}}^{2}-6\,N{\lambda_{{2}}}^{2}+3\,
{N}^{2}+2\,{\lambda_{{1}}}^{2}+6\,{\lambda_{{2}}}^{2}-8\,N+4
\\
&=\partial_3F_A^{\DW}(\lambda_1,\lambda_2,N).
\end{align*}

\snewpage
Then one finds that the enveloping construction should lie on the surface
\\

$p_{\mathrm K}(x,y,z)\equiv
128\,{x}^{10}+3712\,{x}^{8}{y}^{2}-1035\,{x}^{8}{z}^{2}+38144\,{x}^{6}
{y}^{4}+9108\,{x}^{6}{y}^{2}{z}^{2}+2160\,{x}^{6}{z}^{4}+158976\,{x}^{
4}{y}^{6}+66814\,{x}^{4}{y}^{4}{z}^{2}+8720\,{x}^{4}{y}^{2}{z}^{4}-368
\,{x}^{4}{z}^{6}+217728\,{x}^{2}{y}^{8}+118548\,{x}^{2}{y}^{6}{z}^{2}+
25296\,{x}^{2}{y}^{4}{z}^{4}+3360\,{x}^{2}{y}^{2}{z}^{6}+384\,{x}^{2}{
z}^{8}+93312\,{y}^{10}+127413\,{y}^{8}{z}^{2}+63792\,{y}^{6}{z}^{4}+
13968\,{y}^{4}{z}^{6}+1152\,{y}^{2}{z}^{8}+2948\,{x}^{8}z-45488\,{x}^{
6}{y}^{2}z-11596\,{x}^{6}{z}^{3}-247336\,{x}^{4}{y}^{4}z-56420\,{x}^{4
}{y}^{2}{z}^{3}-64\,{x}^{4}{z}^{5}-395568\,{x}^{2}{y}^{6}z-158436\,{x}
^{2}{y}^{4}{z}^{3}-27648\,{x}^{2}{y}^{2}{z}^{5}-2880\,{x}^{2}{z}^{7}-
196668\,{y}^{8}z-203724\,{y}^{6}{z}^{3}-68544\,{y}^{4}{z}^{5}-7488\,{y
}^{2}{z}^{7}-2300\,{x}^{8}+42128\,{x}^{6}{y}^{2}+24856\,{x}^{6}{z}^{2}
+138136\,{x}^{4}{y}^{4}+139664\,{x}^{4}{y}^{2}{z}^{2}+8152\,{x}^{4}{z}
^{4}+140688\,{x}^{2}{y}^{6}+318168\,{x}^{2}{y}^{4}{z}^{2}+81648\,{x}^{
2}{y}^{2}{z}^{4}+8928\,{x}^{2}{z}^{6}+46980\,{y}^{8}+201312\,{y}^{6}{z
}^{2}+122904\,{y}^{4}{z}^{4}+19872\,{y}^{2}{z}^{6}-25323\,{x}^{6}z-
129681\,{x}^{4}{y}^{2}z-23688\,{x}^{4}{z}^{3}-181857\,{x}^{2}{y}^{4}z-
113568\,{x}^{2}{y}^{2}{z}^{3}-14736\,{x}^{2}{z}^{5}-77499\,{y}^{6}z-
106776\,{y}^{4}{z}^{3}-27792\,{y}^{2}{z}^{5}+10296\,{x}^{6}+27576\,{x}
^{4}{y}^{2}+26361\,{x}^{4}{z}^{2}+24264\,{x}^{2}{y}^{4}+75474\,{x}^{2}
{y}^{2}{z}^{2}+13896\,{x}^{2}{z}^{4}+6984\,{y}^{6}+45801\,{y}^{4}{z}^{
2}+22104\,{y}^{2}{z}^{4}-10854\,{x}^{4}z-19224\,{x}^{2}{y}^{2}z-7452\,
{x}^{2}{z}^{3}-8370\,{y}^{4}z-10044\,{y}^{2}{z}^{3}+324\,{x}^{4}+648\,
{x}^{2}{y}^{2}+2106\,{x}^{2}{z}^{2}+324\,{y}^{4}+2430\,{y}^{2}{z}^{2}-
243\,{x}^{2}z-243\,{y}^{2}z=0
$
~\\~\\
(the indices $\mathrm{CKB(P)}$ are omitted from the coordinates from now here.)

[The polynomial $p_{\mathrm K}(x,y,z)$ can be obtained by any sufficiently
powerful algebraic solver.
However, in the present setting, formally, we have to find the
the multivariable discriminant of
\[P:=
F^{\DW}_A(x,y,(\lambda_1)^2+(\lambda_2)^2-2\lambda_1x-2\lambda_2y+z )=
K^{\DW}_A(2\lambda_1,2\lambda_2,-1,-2\lambda_1-2\lambda_2x+z)\]
in $\lambda_1$ and $\lambda_2$.
In order to obtain this, it is a good enough idea to take iterated discriminants,
cf. Sharipov \cite{Sh}.
Indeed, we find
\[\gcd( \Discr_{\lambda_1}(\Discr_{\lambda_2}(P)), \Discr_{\lambda_2}(\Discr_{\lambda_1}(P))   )\sim (z-2)p_{\mathrm K}(x,y,z). \]
This checks out; in fact, $z-2$ can be omitted.]

As there is a symmetry $y\leftrightarrow-y$, we see that  $\DW^{\mathbb R}_{\mathrm{CKB(P)}}(A)$
which is the projection of the shell to the $y=0$ plane, is also the intersection with the $y=0$
plane.
When we restrict to $y=0$, we find
\begin{align*}
p_{\mathrm K}(x,0,z)=&{x}^{2} \left( x-3+2\,z\right) ^{2} \left( x+3-2\,z \right) ^{2} \\
&\cdot\left( 128\,{x}^{4}-11\,{x}^{2}{z}^{2}+24\,{z}^{4}-124\,{x}^{2
}z-36\,{z}^{3}+4\,{x}^{2}+18\,{z}^{2}-3\,z \right).
\end{align*}
This, compared to the result of Example \ref{ex:envel2}, includes some ungeometrical components.
This behaviour, however, is ungenerical.
(Similar behaviour is observed by Chien and Nakazato \cite{CN} for joint numerical ranges.)

In  Figure \ref{fig:figB06} we include sections for $y=0$, $y=1/200$, $y=1/20$, $y=1/4$,  $y=1/2$, $y=2/3$.
(We also indicate the asymptotical points of the CKB(P) model in the pictures.)

If we remove the ungeometrical parts, then we obtain the geometric enveloping surface which is like in Figure \ref{fig:figB07}(a).

\snewpage
\begin{figure}[H]
  \centering
  \begin{subfigure}[b]{0.4\linewidth}
    \includegraphics[width=2.5in]{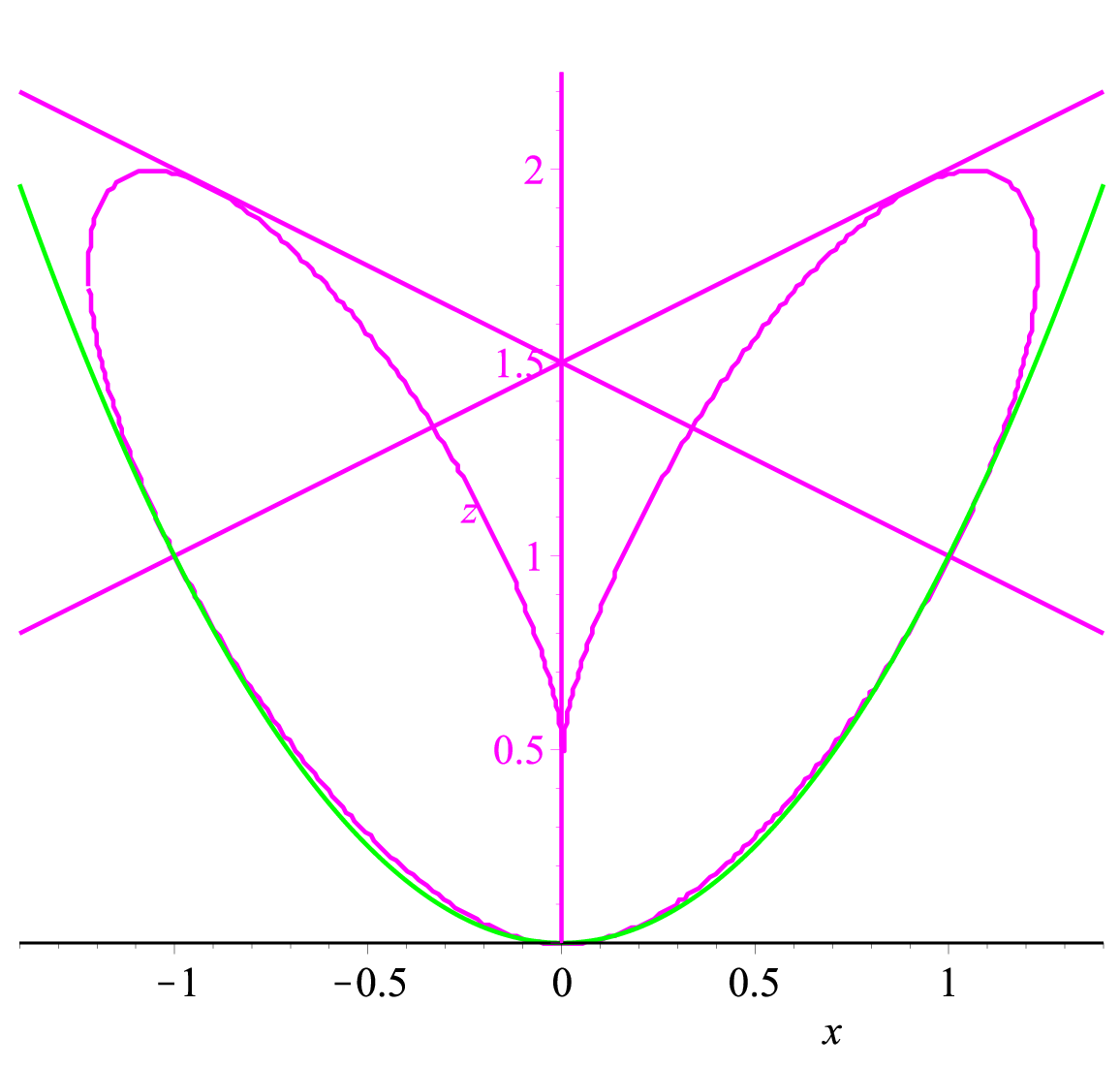}
    \caption*{Fig. \ref{fig:figB06}(a) section $y=0$}
  \end{subfigure}
  \quad
  \begin{subfigure}[b]{0.4\linewidth}
    \includegraphics[width=2.5in]{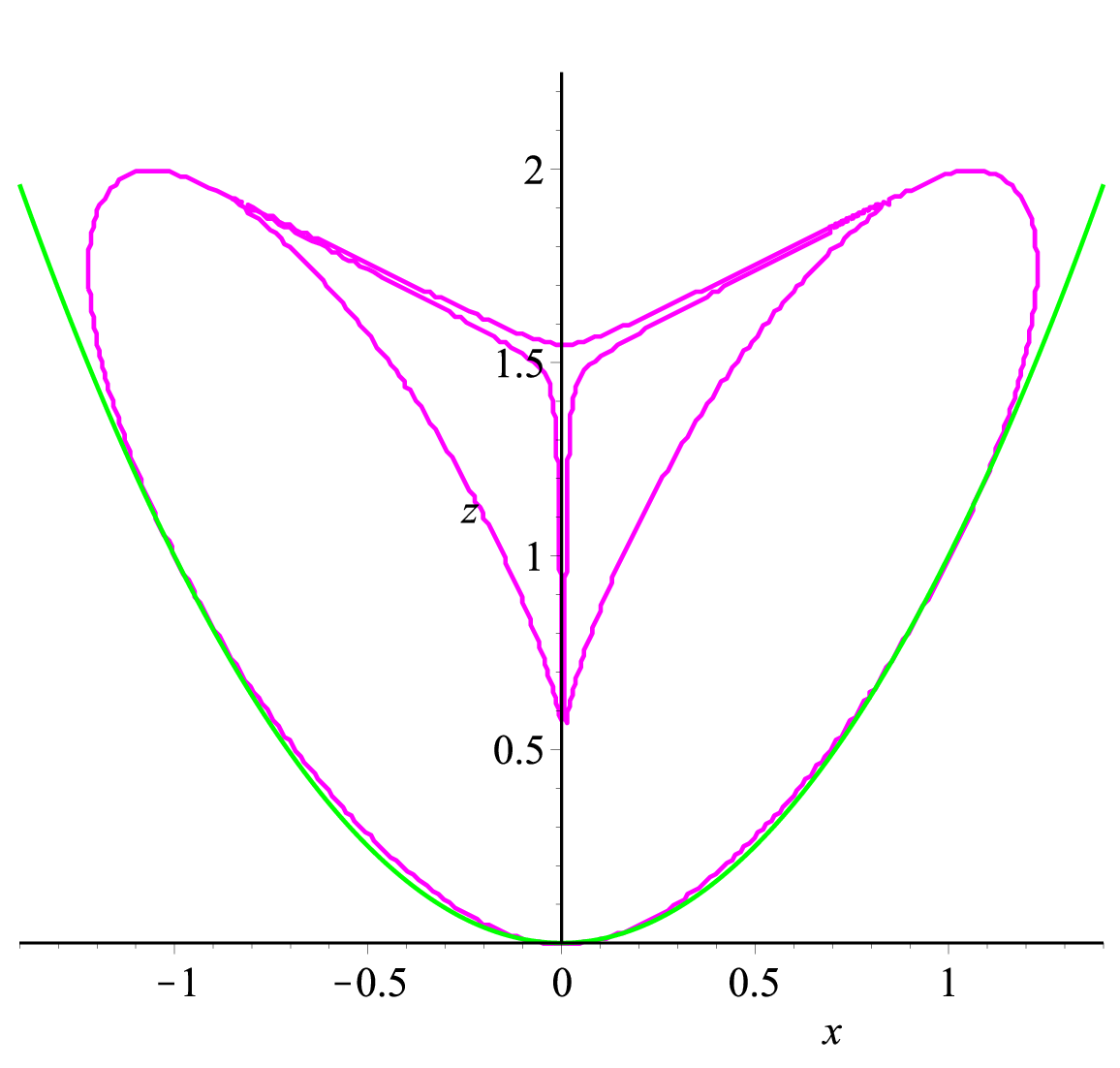}
    \caption*{\ref{fig:figB06}(b) section $y=1/200$  }
  \end{subfigure}
\phantomcaption
\plabel{fig:figB06}
\end{figure}
\begin{figure}[H]
  \ContinuedFloat
  \begin{subfigure}[b]{0.4\linewidth}
    \includegraphics[width=2.5in]{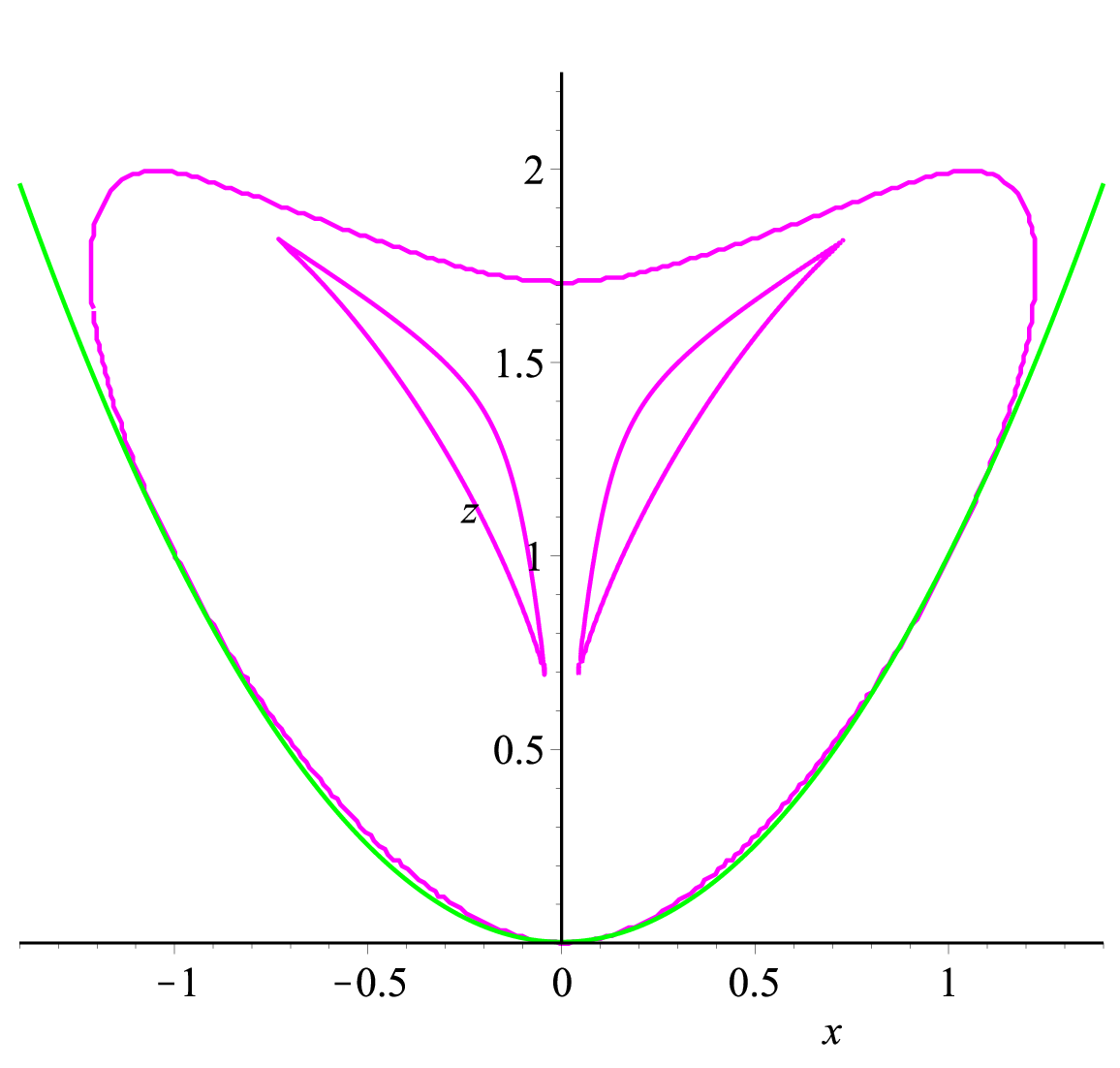}
    \caption*{\ref{fig:figB06}(c) section $y=1/20$  }
  \end{subfigure}
  \qquad
  \begin{subfigure}[b]{0.4\linewidth}
    \includegraphics[width=2.5in]{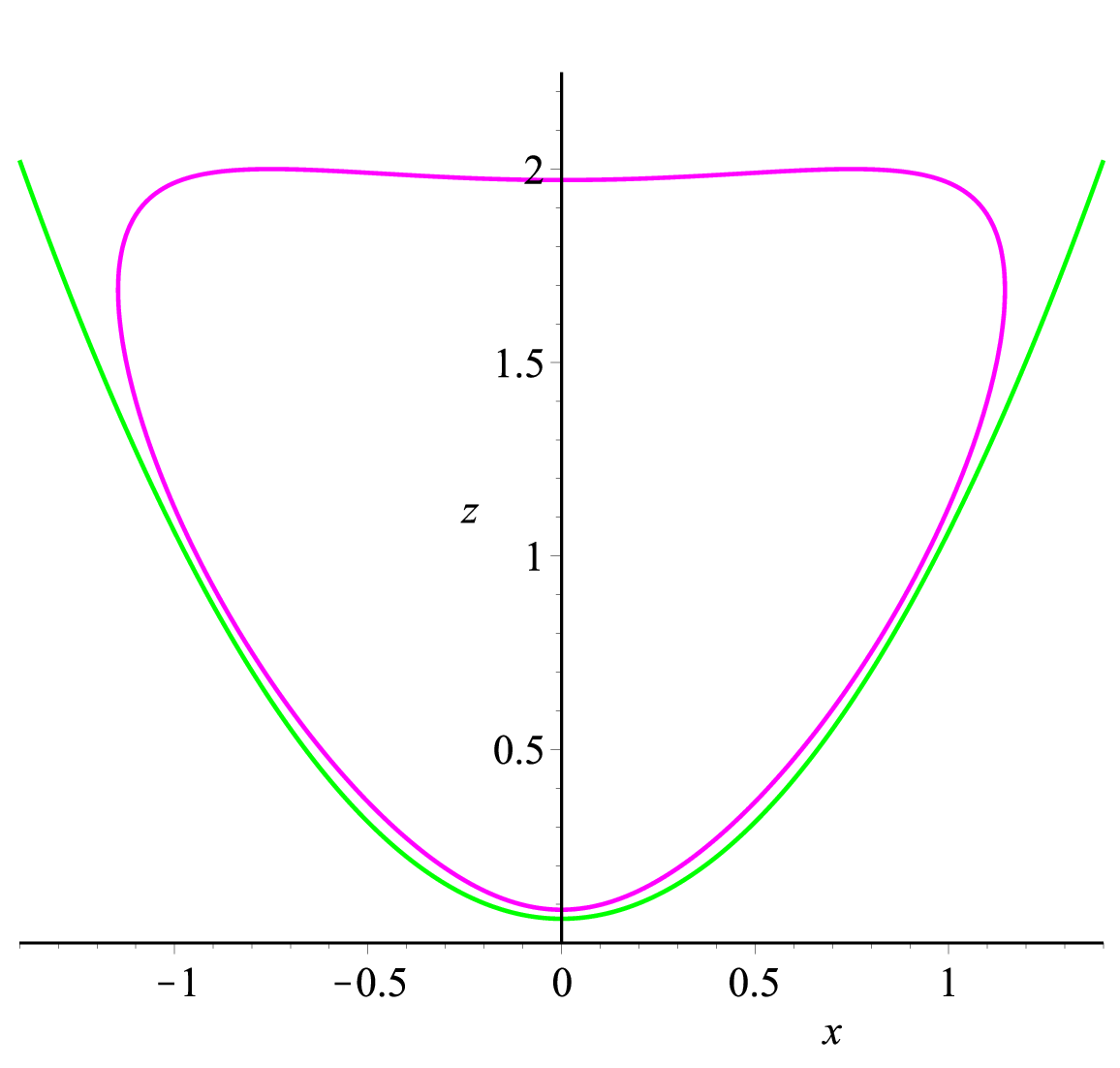}
    \caption*{\ref{fig:figB06}(d) section $y=1/4$  }
  \end{subfigure}
  \phantomcaption
\end{figure}
\vspace{-5mm}
\begin{figure}[H]
   \ContinuedFloat
   \begin{subfigure}[b]{0.4\linewidth}
    \includegraphics[width=2.5in]{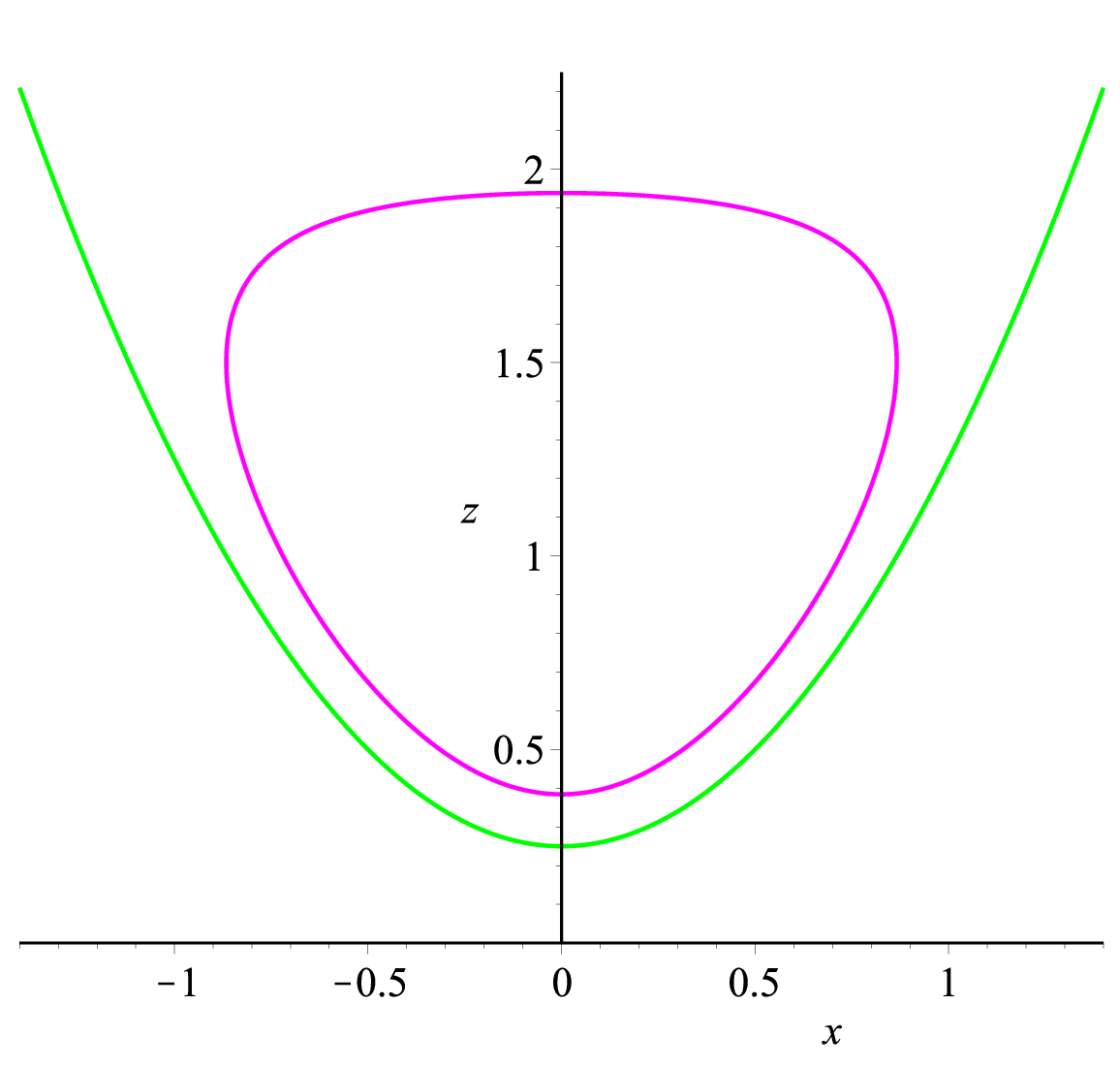}
    \caption*{\ref{fig:figB06}(e) section $y=1/2$}
  \end{subfigure}
  \qquad
  \begin{subfigure}[b]{0.4\linewidth}
    \includegraphics[width=2.5in]{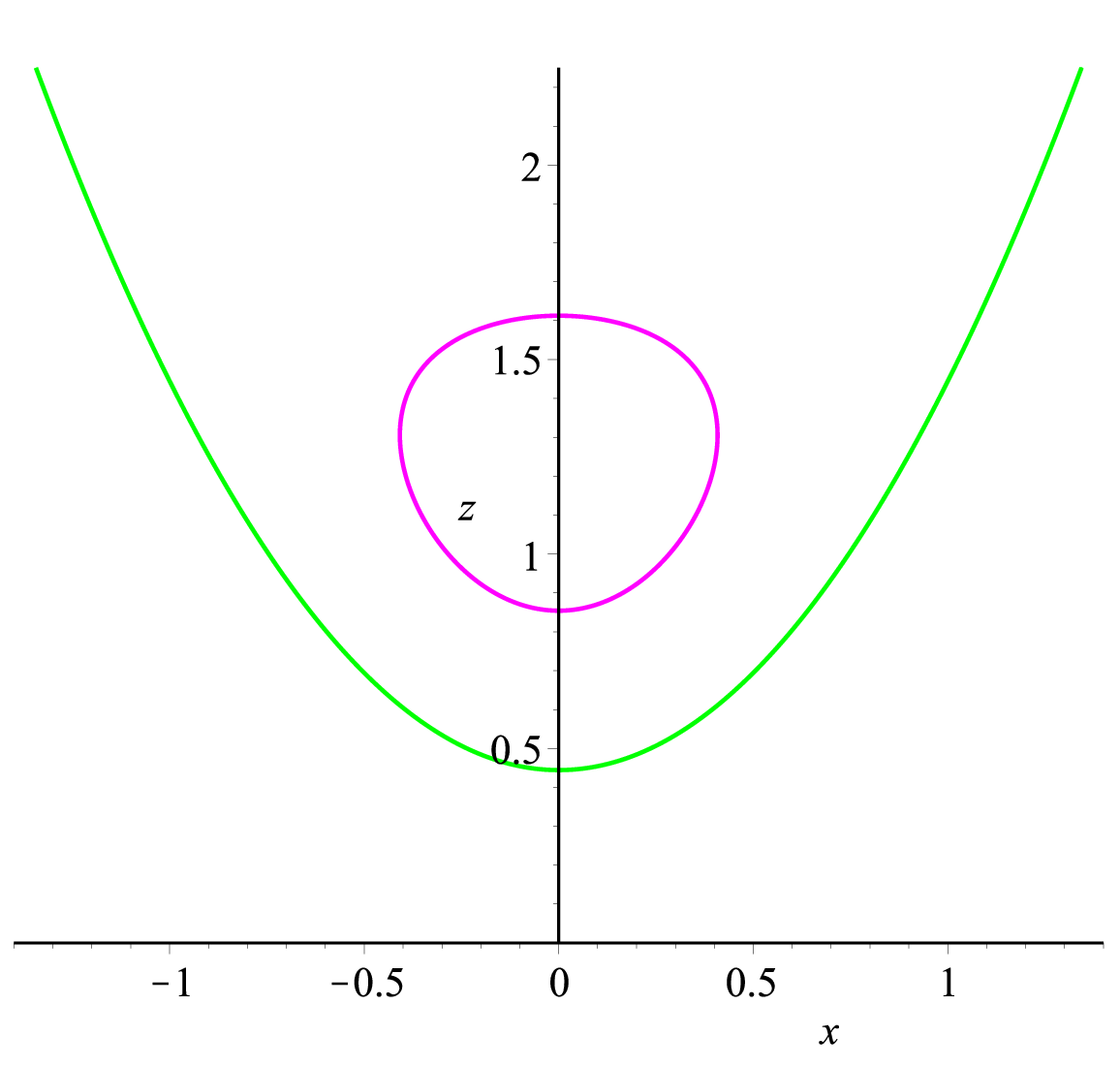}
    \caption*{\ref{fig:figB06}(f) section $y=2/3$}
  \end{subfigure}
\phantomcaption
\end{figure}
\vspace{-5mm}
\begin{figure}[H]
  \begin{subfigure}[b]{2.5in}
    \includegraphics[width=3in]{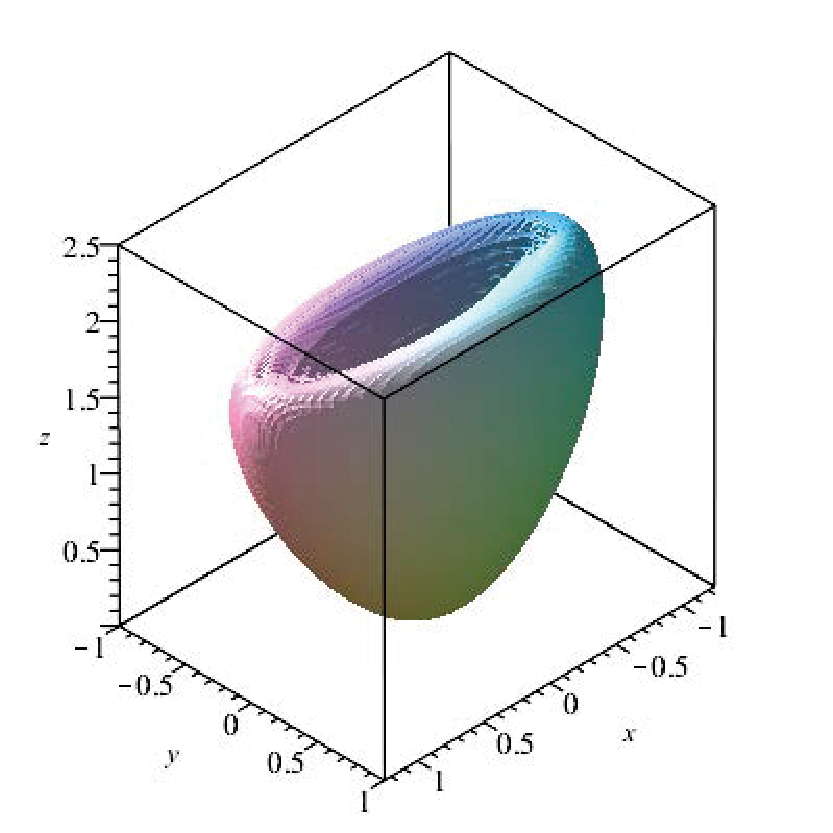}
    \caption*{Fig. \ref{fig:figB07}(a) the geometrical part}
  \end{subfigure}
\qquad
  \begin{subfigure}[b]{2.5in}
    \includegraphics[width=3in]{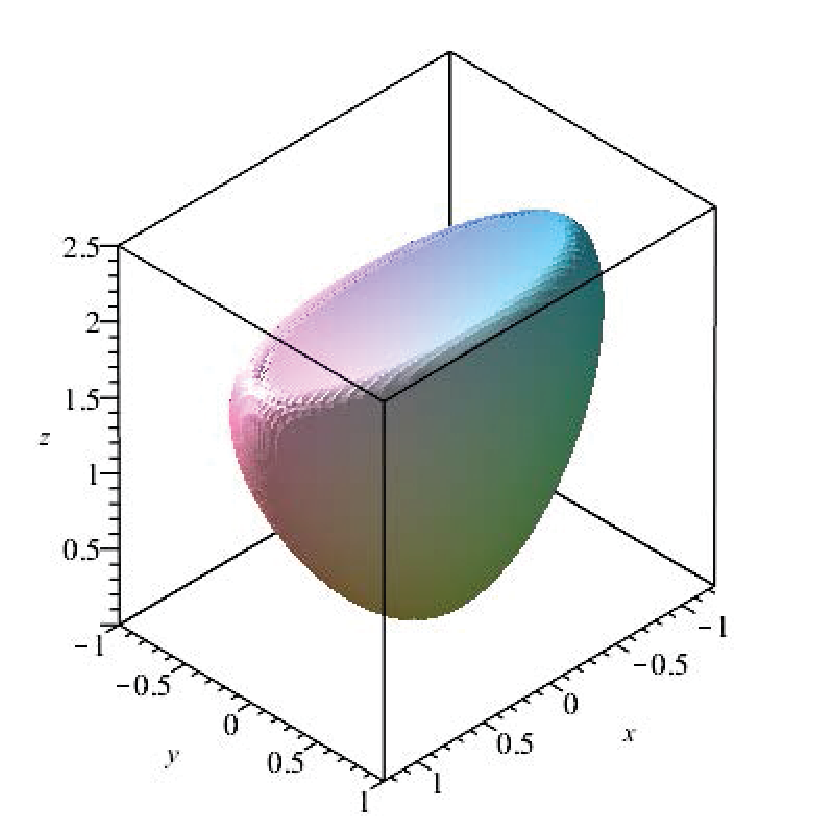}
    \caption*{\ref{fig:figB07}(b) the full shell  }
  \end{subfigure}
\phantomcaption
\plabel{fig:figB07}
\end{figure}

[This means removing a sufficiently large, at most $1$ dimensional, subset from the surface, then taking the
closure.
Nevertheless, it might be interesting to see certain details regarding how the $p_{\mathrm K}(x,y,z)=0$
can be bigger than the closure of the smooth part:
First of all, the lift of the boundary of numerical range does not occur (corresponding to $\lambda_1+\mathrm i\lambda_2\sim\infty$), although it gets restored
by taking the closure.
Another kind of phenomenon occurs when we restrict to say, to the line $x+2z-3=0$ on the plane $y=0$
(cf. Figure \ref{fig:figB06}(a)). Its points are obtained from the enveloping construction as follows:
the lower crossing point $(x,y,z)=(1,0,1)$ comes from $(\lambda_1,\lambda_2,N)=(1,0,0)$,
 $(x,y,z)=\left(\frac13,0,\frac43\right)$ comes from $(\lambda_1,\lambda_2,N)=\left(\frac45,0,\frac{36}{25}\right)$,
the segment between
 $(x,y,z)=\left( -\frac{12}{13}, 0, \frac{51}{26}  \right)$
and
 $(x,y,z)=(0,0,\frac32)$
comes from $(\lambda_1,\lambda_2,N)=\left(-\frac14,\lambda_2,(\lambda_2)^2+\frac{25}{16}\right)$ as
$(x,y,z)=\left(\frac{-12}{ 16\,({\lambda_{{2}}})^{2}+13 },0,\frac32\,{\frac {
16\,({\lambda_{{2}}})^{2}+17}{16\,({\lambda_{{2}}})^{2}+13}}\right)$ with $\lambda_2\sim\infty$ corresponding to the second end point.
But this second end point $(x,y,z)=(0,0,\frac32)$ is also produced by other $\lambda_1+\mathrm i\lambda_2\sim\infty$
(as the ``middle eigenvalue branch'').
There is also a general singularity corresponding to $\lambda_1+\mathrm i\lambda_2=0$, best to be ignored.
These points are perturbable, so they will remain in any topological closure.

Somewhat more systematically: Irregularities in the branches of $N$ are coming from the
discriminant of $F_A(\lambda_1,\lambda_2,\nu)$ in $\nu$.
In the present case this yields only $\lambda_1=\lambda_2=0$.
Thus, no proper singularity analysis is needed, only the planes $z=0$ and $z=2$ may be critical.
Every other singularity plus the limit of $\lambda_1+\mathrm i\lambda_2\rightarrow \infty$ will be removed if the
lift of the discriminant of  $p_{\mathrm K}(x,y,z)$ in $z$ is removed.
Cf. Figure \ref{fig:figB08}(a);
the axes belong to the discriminant set, the boundary of the numerical range is also transparent.
(Regarding the smooth part: The obstruction to computing $(\lambda_1,\lambda_2)$ from $(x,y,z)$
is $\partial_3 p_{\mathrm K}(x,y,z)=1$, and that is taken care by the discriminant.)
But this is an overkill, as  by conformal transformations in $A$ it can be reduced further
(corresponding to removing $\lambda_1+\mathrm i\lambda_2\sim \infty$ from around the ideal line).
Ultimately,
while there are some singularities in the lift of the inner leaves of the discriminant set in $z$
(cf. Figure \ref{fig:figB08}),
it turns out that in order obtain the closure of smooth enveloping construction it is sufficient to remove the three lines
seen in the $y=0$ section, Figure \ref{fig:figB06}(a), and take closure (which restores parts of the lines).
Then we are ready to take convex closure.]

\begin{figure}[H]
  \centering
  \begin{subfigure}[b]{0.4\linewidth}
    \includegraphics[width=2.5in]{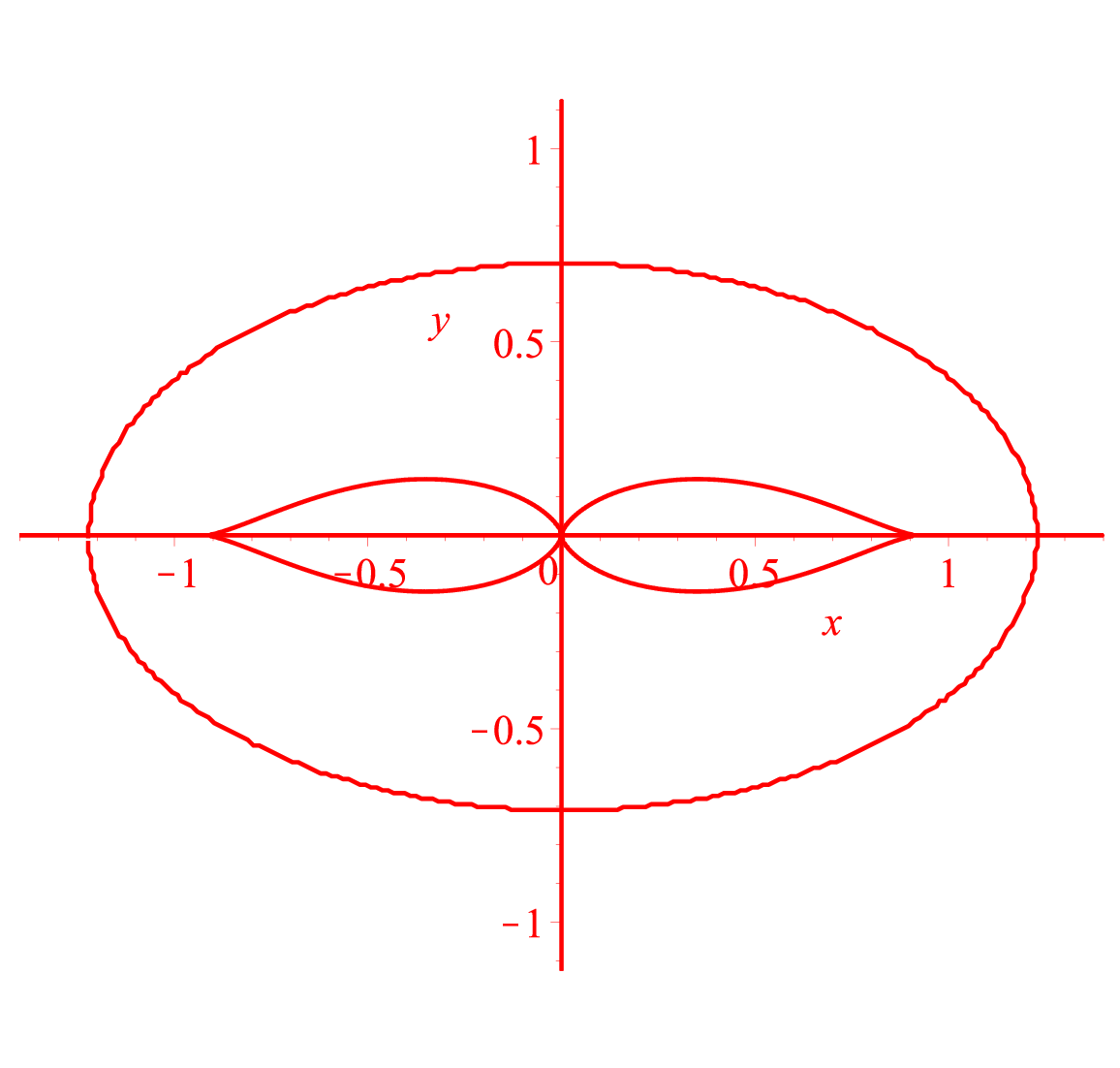}
    \caption*{Fig. \ref{fig:figB08}(a) discriminant in $z$}
  \end{subfigure}
  \qquad
  \begin{subfigure}[b]{0.4\linewidth}
    \includegraphics[width=2.5in]{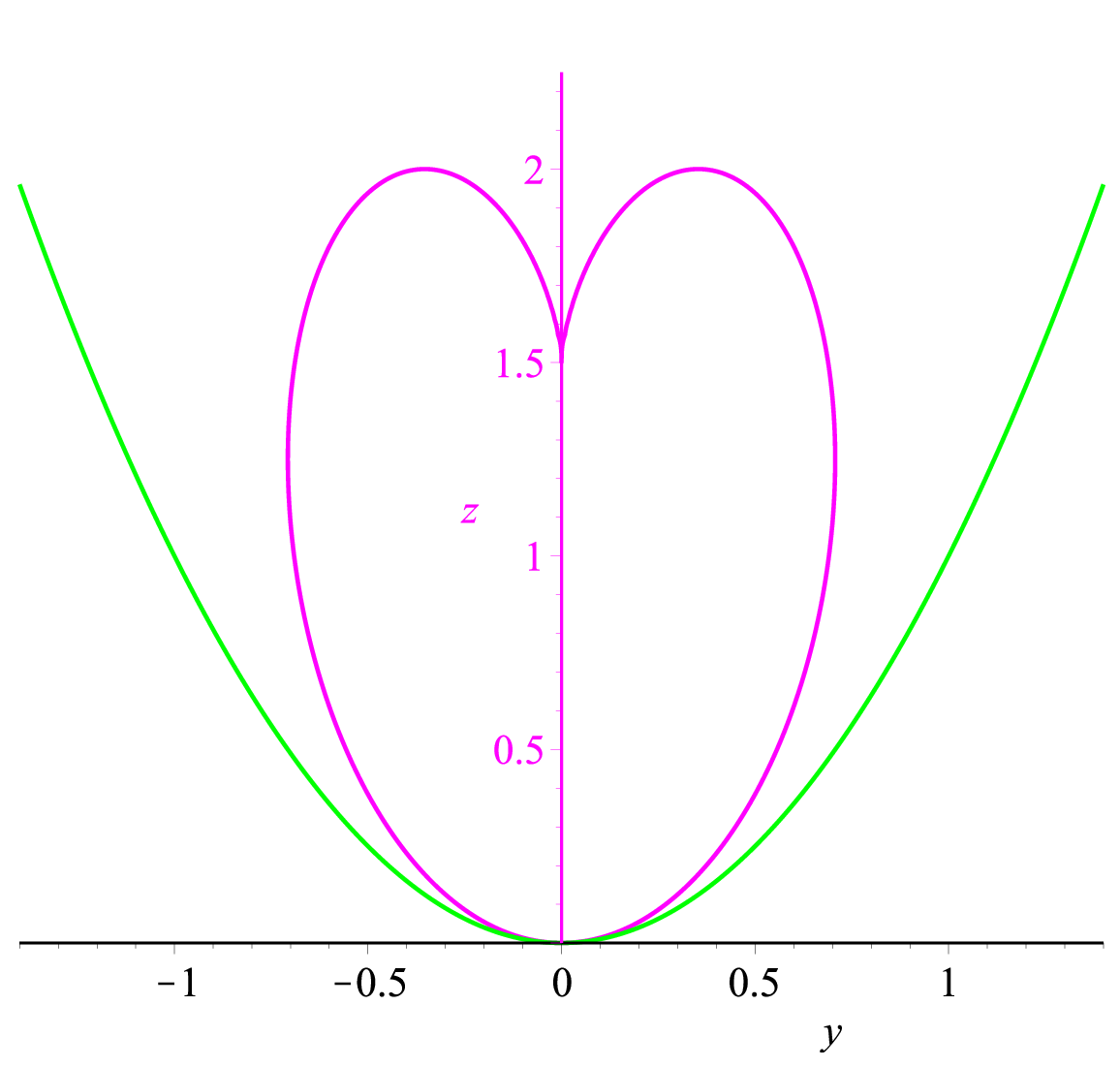}
    \caption*{\ref{fig:figB08}(b) section $x=0$  }
  \end{subfigure}
\phantomcaption
\plabel{fig:figB08}
\end{figure}
\begin{commentx}
\begin{figure}[H]
   \ContinuedFloat
   \begin{subfigure}[b]{0.4\linewidth}
    \includegraphics[width=2.5in]{figB08c}
    \caption*{\ref{fig:figB08}(cc) discriminant in $y$}
  \end{subfigure}
  \qquad
  \begin{subfigure}[b]{0.4\linewidth}
    \includegraphics[width=2.5in]{figB08d}
    \caption*{\ref{fig:figB08}(dd) discriminant in $x$}
  \end{subfigure}
\phantomcaption
\end{figure}
\end{commentx}
\begin{figure}[H]
  \ContinuedFloat
  \centering
  \begin{subfigure}[b]{0.4\linewidth}
    \includegraphics[width=2.5in]{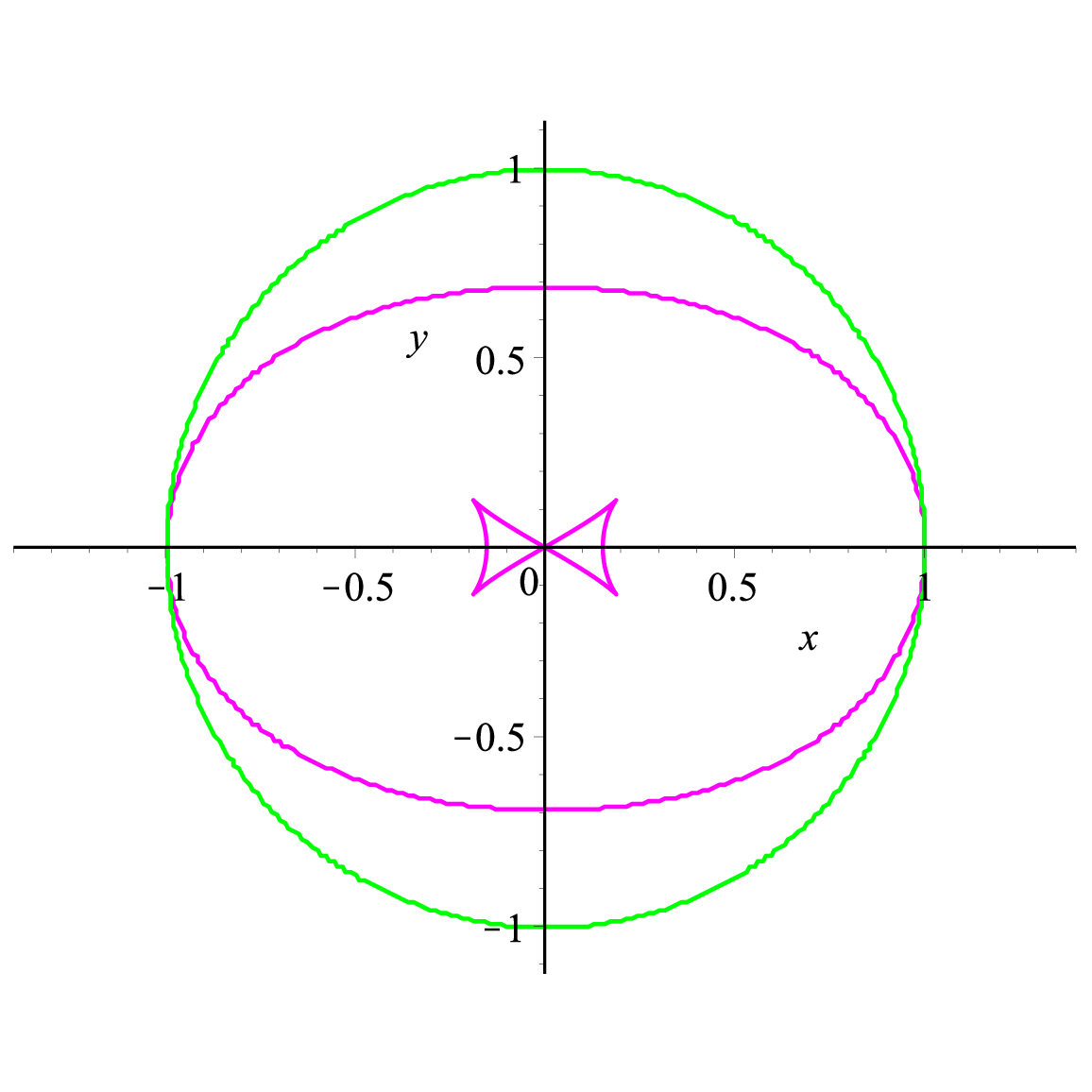}
    \caption*{Fig. \ref{fig:figB08}(c) section $z=1$}
  \end{subfigure}
  \qquad
  \begin{subfigure}[b]{0.4\linewidth}
    \includegraphics[width=2.5in]{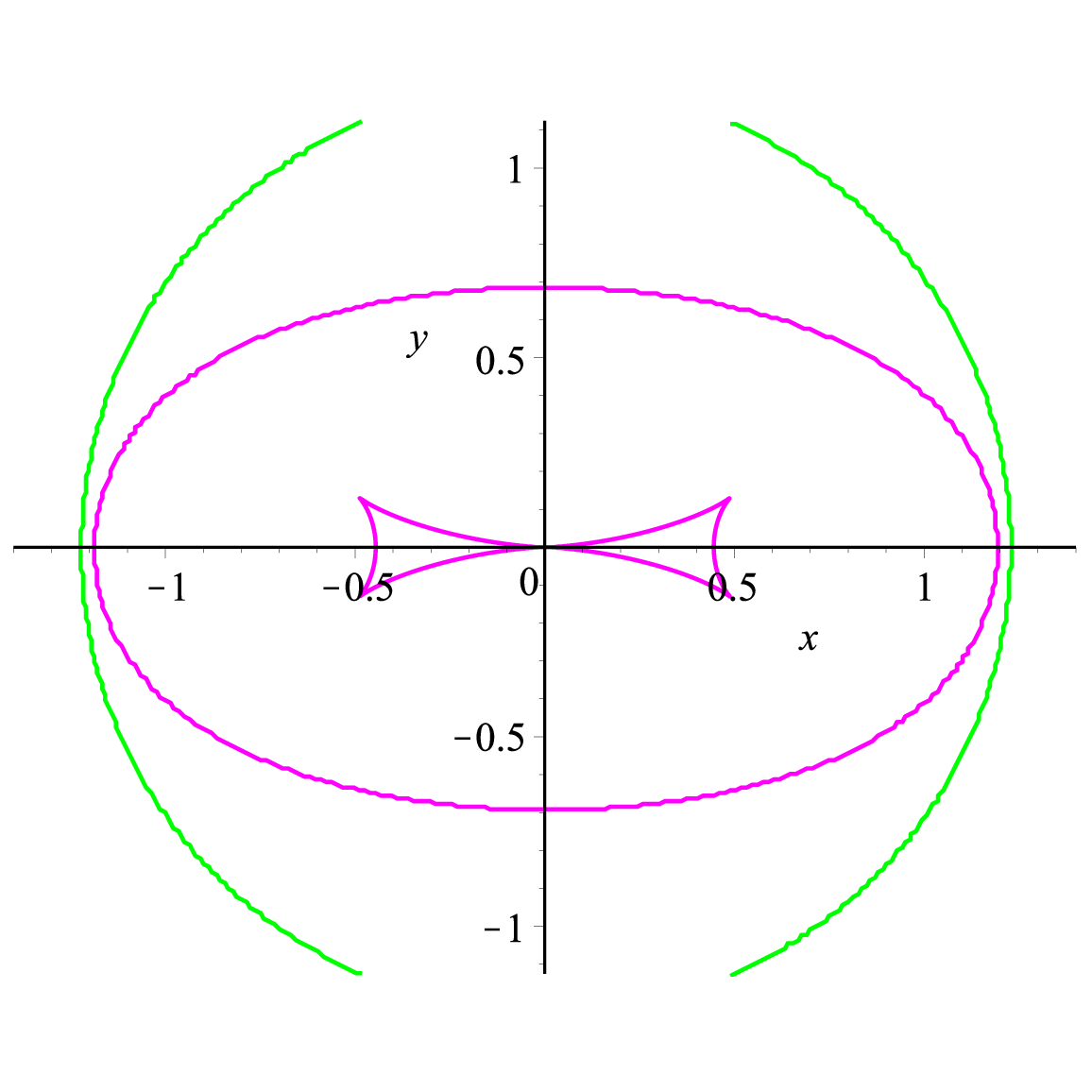}
    \caption*{\ref{fig:figB08}(d) section $z=3/2$  }
  \end{subfigure}
\phantomcaption
\end{figure}
\begin{figure}[H]
   \ContinuedFloat
   \begin{subfigure}[b]{0.4\linewidth}
    \includegraphics[width=2.5in]{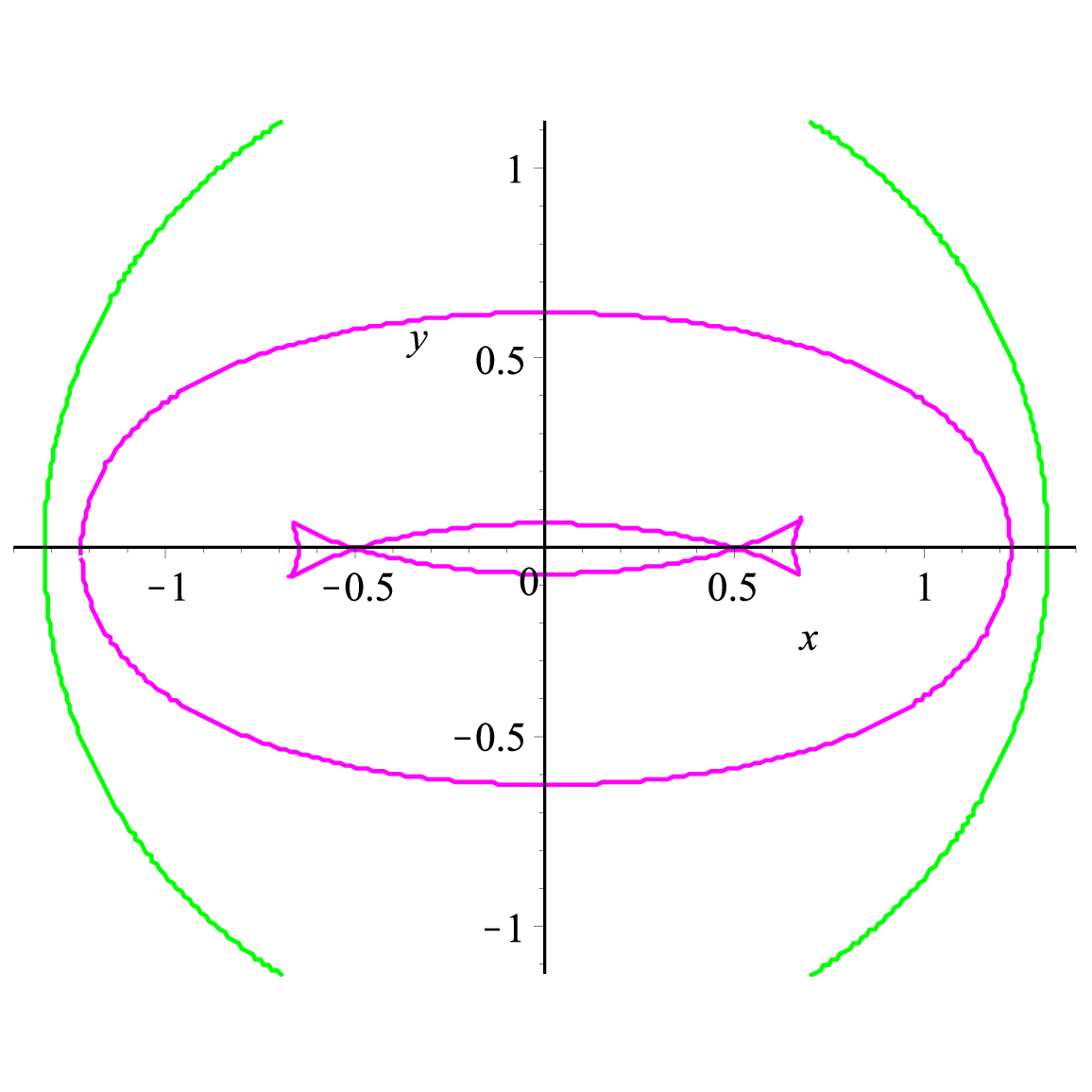}
    \caption*{\ref{fig:figB08}(e) section $z=7/4$}
  \end{subfigure}
  \qquad
  \begin{subfigure}[b]{0.4\linewidth}
    \includegraphics[width=2.5in]{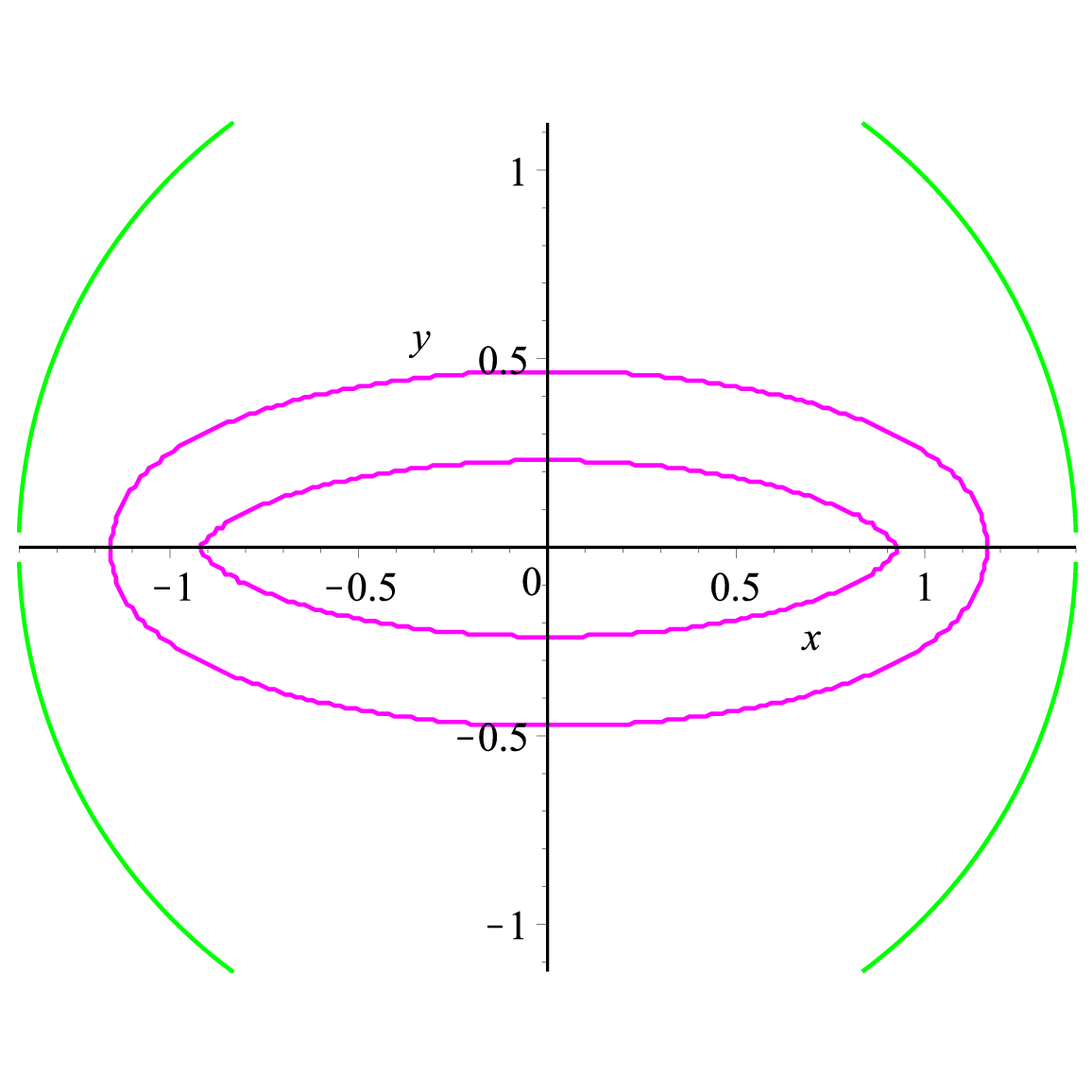}
    \caption*{\ref{fig:figB08}(f) section $z=51/26$}
  \end{subfigure}
\phantomcaption
\end{figure}

Now we take convex closure, in order to obtain $\DW_{\mathrm{CKB(P)}}(A)$.
This may be complicated in general, but in the present case this is geometrically simple:
\[p_{\mathrm K}(x,y,2)=
2\, \left( {x}^{2}+9\,{y}^{2} \right)  \left( {x}^{2}+{y}^{2}+2\,x+1
 \right)  \left( {x}^{2}+{y}^{2}-2\,x+1 \right)  \left( 8\,{x}^{2}+72
\,{y}^{2}-9 \right) ^{2}\]
indicates that (artificial parts ignored) the surface must by augmented by the
convex closure of the ellipse
\[\left\{x,y,z\,:\, \frac{x^2}{\frac98} +\frac{y^2}{\frac18}=1 ,z=2\right\},\]
yielding the full Davis-Wielandt shell as in Figure \ref{fig:figB07}(b).

In particular, we see that the direct algebraic description of the
shell is cumbersome.
For larger matrices is better proceed otherwise:
If one is interested in the visualization of the full shell (as in  Figure \ref{fig:figB07}(b))
then it is better take several supporting half-spaces using norms an co-norms, and intersect them.
If one is interested in the enveloping construction (as in  Figure \ref{fig:figB07}(a)), then it is
better to proceed parametrizing by $\lambda_1,\lambda_2$; but possibly using multiple patches (up to M\"obius transformation)
in order to obtain the points with vertical tangents spaces more precisely.
\qedexer
\end{example}
\begin{commentx}
\begin{remark}
If one is to work only with the dual of the surface, then using projective
coordinates one can proceed as follows.
Making the correspondence
\[-2\lambda_1 x_{\mathrm{CKB(P)}} -2\lambda_2 y_{\mathrm{CKB(P)}}+z_{\mathrm{CKB(P)}}+\lambda_1^2+\lambda_2^2-N(A-(\lambda_1+\mathrm i\lambda_2)\Id)=0\]
\[\qquad\sim\qquad ux_{\mathrm{CKB(P)}}+vy_{\mathrm{CKB(P)}}+w_{\mathrm{CKB(P)}}z+t=0,\]
via
\[\det \left( \left(\frac{-u}{2w}\right)^2+\left(\frac{-v}{2w}\right)^2-\left(\frac{t}{w}\right) -
\left(A-\left(\frac{-u}{2w}+\frac{-v}{2w}\mathrm i  \right)\right)^*\left(A-\left(\frac{-u}{2w}+\frac{-v}{2w}\mathrm i  \right)\right)
\right)=0\]
one obtains
\[\det\left(u\frac{A+A^*}{2}+v\frac{A-A^*}{2\mathrm i}+wA^*A+t\Id \right)=0, \]
which already contains the supporting planes of the ordinary numerical range (corresponding to $w=0$).
\end{remark}
\end{commentx}
\snewpage

\end{document}